\documentclass[reqno,english]{amsart}
\usepackage{amsmath,amssymb,amsfonts,array,color}
\usepackage{mathtools}
\usepackage[colorlinks=true]{hyperref}
\definecolor{referencecolor}{rgb}{0.48,0.12,0.15}
\definecolor{citationcolor}{rgb}{0.09,0.35,0.19}
\definecolor{urlcolor}{rgb}{0.35,0.14,0.42}
\hypersetup{
 pdftitle={The regular N-gon maximizes the first positive Neumann eigenvalue
  among simple polygons with at most N sides and fixed area},
 pdfauthor={Senping Luo and Bingbing Xie},
 linkcolor=referencecolor,
 citecolor=citationcolor,
 urlcolor=urlcolor
}

\numberwithin{equation}{section}
\let\oldbibliography\thebibliography
\renewcommand{\thebibliography}[1]{%
  \oldbibliography{#1}%
  \setlength{\itemsep}{0pt}%
}

\newtheorem{theorem}{Theorem}[section]
\newtheorem{lemma}[theorem]{Lemma}
\newtheorem{proposition}[theorem]{Proposition}
\newtheorem{numericallemma}[theorem]{Numerical Lemma}
\newtheorem{corollary}[theorem]{Corollary}

\newcommand{\tr}{\operatorname{tr}}
\newcommand{\diag}{\operatorname{diag}}

\title[The first Neumann eigenvalue of simple polygons]{The regular $N$-gon maximizes first positive Neumann eigenvalue among simple polygons with at most $N$ sides and fixed area}

\author[S.~Luo]{Senping Luo}
\author[B.~Xie]{Bingbing Xie}

\address[S.~Luo and B.~Xie]{School of Mathematics and Statistics,
Jiangxi Normal University, Nanchang 330022, China, and Jiangxi Applied
Mathematical Research Center}
\email[S.~Luo]{luosp1989@163.com}
\email[B.~Xie]{xiebb2074@163.com}
\begin{document}
\raggedbottom
\setlength{\emergencystretch}{1em}
\begin{abstract}
Let $\Phi(P)=|P|\mu_1(P)$, where $\mu_1(P)$ is the first positive Neumann
eigenvalue of $P$. For every $N\ge3$, the regular $N$-gon $R_N$ uniquely
maximizes $\Phi$, up to congruence, among simple polygons of fixed area
with at most $N$ sides. The values $\Phi(R_N)$ increase strictly with $N$.
The maximization theorem answers Laugesen's question
\cite[Problem~(1), p.~405]{ABLW} and proves the polygonal
Szeg\H{o}--Weinberger conjecture of Laugesen and Siudeja
\cite[Conjecture~6.36]{LSH}. It also implies that the triangular lattice
maximizes the first positive Neumann eigenvalue of the Voronoi cell at
fixed covolume.

\end{abstract}
\maketitle

\section{Introduction}
Let \(\Omega\subset\mathbb R^2\) be a bounded connected Lipschitz domain
with outward unit normal \(\nu\). The Dirichlet and Neumann spectra are
defined by
\begin{equation}\label{neumann-boundary-value-problem}
 \begin{cases}-\Delta u=\lambda u&\text{in }\Omega,\\u=0&\text{on }\partial\Omega,\end{cases}
 \qquad
 \begin{cases}-\Delta v=\mu v&\text{in }\Omega,\\\partial_\nu v=0&\text{on }\partial\Omega.\end{cases}
\end{equation}
The eigenvalues, counted with multiplicity, satisfy
\[
 0<\lambda_1<\lambda_2\le\cdots,\qquad
 0=\mu_0<\mu_1\le\mu_2\le\cdots .
\]
For the disk \(B\) of area \(|\Omega|\), the Faber--Krahn inequality gives
\(\lambda_1(B)\le\lambda_1(\Omega)\). Szeg\H{o} \cite{Szego}
proved \(\mu_1(\Omega)\le\mu_1(B)\) for simply connected planar domains.
Weinberger \cite{Weinberger} established the corresponding fixed-volume
ball inequality for general domains in Euclidean space. We use the
normalization
\begin{equation}\label{neumann-first-rayleigh-quotient}
 \Phi(\Omega)=|\Omega|\mu_1(\Omega),\qquad
 \mu_1(\Omega)=\min_{\substack{u\in H^1(\Omega)\setminus\{0\}\\
                      \int_\Omega u\,dx=0}}
 \frac{\int_\Omega |\nabla u|^2\,dx}{\int_\Omega |u|^2\,dx}.
\end{equation}
Thus \(\mu_1(t\Omega)=t^{-2}\mu_1(\Omega)\), \(\Phi(t\Omega)=\Phi(\Omega)\)
for \(t>0\). Let \(R_N\) denote a regular \(N\)-gon. A simple polygon is the
bounded domain enclosed by a polygonal Jordan curve.

\begin{theorem}\label{main-polygonal-comparison}
For every integer \(N\ge3\), the regular \(N\)-gon uniquely maximizes
the first positive Neumann eigenvalue among simple polygons of the
same area with at most \(N\) sides, up to congruence.
The values on regular polygons increase strictly with \(N\), that is,
\[
 \Phi(R_N)<\Phi(R_{N+1})\qquad(N\ge3).
\]
\end{theorem}

\begin{corollary}\label{lattice-voronoi}
Among planar lattices of fixed covolume, the triangular lattice uniquely
maximizes the first positive Neumann eigenvalue of the Voronoi cell,
up to orthogonal transformations.
\end{corollary}

Deng and Luo \cite[Theorem~1.1]{DL} proved that, among three-dimensional
lattices of fixed covolume, the body-centered cubic lattice uniquely
maximizes the first positive Neumann eigenvalue of the Voronoi cell, up to
orthogonal transformations.

In the Oberwolfach report of the workshop organized by Ashbaugh, Benguria,
Laugesen, and Weidl
\cite[Problem~(1), p.~405]{ABLW}, Laugesen asked whether the regular
\(N\)-gon maximizes \(\Phi\) among \(N\)-gons for \(N\ge4\).
Laugesen and Siudeja
\cite[Conjecture~6.36]{LSH} formulated this question as the polygonal
Szeg\H{o}--Weinberger conjecture. Antunes and Henrot \cite[Section~3(c),
Figures~14--15]{AH} computed the ranges of the first two positive Neumann
eigenvalues for several unit-area polygon classes. Their computations gave
evidence for maximality of the square among convex quadrilaterals and of
the regular octagon among convex octagons.

Laugesen and Siudeja \cite{LS} proved that the equilateral triangle
uniquely maximizes the first positive Neumann eigenvalue among triangles
of prescribed perimeter and deduced the fixed-area inequality. Endo and
Osting \cite{EO} proved the fixed-area inequality for convex
quadrilaterals, with equality only for the square. Dong and Wang
\cite[Theorem~1]{DW} proved in their preprint that the regular pentagon
uniquely maximizes the first positive Neumann eigenvalue among simple
pentagons of fixed area, including pentagons with reflex or straight
angles.

P\'olya \cite{Polya1948} developed symmetrization arguments for torsional
rigidity, principal frequency, and electrostatic capacity. P\'olya and
Szeg\H{o} \cite[p.~158]{PS} recorded the sharp Dirichlet and torsion
inequalities for triangles and quadrilaterals and posed the corresponding
polygonal problem for general \(N\). Bogosel and Bucur \cite{BB} reduced
the first Dirichlet eigenvalue inequality for fixed \(N\) to finitely many
certified checks. Cheng, Gui, Hu, and Li \cite{CGHL} proved that the
regular \(N\)-gon uniquely minimizes the first Dirichlet eigenvalue among
convex polygons of fixed area with at most \(N\) sides for all
sufficiently large \(N\). Wang \cite{Wang} proved that the regular
\(N\)-gon uniquely minimizes the first Dirichlet eigenvalue and uniquely
maximizes torsional rigidity among simple polygons of fixed area with at
most \(N\) sides for every \(N\ge5\).

Gui, Hu, and Li \cite{GHL} bounded the Dirichlet torsion deficit in terms
of the perimeter ratio for tangential polygons of fixed area and deduced
that the regular polygon uniquely maximizes torsional rigidity in this
class. Gui, Hu, Li, and Yang \cite{GHLY} obtained an explicit deficit
bound for Robin torsional rigidity in the same class and proved that the
regular polygon uniquely maximizes it for each positive Robin parameter.
Solynin and Zalgaller \cite{SZ} proved that the regular \(N\)-gon uniquely
minimizes logarithmic capacity among simple
\(N\)-gons of fixed area.

For regular polygons of fixed area, Dahne, G\'omez-Serrano, and
Pech-Alberich \cite{DGP} proved strict decrease of the first Dirichlet
eigenvalue with the number of sides. Cheng, Gui, Hu, Li, and Yao
\cite{CGHLY} proved strict increase of the first positive Steklov
eigenvalue at fixed perimeter and derived its expansion through order
\(N^{-5}\). Gui, Hu, Li, and Zhang \cite{GHLZ} proved strict increase of
Dirichlet torsional rigidity at fixed area and obtained its expansion
through order \(N^{-4}\).

Berghaus, Georgiev, Monien, and Radchenko
\cite{BGMR} proved an all-orders expansion of the first Dirichlet
eigenvalue of a regular polygon with fixed area, with coefficients
expressed in terms of multiple zeta values. Berghaus, Jones, Monien, and
Radchenko
\cite[Section~8.1]{BJMR} computed Laplacian eigenvalues on dihedrally
symmetric planar domains to high precision. For regular polygons their
computations gave evidence for an expansion of the Neumann branch
converging to \(j_{1,1}^2\), where \(j_{1,1}\) is the first positive zero
of \(J_1\).

Ni and Wang \cite{NiWang} used the intrinsic diameter and examples to
study the dependence of the first positive Neumann eigenvalue on the
geometry of planar domains. Freitas and Kennedy \cite{FK} obtained the
sharp comparison
\[
 \mu_1(\Omega_1)\ge\frac{\pi^2}{4j_{0,1}^2}\mu_1(\Omega_2)
 \qquad (\Omega_1\subset\Omega_2)
\]
for bounded convex planar domains, where \(j_{0,1}\) is the first positive
zero of \(J_0\). Cavallina, Funano, Henrot, Lemenant, Lucardesi, and
Sakaguchi \cite{CFLHLS} studied existence for two convex planar shape
optimization problems, minimizing Neumann eigenvalues inside a prescribed
box and maximizing them over domains containing a prescribed obstacle.
Huang, Li, Li, and Yao \cite{HLLY} showed that the best Rayleigh constant
under the constraint of zero boundary integral is at most the first
positive Neumann eigenvalue. They proved equality for balls, rectangular
boxes, and equilateral triangles. Ashbaugh and Levine \cite{AL} proved
\(\mu_k(\Omega)\le\lambda_k(\Omega)\) for \(k\ge1\) on smooth spherical
domains with nonnegative boundary mean curvature. A hemisphere gives
equality for \(k=1\).

Let \(p(P)\) be the number of interior angles of \(P\) smaller than
\(\pi\), and put \(x=(j'_{1,1})^2\), the square of the first positive
zero of \(J_1'\). Sections~\ref{part:spectrum} and~\ref{sec:nr-nonconvex} give
\[
 \begin{aligned}
 p(P)<N\quad&\Longrightarrow\quad\Phi(P)<\Phi(R_N)\qquad(N\ge5),\\
 \lim_{N\to\infty}\Phi(R_N)&=\pi x,\\
 \Phi(R_N)-\sup_{p(P)<N}\Phi(P)
   &=\frac{12\pi x\zeta(3)}{(x-1)N^4}+O(N^{-5}).
 \end{aligned}
\]
The supremum extends over simple polygons with arbitrarily many sides.
Sections~\ref{part:finite} and~\ref{part:uniform} compare a convex
polygon with \(R_N\) sector by sector. With vertices \(p_1,\ldots,p_N\)
in cyclic order and \(p_{N+1}=p_1\), the sectors are
\[
 S_j=\operatorname{conv}\!\left\{
 \frac1{|P|}\int_Px\,dx,p_j,p_{j+1}\right\}
\]
with centroid angles \(\vartheta_j\). The linear terms cancel by
\[
 \sum_{j=1}^N\left(\frac{N|S_j|}{|P|}-1\right)
 =\sum_{j=1}^N\left(\frac{N\vartheta_j}{2\pi}-1\right)=0.
\]
For \(5\le N\le32\), the finite inequalities use the monomial integrals,
integer matrix recursions and interval bounds in
Sections~\ref{part:spectrum}--\ref{part:finite} and Appendix~\ref{app:remaining}.
For \(N\ge33\), the bounds are uniform in \(\pi/N\).

\section{Regular polygons and conformal Rayleigh bounds}\label{part:spectrum}

We identify the first eigenspace of a regular polygon and give bounds
for the hexagon's low eigenvalues. Conformal trials yield the comparison
estimates and the expansion of \(\Phi(R_N)\) through order \(N^{-5}\).

\subsection{The first eigenspace and the hexagon bounds}

\begin{lemma}\label{affine-all-N}
For \(N\ge4\), the first positive Neumann eigenspace of \(R_N\)
is the standard two-dimensional representation of \(D_N\).
\end{lemma}

\begin{proof}
On \([-1,1]^2\), the first eigenspace is
\(\operatorname{span}\{\sin(\pi x_1/2),\sin(\pi x_2/2)\}\).
For \(N\ge5\) and circumradius one, the radial map on the sector with
outward normal at angle \(\alpha\) is
\[
 h(\theta)=\frac{\cos(\pi/N)}{\cos(\theta-\alpha)},\qquad
 (r,\theta)\longmapsto(rh(\theta),\theta),\qquad 0<h^2\le1.
\]
In the coordinates \((\partial_r,r^{-1}\partial_\theta)\) the energy matrix is
\[
 B=\begin{pmatrix}1+(h'/h)^2&-h'/h\\-h'/h&1\end{pmatrix},\qquad
 \det B=1,\quad \tr B\le2+\tan^2(\pi/5)<0.49+1/0.49 .
\]
Hence \(B\ge0.49I\). By disk separation, angular monotonicity and
\(j'_{1,1}<j_{1,1}<j'_{1,2}\),
\[
 \mu_3(\mathbb D)\ge\min\{j_{1,1}^2,(j'_{2,1})^2\}>9.3.
\]
For the last inequality, use \(s=z^2/4\in[0,93/40]\) and the alternating
bounds
\[
 \begin{aligned}
 2J_1(z)/z&\ge1-s/2+s^2/12-s^3/144>0,\\
 4J'_2(z)/z&\ge1-2s/3+s^2/8-s^3/90+s^4/1728-s^5/50400>0 .
 \end{aligned}
\]
Both polynomials decrease and are positive at \(93/40\).
The trial \((r-r^3/3)\cos\theta\) gives \(\mu_1(\mathbb D)\le112/33\).
Use Szeg\H{o}--Weinberger, \(|R_N|=\tfrac N2\sin(2\pi/N)>47/20\) and min--max:
\[
 \mu_1(R_N)\le\frac{112\pi}{33|R_N|}<4.55
 <0.49(9.3)<0.49\mu_3(\mathbb D)\le\mu_3(R_N).
\]
For a first eigenfunction in the other rotation classes, its
pullback \(v\) satisfies \(\int_{\mathbb D}h^2v=0\) and is
orthogonal to the disk dipoles. Thus
\[
 \int_{\mathbb D}h^2|v|^2
 =\min_{b\in\mathbb C}\int_{\mathbb D}h^2|v-b|^2
 \le\min_{b\in\mathbb C}\int_{\mathbb D}|v-b|^2
 \le\frac1{\mu_3(\mathbb D)}\int_{\mathbb D}|\nabla v|^2.
\]
Thus \(v\ne0\) would give \(\mu_1(R_N)\ge0.49\mu_3(\mathbb D)\).
Complex conjugation and \(\mu_1(R_N)<\mu_3(R_N)\) give one dimension
in each of the classes \(e^{\pm2\pi i/N}\).
\end{proof}

For \(H=\operatorname{conv}\{(0,\pm1),(\pm\sqrt3/2,\pm1/2)\}\),
sector integration, the diameter bound \cite{PW}, and the trial
\(x_1\) give
\begin{equation}\label{rough}
 |H|=\frac{3\sqrt3}{2},\qquad
 \frac{1}{|H|}\int_H xx^T\,dx=\frac5{24}I,\qquad
 \frac{\pi^2}{4}\leq\mu_1(H)\leq\frac{24}{5}.
\end{equation}

By Lemma~\ref{affine-all-N}, the first pair can be chosen with
\[
 |H|^{-1}\int_H UU^{\mathsf T}=I,\qquad U(Ox)=OU(x)\quad(O\in D_6),
 \qquad u_j(x-2x_i e_i)=(-1)^{\delta_{ij}}u_j(x).
\]

On the bisector sector \(S\) of the unit-circumradius \(R_N\),
use \((x,y)=(X\cos(\pi/N),Y\sin(\pi/N))\),
\(0\le X\le1\), \(|Y|\le X\). Then
\begin{equation}\label{certified-sector-monomials}
 \frac{1}{|R_N|}\int_S X^iY^j
 =\frac1N\int_0^1\int_{-X}^X X^iY^j\,dY\,dX
 =\frac{1+(-1)^j}{N(j+1)(i+j+2)}\qquad(i,j\in\mathbb Z_{\ge0}).
\end{equation}
For \(A=A^{\mathsf T}\), \(A_{ij}\le0\) \((i\ne j)\), and diagonal \(M\succ0\),
\begin{equation}\label{discrete-positive-test}
 v>0,\quad(A-\ell M)v>0
 \quad\Longrightarrow\quad \lambda_{\min}(A,M)>\ell.
\end{equation}
Since \(|u|^{\mathsf T}A|u|\le u^{\mathsf T}Au\) and
\(|u|^{\mathsf T}M|u|=u^{\mathsf T}Mu\), take a minimizing eigenvector \(u\ge0\):
\[
 0<u^{\mathsf T}(A-\ell M)v=(\lambda_{\min}-\ell)u^{\mathsf T}Mv.
\]
If \(u\) is a unit trial in a spectral subspace and \(w\) its component
orthogonal to the first eigenspace, then, for
\(0\le\lambda_-\le\lambda_1\) and \(\lambda_-<\Lambda\le\lambda_{\rm next}\),
\begin{equation}\label{spectral-projection-bounds}
 \|w\|^2\leq\frac{\mathfrak a(u,u)-\lambda_-}{\Lambda-\lambda_-},\qquad
 \mathfrak a(w,w)\leq\Lambda\frac{\mathfrak a(u,u)-\lambda_-}{\Lambda-\lambda_-}.
\end{equation}
Indeed
\[
 \mathfrak a(u,u)-\lambda_-\ge\mathfrak a(w,w)-\lambda_-\|w\|^2
 \ge\max\left\{(\Lambda-\lambda_-)\|w\|^2,
       (1-\lambda_-/\Lambda)\mathfrak a(w,w)\right\}.
\]
Decimals denote rationals. Cancel common vanishing factors before
rounding intervals outward. The mixed bounds use the integer recursions
in Appendix~\ref{app:finite-record}.
\begin{numericallemma}[Spectrum of the unit-circumradius hexagon]
\label{reference-numerical}
The regular hexagon \(H\) of circumradius one satisfies
\begin{equation}\label{reference-values}
 \mu_1(H)>4.04320119,\qquad
 \mu_3(H)>10,\qquad
 10.5045448<\Phi(H)<10.50465 .
\end{equation}
\end{numericallemma}

\begin{proof}
\emph{(1) Lower bounds.}
On the equilateral mesh of side \(h=1/N\), let \(I_hu\) be the
Crouzeix--Raviart interpolant and set \(w=u-I_hu\) \cite{CR,Liu}.
For each triangle \(T\) and each edge \(e\), integration by parts gives
\[
 \begin{gathered}
 \int_e w=0\ \Longrightarrow\ \int_T\nabla w=0,\qquad
 2\int_Tw=-\int_T\left(x-\frac1{|T|}\int_T y\,dy\right)\cdot\nabla w\,dx,\\
 \int_T|\nabla u|^2=\int_T|\nabla I_hu|^2+\int_T|\nabla w|^2 .
 \end{gathered}
\]
The vector \(x-|T|^{-1}\int_T y\,dy\) has constant normal component
on each edge. By Cauchy--Schwarz and \cite{McC},
\[
 \begin{aligned}
 \frac1{|T|}\left(\int_T w\right)^2
 &\le\frac1{4|T|}\int_T\left|x-\frac1{|T|}\int_T y\,dy\right|^2dx
       \int_T|\nabla w|^2=\frac{h^2}{48}\int_T|\nabla w|^2,\\
 \int_T w^2
 &=\int_T\!\left(w-\frac1{|T|}\int_T w\right)^2
    +\frac1{|T|}\left(\int_Tw\right)^2
 \le h^2\left(\frac9{16\pi^2}+\frac1{48}\right)\int_T|\nabla w|^2 .
 \end{aligned}
\]
Set \(C_h^2=h^2(9/(16\pi^2)+1/48)\). In the first \(k+1\) continuous
modes choose \(u\ne0\) with \(I_hu\) orthogonal to the first \(k\)
broken modes. With \(\nabla_h\) the elementwise gradient,
\[
 \|u\|_2\le
 \frac{\|\nabla_h I_hu\|_2}{\sqrt{\lambda_{h,k}}}
 +C_h\|\nabla_h w\|_2
 \le(\lambda_{h,k}^{-1}+C_h^2)^{1/2}\|\nabla u\|_2 .
\]
Thus
\begin{equation}\label{broken-lower}
      \mu_k(H)\geq\frac{\lambda_{h,k}}{1+C_h^2\lambda_{h,k}},
           \qquad C_h^2<5h^2/64 .
\end{equation}
For related certified lower bounds, see \cite{CG,NSY}.
The mixed first mode obeys \eqref{broken-lower} because
\(u|_{\Gamma_D}=0\Rightarrow\int_e I_hu=0\) on \(e\subset\Gamma_D\).

For triangle-edge incidence \(B_{Te}=\mathbf1_{\{e\subset\partial T\}}\)
and \(D=\diag(\sum_T B_{Te})\), the basis \(1-2\ell_i\) in
barycentric coordinates gives
\[
 M_h=(|T|/3)D,\qquad
 A_h=\frac{8|T|}{3h^2}(3D-B^{\mathsf T}B).
\]
For the face adjacency Laplacian \(L\), the full mesh gives
\(2(3I-BD^{-1}B^{\mathsf T})=L\). On the upper half-mesh
orient the cut through opposite vertices horizontally and delete its
Dirichlet columns from \(B,D\). Then
\[
 2(3I-BD^{-1}B^{\mathsf T})=L+2\operatorname{diag}(\delta_i),\qquad
 \delta_i=\#\{\text{Dirichlet edges of triangle }i\}.
\]
The positive spectra of \(BD^{-1}B^{\mathsf T}\) and
\(D^{-1/2}B^{\mathsf T}BD^{-1/2}\) coincide. The broken spectrum below
\(24/h^2\) equals \(4/h^2\) times the face spectrum below \(6\), with
multiplicities. Odd reflection gives \(\lambda_1^{\rm mix}=\mu_1(H)\).

Take the lattice points \(\max(|i|,|j|,|i+j|)\le N\), restricted to
\(j\ge0\) for the half-hexagon. The faces are the admissible triples
\[
 T_{ijb}=\{(i+1,j),(i,j+1),(i+b,j+b)\},\quad b\in\{0,1\},
 \qquad T\sim T'\iff \#(T\cap T')=2 .
\]
On the half-mesh \(\delta_{ijb}=\mathbf1_{\{j=b=0\}}\).

For the \(3N^2=172800\) faces of the half-hexagon set
\[
 N=240,\qquad G=L+2\diag(\delta),\qquad \sigma=\frac{404322337}{4N^2\cdot10^8}.
\]
Order a face set \(E\) by lexicographic order for \(\#E\le32\).
Otherwise choose \(s=i\) if \(\operatorname{span}_E i\ge
\operatorname{span}_E j\), and \(s=j\) otherwise. For
\(t=\lfloor(\min_E s+\max_E s)/2\rfloor\), concatenate the recursive
orders of \(E\cap\{s<t\}\), \(E\cap\{s>t\}\) and the lexicographic
order of \(E\cap\{s=t\}\). From \(C^{(0)}=\lfloor 2^{58}(G-\sigma I)\rfloor\)
(entrywise), compute
\begin{equation}\label{hex-rounded-schur}
 c_k=C^{(k-1)}_{kk},\qquad
 C^{(k)}_{ij}=C^{(k-1)}_{ij}-
 \left\lceil\frac{C^{(k-1)}_{ik}C^{(k-1)}_{kj}}{c_k}\right\rceil
 \qquad(i,j>k).
\end{equation}
The integer minima are
\[
 \min_{k<172800}c_k=82430894203534619,\qquad c_{172800}=133204153439.
\]
For the exact Schur complements \(S^{(k)}\) of \(2^{58}(G-\sigma I)\)
the off-diagonal entries are nonpositive. At step \(k\),
\[
 C^{(k-1)}\le S^{(k-1)},\quad 0<c_k\le S^{(k-1)}_{kk}
 \quad\Longrightarrow\quad
 \frac{C_{ik}C_{kj}}{c_k}\ge\frac{S_{ik}S_{kj}}{S_{kk}}\quad(i,j>k).
\]
Here the matrix inequality is entrywise and the entries in the last
quotients are from step \(k-1\). Therefore \(C^{(k)}\le S^{(k)}\), and
\[
 G\succ\sigma I
 \quad\Longrightarrow\quad
 \lambda_{h,1}^{\rm mix}>4N^2\sigma
 \quad\overset{\eqref{broken-lower}}{\Longrightarrow}\quad
 \mu_1(H)>4.04320119,\qquad \Phi(H)>10.5045448.
\]

For the full hexagon set \(N=8\), \(A^{(0)}=25600L-1021 I\),
and retain lexicographic face order. Starting at \(k=1\), use
\[
 d_k=A^{(k-1)}_{kk},\qquad
 A^{(k)}_{ij}=A^{(k-1)}_{ij}
 -\frac{A^{(k-1)}_{ik}A^{(k-1)}_{kj}}{d_k}\quad(i,j>k).
\]
The half-bandwidth is \(31\), and the \(384\) rational pivots satisfy
\[
 \{k:d_k<0\}=\{94,197,307\},\qquad \min_k|d_k|>326 .
\]
Sylvester's inertia theorem and \eqref{broken-lower} give
\[
 \lambda_4(L)>1021/25600\quad\Longrightarrow\quad
 \lambda_{h,3}>10.21\quad\Longrightarrow\quad\mu_3(H)>10,
 \qquad\lambda_1(L)=0.
\]

\emph{(2) Upper bound.}
Order \(f_i=z^{p_i}\bar z^{q_i}\), \(1\le p_i+q_i\le19\),
\(p_i-q_i\equiv1\pmod6\), by \((p_i+q_i,p_i)\).
Their mass and energy matrices are
\[
 \begin{aligned}
 M_{ij}&=\frac1{|H|}\int_H z^{p_i+q_j}\bar z^{q_i+p_j}\,dx,\\
 K_{ij}&=\frac{2(p_ip_j+q_iq_j)}{|H|}
             \int_H z^{p_i+q_j-1}\bar z^{q_i+p_j-1}\,dx .
 \end{aligned}
\]
Set the energy term to zero when its prefactor vanishes.
For \(\int_H z^a\bar z^b\), rotation gives zero if
\(a-b\not\equiv0\pmod6\). Otherwise \(a+b\) is even and reflection
on each bisector kills odd \(Y\)-powers in
\eqref{certified-sector-monomials}. The remaining powers of \(X\) are
even, so all entries are rational.
Elimination of \(K-4.043240668M\) gives
\[
 d_1,\ldots,d_{36}>0,\qquad -31<10^6d_{37}<-30.
\]
The negative direction has zero integral by rotation. One real
component therefore gives \(\Phi(H)<(3\sqrt3/2)(4.043240668)<10.50465\).
\end{proof}

\subsection{Conformal Rayleigh bounds}\label{sec:conformal-rayleigh}

We use the conformal centering of Szeg\H{o} \cite{Szego}. For the
corresponding Neumann, Steklov and Robin trials, see \cite{GP,FL}.
The balancing argument is also given in \cite[Lemma~5]{FL}.

\begin{lemma}[Conformal dipole trials]\label{conformal-dipole}
Let \(P\) be a bounded simply connected Lipschitz domain and let
\(\phi\in C^1([0,1])\) satisfy \(\phi(0)=0<\phi(1)\).
Write \(V(re^{i\theta})=\phi(r)e^{i\theta}\). There is a conformal map
\(f:\mathbb D\to P\) such that
\begin{equation}\label{balance}
 \int_{\mathbb D}V(z)|f'(z)|^2\,dA(z)=0 .
\end{equation}
For any such map put
\[
 \begin{gathered}
 \frac{f'(z)}{f'(0)}=\sum_{j\ge0}b_jz^j,\quad b_0=1,\qquad
 \beta=\frac{\pi|f'(0)|^2}{|P|},\qquad
 S=\sum_{j\ge1}\frac{|b_j|^2}{j+1}=\beta^{-1}-1,\\
 E=\int_0^1\left(r\phi'(r)^2+\frac{\phi(r)^2}{r}\right)\,dr,\qquad
 d_j=2j\int_0^1r^{2j-1}\phi(r)^2\,dr .
 \end{gathered}
\]
Then
\begin{equation}\label{conformalray}
 \Phi(P)\le\frac{2\pi E}{D_*},\qquad
 D_*=\beta\left(\sum_{j\ge0}\frac{d_{j+1}|b_j|^2}{j+1}+|K|\right),
 \qquad K=\sum_{j\ge0}\frac{d_{j+2}}{j+2}b_j\overline{b_{j+2}}.
\end{equation}
In particular, omitting \(|K|\) gives
\begin{equation}\label{nr-ray-dipole}
 \Phi(P)\ge T>0\quad\Longrightarrow\quad
 d_1+\sum_{j\ge1}d_{j+1}\frac{|b_j|^2}{j+1}
 \le\frac{2\pi E}{T}\left(1+\sum_{j\ge1}\frac{|b_j|^2}{j+1}\right).
\end{equation}
\end{lemma}

\begin{proof}
Given a conformal map \(f_0:\mathbb D\to P\), put
\(T_a(z)=(z-a)/(1-\bar a z)\) and
\(B(a)=\int_{\mathbb D}V(T_a(z))|f_0'(z)|^2\,dA(z)\).
Dominated convergence gives \(B(\xi)=-|P|\phi(1)\xi\) for
\(|\xi|=1\), so \(\deg(B,\mathbb D,0)=1\). At a zero \(a\) take
\(f=f_0\circ T_a^{-1}\). A change of variables yields \eqref{balance}.

Each \(\Re(e^{it}V\circ f^{-1})\) has integral zero and energy \(\pi E\).
Maximizing its squared \(L^2(P)\) norm over \(t\) gives
\[
 \frac{1}{2}\left(\int_{\mathbb D}|V|^2|f'|^2\,dA+
              \left|\int_{\mathbb D}V^2|f'|^2\,dA\right|\right).
\]
Angular orthogonality gives the two normalized integrals
\[
 \frac1{|P|}\int_{\mathbb D}|V|^2|f'|^2\,dA
 =\beta\sum_{j\ge0}\frac{d_{j+1}|b_j|^2}{j+1},\qquad
 \frac1{|P|}\int_{\mathbb D}V^2|f'|^2\,dA=\beta K .
\]
The angular products converge in \(L^1\) since \(f'\in L^2\).
The quotient is bounded by \eqref{conformalray}.
\end{proof}

For the first two applications take
\begin{equation}\label{dipole-polynomial-profile}
 \phi(r)=r-\frac{39}{100}r^3+\frac{17}{500}r^5,\qquad
 \phi'(r)=(1-r^2)(1-17r^2/100)\ge0 .
\end{equation}

\begin{theorem}\label{pentagon}
For a convex polygon \(P\) with at most five sides, \(\Phi(P)<10.5\).
For a bounded simply connected Lipschitz domain \(P\) and a map
\(f\) satisfying \eqref{balance} for \eqref{dipole-polynomial-profile},
\[
 \frac{\pi|f'(0)|^2}{|P|}\le0.9585
 \quad\Longrightarrow\quad\Phi(P)<10.5013 .
\]
\end{theorem}

\begin{proof}
\emph{(1) Estimates for the mass.}
Take \(h_0=\int_0^1r^2\phi(r)\,dr\).
By \eqref{balance}, the formula for \(K\), and Cauchy--Schwarz,
\begin{equation}\label{harmonicbounds}
\begin{aligned}
 h_0|b_1|&\le \phi(1)\sum_{j\ge1}\frac{|b_j b_{j+1}|}{2j+3}
 \le\tfrac12\phi(1)S,\\
 \big|K-\tfrac12 d_2\overline b_2\big|
 &\le\phi(1)^2\sum_{j\ge1}\frac{|b_j b_{j+2}|}{j+2}\le\phi(1)^2 S.
\end{aligned}
\end{equation}

By \(d_n=\int_0^1\phi(s^{1/(2n)})^2\,ds\uparrow\),
\begin{equation}\label{conformaltrace}
             D_*\geq d_1+(d_2-d_1)\frac{S}{1+S}.
\end{equation}
For \(S\leq7/100\), use \(d_{j+1}\geq d_4\) for \(j\geq3\)
and retain \(|K|\). For \(0\le y=|b_2|\le\sqrt{3S}\),
\[
 \frac{d_4-d_3}{3}y^2-\left(\frac{d_2}2y-\phi(1)^2S\right)_+
 \le\frac{d_4-d_3}{3}\left(\frac{2\phi(1)^2S}{d_2}\right)^2.
\]
The maximum lies at the breakpoint \(2\phi(1)^2S/d_2\), or is smaller
if the breakpoint lies beyond the interval, since
\(2(d_4-d_3)\sqrt{0.21}/3-d_2/2<0\).
With \eqref{harmonicbounds} this gives
\begin{equation}\label{conformalgain}
  D_*\geq d_1+\frac{(d_4-d_1)S-CS^2}{1+S},\quad
 C=\frac{(d_4-d_2)\phi(1)^2}{8h_0^2}
             +\frac{d_4-d_3}{3}\left(\frac{2\phi(1)^2}{d_2}\right)^2 .
\end{equation}
The lower bound is increasing on \([0,0.07]\) since
\(d_4-d_1-C(2S+S^2)>0\).

\emph{(2) The polygonal bound.}
Solynin's theorem \cite{Solynin} and the regular derivative
proportional to \((1-z^n)^{-2/n}\) give
\[
 \beta\le\beta_n
 =\left[\sum_{j\ge0}\frac{[(2/n)_j/j!]^2}{1+nj}\right]^{-1}
 \le\beta_5<0.95795\qquad(n\le5).
\]
Here \((x)_j=x(x+1)\cdots(x+j-1)\), the summands decrease with
\(n\), and the final bound is Numerical Lemma~\ref{rationalconformal}.
That lemma applied to \eqref{conformalray} proves both assertions.
\end{proof}

\begin{proposition}\label{improved-conformal}
Let \(P\) be a convex polygon with \(\Phi(P)\geq10.504435\).
Every map satisfying \eqref{balance} for \eqref{dipole-polynomial-profile} has \(\beta>0.964\).
\end{proposition}

\begin{proof}
Set \(T_0=10.504435\). By Numerical Lemma~\ref{rationalconformal},
it suffices to exclude \(\beta\in[0.9585,0.964]\).
On \([\beta_a,\beta_b]\) set \(S_i=\beta_i^{-1}-1\) for \(i=a,b\).

\emph{(1) Auxiliary modes.} Put \(Q_m(re^{i\theta})=\psi_m(r)e^{im\theta}\),
\[
 \psi_m(r)=r^m\bigl[1+(m/2)(1-r^2)\bigr],\qquad
 e_m=m+\frac{m^2}{2(m+2)},\qquad 2\leq m\leq20 .
\]
Here \(\psi_m'(r)=m(1+m/2)r^{m-1}(1-r^2)\ge0\) and
\(0\le\psi_m\le1\). The energy blocks of \(V,Q_2,\ldots,Q_{20}\)
are \(\pi E I,\pi e_2 I,\ldots,\pi e_{20}I\).
Schwarz--Christoffel gives
\[
 \frac{|f'(re^{i\theta})|}{|f'(0)|}
   =\prod_s|1-re^{i\theta}\overline{\xi_s}|^{-\alpha_s}
   \geq(1+r)^{-2},
 \qquad \sum_s\alpha_s=2,\quad\alpha_s\geq0 ,
\]
where \(\xi_s\) are the prevertices and \(\pi\alpha_s\) the exterior
angles. Take
\[
 L(r)=\sum_{k=0}^3\binom{k+3}{3}\frac{(1-r)^k}{2^{k+4}}
 \le(1+r)^{-4},\qquad
 \tau_m=\int_0^1 rL(r)\psi_m(r)^2\,dr .
\]
For \(W=\sum_{m=2}^{20}\Re(\bar\eta_m Q_m)\), \(\eta_m\in\mathbb C\),
the angular integral of \(WL\) is zero. Hence
\[
 \frac1{|P|}\min_{c\in\mathbb R}\int_{\mathbb D}(W-c)^2 |f'|^2\,dA
 \ge\frac{\beta_a}{\pi}\int_{\mathbb D}W^2 L\,dA
 =\beta_a\sum_{m=2}^{20}\tau_m|\eta_m|^2 .
\]
For \(U=\Re V,\Im V\), apply Rayleigh to the centered trial
\(U+W-|P|^{-1}\int_{\mathbb D}W|f'|^2\,dA\). By \eqref{balance},
\[
 \begin{aligned}
 &\frac{T_0}{|P|}\int_{\mathbb D}U^2|f'|^2\,dA+
 \sum_{m=2}^{20}\left[
 \frac{2T_0}{|P|}\Re\!\left(\bar\eta_m\int_{\mathbb D} UQ_m |f'|^2\,dA\right)
 -(\pi e_m-T_0\beta_a\tau_m)|\eta_m|^2\right]\le\pi E .
 \end{aligned}
\]
Use \(e_m\ge2\) and \(0\le L,\psi_m^2\le1\) to obtain
\(\pi e_m-T_0\beta_a\tau_m\ge2\pi-T_0/2>0\).
Maximize each quadratic in \(\eta_m\) and add the two bounds:
\begin{equation}\label{improved-schur}
 \begin{aligned}
 &\frac1{|P|}\int_{\mathbb D}|V|^2|f'|^2\,dA\\
 &\quad+\sum_{m=2}^{20}
 \frac{T_0}{2|P|^2(\pi e_m-T_0\beta_a\tau_m)}
 \sum_{Q\in\{Q_m,\bar Q_m\}}\left|\int_{\mathbb D}VQ |f'|^2\,dA\right|^2
 \le\frac{2\pi E}{T_0}.
 \end{aligned}
\end{equation}

\emph{(2) Fourier coefficients.} For \(j\ge1\) put
\[
 x_j=\frac{|b_j|^2}{(j+1)S},\qquad s_j=1-\sum_{k=1}^j x_k .
\]
For \((j,Q)=(m+1,Q_m),(m-1,\bar Q_m)\) and \(j\ge2\), the leading
Fourier term has modulus \(\beta c_{jm}|b_j|\), where
\[
 c_{jm}=2\int_0^1r^{j+1}\phi(r)\psi_m(r)\,dr .
\]
By Cauchy--Schwarz and \(S\le S_a\),
\[
 \begin{aligned}
 \frac1{\beta|P|}\left|\int_{\mathbb D}VQ|f'|^2\,dA\right|
 &\ge\left(c_{jm}|b_j|-\phi(1)\sum_{k\ge1}
                   \frac{|b_kb_{k+j}|}{k+1+j/2}\right)_+\\
 &\ge\sqrt S\left(c_{jm}\sqrt{(j+1)x_j}
                         -\phi(1)\sqrt{S_as_j}\right)_+ .
 \end{aligned}
\]
For
\(K_0\ge\phi(1)\sqrt{S_a}\), the left side of
\eqref{improved-schur} is at least
\[
 d_1+(1-\beta)\left[d_{21}-d_1-\sum_{j=1}^{19}\ell_j\right],
\]
where
\[
 \begin{aligned}
 \ell_j&=(d_{21}-d_{j+1})x_j-
 \sum_{\substack{m=j\pm1\\2\leq m\leq20,\ j\geq2}}
 A_m\bigl(c_{jm}\sqrt{(j+1)x_j}-K_0\sqrt{s_j}\bigr)_+^2,\\
 A_m&=\frac{\beta_a T_0}{2(\pi e_m-T_0\beta_a\tau_m)} .
 \end{aligned}
\]
We have \(x_1\leq u_1=\phi(1)^2 S_a/(8h_0^2)\).
For the bound on \(x_2\), write \(y=|b_2|\) and put
\[
 R(y)=(d_4-d_1)S_b-(d_4-d_2)u_1S_a
        -(d_4-d_3)y^2/3+(d_2/2)y-\phi(1)^2 S_a .
\]
The mass inequality \eqref{conformalray} gives
\(d_1+\beta_a R(y)\leq 2\pi E/T_0\) whenever \(R(y)>0\).
On \(0\le y\le\sqrt{3S_a}\) its derivative is
\(R'(y)=d_2/2-2(d_4-d_3)y/3>0\). Thus any rational \(Y>0\) with
\(d_1+\beta_a R(Y)>2\pi E/T_0\) gives \(x_2\leq u_2=Y^2/(3S_b)\).

\emph{(3) Recursion.}
With \(L_{20}=0\) and \(y^2=(j+1)x_j/s_j\), use
\[
 \frac{L_{j+1}s_j+\ell_j}{s_{j-1}}\le L_j
 \quad\Longrightarrow\quad\sum_{i=j}^{19}\ell_i\le L_js_{j-1}.
\]
The quotient is
\begin{equation}\label{schur-recursion}
 \frac{L_{j+1}+(d_{21}-d_{j+1})y^2/(j+1)
       -\sum_{m=j\pm1,\ m\geq2} A_m(c_{jm}y-K_0)_+^2}
      {1+y^2/(j+1)} .
\end{equation}
Use the limit if \(s_j=0\). For \(j\ge3\) maximize over \(y\ge0\),
and for \(j=2\) use \(y^2\le3u_2/(1-u_1-u_2)\).
After clearing the denominator the tests are quadratic between
the breakpoints \(K_0/c_{jm}\). In particular,
\[
 \sum_{j=1}^{19}\ell_j
 \le (d_{21}-d_2)x_1+L_2(1-x_1)
 \le L_2+(d_{21}-d_2-L_2)_+u_1,
\]
which contradicts \eqref{improved-schur} by Numerical
Lemma~\ref{schur-numerical}.
\end{proof}

\begin{numericallemma}[Basic conformal bounds]\label{rationalconformal}
With \eqref{conformaltrace} for \(S\ge0.07\) and \eqref{conformalgain}
for \(S\le0.07\), the bound in \eqref{conformalray} satisfies
\[
 \beta\le0.95795\ \Longrightarrow\ 2\pi E/D_*<10.5,\qquad
 \beta\le0.9585\ \Longrightarrow\ 2\pi E/D_*<10.5013,
 \qquad \beta_5<0.95795 .
\]
\end{numericallemma}

\begin{proof}
Integrating the polynomial \(\phi\) gives
\begin{gather*}
 E=\frac{153583}{312500},\qquad \phi(1)^2=\frac{25921}{62500},
 \qquad h_0=\frac{757}{4000},\\
 (d_1,d_2,d_3,d_4)
 =\left(\frac{869741}{3000000},\frac{935017}{2625000},
         \frac{5335649}{14000000},\frac{6190147}{15750000}\right).
\end{gather*}
In \eqref{conformaltrace} and \eqref{conformalgain} the denominators
increase on their respective ranges. Using \(\pi<3.141592654\) gives
\begin{align*}
 \left.\frac{2\pi E(1+S)}{d_1+d_2S}\right|_{S=0.07}&<10.4944,\\
 \left.\frac{2\pi E(1+S)}{d_1+d_4S-CS^2}\right|_{S=1/\beta-1}
 &<\begin{cases}10.4994,&\beta=0.95795,\\10.5013,&\beta=0.9585.\end{cases}
\end{align*}

For the conformal radius of \(R_5\), take \(p_0=10^{12}\) and
\[
 p_j=\left\lfloor\frac{(5j-3)p_{j-1}}{5j}\right\rfloor
 \le10^{12}\frac{(2/5)_j}{j!}.
\]
The integer sum yields
\[
 \beta_5^{-1}\ge
 10^{-12}\sum_{j=0}^{512}\left\lfloor\frac{p_j^2}{10^{12}(1+5j)}\right\rfloor
 =\frac{1043897249658}{10^{12}}>\frac1{0.95795}.
\]
\end{proof}

\begin{numericallemma}[Conformal Rayleigh bounds]\label{schur-numerical}
On each of the intervals with consecutive endpoints
\[
 0.9585,\quad 0.962,\quad 0.963,\quad 0.9635,\quad 0.96375,\quad 0.9639,\quad 0.964,
\]
the recursion \eqref{schur-recursion} gives choices of \(K_0,Y,L_2\)
for which
\[
 d_1+(1-\beta_b)
 \bigl[d_{21}-d_1-L_2-(d_{21}-d_2-L_2)_+u_1\bigr]
       > 2\pi E/T_0+1.1\cdot10^{-7}.
\]
The bracket on the left is positive.
\end{numericallemma}

\begin{proof}
Write \(p=3.141592654\),
\(s(x)=10^{-10}(1+\lfloor\sqrt{10^{20}x}\rfloor)\). Set
\[
 K_0=\phi(1)s(S_a),\qquad
 Y=10^{-4}\min\{k\in\mathbb N:d_1+\beta_aR(k/10^4)>2pE/T_0\},
 \qquad u_2=Y^2/(3S_b).
\]
Use \(p\) in place of \(\pi\) in \(A_m\), set \(L_{20}=0\), and take
the least positive \(L_j\in10^{-10}\mathbb Z\) such that
\[
 \begin{aligned}
 (1+y^2/(j+1))L_j &\ge L_{j+1}+(d_{21}-d_{j+1})y^2/(j+1)
  -\!\!\sum_{\substack{m=j\pm1\\m\ge2}}\! A_m(c_{jm}y-K_0)_+^2,\\
 &\quad y\in
 \begin{cases}
 [0,\infty),&j=19,\ldots,3,\\
 [0,s(3u_2/(1-u_1-u_2))],&j=2.
 \end{cases}
 \end{aligned}
\]
For each quadratic between consecutive \(K_0/c_{jm}\), test the
endpoints, any interior minimum, and the sign at infinity on an
unbounded interval. The lower bounds for the excess over \(2pE/T_0\) are
\[
 10^{-8}(17400,9600,4800,2500,1000,11).
\]
\end{proof}

\subsection{The regular-polygon expansion and conformal residuals}\label{sec:regular-sc-certificate}

Cureton--Kuttler \cite{CK} used conformal transfer and symmetry for
spectral bounds on regular polygons. Berghaus, Jones, Monien and
Radchenko \cite[Section~8.1]{BJMR} computed asymptotic coefficients
numerically for the Neumann branch tending to \(j_{1,1}^2\).
Here we treat the first positive branch, which tends to \((j'_{1,1})^2\).
Dong and Wang list a related note on the large-side limit in
\cite[reference~6]{DW}.

On the unit disk \(\mathbb D\), put
\[
 x=(j'_{1,1})^2=\mu_1(\mathbb D),\quad
 \varphi(r)=\frac{J_1(\sqrt{x}\,r)}{J_1(\sqrt{x})},\quad
 d_k=2k\int_0^1r^{2k-1}\varphi^2\,dr .
\]
The Bessel mass identity and \(J'_1(\sqrt x)=0\) give
\[
 d_1=1-\frac{J_0(\sqrt x)J_2(\sqrt x)}{J_1(\sqrt x)^2}=1-x^{-1},
 \qquad
 \varphi(0)=0,\quad \varphi(1)=1,\quad
 \varphi'>0\ \text{on }(0,1),\quad\varphi'(1)=0 .
\]
\begin{theorem}\label{neumann-fifth-order}
As \(N\to\infty\), the regular-polygon value satisfies
\begin{equation}\label{regular-neumann-asymptotic-formula}
 \Phi(R_N)=\pi x-\frac{4\pi x\zeta(3)}{(x-1)N^3}
 -\frac{\pi x[12+4x(x+1/2)]\zeta(5)}{(x-1)N^5}
 +O(N^{-6}).
\end{equation}
\end{theorem}

\begin{proof}
\emph{(1) Conformal weight.}
Take \(|R_N|=\pi\) and \(c_k=(2/N)_k/k!\). By
Schwarz--Christoffel \cite{DT} and angular orthogonality,
\[
 \begin{gathered}
 w_N=\beta_N|1-z^N|^{-4/N},\qquad
 \beta_N^{-1}=1+\sum_{k\ge1}\frac{c_k^2}{Nk+1},\\
 \frac1{2\pi}\int_0^{2\pi}|1-r^Ne^{iN\theta}|^{-4/N}\,d\theta
 =\sum_{k\ge0}c_k^2r^{2Nk}.
 \end{gathered}
\]
The area formula with circumradius
\(\sqrt{\beta_N}B(1/N,1-2/N)/N\) and Euler's reflection identity give
\begin{align*}
 \beta_N^{-1}
 &=\frac{\sin(2\pi/N)}{2\pi N}B(1/N,1-2/N)^2
 =\frac{\Gamma(1+1/N)^2\Gamma(1-2/N)}
 {\Gamma(1-1/N)^2\Gamma(1+2/N)},\\
 -\log\beta_N
 &=\sum_{\substack{k\ge3\\k\ {\rm odd}}}\frac{(2^{k+1}-4)\zeta(k)}{kN^k},\qquad
 1-\beta_N=4\zeta(3)N^{-3}+12\zeta(5)N^{-5}+O(N^{-6}).
\end{align*}
On \(\{u\in H^1(\mathbb D):\int u=0\}\), conformal invariance gives
\[
 a(u,v)=\int_{\mathbb D}\nabla u\cdot\overline{\nabla v},\quad
 a(K_Nu,v)
 =\int_{\mathbb D}w_Nu\bar v-\frac1\pi\int_{\mathbb D}w_Nu
                                  \int_{\mathbb D}w_N\bar v ,
\]
with \(\lambda_{\max}(K_N)=\mu_1(R_N)^{-1}\). Write \(K_0\) for the
compact operator for unit weight and \(\|v\|_a^2=a(v,v)\).
For \(0<t<1\), \(k\ge1\), use
\begin{equation}\label{pochhammer-power-bound}
 \frac{(t)_k}{k!}\le\frac{t}{k}\exp\!\left(t\sum_{i=1}^{k-1}i^{-1}\right)
 \le t e^t k^{t-1}\le\frac{t}{1-t} k^{t-1} .
\end{equation}

For \(N\ge16\) angular integration gives
\begin{align*}
 \|w_N-1\|_2^2
 &=\pi\left(\beta_N^2\sum_{k\ge0}
       \frac{[(4/N)_k/k!]^2}{Nk+1}-1\right)
 \le\pi\sum_{k\ge1}\frac{[(4/N)_k/k!]^2}{Nk+1}=O(N^{-3}),\\
 \left\|\frac1{2\pi}\int_0^{2\pi}w_N(re^{i\theta})\,d\theta-1\right\|_2
 &\le\sqrt\pi(1-\beta_N)+\beta_N\sum_{k\ge1}c_k^2
                          \sqrt{\frac\pi{2Nk+1}}=O(N^{-5/2}).
\end{align*}
For \(\int u=0\), Sobolev embedding and
\(\int w_Nu=\int(w_N-1)u\) give
\[
 \|K_N-K_0\|\le C\bigl(\|w_N-1\|_2+\|w_N-1\|_2^2\bigr)=O(N^{-3/2}).
\]

\emph{(2) Compression to frequency one.}
By Lemma~\ref{affine-all-N} and reflection, it suffices to use
\(m\equiv1\pmod N\). Let \(Q\) be the energy projection onto the
orthogonal complement of \(U=\varphi(r)e^{i\theta}\) in that class.
Pair the \(m=1\) component with the angular mean of \(w_N-1\).
For \(|m|\ge N-1\), \(\|v_m\|_2\le |m|^{-1}\|v_m\|_a\), so
\[
 \|QK_NU\|_a\le C\left(
 \left\|\frac1{2\pi}\int_0^{2\pi}(w_N(re^{i\theta})-1)\,d\theta\right\|_2
   +\frac{\|w_N-1\|_2}{N-1}\right)=O(N^{-5/2}).
\]
On \(\operatorname{ran}Q\), the disk gap and
\(|\mu_1(R_N)^{-1}-x^{-1}|\le\|K_N-K_0\|\) give
\[
 (\mu_1(R_N)^{-1}I-QK_NQ)^{-1}
 =(x^{-1}I-QK_0Q)^{-1}+O(N^{-3/2}).
\]
Since \(a(U,U)=\pi x d_1\), the Schur complement on \(\mathbb C U\) yields
\[
 \frac{\pi x d_1}{\mu_1(R_N)}
 =a(K_NU,U)+a\!\left(QK_NU,(x^{-1}I-QK_0Q)^{-1}QK_NU\right)
       +O(N^{-13/2}).
\]
Write \(H(z)=\sum_{k\ge1}c_kz^{Nk}\), so
\(w_N=\beta_N(1+H+\overline H+|H|^2)\). The constant-weight
part gives \(Q(\beta_N K_0U)=0\). For the quadratic source
\(\beta_NU|H|^2\), a test of frequency one satisfies
\[
 |v_1(r)|^2\le\int_0^1(s|v'_1|^2+s^{-1}|v_1|^2)\,ds,\qquad
 \frac1\pi\int_{\mathbb D}|H|^2\,dA=\beta_N^{-1}-1=O(N^{-3}).
\]
Thus the frequency-one part has dual energy norm \(O(N^{-3})\).
For higher frequencies, the convolution bound
\[
 |[z^{Nk}]H^2|\le C N^{-2}\sum_{l=1}^{k-1}
                l^{-3/4}(k-l)^{-3/4}\le C N^{-2}k^{-1/2}
\]
gives \(\||H|^2\|_2=\|H^2\|_2=O(N^{-5/2})\) and dual norm
\(O(N^{-7/2})\) after division by \(N-1\). The disk resolvent preserves frequencies. The quadratic source costs
\(O(N^{-3})^2+O(N^{-7/2})O(N^{-5/2})=O(N^{-6})\).
For the linear source \(\beta_NU(H+\bar H)\), every \(|m|\ge N-1\), so
\[
 0\le K_0\le(N-1)^{-2}I,\qquad
 (x^{-1}I-QK_0Q)^{-1}=xI+O(N^{-2}).
\]

\emph{(3) The coefficient of \(N^{-5}\).}
For \(p=Nk\), \(m=p\pm1\), use the Neumann inverse of
\(L_m=-\partial_r^2-r^{-1}\partial_r+m^2r^{-2}\) on
\(L^2((0,1),r\,dr)\). The squared dual energy norm of the linear
source and the diagonal mass are
\[
 2\pi\beta_N^2\sum_{k\ge1}c_k^2\sum_{m=p\pm1}
 \int_0^1\varphi(r)r^{p+1}[L_m^{-1}(\varphi r^p)](r)\,dr,\qquad
 \frac{a(K_NU,U)}{\pi}=\beta_N\left(d_1+\sum_{k\ge1}
                         \frac{c_k^2d_{Nk+1}}{Nk+1}\right).
\]
From \(\varphi''(1)=1-x\) and the beta integral,
\[
 1-\varphi(r)^2=(x-1)(1-r)^2+O((1-r)^3),\quad
 1-d_{p+1}=\frac{2(x-1)}{(2p+3)(2p+4)}+O(p^{-3}).
\]
In \(L^2((0,1),r\,dr)\),
\[
 \|L_m^{-1}\|\le m^{-2},\qquad \|r^p\|=O(p^{-1/2}),\qquad
 \|(1-\varphi)r^p\|=O(p^{-5/2}).
\]
Thus the two replacements \(\varphi\mapsto1\) cost \(O(p^{-5})\), and
\begin{align*}
 L_m^{-1}(r^p)&=\frac{(p+2)r^m/m-r^{p+2}}{(p+2)^2-m^2},\\
 \int_0^1\varphi(r)r^{p+1}[L_m^{-1}(\varphi r^p)](r)\,dr
 &=\frac{2p+4+m}{m(p+m+2)^2(2p+4)}+O(p^{-5})
   =\frac{3}{8p^3}+O(p^{-4}).
\end{align*}
Insert these expansions and
\(\beta_N\sum_{k\ge1}c_k^2/(Nk+1)=1-\beta_N\) into the Schur complement.
Using \(1-d_1=1/x\) gives
\[
 \frac{x d_1}{\mu_1(R_N)}-d_1
 =\frac{1-\beta_N}{x}
 +(x+\tfrac12)\sum_{k\ge1}\frac{c_k^2}{(Nk)^3}+O(N^{-6}).
\]
Here \(\sum_{k\ge1}c_k^2/(Nk)^4=O(N^{-6})\), and
\(\beta_N=1+O(N^{-3})\) in the terms of order \(N^{-5}\). Since
\[
 c_k=\frac{2}{Nk}\prod_{j=1}^{k-1}\left(1+\frac2{Nj}\right)
 =\frac2{Nk}+O\!\left(\frac{1+\log k}{N^2k^{3/4}}\right),\qquad
 \sum_{k\ge1}\frac{c_k^2}{(Nk)^3}=\frac{4\zeta(5)}{N^5}+O(N^{-6}),
\]
inversion gives \eqref{regular-neumann-asymptotic-formula}.
\end{proof}

\begin{lemma}\label{large-N-conformal-trials}
For a bounded simply connected Lipschitz domain \(P\), take the map
\(f\) and \(b_j,\beta,S,d_j\) of Lemma~\ref{conformal-dipole} with
\(V(re^{i\theta})=\varphi(r)e^{i\theta}\). Then
\begin{equation}\label{large-N-trace}
 \begin{aligned}
 \frac1{|P|}\left(\int_{\mathbb D}|V|^2|f'|^2\,dA+
             \left|\int_{\mathbb D}V^2|f'|^2\,dA\right|\right)
 &\le\frac{\pi x d_1}{\Phi(P)},\\
 \frac1{|P|}\int_{\mathbb D}|V|^2|f'|^2\,dA
 &=d_1+\beta\sum_{j\ge1}(d_{j+1}-d_1)\frac{|b_j|^2}{j+1}.
 \end{aligned}
\end{equation}
For every integer \(J\ge2\),
\begin{equation}\label{large-N-schur}
 \begin{aligned}
 &\frac1{|P|}\int_{\mathbb D}|V|^2|f'|^2\,dA\\
 &\qquad+\sum_{m=2}^J\frac{\Phi(P)}{2\pi m|P|^2}
 \sum_{Q\in\{z^m,\bar z^m\}}\left|\int_{\mathbb D}VQ |f'|^2\,dA\right|^2
 \le\frac{\pi x d_1}{\Phi(P)} .
 \end{aligned}
\end{equation}
For \(j\ge3\),
\begin{align*}
 \left|\frac1{|P|}\int_{\mathbb D}Vz^{j-1}|f'|^2\,dA
              -2\beta\bar b_j\int_0^1r^{2j}\varphi(r)\,dr\right|&\le\beta S,\\
 \left|\frac1{|P|}\int_{\mathbb D}V\bar z^{j+1}|f'|^2\,dA
              -2\beta b_j\int_0^1r^{2j+2}\varphi(r)\,dr\right|&\le\beta S ,
\end{align*}
and
\begin{equation}\label{large-N-radial-moment}
 1-d_{j+1}\le \frac{8}{(2j+3)(2j+4)},\qquad
 2\int_0^1r^{2j}\varphi\,dr\ge\frac{1}{j+1},\qquad
 2\int_0^1r^{2j+2}\varphi\,dr\ge\frac{1}{j+2}.
\end{equation}
\end{lemma}

\begin{proof}
Apply \eqref{conformalray} with \(2E=xd_1\). Since
\(\beta\sum_{j\ge0}|b_j|^2/(j+1)=1\), the radial part of its
denominator is
\[
 \beta\sum_{j\ge0}\frac{d_{j+1}|b_j|^2}{j+1}
 =d_1+\beta\sum_{j\ge1}(d_{j+1}-d_1)\frac{|b_j|^2}{j+1}.
\]
For \eqref{large-N-schur}, complete the squares as in
\eqref{improved-schur} with \(T_0=\Phi(P)\), \(E=xd_1/2\),
\(Q_m=z^m\), \(e_m=m\), \(\tau_m=0\), and \(2\le m\le J\).
For \(j=m+1,m-1\), the terms with \(k\ge1\) contribute at most
\[
 \beta\sum_{k\ge1}\frac{|b_kb_{k+j}|}{k+1+(j+m)/2}
 \le \beta\sum_{k\ge1}\frac{|b_kb_{k+j}|}{\sqrt{(k+1)(k+j+1)}}\le\beta S.
\]
The first two modes give
\[
 \begin{aligned}
 \left|\frac1{|P|}\int_{\mathbb D}V |f'|^2\,dA
           -2\beta\overline{b_1}\int_0^1r^2\varphi\,dr\right|&\le\beta S,\\
 \left|\frac1{|P|}\int_{\mathbb D}V^2 |f'|^2\,dA
           -\tfrac12\beta d_2\overline{b_2}\right|&\le\beta S .
 \end{aligned}
\]
The radial equation gives
\((r^3(\varphi/r)')'=-xr^2\varphi\), so
\((\varphi/r)'\le0\), \(\varphi(r)\ge r\), and
\(\varphi''=-x\varphi-(\varphi/r)'\ge-x>-4\).
Integrating at \(r=1\) gives \(1-\varphi^2\le4(1-r)^2\), whence
\begin{align*}
 1-d_{j+1}&\le8(j+1)\int_0^1r^{2j+1}(1-r)^2\,dr
      =\frac8{(2j+3)(2j+4)},\\
 2\int_0^1r^{j+m+1}\varphi\,dr&\ge2\int_0^1r^{j+m+2}\,dr
      =\frac2{j+m+3}.
\end{align*}
\end{proof}

For a convex \(m\)-gon, \(m\le N\), Solynin's inequality
\cite{Solynin,BHS} and the gamma series give
\(\beta(P)\le\beta_m\le\beta_N\).

\smallskip\noindent\emph{Residual of a conformal trial.}
Take \(n\ge17\) and \(R_n\) of circumradius one. The map
with derivative \((1-z^n)^{-2/n}\) has an image similar to \(R_n\).
Its mass weight and Fourier coefficients at frequencies \(mn\) are
\[
 w(r,\theta)=|1-r^n e^{in\theta}|^{-4/n},\qquad
 a_l=\frac{(2/n)_l}{l!},\qquad
 \widehat w_m(r)=\sum_{l\ge0}a_l a_{l+|m|}
                      r^{n(2l+|m|)} .
\]

In the class \(k\equiv1\pmod n\), set
\[
 \langle Ku,v\rangle_E=\frac1{2\pi}\int_{\mathbb D} w\,u\bar v\,dA ,
 \qquad
 \langle u,v\rangle_E=\frac1{2\pi}\int_{\mathbb D}
                      \nabla u\cdot\overline{\nabla v}\,dA .
\]

By Lemma~\ref{affine-all-N}, \(\lambda_{\max}(K)\) is the reciprocal
first eigenvalue of the unnormalized image.

\emph{(1) Separation of the dipole.}
The second largest eigenvalue \(\gamma_2\) of \(K\) satisfies
\(\gamma_2\le\overline\gamma_2\), where
\[
 \overline\gamma_2=\frac{100}{2809}+\kappa,\qquad
 \kappa=\frac{8(n-2)}{n(n-4)^2}.
\]
For an angular mode \(k\), set
\(A_k(v)=\int_0^1(r|v'|^2+k^2|v|^2/r)\,dr\). For \(k\ne0\),
\[
 |k|\,|v(r)|^2\le2|k|\int_0^r|v'||v|\,ds\le A_k(v).
\]
By angle integration and \(w-1\le(1-r^n)^{-4/n}-1\), the excess
mass in the dipole class is at most \(\kappa\) times the energy:
\[
 \int_0^1[(1-r^n)^{-4/n}-1]r\,dr
 =\sum_{l\ge1}\frac{(4/n)_l}{l!(nl+2)}
 \le\frac4{n(n-4)}\left(1+\int_1^\infty t^{4/n-2}\,dt\right)=\kappa,
\]
by \eqref{pochhammer-power-bound}.
For the disk modes \(|k|\ge n-1\),
\(A_k(v)\ge k^2\int_0^1r|v|^2\,dr\), while
\(j'_{1,2}\ge j_{1,1}>\sqrt{9.3}>2\).
On \(1\le s\le 2809/400\) the alternating Bessel series gives
\[
 2J'_1(2\sqrt s)\le P(s)=\sum_{j=0}^{10}
 \frac{(-1)^j(2j+1)s^j}{j!(j+1)!}.
\]
Its Bernstein coefficients on this interval satisfy
\[
 \max_{0\le i\le10}\sum_{j=0}^{i}
 \frac{\binom ij}{\binom{10}j}\frac{(2409/400)^jP^{(j)}(1)}{j!}
 <-\frac1{50}.
\]

Thus \(j'_{1,2}>5.3\), and min--max gives
\(\gamma_2\le(5.3)^{-2}+\kappa=\overline\gamma_2\).
In residue zero take \(\int_0^1 v(s)s\,ds=0\).
Cauchy--Schwarz applied to \(v(r)=2\int_0^1[v(r)-v(s)]s\,ds\) gives
\[
 |v(r)|^2\le 2 A_0(v)\int_0^1 |\log(r/s)| s\,ds
            =A_0(v)\bigl(-\log r+r^2-\tfrac12\bigr).
\]
Radial integration gives
\[
 \int_0^1[(1-r^n)^{-4/n}-1](-\log r)r\,dr
 =\sum_{l\ge1}\frac{(4/n)_l}{l!(nl+2)^2}\le\frac{\kappa}{n+2}.
\]

Use \(r^2-1/2\le1\) and \(1/\mu_3(\mathbb D)<1/6\)
(Lemma~\ref{affine-all-N}). For
\(u(e^{2\pi i/n}z)=e^{2\pi ij/n}u(z)\), \(j\not\equiv\pm1\pmod n\),
\[
 \inf_{b\in\mathbb C}\int_{\mathbb D}w|u-b|^2
 \le \left[\frac16+\kappa\left(1+\frac1{n+2}\right)\right]
                       \int_{\mathbb D}|\nabla u|^2.
\]
\emph{(2) Mass and residual.}
For a finite rational trial \(u\) in residue one use
\[
 r^{|k|}\bigl(1+\tfrac{|k|}{2}(1-r^2)\bigr),\qquad
 r^{|k|}|k|^q(1-r^2)^q\quad(q\ge2).
\]
Set \(\rho=\langle Ku,u\rangle_E/\|u\|_E^2\).
For \(p=|k|+|k'|\), \(m=(k-k')/n\), each mass entry is a
nonnegative combination of sums
\[
 \frac12\sum_{l\ge0} a_la_{l+|m|}
  \frac{q!}{\prod_{j=1}^{q+1}(j+(p+n(2l+|m|))/2)},
 \qquad
 \int_0^1t^b(1-t)^q\,dt=\frac{q!}{\prod_{j=1}^{q+1}(b+j)} .
\]
For \(l>L\ge1\), \((p+n(2l+|m|))/2\ge nl\) and
\(a_l a_{l+|m|}\le 4l^{4/n-2}/(n-2)^2\). Since \(4/n\le1/4\),
the tail lies between zero and
\[
 \frac{2q!}{(n-2)^2 n^{q+1}}\int_L^\infty t^{-q-11/4}\,dt
 =\frac{2q!}{(n-2)^2 n^{q+1}(q+7/4)L^{q+7/4}} .
\]
For the area \(A_w=2\pi\int_0^1\widehat w_0(r)r\,dr\), the case
\(p=q=m=0\) gives this same bound for the tail of \(A_w/(2\pi)\).
For an enclosure \(\rho\in[\rho_-,\rho_+]\), set
\(\rho_0=(\rho_-+\rho_+)/2\) and \(F=wu+\rho_0\Delta u\).
Green's identity gives
\[
 \langle(K-\rho_0)u,v\rangle_E=\frac1{2\pi}\int_{\mathbb D}F\bar v,
 \qquad \|(K-\rho_0)u\|_E^2
 =\sum_{k=1+jn}\int_0^1\overline{F_k}(L_{|k|}^{-1}F_k)\,r\,dr .
\]
The Neumann radial equation gives
\[
 \int_0^1 r^p L_k^{-1}(r^q)\,r\,dr
 =\frac{p+q+4+k}{k(p+k+2)(q+k+2)(p+q+4)}\qquad(k\ge1),
\]
including \(q+2=k\) by continuity.
Retain the output modes \(k=1+mn\) with \(|m|\le M\), using
\(0\le l\le L_r\) in each Fourier coefficient of \(w\). For
\(u=\sum_{|j|\le J}u_j(r)e^{i(1+jn)\theta}\), \(J<M\), the
respective radial and angular tail norms, divided by
\(\sum_j\|u_j\|_\infty\), are bounded by
\[
 \frac{4nL_r^{-3/4}}{(n-2)^2(n-4)}
 \left(\frac{1+4/(n-1)^2}{2(n(2L_r+2)+1)}\right)^{1/2},
 \qquad
 \frac{2n}{(n-2)(n-4)\sqrt{(7/2)n(n-1)^2(M-J)^{7/2}}}.
\]
For the first, use \(\|L_k^{-1/2}\|\le1/|k|\) and
\[
 \sum_{l>L_r}a_la_{l+|m|}r^{n(2l+|m|)}
 \le\frac{4nL_r^{-3/4}}{(n-2)^2(n-4)}r^{n(2L_r+2)}.
\]
For the second, split the coefficient sum at \(l=0,|m|\). As
\[
 \sum_{l=1}^{|m|}l^{2/n-1}\le1+\frac n2(|m|^{2/n}-1),\qquad
 \sum_{l>|m|}l^{4/n-2}\le\frac{n|m|^{4/n-1}}{n-4},
 \qquad 1+\frac2{n-2}-\frac{n}{n-2}=0 ,
\]
one obtains
\[
 |\widehat w_m(r)|\le\frac{2n}{(n-2)(n-4)}|m|^{4/n-1}r^{n|m|}
 \qquad(m\ne0).
\]
Using \(\|L_k^{-1/2}\|\le1/|k|\), the radial integral of
\(r^{2n(|m|-J)}\) and \(8/n\le1/2\), the squared tail after the same division is
at most
\[
 \frac{4n}{(n-2)^2(n-4)^2(n-1)^2}\int_{M-J}^\infty s^{-9/2}\,ds .
\]
For \(q\ge1\),
\[
 \sup_{0\le r\le1}r^k(1-r^2)^q
 =\left(\frac{k}{k+2q}\right)^{k/2}
  \left(\frac{2q}{k+2q}\right)^q
 \le\left(\frac{2q}{k+2q}\right)^q .
\]
\emph{(3) Eigenvalue enclosure.}
Use \(\mathcal E\ge\|(K-\rho_0)u\|_E/\|u\|_E+(\rho_+-\rho_-)/2\).
For \(\gamma=\lambda_{\max}(K)\), the spectral gap gives
\[
 \begin{aligned}
 0&\le\frac{\langle(K-\gamma I)(K-\overline\gamma_2 I)u,u\rangle_E}
              {\|u\|_E^2}\\
 &=\frac{\|(K-\rho I)u\|_E^2}{\|u\|_E^2}
                 -(\gamma-\rho)(\rho-\overline\gamma_2).
 \end{aligned}
\]
By Temple's bound \cite{Temple} and scaling,
\[
 \rho\le\gamma\le\rho+\frac{\mathcal E^2}{\rho-\overline\gamma_2}
 \quad(\rho>\overline\gamma_2),
 \qquad
 \mu_1(R_n)=\frac{2 A_w}{n\sin(2\pi/n)\,\gamma},
 \qquad \Phi(R_n)=\frac{A_w}{\gamma}.
\]

Use outward rational endpoints throughout and Machin's identity
\(\pi/4=4\arctan(1/5)-\arctan(1/239)\) with alternating remainders.

\section{Centroid sectors and the finite convex comparisons}\label{part:finite}

For convex \(N\)-gons with \(5\le N\le32\), we reduce the comparison
to triangles with a common vertex at the centroid. Summing the sector
inequalities gives \(\Phi(P)<\Phi(R_N)\) for every nonregular competitor.

\subsection{Localization and sector comparison}\label{sec:centroid-sector-comparison}

\begin{lemma}\label{diam}
For a planar convex body \(P\),
\[
 \mu_1(P)\le
 \frac{|P|}{\displaystyle\max_{|e|=1}
 \int_P\left[e\cdot\left(x-\frac1{|P|}\int_P y\,dy\right)\right]^2\,dx}
 \le\frac{24}{\operatorname{diam}(P)^2}.
\]
\end{lemma}

\begin{proof}
Translate to \(\int_P x\,dx=0\). For \(|e|=1\), let \([a,a+D]\) be the
projection of \(\overline P\) on \(e\). The Rayleigh principle and the
slice identity give
\[
 \mu_1(P)\int_P(e\cdot x)^2\,dx\le |P|,
 \qquad
 L(t)=\frac{D}{|P|}\mathcal H^1\{x\in P:e\cdot x=a+Dt\},
 \quad \int_0^1L(t)\,dt=1 .
\]
The slice length \(L\) is concave, with one-sided limits at \(0,1\).
For its piecewise linear interpolants with knots \(0=s_0<\cdots<s_m=1\),
\[
 \begin{gathered}
 L=\sum_{j=0}^m c_j H_{s_j},\quad
 H_s(t)=2\min\{t/s,(1-t)/(1-s)\},\quad
 (H_0(t),H_1(t))=2(1-t,t),\\
 c_0=L(0)/2,\quad c_m=L(1)/2,\quad
 c_j=\tfrac12 s_j(1-s_j)(L'_-(s_j)-L'_+(s_j))\quad(0<j<m),\\
 c_j\ge0,\quad \sum_jc_j=1.
 \end{gathered}
\]
Normalize the interpolants to have integral one and pass to the
\(L^1(0,1)\) limit after using
\[
 \int_0^1 t^k H_s(t)\,dt=\frac{2\sum_{i=0}^k s^i}{(k+1)(k+2)}
 \qquad (k\ge0).
\]
In the diameter direction \(D=\operatorname{diam}(P)\). The moments for
\(k=1,2\) give
\[
 \begin{aligned}
 \frac{1}{|P|D^2}\int_P(e\cdot x)^2\,dx
 &=\min_b\int_0^1(t-b)^2L(t)\,dt
 \ge\sum_j c_j\min_b\int_0^1(t-b)^2H_{s_j}(t)\,dt\\
 &=\frac1{24}+\frac1{18}\sum_j c_j(s_j-1/2)^2
 \ge\frac1{24}.
 \end{aligned}
\]
\end{proof}

\begin{lemma}\label{cubic-directional-integral}
For a planar convex body \(P\),
\[
 \begin{aligned}
 \frac{\mu_1(P)}{|P|}\max_{|e|=1}
 \int_P\left[e\cdot\left(x-\frac1{|P|}\int_P y\,dy\right)\right]^2\,dx
 &\le\frac{35}{36},\\
 \mu_1(P)\operatorname{diam}(P)^2&\le\frac{70}{3}.
 \end{aligned}
\]
\end{lemma}

\begin{proof}
Translate to \(\int_P x\,dx=0\). For the slices of Lemma~\ref{diam},
the four monomial integrals give
\[
 \begin{aligned}
 &36\left(\min_b\int_0^1(t-b)^2L(t)\,dt\right)^2
 -15\int_0^1\left(t-\int_0^1 yL(y)\,dy\right)^4L(t)\,dt\\
 &=\left[\sum_jc_j\left(s_j-\sum_i c_i s_i\right)^2\right]^2
 +\sum_jc_j\left(s_j-\sum_i c_i s_i\right)^2s_j(1-s_j)\ge0 .
 \end{aligned}
\]
Pass to the limit in \(L^1(0,1)\) and set
\[
 \begin{gathered}
 Z(x)=\frac{e\cdot x}{\sqrt{|P|^{-1}\int_P(e\cdot y)^2\,dy}},
 \qquad k=\frac1{|P|}\int_P Z^4\,dx\in[1,12/5],\\
 u=Z+\frac{k-3}{6k}\left(Z^3-\frac1{|P|}\int_P Z^3\,dx\right).
 \end{gathered}
\]
Since \(\int_P Z=\int_P u=0\), \(\int_P Z^2=|P|\), Cauchy--Schwarz and
differentiation give
\[
 \begin{aligned}
 \frac{1}{|P|}\int_P u^2
 &\ge\left(\frac{1}{|P|}\int_P uZ\right)^2=\frac{(k+3)^2}{36},\\
 \frac{\int_P(e\cdot y)^2\,dy}{|P|^2}\int_P|\nabla u|^2
 &=\frac{1}{|P|}\int_P\left(1+\frac{k-3}{2k}Z^2\right)^2
 =\frac{(k-1)(k+3)}{4k},\\
 \frac{\mu_1(P)}{|P|}\int_P(e\cdot x)^2\,dx
 &\le\frac{9(k-1)}{k(k+3)}\le\frac{35}{36}\qquad(1\le k\le12/5).
 \end{aligned}
\]
Maximize over \(|e|=1\) and apply Lemma~\ref{diam}.
\end{proof}

Kr\"oger's sharp diameter bound \cite{Kroger} gives
\(\mu_1(P)\operatorname{diam}(P)^2\le4j_{0,1}^2<70/3\).
We use the directional integral estimate above.

\begin{proposition}\label{balanced-centroid-offset}
Let \(P\) be a convex polygon of area one, and let \(z_0=f(0)\) for
a map satisfying \eqref{balance} for the profile
\eqref{dipole-polynomial-profile}. For \(\beta=\pi|f'(0)|^2\),
\[
 \left|\int_P x\,dx-z_0\right|
 \le\frac{1-\beta}{\beta}\left[
 0.301|f'(0)|+\sqrt{8\max_{|e|=1}\int_P
 \left(e\cdot x-\int_P e\cdot y\,dy\right)^2\,dx}\right].
\]
In particular, \(\Phi(P)\ge10.504435\) implies
\(\left|\int_P x\,dx-z_0\right|<0.039\).
\end{proposition}

\begin{proof}
Normalize the map by
\[
 F(w)=\frac{f(w)-f(0)}{f'(0)}=w+g(w),\qquad
 g'(w)=\sum_{j\ge1}b_jw^j,\qquad
 \frac1\pi\int_{\mathbb D}|g'|^2\,dA
 =\sum_{j\ge1}\frac{|b_j|^2}{j+1}=\beta^{-1}-1 .
\]
Write the balancing trial as \(V(re^{i\theta})=\psi(r)e^{i\theta}\).
By \eqref{balance} and angular integration,
\begin{equation}\label{centroid-moment-expansion}
 \begin{aligned}
 \frac{\int_P x\,dx-f(0)}{\beta f'(0)}
 &=\frac{1}{\pi}\int_{\mathbb D}
    \left(F-\frac{V}{4\int_0^1 r^2\psi(r)\,dr}\right)|1+g'|^2\,dA\\
 &=\sum_{j\ge1}\frac{b_j\overline b_{j+1}}{j+1}
 \left(1-\frac{(j+1)\int_0^1 r^{2j+2}\psi(r)\,dr}
                   {2\int_0^1 r^2\psi(r)\,dr}\right)
 +\frac{1}{\pi}\int_{\mathbb D}g|g'|^2\,dA .
 \end{aligned}
\end{equation}
The coefficient from \(w|g'|^2+g\overline{g'}\) is
\(1/(j+2)+1/[(j+1)(j+2)]=1/(j+1)\). The terms linear in \(b_1\)
cancel. Pass to the limit by \(L^2\) truncation of \(g'\) and radial
dilation of bounded \(g\).

For \(\psi(r)=r-0.39r^3+0.034r^5\), the coefficient of
\(b_j\overline b_{j+1}/\sqrt{(j+1)(j+2)}\) is \(a_{j+1}\), where
\begin{equation}\label{centroid-jacobi}
 a_x=\sqrt{\frac{x+1}{x}}
 \left[1-\frac{1}{0.757}
       \left(\frac{x}{x+1}
        -\frac{0.39 x}{x+2}+\frac{0.034 x}{x+3}\right)\right].
\end{equation}
For \(x\ge2\) the inner sum increases, by
\((x+1)^{-2}-0.78(x+2)^{-2}>0\), to \(1-0.39+0.034<0.757\).
Hence \(a_x\) is positive and decreasing. Numerical
Lemma~\ref{localization-numerical} gives
\[
 \left|\sum_{j\ge1}\frac{a_{j+1}b_j\overline b_{j+1}}{\sqrt{(j+1)(j+2)}}\right|
 \le0.301\sum_{j\ge1}\frac{|b_j|^2}{j+1}.
\]

By the Marx--Strohh\"acker theorem \cite[Theorem~2.6(a), p.~57]{MM},
the normalized convex map \(F\) satisfies \(\Re(F(w)/w)>1/2\).
Thus
\begin{equation}\label{centroid-convex-integral}
 \frac1\pi\left|\int_{\mathbb D} g|g'|^2\,dA\right|
 \le (\beta^{-1}-1)\|F\|_\infty,
 \qquad |g(w)|=|F(w)-w|\le |F(w)|.
\end{equation}

For \(\int_P y\,dy=0\), the first two moments of the slices of
Lemma~\ref{diam} also give
\[
 \begin{aligned}
 \left(1-\int_0^1 tL(t)\,dt\right)^2
 &=\tfrac19(2-\textstyle\sum_j c_j s_j)^2\
 \le \tfrac49\left[1-\textstyle\sum_jc_js_j+
                              \left(\textstyle\sum_j c_js_j\right)^2\right]\\
 &\le8\min_b\int_0^1(t-b)^2L(t)\,dt ,
 \end{aligned}
\]
where \((2-t)^2\le4(1-t+t^2)\). Passing to the limit, rescaling and
maximizing over all directions gives
\begin{equation}\label{centroid-radii}
 \sup_{x\in P}|x|^2\le\frac8{|P|}\max_{|e|=1}\int_P(e\cdot y)^2\,dy .
\end{equation}
Now \eqref{centroid-moment-expansion}--\eqref{centroid-convex-integral} give
\[
 \begin{aligned}
 \left|\int_P x\,dx-f(0)\right|
 &\le(1-\beta)\left(0.301|f'(0)|+\sup_{x\in P}|x-f(0)|\right)\\
 &\le(1-\beta)\left(0.301|f'(0)|+\sup_{x\in P}\left|x-\int_P y\,dy\right|
                         +\left|\int_P x\,dx-f(0)\right|\right).
 \end{aligned}
\]
Use \eqref{centroid-radii} after translation. The numerical bound follows
from Proposition~\ref{improved-conformal} and Numerical
Lemma~\ref{localization-numerical}.
\end{proof}

\begin{lemma}[Conformal radius and the sectors at the pole]\label{fan-shape-bound}
Let \(z_0=f(0)\) for a conformal map \(f:\mathbb D\to P\) onto a
convex \(m\)-gon. Put \(\beta=\pi|f'(0)|^2/|P|\), and divide \(P\) into
the triangles \(T_j=(z_0,p_j,p_{j+1})\). Write \(2\pi a_j\) for the angle
at \(z_0\), \(\pi(y_j\pm d_j)\) for the base angles, and
\(b_j=|T_j|/|P|\). Thus \(y_j=1/2-a_j\) and
\(\sum_j a_j=\sum_j b_j=1\). For \(n\ge m\) let
\[
 F(a)=\log\frac{\Gamma(1+a)^2\Gamma(1-2a)}
                   {\Gamma(1-a)^2\Gamma(1+2a)}
 =\sum_{\substack{k\ge3\\k\ {\rm odd}}}
   \frac{(2^{k+1}-4)\zeta(k)}{k}\,a^k,
 \qquad \beta_n=e^{-F(1/n)} .
\]
Then \(\beta\le\beta_m\le\beta_n\), and
\begin{equation}\label{fan-deficit}
 \begin{aligned}
 \log\frac{\beta_n}{\beta}\ \ge{}&
  \sum_j a_j\log\frac{a_j}{b_j}
  +\sum_j a_j F(a_j)-F(1/n)\\
 &+\sum_j a_j\sum_{h\ge0}\log
   \frac{1-d_j^2/(h+1-y_j)^2}{1-d_j^2/(h+y_j)^2}.
 \end{aligned}
\end{equation}
All three terms on the right are nonnegative. Moreover,
\begin{equation}\label{fan-angular-remainders}
 \begin{aligned}
 \sum_j a_jF(a_j)-F(1/n)
 &\ge4\zeta(3)\left(\sum_j a_j^4-n^{-3}\right),\\
 \sum_{h\ge0}\log\frac{1-d^2/(h+1-y)^2}{1-d^2/(h+y)^2}
 &\ge d^2\big(y^{-2}-(1-y)^{-2}\big)
 \qquad (|d|<y<1/2).
 \end{aligned}
\end{equation}
\end{lemma}

\begin{proof}
Reduced-modulus splitting \cite{Solynin,BHS} and the triangle
Schwarz--Christoffel formula with prevertices \((0,1,\infty)\) \cite{DT} give
\[
 \frac{\log|f'(0)|}{2\pi}\le\sum_j a_j^2m(T_j),\qquad
 m(T_j)=\frac{\log4}{\pi}
   +\frac{1}{2\pi a_j}\log\frac{|p_j-z_0|}{2a_j B(2a_j,y_j+d_j)} .
\]
Average over the two orders of the base vertices and substitute
\(|p_j-z_0||p_{j+1}-z_0|=2|T_j|/\sin(2\pi a_j)\).
The reflection and duplication formulas give
\[
 -\log\beta\ge\sum_j a_j\log(a_j/b_j)+\sum_j a_jF(a_j)
 +\sum_j a_j\log\frac{\Gamma(y_j+d_j)\Gamma(y_j-d_j)\Gamma(1-y_j)^2}
 {\Gamma(1-y_j-d_j)\Gamma(1-y_j+d_j)\Gamma(y_j)^2}.
\]
Expand the last term by the Euler product. In the regular case,
\[
 \frac{f'(w)}{f'(0)}=(1-w^n)^{-2/n},\qquad
 \frac{\pi|f'(0)|^2}{|R_n|}
 =\frac{2\pi n}{\sin(2\pi/n)B(1/n,1-2/n)^2}=e^{-F(1/n)}.
\]
The first two terms are nonnegative by convexity of \(t\log t\)
and \(\sum_j a_j^{k+1}\ge m^{-k}\ge n^{-k}\). For the last,
\[
 \log\frac{1-d^2/(1-y)^2}{1-d^2/y^2}
 =\sum_{i\ge1}\frac{d^{2i}}{i}\left(y^{-2i}-(1-y)^{-2i}\right)
 \ge d^2\left(y^{-2}-(1-y)^{-2}\right).
\]
\end{proof}

\begin{proposition}[Distance from the balancing pole to the centroid]
\label{bessel-centroid}
Let \(|P|=1\), \(P\) convex, and let \(f\) satisfy \eqref{balance}
with \(V(re^{i\theta})=\varphi(r)e^{i\theta}\) and
\(r\le\varphi(r)\le1\). Then
\[
 \left|\int_P x\,dx-f(0)\right|
 \le(\beta^{-1}-1)\left(\sqrt{3/2}\,|f'(0)|
       +\sup_{x\in P}\left|x-\int_P y\,dy\right|\right).
\]
In particular, \eqref{centroid-radii} and
Lemma~\ref{cubic-directional-integral} bound the supremum by
\(\sqrt{8(35/36)/\Phi(P)}\).
\end{proposition}

\begin{proof}
Apply \eqref{centroid-moment-expansion} with \(\psi=\varphi\). For
\(j\ge1\),
\[
 0\le
 \frac{(j+1)\int_0^1r^{2j+2}\varphi(r)\,dr}{2\int_0^1r^2\varphi(r)\,dr}
 \le\frac{2(j+1)}{2j+3}<1,
 \qquad \sqrt{\frac{j+2}{j+1}}\le\sqrt{\frac32}.
\]
Cauchy--Schwarz and \eqref{centroid-convex-integral} give
\[
 \left|\int_P x\,dx-f(0)\right|
 \le(1-\beta)\left(\sqrt{\tfrac32}|f'(0)|+\sup_{x\in P}|x-f(0)|\right).
\]
Apply the triangle inequality as in Proposition~\ref{balanced-centroid-offset}
and use \eqref{centroid-radii}.
\end{proof}

For \(0<a,b<1\), define
\begin{equation}\label{sector-log-function}
 \mathcal L(a,b)
 =a\log\frac ab+(1-a)\log\frac{1-a}{1-b}
 =\int_a^b\frac{t-a}{t(1-t)}\,dt .
\end{equation}
For \(n\ge2\), \(a_i,b_i>0\), \(\sum_i a_i=\sum_i b_i=1\), convexity of
\(t\log t\) gives
\begin{equation}\label{sector-log-sum-bound}
 \begin{aligned}
 \sum_i a_i\log\frac{a_i}{b_i}
 &\ge a_k\log\frac{a_k}{b_k}
     +(1-a_k)\log\frac{1-a_k}{1-b_k}\\
 &=\int_{a_k}^{b_k}\frac{t-a_k}{t(1-t)}\,dt
 =\mathcal L(a_k,b_k)\ge0 .
 \end{aligned}
\end{equation}

\begin{lemma}[Translation of the sector pole]\label{sector-pole-translation}
Let \(p,p'\) be counterclockwise successive vertices of a convex
polygon, and let \(z,z'\) be interior poles. At \(z\), write \(\theta=\angle(p-z,p'-z)\),
\(\rho=\log(|p'-z|/|p-z|)\), and \(b=|\operatorname{conv}(z,p,p')|\).
Use \(\theta',\rho'\) at \(z'\). If
\[
 |z'-z|\le\Delta<\min(|p-z|,|p'-z|),
\]
then
\begin{equation}\label{center-translation}
 |(\rho'-\rho)+i(\theta'-\theta)|
 \le\frac{\Delta\sqrt{2\cosh\rho-2\cos\theta}}
           {\sqrt{2b/\sin\theta}-\Delta e^{|\rho|/2}} .
\end{equation}
\end{lemma}

\begin{proof}
Along \([z,z']\) the argument of \((p'-w)/(p-w)\) lies in \((0,\pi)\).
For \(w=z+s(z'-z)\), \(0\le s\le1\), its logarithmic derivative is bounded by
\[
 \left|\frac{d}{dw}\log\frac{p'-w}{p-w}\right|
 =\frac{|p'-p|}{|p-w||p'-w|}
 \le\frac{|p'-p|}{|p-z||p'-z|
       (1-s\Delta/\min(|p-z|,|p'-z|))^2}.
\]
Integrate in \(s\) and use
\[
 \begin{gathered}
 |p-z||p'-z|=\frac{2b}{\sin\theta},\qquad
 \min(|p-z|,|p'-z|)=\sqrt{\frac{2b}{\sin\theta}}e^{-|\rho|/2},\\
 \frac{|p'-p|^2}{|p-z||p'-z|}=2\cosh\rho-2\cos\theta .
 \end{gathered}\qedhere
\]
\end{proof}

\begin{proposition}\label{prop:fan-localization}
Let \(|P|=1\), \(P\) a convex hexagon, \(\Phi(P)\ge10.504435\), and
\(z=f(0)\) as in \eqref{balance}. On \(T_j=\operatorname{conv}(z,p_j,p_{j+1})\)
let the three angles be \(2\pi a_j,\pi\beta_j,\pi\gamma_j\) in that order.
Set \(b_j=|T_j|/|P|\), \(d_j=(\beta_j-\gamma_j)/2\). Then
\[
 0.036<a_j<0.254,\qquad
 Q(a_j)+\mathcal L(a_j,b_j)
       +0.0948\,\rho_j^2<0.01268,
\]
where \(\rho_j=\log(|p_{j+1}-z|/|p_j-z|)\) and \(Q\) is given by \eqref{fan-angular}.
\end{proposition}

\begin{proof}
The base-angle term in \eqref{fan-deficit}, with \(y=1/2-a\), is
\begin{equation}\label{hex-gamma-second-difference}
 \log\frac{\Gamma(y+d)\Gamma(y-d)\Gamma(1-y)^2}
          {\Gamma(1-y-d)\Gamma(1-y+d)\Gamma(y)^2}
 =\sum_{k\ge0}\log\frac{1-d^2/(k+1-y)^2}{1-d^2/(k+y)^2}.
\end{equation}
By \eqref{fan-deficit} at \(n=6\), Proposition~\ref{improved-conformal}
and Numerical Lemma~\ref{localization-numerical},
\[
 0\le\log(\beta_6/\beta)<\log(0.9763/0.964)<0.01268.
\]
For each individual \(a_j=x\), Jensen bounds the angular term below
by
\begin{equation}\label{fan-angular}
 Q(x)=\sum_{\substack{3\leq k\leq13\\k\ {\rm odd}}}
 \frac{2^{k+1}-4}{k}\left(\sum_{i=1}^{30}i^{-k}\right)
 \left[x^{k+1}+5((1-x)/5)^{k+1}-6^{-k}\right].
\end{equation}
By convexity, \(Q'(1/6)=0\) and
\(\min\{Q(0.017),Q(0.263)\}>0.0158\)
(Numerical Lemma~\ref{localization-numerical}), so \(0.017<a_j<0.263\).
Apply \eqref{sector-log-sum-bound} with \(\sum_i b_i=\sum_i a_i=1\).
By the sine rule,
\begin{equation}\label{fan-radius-angle-identity}
 \tanh(|\rho_j|/2)
 =\frac{|\sin(\pi\beta_j)-\sin(\pi\gamma_j)|}
        {\sin(\pi\beta_j)+\sin(\pi\gamma_j)}
 =\tan(\pi a_j)\tan(\pi |d_j|).
\end{equation}
Hence \(e^{|\rho_j|}\ge8/5\) implies
\[
 |d_j|\ge d_0(a_j),\qquad
 d_0(a)=\pi^{-1}\arctan[(3/13)\cot(\pi a)].
\]
The \(k=0\) summand of \eqref{hex-gamma-second-difference} increases
in \(|d|\) and decreases in \(y\). By \eqref{radial-exclusion-number},
\[
 e^{|\rho_j|}\ge8/5\quad\Longrightarrow\quad
 \log(\beta_6/\beta)>0.01761>0.01268.
\]
Hence \(\tanh(|\rho_j|/2)<3/13\).
Numerical Lemma~\ref{radial-multiplier} yields \(0.036<a_j<0.254\)
and total base-angle contribution \(\ge0.079\sum_i\rho_i^2\).
Since \(\sum_i\rho_i=0\),
\[
 \sum_i\rho_i^2\ge\rho_j^2+\tfrac15\Big(\sum_{i\ne j}\rho_i\Big)^2
 =\tfrac65\rho_j^2,\qquad 0.079\cdot\tfrac65=0.0948 .
\]
Substitution in \eqref{fan-deficit} proves the claim.
\end{proof}

\begin{numericallemma}[Localization bounds]\label{localization-numerical}
For the Jacobi matrix with off-diagonal entries \(a_{j+1}/2\)
from \eqref{centroid-jacobi}, the upper bound is \(0.301 I\).
The centroid estimate in Proposition~\ref{balanced-centroid-offset}
then gives \(\left|\int_P x\,dx-z_0\right|<0.039\) when \(\Phi(P)\geq10.504435\).
For the functions in \eqref{fan-deficit} and
\eqref{fan-angular},
\[
 F(1/6)>-\log(0.9763),\quad
 \log(0.9763/0.964)<0.01268,\quad
 Q(0.017),\,Q(0.263)>0.0158.
\]
With \(d_0(a)=\pi^{-1}\arctan[(3/13)\cot(\pi a)]\), one also has
\begin{equation}\label{radial-exclusion-number}
 Q(a)+a\log
 \frac{1-d_0(a)^2/(1/2+a)^2}{1-d_0(a)^2/(1/2-a)^2}
 >0.01761,\qquad 0.017\leq a\leq0.263 .
\end{equation}
\end{numericallemma}

\begin{proof}
For \(0.301I\) minus the Jacobi matrix, symmetric elimination gives
\[
 p_1=0.301,\qquad p_{j+1}=0.301-\frac{a_{j+1}^2}{4p_j}.
\]
Rational substitution gives \(p_1,\ldots,p_8>0\), \(p_8>0.208\).
For \(j\ge8\), \(a_{j+1}\le a_9<0.21\), so
\(p_j\ge0.15\Rightarrow p_{j+1}\ge0.301-0.21^2/(4\cdot0.15)>0.15\).
Positivity of all compressions gives the \(\ell^2\) bound. By
Lemma~\ref{cubic-directional-integral},
\[
 \max_{|e|=1}\int_P\left(e\cdot x-\int_P e\cdot y\,dy\right)^2\,dx
 \le 35/(36\cdot10.504435)<0.09258 .
\]
Since \(0.964<\beta\le1\), \(|f'(0)|=\sqrt{\beta/\pi}<0.565\) and
\[
 \left|\int_P x\,dx-z_0\right|<\frac{1-0.964}{0.964}\bigl(0.301\cdot0.565+0.861\bigr)<0.039 .
\]

For \(F(1/6)\) use \(3\le k\le13\), \(\zeta(k)\ge\sum_{i=1}^{30}i^{-k}\)
and \(\log(1+x)\le x-x^2/2+x^3/3\), \(0<x<1\). For
\eqref{radial-exclusion-number} put \(I=[l,u]=[j,j+1]/1000\),
\(17\le j<263\). Choose rational
\(0\le q<1/2\) with \(\tan(\pi q)\tan(\pi u)\le3/13\) by
alternating sine and cosine bounds. The left side on \(I\) is at least
\[
 Q(\operatorname{proj}_I(1/6))+
 l\log\frac{1-q^2/(1/2+l)^2}{1-q^2/(1/2-l)^2}.
\]
All \(246\) lower bounds exceed \(0.01761\), and \(Q(0.017),Q(0.263)>0.0158\).
\end{proof}

\begin{numericallemma}[Base-angle term]\label{radial-multiplier}
Let \(Q\) be the angular lower bound in the proof of
Proposition~\ref{prop:fan-localization}. For \(0.017\leq a\leq0.263\)
and \(0\leq t\leq3/13\), set \(y=1/2-a\) and
\[
 d=\frac{\arctan(t\cot(\pi a))}{\pi},\qquad
 B(a,t)=a\sum_{k\ge0}\log\frac{1-d^2/(k+1-y)^2}{1-d^2/(k+y)^2}.
\]
If \(Q(a)+B(a,t)\leq0.01268\), then
\[
 0.036<a<0.254,\qquad
 B(a,t)\geq0.079\,\bigl(2\operatorname{arctanh}t\bigr)^2 .
\]
\end{numericallemma}

\begin{proof}
Convexity of \(Q\) and
\(\min\{Q(0.036),Q(0.254)\}>0.01268\) give the angular bounds.
For \(t>0\) expand the logarithms in the nonnegative series
\eqref{hex-gamma-second-difference}. Truncation gives
\[
 \frac{B(a,t)}{a\,d^2}\ \geq\
 \sum_{n=0}^{2}\sum_{i=1}^{11}
 \frac{(n+y)^{-2i}-(n+1-y)^{-2i}}{i}\,d^{\,2i-2}.
\]
Each summand increases in \(a,d\). On
\([a_0,a_1]\times[t_0,t_1]\) choose \(0\le d_0\le d(a_1,t_0)\).
Write \(L\) for the displayed sum at \((a_0,d_0)\). Since \(d\ge d_0\),
\[
 Q(a)+B(a,t)\ge
 \min_{a_0\leq a\leq a_1}Q(a)+a_0d_0^2L .
\]
Also \(d/t\) decreases in \(t\) for fixed \(a\), and decreases in \(a\)
for fixed \(t\). Since
\[
 \frac{\operatorname{arctanh}t}{t}
 \leq 1+\frac{t^2}{3}+\frac{t^4}{5(1-t^2)},
\]
the quotient \(B(a,t)/(2\operatorname{arctanh}t)^2\) on the
rectangle is at least
\[
  a_0 L\left[
      \frac{d_1}
      {2t_1(1+t_1^2/3+t_1^4/[5(1-t_1^2)])}
      \right]^2 ,
\]
where \(d_1\) is a rational lower enclosure of \(d(a_1,t_1)\).
Alternating sine and cosine bounds check
\(\tan(\pi a_1)\tan(\pi d_i)\le t_i\) \((i=0,1)\).
A dyadic partition of \([0.036,0.254]\times[0,3/13]\) has depth
\(\le19\) and \(6146\) leaves: \(792\) have \(Q+B>0.01268\);
the other \(5354\) have \(B/(2\operatorname{arctanh}t)^2>0.079\)
for \(t>0\). At \(t=0\) both sides vanish.
\end{proof}

Give \(R_N\) circumradius one and translate and scale \(P\) so that \(\int_P x\,dx=0,\ |P|=|R_N|\). Let \(v_j,p_j\) be their counterclockwise vertices (indices
modulo \(N\)), \(S_j=\operatorname{conv}(0,v_j,v_{j+1})\), and
\(L_j(v_j,v_{j+1})=(p_j,p_{j+1})\), as an identity of column matrices.
For \(\vartheta_j=\angle(p_j,p_{j+1})\), put
\[
 \begin{gathered}
 \alpha=\pi/N,\qquad w_j=\det L_j,\\
 r_j=\frac{\tan(\vartheta_j/2)}{\tan\alpha},\quad
 \rho_j=\log\frac{|p_{j+1}|}{|p_j|},\quad d_j=2\cosh\rho_j-2 .
 \end{gathered}
\]
In the sector bisector frame, suppress \(j\) and write
\(M=L^{\mathsf T}L/w=\left(\begin{smallmatrix}x&z\\z&y\end{smallmatrix}\right)\).
Applying \(L\) to the reference endpoints gives
\begin{equation}\label{common-fan-coordinates}
 \begin{gathered}
 x=r^{-1}+\frac{d\tan\alpha}{2\sin\vartheta},\qquad
 y=r+\frac{d\cot\alpha}{2\sin\vartheta},\qquad
 q=z^2=\frac{d(d+4)}{4\sin^2\vartheta},\\
 \vartheta=2\arctan(r\tan\alpha),\qquad
 \mathcal A(r)=\frac{\vartheta-2\alpha}{\sin(2\alpha)} .
 \end{gathered}
\end{equation}
Thus
\[
 \operatorname{sgn} z=\operatorname{sgn}\rho,\qquad
 \sum_j w_j=N,\qquad \sum_j\mathcal A(r_j)=\sum_j\rho_j=0 .
\]
For a dihedrally equivariant pair \(U_0\) choose sector trials \(U_j\)
whose transports \(W|_{L_jS_j}=L_jU_j\circ L_j^{-1}\) have matching
traces. In particular, \((U_j-U_0)|_{[0,v_j]\cup[0,v_{j+1}]}=0\)
suffices by equivariance. For a single scalar \(\gamma\) and all \(j\),
reflection and rotation give
\(\int_{S_j}U_0\,dx=\gamma\int_{S_j}x\,dx\).
Since \(\sum_jw_j L_j\int_{S_j}x\,dx=\int_Py\,dy=0\),
\begin{equation}\label{fan-sector-integral}
 \begin{aligned}
 \frac1{|P|}\int_P W\,dy
 &=\frac1{|R_N|}\sum_j w_j L_j\int_{S_j}(U_j-U_0)\,dx,\\
 \left|\frac1{|P|}\int_P W\,dy\right|^2
 &\le\frac{N}{|R_N|^2}\sum_jw_j
 \left|L_j\int_{S_j}(U_j-U_0)\,dx\right|^2 .
 \end{aligned}
\end{equation}
Hence \(M_W=|P|^{-1}\int_P|W|^2\,dy-|P|^{-2}|\int_P W\,dy|^2
\ge\sum_jw_j^2F_j\), and \(|P|^{-1}\int_P|DW|^2\,dy=\sum_jw_jG_j\), with
\begin{equation}\label{finite-sector-mass-energy}
 \begin{aligned}
 F&=\frac1{|R_N|}\int_{S_j}U^{\mathsf T}MU\,dx
 -\frac{N}{|R_N|^2}
 \left(\int_{S_j}(U-U_0)\,dx\right)^{\mathsf T}M\int_{S_j}(U-U_0)\,dx,\\
 G&=\frac1{|R_N|}\int_{S_j}
 \operatorname{tr}\big((DU)^{\mathsf T}M(DU)M^{-1}\big)\,dx .
 \end{aligned}
\end{equation}
For \(U=(u,v),\,U_0=(u_0,v_0)\), set
\[
 \begin{aligned}
 (A,B,C,D,E)&=\frac1{|R_N|}\int_{S_j}
       (u_x^2,u_y^2,v_x^2,v_y^2,u_xv_y+u_yv_x)\,dA,\\
 m&=\frac1{|R_N|}\int_{S_j}u^2\,dA
 -\frac{N}{|R_N|^2}\left(\int_{S_j}(u-u_0)\,dA\right)^2,\\
 n&=\frac1{|R_N|}\int_{S_j}v^2\,dA
 -\frac{N}{|R_N|^2}\left(\int_{S_j}(v-v_0)\,dA\right)^2 .
 \end{aligned}
\]
Since \(xy-z^2=1\),
\begin{equation}\label{finite-sector-signed-expansion}
 \begin{aligned}
 F={}&mx+ny+\frac{2z}{|R_N|}\int_{S_j}uv\,dA
 -\frac{2Nz}{|R_N|^2}\int_{S_j}(u-u_0)\,dA\int_{S_j}(v-v_0)\,dA,\\
 G={}&A+D+Bx^2+Cy^2+(A+D-2E)z^2\\
 &+\frac{2z}{|R_N|}\int_{S_j}(yv_x-xu_y)(u_x-v_y)\,dA .
 \end{aligned}
\end{equation}
For reflection-equivariant \(U\), both terms linear in \(z\) vanish.
Reflecting a zero-trace correction carries the estimate for \(z\ge0\) to \(z<0\).

With \(\nu>0\) and \(c\in\mathbb R\) we have
\[
 \mu_1(P)M_W\le\sum_jw_jG_j,\qquad
 \sum_j\left[\frac{2\nu}{N}(w_j-1)-c\mathcal A(r_j)\right]=0 .
\]
Together with \(M_W\ge\sum_jw_j^2F_j\), these give
\begin{equation}\label{finite-sector-direct-comparison}
 (\nu-\mu_1(P))M_W\geq-\sum_j
 \bigl[w_j(2\nu/N+G_j)-\nu w_j^2F_j-2\nu/N-c\mathcal A(r_j)\bigr],
 \qquad c\in\mathbb R .
\end{equation}
For \(F_j>0\), the identity
\[
 w_j(2\nu/N+G_j)-\nu F_jw_j^2
 =\frac{(2\nu/N+G_j)^2}{4\nu F_j}
 -\nu F_j\left(w_j-\frac{2\nu/N+G_j}{2\nu F_j}\right)^2
\]
yields
\begin{equation}\label{finite-sector-square-completion}
 (\nu-\mu_1(P))M_W\geq
     -\sum_j\left[\frac{(2\nu/N+G_j)^2}{4\nu F_j}-\frac{2\nu}{N}-c\mathcal A(r_j)\right].
\end{equation}
Hold \(r\) and the reflection-equivariant trial \(U\) fixed.
Then \(F\) is affine in \(d\) and \(G\) is quadratic. Positive radii realize every \(d\ge0\),
so \(G(d)\ge0\) and \(\partial_d^2 G\ge0\).
For \(F>0\) and \(a>0\),
\begin{equation}\label{sector-square-convexity}
 \partial_d^2\!\left(\frac{(a+G)^2}{F}\right)
 =\frac2F\left(\partial_dG-\frac{(a+G)\partial_d F}{F}\right)^2
  +\frac{2(a+G)\partial_d^2G}{F}\ge0 .
\end{equation}

For two positive sector masses at \((r,0)\) and \((r,d)\), write
\(F_1=(F(r,d)-F(r,0))/d\), \(G_1=(G(r,d)-G(r,0))/d\).
With \(a=2\nu/N\) and \(T=(a+G)^2/(4\nu F)\), subtraction gives
\begin{equation}\label{finite-fan-radial-quotient}
 \frac{T(r,d)-T(r,0)}d
 =\frac{[2(a+G(r,0))+dG_1]G_1-(a+G(r,0))^2F_1/F(r,0)}
              {4\nu F(r,d)} .
\end{equation}
Subtracting an angular term leaves this quotient fixed. At \(d=0\)
it extends by differentiation whenever \(F,G\) are differentiable in \(d\).

For an exact trial at \(d=0\), write \(\nu=\lambda=\mu_1(R_N)\),
\(f_0=2/N\), \(\xi=r-1\), and
\[
 F=f_0+\xi f_1,\qquad G=a+\xi g_1,\qquad
 h_1=g_1(1)-\lambda f_1(1),\qquad
 h_2=\frac{g_1-\lambda f_1-h_1}{\xi}.
\]
Here \(F(1,0)=f_0\), \(G(1,0)=a\) and \(h_1=T_r(1,0)\).
Expanding \(T-a=[f_0\xi(g_1-\lambda f_1)+\xi^2g_1^2/(4\lambda)]/F\)
gives, for \(B_0=T-a-h_1\mathcal A\),
\begin{equation}\label{finite-fan-angular-quotient}
 \frac{B_0(r,0)}{\xi^2}
 =\frac{f_0h_2-h_1f_1+g_1^2/(4\lambda)}{F(r,0)}
 -h_1\frac{\mathcal A(r)-\xi}{\xi^2}.
\end{equation}
For \(t=\tan\alpha\) and \(v=t\xi/(1+t^2r)\), the last quotient is
\[
 \frac{\mathcal A(r)-\xi}{\xi^2}
 =-\frac{t^2}{1+t^2r}-\frac{(1+t^2)t^2\xi}{(1+t^2r)^3}
                           \frac{1-\arctan(v)/v}{v^2}.
\]
At \(r=1\) this equals \(-\sin^2\alpha\). In particular, for the
uncorrected exact pair, \eqref{finite-sector-signed-expansion} gives
\begin{equation}\label{uncorrected-exact-angular}
 \begin{aligned}
 f_1&=n-m/r,&
 g_1&=(r+1)(C-B/r^2),\\
 h_1&=2(C-B)-\lambda(n-m),&
 h_2&=C+B(2r+1)/r^2-\lambda m/r .
 \end{aligned}
\end{equation}

\begin{lemma}[Equivariant spectral projection]\label{equivariant-sector-projection}
For \(N\ge5\), normalize the first equivariant pair \(U\) and a trial pair
\(V\) in the same equivariant class by
\(\int_{R_N} U_i^2=\int_{R_N}V_i^2=|R_N|\), \(i=1,2\).
Let
\[
 \lambda_-<\mu_1(R_N)\le\mathcal R(V)\le\lambda_+<\Lambda,
 \qquad 2(\lambda_+-\lambda_-)<\Lambda-\lambda_- ,
\]
with \(\Lambda\) a lower bound for the next eigenvalue in the same signed
rotation class. Put
\begin{equation}\label{finite-trial-pair-constants}
 e=\frac{2(\lambda_+-\lambda_-)}{N(\Lambda-\lambda_+)},\qquad
 \eta=\sqrt{1-Ne/2},\qquad
 \delta^2=\frac{2e}{1+\eta},\quad
 \epsilon^2=\Lambda e+\frac{2\lambda_+}N(1-\eta)^2 .
\end{equation}
For the norm on one sector \(S\) given by \(|R_N|^{-1}\int_S|\cdot|^2\),
the sign of \(U\) may be chosen so that
\begin{equation}\label{finite-trial-pair-errors}
 \|U-V\|_S\le\delta,\qquad \|DU-DV\|_S\le\epsilon .
\end{equation}
\end{lemma}

\begin{proof}
The projection of \(V\) onto the first pair commutes with rotation and
reflection, so it equals \(\gamma U\). Orient \(U\) so \(\gamma\ge0\).
By \eqref{spectral-projection-bounds} and rotation, the squared sector
mass and energy of \(V-\gamma U\) are at most \(e,\Lambda e\).
Orthogonality for both forms restricts to each sector by rotation. Thus
\[
 \begin{gathered}
 \gamma^2=1-\tfrac N2\|V-\gamma U\|_S^2\ge\eta^2,\qquad
 \|U-V\|_S^2=\frac{2\|V-\gamma U\|_S^2}{1+\gamma},\\
 \|DU-DV\|_S^2=\|D(V-\gamma U)\|_S^2+\frac{2\mu_1(R_N)}N(1-\gamma)^2.
 \end{gathered}\qedhere
\]
\end{proof}

\begin{lemma}[Central and exterior sector comparison]\label{central-exterior-comparison}
Normalize \(|P|=|R_N|\), and write \(\lambda=\mu_1(R_N)\).
Classify each sector as central or exterior. For a piecewise trial with
matching traces let \(R_{\nu,j}\) satisfy
\[
 (\nu-\mu_1(P))M_W\ge-\sum_{j=1}^N R_{\nu,j}.
\]
One may use either \eqref{finite-sector-direct-comparison} or, when
\(F_j>0\), \eqref{finite-sector-square-completion}.
\begin{enumerate}
\item If all sectors are central and a trial has
\(R_{\lambda,j}\le-\epsilon[(r_j-1)^2+d_j]\) with \(\epsilon>0\), then
\(\mu_1(P)\le\lambda\), with equality only if \(P\cong R_N\).
\item If some sector is exterior and a trial has
\(R_{\nu,j}\le u\) on central sectors and \(R_{\nu,j}\le-v\) on
exterior sectors, where \(0\le(N-1)u<v\) and \(\nu<\lambda\), then
\(\mu_1(P)<\nu\). The same conclusion follows whenever
\(\sum_j R_{\nu,j}<0\).
\end{enumerate}
\end{lemma}

\begin{proof}
In (1), for \(\mu_1(P)\ge\lambda\) we have
\[
 0\ge (\lambda-\mu_1(P))M_W\ge\epsilon\sum_j[(r_j-1)^2+d_j]
 \quad\Longrightarrow\quad
 \vartheta_j=2\pi/N,\quad |p_{j+1}|=|p_j|\quad(1\le j\le N).
\]
The centroid angles and the common radii give \(P\cong R_N\).
In (2),
\[
 (\nu-\mu_1(P))M_W\ge-\sum_jR_{\nu,j}
 \ge v-(N-1)u>0 .
\]
Any strictly negative residual sum also gives \((\nu-\mu_1(P))M_W>0\).
\end{proof}

\begin{theorem}[A comparison from the centroid fan]\label{covariantfan}
Translate the convex hexagon \(P\) so that \(\int_P x\,dx=0\).
For its cyclic vertices \(p_j\), set
\(\vartheta_j=\angle(p_j,p_{j+1})\in(0,\pi)\),
\[
 r_j=\sqrt3\tan(\vartheta_j/2),\qquad
 d_j=\frac{|p_j|}{|p_{j+1}|}+\frac{|p_{j+1}|}{|p_j|}-2.
\]
Suppose \(0.75\leq r_j\leq4\) for every \(j\) and
\[
 d_j\leq 0.10
  +0.15\,\mathbf 1_{\{r_j\geq0.78\}}
  +0.125\,\mathbf 1_{\{r_j\geq0.80\}}
  +0.125\,\mathbf 1_{\{r_j\geq1\}} .
\]
Then \(\Phi(P)\leq\Phi(H)\), with equality precisely when \(P\)
is similar to \(H\).
\end{theorem}

\begin{proof}
Normalize \(|P|=|H|\) and set \(U_j=U_0\), the exact equivariant pair with
unit component masses. By \eqref{fan-sector-integral},
\(\int_P W\,dy=0\). Write
\[
 m+n=1/3,\qquad a=A+B+C+D=\lambda/3,\qquad
 K=A+D-2E,\qquad \lambda=\mu_1(H).
\]
With \(\tau=(3+r^2)/(4r)\), \eqref{common-fan-coordinates} gives
\[
 x=r^{-1}+\tau d/3,\qquad y=r+\tau d,\qquad z^2=\tau^2d(d+4)/3 .
\]
Hence \eqref{finite-sector-signed-expansion} reduces to
\begin{equation}\label{hexagonal-sector-forms}
 \begin{aligned}
 F(r,d)&=m/r+nr+\tau(m/3+n)d,\\
 G(r,d)&=A+D+B/r^2+Cr^2+2\tau(B/(3r)+Cr)d
          +\tau^2d\left[(B/9+C)d+K(d+4)/3\right].
 \end{aligned}
\end{equation}
For \(T(r,d)=(a+G)^2/(4\lambda F)\), Numerical
Lemma~\ref{centroid-numerical} gives
\begin{equation}\label{covariantsector}
 T(r,d)-a-T_r(1,0)\mathcal A(r)
       \le-\frac{a}{1000}\big[(r-1)^2+d\big].
\end{equation}
Apply Lemma~\ref{central-exterior-comparison}(1).
\end{proof}

\subsection{The hexagon and the pentagon}\label{sec:pentagon-comparison}

\begin{proposition}\label{hexagon-comparison}
For every convex polygon \(P\) with at most six sides,
\(\Phi(P)\leq\Phi(H)\), with equality precisely when \(P\)
is similar to \(H\).
\end{proposition}

\begin{proof}
By \eqref{reference-values} the comparison is strict for
\(\Phi(P)<10.504435\). Otherwise \(\Phi(P)>10.5\), hence
\(P\) has six genuine vertices by Theorem~\ref{pentagon}.
Translate and scale so that \(\int_P x\,dx=0,\ |P|=1\). For \(z,a_j,b_j,\rho_j\)
from Proposition~\ref{prop:fan-localization}, set
\(\vartheta_j=\angle(p_j,p_{j+1})\),
\(\widetilde\rho_j=\log(|p_{j+1}|/|p_j|)\).

Proposition~\ref{balanced-centroid-offset} gives \(|z|<0.039\).
Apply Lemma~\ref{sector-pole-translation} with \(z'=0\),
\(\Delta=0.039\) and \(\theta=2\pi a_j\). Put \(q=|\rho_j|\).
Numerical
Lemma~\ref{global-mixed-fan-check} checks that the denominator
of \eqref{center-translation} is positive.
For \(0<b\le a<1/2\),
\[
 a\log\frac ab+(1-a)\log\frac{1-a}{1-b}
 =\int_b^a\frac{a-t}{t(1-t)}\,dt\ge\frac{(a-b)^2}{2a(1-a)} .
\]
Proposition~\ref{prop:fan-localization}
therefore yields
\begin{equation}\label{bridge-area}
 b_j\geq a_j-\sqrt{2a_j(1-a_j)(0.01268-Q(a_j)-0.0948\,q^2)}.
\end{equation}
The same lower bound holds for \(b_j\geq a_j\).

Dilate to \(|P|=|H|\), with \(\int_P x\,dx=0\), preserving
\(r_j=\sqrt3\tan(\vartheta_j/2)\) and \(d_j=2\cosh\widetilde\rho_j-2\).
If all sectors satisfy Theorem~\ref{covariantfan}, apply it. Otherwise
at least one exterior sector enters the following comparison.
On the reference triangle \(S_0=(0,\sqrt3/2+i/2,\sqrt3/2-i/2)\)
use the bisector coordinate \(\xi=x+iy\) and the pair
\begin{equation}\label{hexagonal-polynomial-pair}
 \begin{aligned}
 U_*&=\xi-0.498\xi^2\bar\xi-0.047\bar\xi^5
                     +0.072\xi^3\bar\xi^2+0.014\xi^7,\\
 U_0&=\xi\left[
  1-0.498(2x/\sqrt3)^2+0.119(2x/\sqrt3)^4-0.014(2x/\sqrt3)^6
  \right].
 \end{aligned}
\end{equation}
On either ray \(\xi=s e^{\pm i\pi/6}\), \(\bar\xi^5=-s^4\xi\) and
\(\xi^7=-s^6\xi\), so \(U_*=U_0\parallel\xi\).
Use the same scalar for both pairs, so that
\(\int_H(U_*)_1^2\,dA=\int_H(U_*)_2^2\,dA=|H|\). On \(L_j S_j\) put
\[
 W=L_j[(1-t_j)U_*+t_jU_0]\circ L_j^{-1},\qquad 0\le t_j\le1 .
\]
The vector \(\int_{S_j}(U_0-U_*)\,dA\) lies on the sector bisector.
Thus \eqref{finite-sector-mass-energy} subtracts
\(6t_j^2|H|^{-2}|\int_{S_j}(U_0-U_*)\,dA|^2\)
from the first diagonal sector mass in \eqref{finite-sector-signed-expansion}.
Set
\[
 \lambda_*=4.0431978,\qquad a_*=\lambda_*/3,\qquad
 T_*(r,d,t)=\frac{(a_*+G)^2}{4\lambda_*F}.
\]
Numerical Lemma~\ref{global-mixed-fan-check} gives for each centroid sector a
choice of \(t\) for which
\begin{equation}\label{hybrid-sector-test}
 T_*(r,d,t)-a_*-T_{*,r}(1,0,0)\mathcal A(r)
 < \begin{cases}
 -0.002,& (r,d)\text{ outside the ranges of Theorem~\ref{covariantfan}},\\
 0.00023,& \text{otherwise}.
 \end{cases}
\end{equation}
Apply Lemma~\ref{central-exterior-comparison}(2) with
\(\nu=\lambda_*<\mu_1(H)\) by \eqref{reference-values},
and \(0.002-5(0.00023)=0.00085>0\).
\end{proof}

\begin{numericallemma}[Exact eigenfunction trials]\label{centroid-numerical}
The sector inequality \eqref{covariantsector} holds on each range of
Theorem~\ref{covariantfan}. The sector masses satisfy \(m+n=1/3\) and \(n<1/25\).
\end{numericallemma}

\begin{proof}
By \eqref{sector-square-convexity}, it suffices to test \(d=0\)
and \(d=\kappa\), the upper endpoint in Theorem~\ref{covariantfan}.
For \(d=0\) put
\[
 L=1-6n+2(C-B)/a,\qquad
 D=1+3n(r^2-1),\qquad t=\frac{3(r-1)^2}{(r+3)^2}.
\]
Direct differentiation at \(r=1\) gives \(\partial_r(T/a)=L>0\).
Moreover, with the value at zero taken by continuity,
\[
  \vartheta-\frac\pi3
 =\frac{2\sqrt3(r-1)}{r+3}\frac{\arctan\sqrt t}{\sqrt t}.
\]
Factor \((r-1)^2\) after substituting
\(T(r,0)/a=r[2a+(r^2-1)(C-B/r^2)]^2/(4a^2D)\).
The bounds
\(1-t/3\le\arctan(\sqrt t)/\sqrt t\le1-t/3+t^2/5\)
from integration of \(1/(1+s^2)\) give
\[
 \frac{T(r,0)-a-T_r(1,0)\mathcal A(r)}a\le(r-1)^2J(r),
\]
where
\begin{align*}
 J(r)&=\frac{-3an+C(r+2)+B/r}{aD}
      +\frac{r(r+1)^2(C-B/r^2)^2}{4a^2D}\\
 &\quad+\frac{L[1-12n(r+1)]}{D(r+3)}
      +\frac{4L(r-1)}{(r+3)^3}\left(1-\frac{3t}{5}\mathbf1_{\{r<1\}}\right).
\end{align*}
For \(d=\kappa\), use \eqref{hexagonal-sector-forms} in
\[
 \begin{aligned}
 12a F(r,0)F(r,\kappa)\bigl[T(r,\kappa)-T(r,0)\bigr]
 &=F(r,0)\bigl[G(r,\kappa)-G(r,0)\bigr]
            \bigl[2a+G(r,\kappa)+G(r,0)\bigr]\\
 &\quad-[a+G(r,0)]^2\bigl[F(r,\kappa)-F(r,0)\bigr].
 \end{aligned}
\]
Let \(I_\kappa\) denote the \(r\)-interval corresponding to \(\kappa\)
in Theorem~\ref{covariantfan}. On cells of width \(1/1000\),
outward evaluation gives \(\sup_{I_\kappa}J\le-0.0048\) and
\[
 \sup_{I_\kappa}\!\left((r-1)^2J(r)
     +\frac{T(r,\kappa)-T(r,0)}a\right)
 +0.001\left(\max_{I_\kappa}(r-1)^2+\kappa\right)\le-0.0005.
\]
Now \eqref{sector-square-convexity} gives \eqref{covariantsector} on
\([0,\kappa]\). It remains to verify the enclosures used above. Normalize \(V=U_*\) from \eqref{hexagonal-polynomial-pair} by
\(\int_H V_1^2=\int_H V_2^2=|H|\).
The monomial rule \eqref{certified-sector-monomials} gives
\begin{equation}\label{sectorconstants}
\begin{gathered}
 (A_0,B_0,C_0,D_0)=
 \frac{(1114479983884,16742221548,10252618428,2315149218844)}
      {2564381509059},\\
 E_0=\frac{1478365446656}{2564381509059},\\
 (m_0,n_0)=\frac{(1528383804113,181203868593)}{5128763018118}.
\end{gathered}
\end{equation}
Indeed \eqref{certified-sector-monomials} gives
\[
 \frac{1}{|H|}\int_S x^iy^j
 =\frac{(1+(-1)^j)(\sqrt3/2)^i}{6\cdot2^j(j+1)(i+j+2)} .
\]

By rotation through \(\pi/3\), the component quotient of \(V\) is
\(\rho=3(A_0+B_0+C_0+D_0)=3456624042704/854793836353<4.043811\).
Use Lemma~\ref{equivariant-sector-projection} at
\((\lambda_-,\lambda_+,\Lambda)=(4.04206,\rho,10)\), since \(\mu_3(H)>10\).
For \(e=(\rho-4.04206)/(10-\rho)<0.000294\) in this calculation,
\(1+\sqrt{1-e}\ge2-e\) gives the bounds used for the sector norms,
\[
 \|DU-DV\|_S^2\le\delta^2=\frac{10e+\rho(e/(2-e))^2}{3},\qquad
 \|U-V\|_S^2\le\delta_m^2=\frac{2e}{3(2-e)} .
\]
Thus
\[
 |\sqrt B-\sqrt{B_0}|,\ |\sqrt C-\sqrt{C_0}|\le\delta,\qquad
 |\sqrt n-\sqrt{n_0}|\le\delta_m,\qquad m=1/3-n .
\]
Write \(V=(f,g)\), with the same normalization. The monomial rule gives
\[
 \frac1{|H|}\int_S(f_x-g_y)^2\,dA=\frac{496907402992}{2564381509059}.
\]
Expansion of \((u_x-v_y)^2-2u_yv_x\) about \(V\) and Cauchy--Schwarz give
\[
 \begin{aligned}
 \left|K-\frac1{|H|}\int_S[(f_x-g_y)^2-2f_yg_x]\,dA\right|
 &\le 2\delta
 \sqrt{\frac{2}{|H|}\int_S(f_x-g_y)^2\,dA+B_0+C_0}+2\delta^2,\\
 n&\le(\sqrt{n_0}+\delta_m)^2<0.03916<1/25 .
 \end{aligned}
\]
\end{proof}

\begin{numericallemma}[Mixed polynomial trials]\label{global-mixed-fan-check}
For the parameters of Proposition~\ref{prop:fan-localization},
transferred to the centroid by \eqref{center-translation}, there is a
parameter \(t\in[0,1]\) for which \eqref{hybrid-sector-test} holds.
The uncentered sector masses for every \(t\in[0,1]\) satisfy \(m+n\leq1/3\) and \(n<1/25\).
\end{numericallemma}

\begin{proof}
Proposition~\ref{prop:fan-localization} gives \(0.036<a_j<0.254\),
\(q=|\rho_j|\le\sqrt{0.01268/0.0948}<0.366\). On
\(I\times J=[a_-,a_+]\times[q_-,q_+]\), formula \eqref{fan-angular} gives
\[
 v=0.01268-\min_{a\in I}Q(a)-0.0948q_-^2,\qquad
 b_- = a_- -\sqrt{2a_+(1-a_+)v}.
\]
Discard the rectangle when \(v<0\). Otherwise \eqref{bridge-area}
gives \(b_j\geq b_-\). Use
\[
 \Delta =
 \frac{0.039\sqrt{2\cosh q_+-2+4\sin^2(\pi a_+)}}
 {\sqrt{2b_-/\sup_{[a_-,a_+]}\sin(2\pi a)}
                        -0.039e^{q_+/2}}
\]
as an upper bound on the right of \eqref{center-translation},
with an upper enclosure for the sine supremum. These bounds are used
on rectangles where \(b_-\) and the denominator are positive.
As \(|\vartheta-2\pi a_j|,|\widetilde\rho_j-\rho_j|\le\Delta\),
the projection gives
\[
 \sqrt3\tan(\pi a_--\Delta/2)\le r\le\sqrt3\tan(\pi a_++\Delta/2),
 \qquad 0\le d\le2\cosh(q_++\Delta)-2 ,
\]
provided both tangent arguments lie in \((0,\pi/2)\).

Use outward endpoints in \(10^{-32}\mathbb Z\), integer square roots,
ten alternating terms for sine and cosine, and eighteen, respectively
nine, terms for the exponential and hyperbolic cosine. On \([0,2)\)
their positive remainders are bounded by twice the first omitted term.
Set \(d_*(r)=\max\{0.023,0.36r-0.05\}\) and
\[
 \mathcal D=\{(r,d):0.15\le r\le 0.75,\ 0\le d\le d_*(r)\}
            \cup\bigl([0.75,4]\times[0,0.25]\bigr).
\]
Bisect \(I\) if \(6|I|>|J|\) and \(J\) otherwise. Stop for \(v<0\) or for
a projected box in \(\mathcal D\), testing the increasing \(d_*\) at the
left \(r\)-endpoint. This yields \(57\) exclusions and \(112\) inclusions
of depth \(\le14\).

For the spectral inequalities, the sector integrals for the
mixed pair are polynomials of degree at most two in \(t\).
In the coordinates \(x=\sqrt3 X,\ y=Y\), the flattened sector
\(\Delta=\{(X,Y):0<X<1/2,\ |Y|<X\}\) satisfies
\[
 \int_\Delta X^iY^j\,dX\,dY
 =\frac{1+(-1)^j}{2^{i+j+2}(j+1)(i+j+2)} .
\]
Use \(U_*,U_0\) before normalization in the following integral over
\(\Delta\), evaluating them at \(\sqrt3 X+iY\). The centering term
in \eqref{finite-sector-mass-energy} contributes
\[
 -\frac{4t^2}{3}
 \frac{\left[\int_\Delta\Re(U_0-U_*)\,dX\,dY\right]^2}
      {\int_\Delta |U_*|^2\,dX\,dY}
\]
to \(m\).
The normalized sector mass of \(U_0\) and its tangential component mass are
\[
       \frac{794079886600}{2564381509059}<\frac13,\qquad
       \frac{79407988660}{2564381509059}<\frac1{25}.
\]
The corresponding quantities for \(U_*\) are \(1/3\) and \(n_0<1/25\)
from \eqref{sectorconstants}. Convexity of the squared norm gives
the bounds for \(|H|^{-1}\int_S|U|^2\,dA\) and \(|H|^{-1}\int_S v^2\,dA\).

The sector integrals give \(m,n,A+D-2E>0\). Thus \(F>0\) and
\eqref{sector-square-convexity} applies. Denote the left side of
\eqref{hybrid-sector-test} by \(H(r,d,t)\). The endpoint bounds are
\[
 \begin{aligned}
 \sup_{0.15\le r\le 0.75}\max_{d\in\{0,d_*(r)\}}H(r,d,1-4r/3)&<-0.002,\\
 \sup_{0.75\le r\le4}\max_{d\in\{0,0.25\}}H(r,d,(1-r)_+)&<0.00023,\\
 \sup_{0.75\le r\le0.78}\max_{d\in\{0.10,0.25\}}H(r,d,1-r-0.10)&<-0.002 .
 \end{aligned}
\]
For the first line perform three initial bisections; bisect again
while the upper bound is \(\ge-0.002\) or positivity is undecided.
On \([l,u]\) use \([d_*(l),d_*(u)]\) and \([1-4u/3,1-4l/3]\).
There are \(394\) leaves of depth \(\le10\). Split the second line
at \(r=1\); the exterior cases in this range lie in the rectangle
of the last line. Convexity covers the intermediate \(d\). Use
\[
 \vartheta-\pi/3
 =\frac{2\sqrt3(r-1)}{r+3}
    \frac{\arctan\sqrt{v}}{\sqrt{v}},
 \qquad v=\frac{3(r-1)^2}{(r+3)^2}.
\]
The quotient is enclosed by consecutive alternating Taylor sums
of degrees \(22\) and \(23\) in \(v\). The parameter \(t\) lies in
\([0,1]\) on every tested interval.
\end{proof}

Give \(R_5\) circumradius one and let \(\alpha=\pi/5\).
\begin{theorem}\label{local-pentagon}
For every convex polygon \(P\) with at most five sides,
\(\Phi(P)\leq\Phi(R_5)\), with equality precisely when \(P\)
is similar to \(R_5\).
\end{theorem}

\begin{lemma}\label{pentagon-centroid-sector-geometry}
Let \(P\) be a convex polygon with \(\Phi(P)\geq10.35\).
Translate so that \(\int_P x\,dx=0\). In polar coordinates at the origin let \(\vartheta\)
be the angular gap between two adjacent vertices of radii \(R,S\).
With \(\alpha=\pi/5\), \(t=\tan\alpha\) and
\[
 r=\tan(\vartheta/2)/t,\qquad
 d=R/S+S/R-2,\qquad
 D(r)=\frac{4t^2r^2(1.71-tr)^2}
 {(1+t^2r^2)(1-t^2r^2+2(1.71)tr)},
\]
one has \(0<r<2.354\) and \(0\leq d\leq D(r)\).
\end{lemma}

\begin{proof}
Normalize \(|P|=1\) and fix
a support line at distance \(h\) from the origin. For the parallel
section of length \(\ell(x)\) at inward distance \(x\),
\[
 \int_0^\infty\ell(x)\,dx=1,\qquad
 \int_0^\infty x\ell(x)\,dx=h,\qquad
 \max_{|e|=1}\int_P(e\cdot y)^2\,dy\ge \frac1{12}\int_0^\infty\ell(x)^3\,dx .
\]

Here \(\min_b\int_I(y-b)^2\,dy=|I|^3/12\). For
\(g(x)=\frac3{5h}(1-2x/(5h))_+^{1/2}\),
\[
 \int_0^\infty g=1,\qquad \int_0^\infty xg=h,
 \qquad \int_0^\infty g^3=\frac{27}{125h^2}.
\]
By convexity,
\[
 \int_0^\infty(\ell^3-g^3)\,dx
 \geq \int_0^\infty 3g^2(\ell-g)\,dx
 \geq \frac{27}{25h^2}\int_0^\infty(1-2x/(5h))(\ell-g)\,dx
 =0 .
\]

The second inequality uses \(\ell-g\ge0\) on \([5h/2,\infty)\).
Thus \(\max_{|e|=1}\int_P(e\cdot y)^2\,dy\ge9/(500h^2)\), and
\eqref{centroid-radii} with Lemma~\ref{cubic-directional-integral} gives
\[
 \frac{\max(R,S)^2}{h^2}
 \le\frac{4000}{9}\left(\max_{|e|=1}\int_P(e\cdot y)^2\,dy\right)^2
 \le\frac{4000}{9}\left(\frac{35}{36\cdot10.35}\right)^2<1+1.71^2 .
\]

For \(\omega=\arctan(1.71)\), let \(\psi,\psi+\vartheta\) be the
endpoint directions from the side normal. Then
\(R=h/\cos\psi\), \(S=h/\cos(\psi+\vartheta)\), so both directions
belong to \([-\omega,\omega]\). Since \(-\log\cos\) is even and convex,
\[
 \max(S/R,R/S)\le\frac{\cos(\omega-\vartheta)}{\cos\omega}
  =\cos\vartheta+1.71\sin\vartheta
  =\frac{1-t^2r^2+2(1.71)tr}{1+t^2r^2}.
\]
Now \(d=(\max(S/R,R/S)-1)^2/\max(S/R,R/S)\le D(r)\), and
\(\vartheta\le2\omega\) gives \(r\le1.71/t<2.354\).
\end{proof}

\begin{proof}[Proof of Theorem~\ref{local-pentagon}]
By \cite{LS,EO} and Numerical Lemma~\ref{pentagon-local-signs},
every convex polygon with at most four sides satisfies
\[
 \Phi(P)\leq\pi^2<10.35<10.35319<\Phi(R_5).
\]

Hence take a pentagon with \(\Phi(P)\ge10.35\), area \(|R_5|\) and
\(\int_Px\,dx=0\). For the normalized first pair, put
\[
 \lambda=\mu_1(R_5),\quad m+n=2/5,\quad a=2\lambda/5=A+B+C+D,
 \qquad T=\frac{(a+G)^2}{4\lambda F},
\]
where \(F,G\) are given by \eqref{finite-sector-signed-expansion}.
For \(r\ge0.600\) in Lemma~\ref{pentagon-centroid-sector-geometry},
Numerical Lemma~\ref{pentagon-global-sectors} gives
\[
 B_0(r,d):=T(r,d)-\frac{2\lambda}{5}-T_r(1,0)\mathcal A(r)
 \le-0.01\bigl((r-1)^2+d\bigr).
\]

Lemma~\ref{central-exterior-comparison}(1) covers \(\min_jr_j\ge0.600\).
Otherwise, take \(V,V_0\) from Numerical Lemma~\ref{pentagon-global-sectors}
with unit component masses for \(V\) and use
\[
 U_j=V+(1-1.1r_j)\mathbf1_{(0,0.600)}(r_j)(V_0-V).
\]
The transports agree on common rays since \(V_0=V\) there. Use
\[
 \lambda_-=4.3543969<\lambda,\qquad
 T_*=\frac{(2\lambda_-/5+G)^2}{4\lambda_-F}
\]
with \(F,G\) from \eqref{finite-sector-mass-energy}.
Numerical Lemma~\ref{pentagon-global-sectors} and
Lemma~\ref{central-exterior-comparison}(2) give the strict inequality,
since \(0.02-4(0.00016)=0.01936>0\).
\end{proof}

\begin{numericallemma}[Reference spectrum and sector derivatives of the pentagon]
\label{pentagon-local-signs}
For the regular pentagon of circumradius one,
\[
 4.3543969<\mu_1(R_5)<4.35448654,\qquad \Phi(R_5)>10.35319 .
\]
For the normalized first pair the reference sector integrals satisfy
\[
 \begin{gathered}
 0.0038<B<0.00665,\quad 0.00190<C<0.00403,\quad 0.2412<K<0.2909,\\
 0.06389<n<0.06779,\qquad m=2/5-n .
 \end{gathered}
\]
The derivatives at \((r,d)=(1,0)\) of \(B_0\) in the proof of
Theorem~\ref{local-pentagon} satisfy
\[
 -0.61<(B_0)_{rr}<-0.41,\qquad -0.51<(B_0)_d<-0.42 .
\]
\end{numericallemma}

\begin{proof}
For \(N=5,\ n=220,\ \sigma=4.354440\), \eqref{mesh-disk-lower} gives
\[
 \mu_1(R_5)>\frac{\sigma}{1+\sigma[0.1893(2\sin(\pi/5))/220]^2}
 >4.3543969=:\lambda_-.
\]
For the upper bound, take \(V=\sum c_{pq}z^p\bar z^q\) in the
bisector coordinates of \eqref{certified-sector-monomials}. In order
of increasing \((p+q,p)\) on \(1\le p+q\le11\), \(p-q\equiv1\pmod5\),
the coefficient vector is
\[
 \begin{aligned}
 (&10^9,-547611997,-89776004,120600855,835307,20759637,\\
  &-67705861,28713270,-14842058,-1085148,69457269,\\
  &-14244427,9286139,6966066,-28991941,-2698633).
 \end{aligned}
\]

The component masses agree by rotation. Substitution of
\(\cos\alpha=(1+\sqrt5)/4\) in \eqref{certified-sector-monomials} gives
\[
 \frac{\int_{R_5}|\nabla\Re V|^2}{\int_{R_5}|\Re V|^2}<4.35448654=:\lambda_+,\qquad
 |R_5|\lambda_-=\tfrac52\sin(2\alpha)\lambda_->10.35319.
\]

In the disk classes \(k\equiv\pm1\pmod5\), modes \(|k|\ge4\) have
quotient \(\ge16\). In mode one the Dirichlet trace constraint gives
\((j'_{1,2})^2\ge j_{1,1}^2>14\). Indeed, for \(0\le s\le7/2\),
\[
 1-s/2+s^2/12-s^3/144+s^4/2880-s^5/86400
 \ \leq\ \frac{2J_1(2\sqrt{s})}{2\sqrt{s}}.
\]

The left polynomial has value \(52853/2764800\) at \(7/2\).
Its derivative is at most
\[
 -\tfrac12+s/6-s^2(1/48-s/720)<0\qquad(0\le s\le7/2),
\]
by splitting at \(s=3\). The radial comparison of Lemma~\ref{affine-all-N} gives the
next distinct eigenvalue in this class \(>0.49\cdot14=6.86\).

Normalize \(V\) to unit component masses. With \(\Lambda=6.86\),
\eqref{finite-trial-pair-constants}--\eqref{finite-trial-pair-errors} give
\[
 \|V-U\|_S<0.003784,\qquad \|DV-DU\|_S<0.009909 .
\]

The reverse triangle inequality bounds the errors in component mass
norms by \(0.003784\), in derivative norms by \(0.009909\), and in
\(\|u_x-v_y\|_S\) by \(\sqrt2(0.009909)\). Here
\[
 K=\frac1{|R_5|}\int_S[(u_x-v_y)^2-2u_yv_x]\,dA,\qquad
 \frac1{|R_5|}\left|\int_S u_yv_x\,dA\right|\le\sqrt{BC}.
\]
At \((1,0)\), \(F=2/5\) and differentiation gives
\begin{align*}
 T_r&=2(C-B)-\lambda(n-m),\\
 (B_0)_{rr}&=6B+2C-2\lambda m
       +\frac5\lambda[C-B-\lambda(n-m)]^2+2\sin^2\alpha\,T_r,\\
 (B_0)_d&=\frac{Bt+C/t}{\sin(2\alpha)}
       +\frac{K}{\sin^2(2\alpha)}-\lambda\frac{mt+n/t}{2\sin(2\alpha)} .
\end{align*}
Substitution with outward endpoints in \(10^{-32}\mathbb Z\) gives
\(-0.61<(B_0)_{rr}<-0.41\) and \(-0.51<(B_0)_d<-0.42\).
\end{proof}

\begin{numericallemma}[Centroid sector comparison for the pentagon]
\label{pentagon-global-sectors}

With \(D(r)\) from Lemma~\ref{pentagon-centroid-sector-geometry}, the exact pair satisfies
\[
 B_0(r,d)\leq -0.01\bigl((r-1)^2+d\bigr)
 \quad\text{for } 0.600\leq r\leq2.354,\quad 0\leq d\leq D(r).
\]
In bisector coordinates take
\[
\begin{aligned}
 V={}&z-0.543954133\,z^2\bar z-0.087987484\,\bar z^4
     +0.099144730\,z^3\bar z^2\\
    &+0.009141043\,z\bar z^5+0.016498801\,z^6
     -0.010094652\,z^4\bar z^3 .
\end{aligned}
\]
For its coefficients \(c_{pq}\), define
\[
 V_0=z\sum_{p,q} c_{pq}(-1)^{(p-q-1)/5}
                \left(\frac{\operatorname{Re}z}{\cos\alpha}
                \right)^{p+q-1}.
\]

Since \((\rho e^{\pm i\alpha})^p(\rho e^{\mp i\alpha})^q
=\rho^{p+q}e^{\pm i\alpha}(-1)^{(p-q-1)/5}\),
\(V_0(\rho e^{\pm i\alpha})=V(\rho e^{\pm i\alpha})\).
For \(s=(1-1.1r)\mathbf1_{(0,0.600)}(r)\), the raw sector trace
mass is \(<0.5\). With \(T_*\) as in the proof of
Theorem~\ref{local-pentagon} and \(T_{*,r}(1,0)\) taken at fixed
\(s=0\), every \(0\le d\le D(r)\) satisfies
\[
 T_*(r,d)-2\lambda_-/5-
       T_{*,r}(1,0)\frac{\vartheta(r)-2\alpha}{\sin(2\alpha)}
 < \begin{cases}
     -0.02,&0<r<0.600,\\
     0.00016,&0.600\leq r\leq2.354.
    \end{cases}
\]
\end{numericallemma}

\begin{proof}
Reflection gives \(F=mx+ny\), where \(x,y>0\) are affine in \(d\).
For the exact pair \(m,n>0\). For the mixed pair the corrected
sector integrals are quadratic in \(s\). Their Bernstein coefficients
\((a,a+b/2,a+b+c)\) for \(a+bs+cs^2\), evaluated by
\eqref{certified-sector-monomials}, give
\[
 \min_{0\le s\le1}m(s)>0.3125,\qquad \min_{0\le s\le1}n(s)>0.055 .
\]

Thus \(F>0\), and \eqref{sector-square-convexity} reduces the comparison
to \(d=0,D(r)\). For the exact pair use Numerical
Lemma~\ref{pentagon-local-signs}. For \(d=0\), use \eqref{finite-fan-angular-quotient} and
\eqref{uncorrected-exact-angular}, keeping the square term in
\eqref{finite-fan-angular-quotient} separate. At the other endpoint use
\eqref{finite-fan-radial-quotient} with
\[
 G(r,0)=a+(r^2-1)(C-B/r^2),\qquad
 G_1=G_d(r,0)+\tfrac d2G_{dd}(r,0),\qquad F_1=F_d(r,0).
\]
Use
\[
 \vartheta(r)-2\alpha=2\arctan\frac{t(r-1)}{1+t^2r},\qquad
 \left|\frac{1-\arctan(v)/v}{v^2}-
       \sum_{j=0}^{18}\frac{(-1)^jv^{2j}}{2j+3}\right|\le\frac{|v|^{38}}{41}
\]
at \(v=t(r-1)/(1+t^2r)\). For the mixed trial, the monomial formula
\eqref{certified-sector-monomials} gives sector integrals quadratic in \(s\),
\[
 m=m_{\rm raw}-\frac{s^2}{5|S|^2}\left(\int_S(\Re V_0-\Re V)\right)^2,\qquad
 \|V+s(V_0-V)\|_S^2\le(1-s)\|V\|_S^2+s\|V_0\|_S^2<0.5 .
\]
On the small sectors \((\Re V_0)_y=0\), hence
\(B=(1-s)^2B(V)=(1.1r)^2 B(V)\).
At the left endpoint use \(rx,y/r,q\), \(B/r^2=(1.1)^2B(V)\), and
\[
 \frac d{r^2}\in\left\{0, \frac{4t^2(1.71-tr)^2}{(1+t^2r^2)(1-t^2r^2+2(1.71)tr)}\right\}.
\]
The bounds extend to \(r=0\). Bisection until both endpoint tests hold
gives the pairs (terminal interval count, maximal depth)
\[
 \begin{aligned}
 [0.600,2.354]_{\rm exact}&:\ (18,7),\qquad
 [0,0.600]_{\rm mixed}:\ (260,10),\qquad
 [0.600,2.354]_{\rm mixed}:\ (5885,19).
 \end{aligned}
\]
\end{proof}

\subsection{Seven, eight and nine sides}
Take \(R_N\) of circumradius one. Use \(|P|=1\) in localization and \(|P|=|R_N|\)
in the Rayleigh comparison. For \(\lambda=\mu_1(R_N)\) and a regular lower endpoint
\(\lambda_*<\lambda\), define
\[
 H_{\nu,j}=\frac{(2\nu/N+G_j)^2}{4\nu F_j}-\frac{2\nu}N
 -\left.\partial_r\frac{(2\nu/N+G)^2}{4\nu F}\right|_{(1,0)}
   \mathcal A(r_j),\qquad \nu\in\{\lambda,\lambda_*\}.
\]
The derivative for a polynomial trial is taken on its uncorrected reference pair.
Write \(\alpha=\pi/N\), \(t=\tan\alpha\), \(a=2\lambda/N\),
\(T=(a+G)^2/(4\lambda F)\), and \(\|w\|_S^2=|R_N|^{-1}\int_S|w|^2\).
\begin{numericallemma}[Regular spectral enclosures]\label{nonagon-local-signs}
For circumradius one,
\[
 \begin{aligned}
 3.861641&<\mu_1(R_7)<3.86174338,&
 3.74694255&<\mu_1(R_8)<3.74695936,\\
 3.66977986&<\mu_1(R_9)<3.66984730,&
 3.4768738&<\mu_1(R_{16})<3.4768795.
 \end{aligned}
\]
In particular
\(\Phi(R_9)>(9/2)\sin(2\pi/9)(3.66977986)>10.61500\).
\end{numericallemma}

\begin{proof}
By \eqref{mesh-disk-lower} and Lemma~\ref{lem:standalone-ritz} we obtain the stated
lower and upper endpoints, respectively. The sector trials use
\[
 (N,d,s)=(7,19,9),(8,39,20),(9,25,10),(9,55,20).
\]
Order \(f_j=z^{p_j}\bar z^{q_j}\) by \((p_j+q_j,p_j)\), with
\(1\le p_j+q_j\le d\) and \(p_j-q_j\equiv1\pmod N\). In the full polygon put
\[
 M_{ij}=\frac1{|R_N|}\int_{R_N}\overline{f_i}f_j,\qquad
 K_{ij}=\frac1{|R_N|}\int_{R_N}\overline{\nabla f_i}\cdot\nabla f_j.
\]
Reflection makes both matrices real. From \(c=e_1\), apply four steps
\[
 y=(K-\lambda_-M)^{-1}Mc,\qquad c\gets y/y_1,\qquad
 (\lambda_-)_{N=7,8,9}=(3.861641,3.74694255,3.66977986).
\]
Here \(K-\lambda_-M>0\) and \(y_1>0\) in each step. After rescaling
\(z\mapsto z/\cos(\pi/N)\) in interval elimination at precision \(2^{-2048}\),
form the trial in physical coordinates:
\[
 10^{-s}V_N(z)=10^{-s}\sum_j\lfloor10^sc_j\rfloor f_j(z).
\]
Use the largest root of \(U_{N-1}\) for \(\cos(\pi/N)\) as in
Lemma~\ref{lem:standalone-ritz}. Integration gives
\[
 \mathcal R(V_N)<(3.8617433681,3.7469593531,3.669847291,3.669794)
\]
in the order \((N,d,s)\) above.
\end{proof}

\begin{numericallemma}[Sector integrals for seven and eight sides]
\label{heptagon-local-signs}\label{octagon-global-reference}
For the first equivariant pair normalized by \(\int_{R_N}U_i^2=|R_N|\),
\(m=2/N-n\). On a bisector sector, the following enclosures hold:
\[
\begin{aligned}
N=7:\quad &0.00454299<B<0.00671092,&0.00270267&<C<0.00442301,\\
 &0.13954961<K<0.17769154,&0.02063726&<n<0.02207740;\\
N=8:\quad &0.00409014<B<0.00479969,&0.00261495&<C<0.00318797,\\
 &0.12585176<K<0.13270359,&0.01372971&<n<0.01414335.
\end{aligned}
\]
At \((r,d)=(1,0)\) for \(B_0=T-a-T_r(1,0)\mathcal A(r)\),
\[
\begin{aligned}
 N=7:\quad &-0.103<(B_0)_{rr}<-0.016,&-0.191&<(B_0)_d<-0.112,\\
 N=8:\quad &-0.032561<(B_0)_{rr}<-0.007621,&-0.092427&<(B_0)_d<-0.034740.
\end{aligned}
\]
\end{numericallemma}
\begin{proof}
Write \(V\) for the normalized degree-\(19\) or degree-\(39\) polynomial
of Numerical Lemma~\ref{nonagon-local-signs}, and
\(\langle X,Y\rangle_S=|R_N|^{-1}\int_S X\cdot Y\).
The second disk eigenvalue in the rotation class \(1\bmod N\)
exceeds \(\min\{(N-1)^2,j_{1,1}^2\}>14\) by the radial interlacing
used in Numerical Lemma~\ref{pentagon-local-signs}. From the matrix
in Lemma~\ref{affine-all-N} the signed polygon gaps exceed
\(\Lambda_7=0.62(14)=8.68\), \(\Lambda_8=0.66(14)=9.24\).
Apply Lemma~\ref{equivariant-sector-projection} with the endpoints of
Numerical Lemma~\ref{nonagon-local-signs} and \(\Lambda=\Lambda_N\).
In particular,
\(
 (\delta_7,\epsilon_7)<(0.002464,0.007260),\
 (\delta_8,\epsilon_8)<(0.000876,0.002663).
\)
For \(U=(u,v)\), \(V=(f,g)\), use \eqref{certified-sector-monomials} and
\eqref{finite-trial-pair-errors}:
\[
\begin{gathered}
 |\sqrt B-\|f_y\|_S|,\ |\sqrt C-\|g_x\|_S|\le\epsilon,
 \qquad |\sqrt n-\|g\|_S|\le\delta,\\
 \big|\|u_x-v_y\|_S-\|f_x-g_y\|_S\big|\le\sqrt2\,\epsilon,
 \qquad K=\|u_x-v_y\|_S^2-2\langle u_y,v_x\rangle_S.
\end{gathered}
\]
Cauchy--Schwarz gives
\[
 \begin{aligned}
 |\langle u_y,v_x\rangle_S|&\le\sqrt{BC} && (N=7),\\
 |\langle u_y,v_x\rangle_S-\langle f_y,g_x\rangle_S|
 &\le\epsilon(\|f_y\|_S+\|g_x\|_S)+\epsilon^2/2 && (N=8).
 \end{aligned}
\]
Substitute the unrounded \(\delta,\epsilon\). At \((r,d)=(1,0)\),
\(T_r=2(C-B)-\lambda(n-m)\), and
\[
\begin{aligned}
 (B_0)_{rr}&=6B+2C-2\lambda m+\frac N\lambda[C-B-\lambda(n-m)]^2
                  +2\sin^2\alpha\,T_r,\\
 (B_0)_d&=\frac{Bt+C/t}{\sin(2\alpha)}+\frac K{\sin^2(2\alpha)}
                   -\lambda\frac{mt+n/t}{2\sin(2\alpha)}.
\end{aligned}
\]
For the second derivative collect coefficients of \(B,C,(B-C)^2\)
using \(m=2/N-n\) before substituting the intervals.
\end{proof}

For \(7\le N\le15\), use \((r,d)\) from \eqref{common-fan-coordinates}.
The localized sets are
\[
 I_k=[k/100,(k+1)/100],\qquad
 \mathcal G_N=\bigcup_k I_k\times[0,B_{N,k}].
\]
The index ranges are specified in the localization statements below.
For \(N=7,8,9\) take \(x=(2k+1)/200-1\) and
\[
 10^4 B_{N,k}=
 \begin{cases}
 \min\{787,\max\{1,\lceil186x^4-183x^3-1604x^2+58x+791\rceil\}\},
       &N=7,\\
 \min\{370,\max\{9,\lceil43x^4-227x^3-909x^2+29x+372\rceil\}\},
       &N=8,\\
 \min\{207,\max\{9,\lceil-191x^3-558x^2+15x+209\rceil\}\},
       &N=9 .
 \end{cases}
\]
\begin{numericallemma}[Centroid localization for seven to nine sides]
\label{heptagon-global-geometry}
For \(N=(7,8,9)\), respectively, put
\[
 \begin{gathered}
 (b_N)=(10.56703,10.59789,10.615),\qquad (D_N)=(0.0787,0.0370,0.0207),\\
 (J_N)=([0.25,1.72],[0.30,1.65],[0.32,1.59]).
 \end{gathered}
\]
If \(P\) is a convex \(N\)-gon with \(|P|=1\) and \(\Phi(P)\ge b_N\),
then \((r,d)\in\mathcal G_N\cap(J_N\times[0,D_N])\) on every sector,
where \(k\) runs over the indices with \(I_k\subseteq J_N\).
\end{numericallemma}

\begin{proof}
\emph{Step 1. The conformal ratio.}
Apply \eqref{schur-recursion} with the balancing profile \(\phi\) of
Section~\ref{sec:conformal-rayleigh}, \(T_0=b_N\). For \(N=7\) the endpoints are
\[
 0.9585,\ 0.964,\ 0.969,\ 0.974,\ 0.977,\ 0.979,\ 0.9795,\ 0.9798,\ 0.9799,\
 0.98,\ 0.98009,\ 0.98015,\ 0.98017 .
\]
For \(N=8,9\), respectively, the consecutive endpoints after \(0.9585\) are
\[
\begin{aligned}
8:\ &0.964,0.969,0.974,0.977,0.979,0.98017,0.983,0.985,0.986,0.9865,\\
&0.9869,0.987,0.9871,0.9872,0.9873,0.9874,0.9875,0.9876,0.9877;\\
9:\ &0.964,0.969,0.974,0.977,0.979,0.983,0.986,0.988,0.9895,0.9903,\\
&0.9907,0.991,0.99125,0.9914,0.9915,0.9916,0.99166,0.9917,0.99172,0.99174 .
\end{aligned}
\]
The lower bound in \eqref{improved-schur} exceeds its right side on every interval.
By Numerical Lemma~\ref{rationalconformal}, \(\beta>0.9585\). Proposition
\ref{balanced-centroid-offset}, Lemma~\ref{cubic-directional-integral}, and the
series for \(F\) yield the bounds at \(z=f(0)\):
\[
 \begin{gathered}
 \beta>\beta_-,\qquad \beta_N<\beta_+,\qquad
 \left|\int_P y\,dy-z\right|<\Delta,\qquad \log(\beta_N/\beta)<\log(\beta_+/\beta_-)<\varepsilon,\\
 (\beta_-)=(0.98017,0.9877,0.99174),\quad(\beta_+)=(0.985304,0.99026,0.993209),\\
 (\Delta)=(0.02076,0.01278,0.00855),\qquad
 (\varepsilon)=(0.005225,0.00259,0.001485).
 \end{gathered}
\]
Indeed \(\beta^{-1/2}-\beta^{1/2}\) and \(\beta^{-1}-1\) decrease.
\emph{Step 2. Localization at the conformal pole.}
Let \(2\pi a,b\) be the sector angle and area at \(z\), and set
\(q=|\log(|p'-z|/|p-z|)|\) for its consecutive vertices \(p,p'\). Put
\[
 Q_N(a)=aF(a)+(1-a)F\bigl((1-a)/(N-1)\bigr)-F(1/N),\qquad
 \delta=\frac1\pi\arctan\frac{\tanh(q/2)}{\tan(\pi a)}.
\]
In \eqref{fan-deficit}, apply Jensen to the other \(N-1\) angles,
\eqref{sector-log-sum-bound} to the area term, and
\eqref{fan-radius-angle-identity} to the radius ratio. This gives
\begin{equation}\label{finite-fan-budget}
 \begin{aligned}
 Q_N(a)+a\log\frac ab+(1-a)\log\frac{1-a}{1-b}
 &+a\sum_{j\ge0}
 \log\frac{1-\delta^2/(j+\tfrac12+a)^2}{1-\delta^2/(j+\tfrac12-a)^2}\\
 &\le\log\frac{\beta_N}{\beta}.
 \end{aligned}
\end{equation}
Retain odd powers \(3,\ldots,15\) in \(Q_N\), \(80\) terms per zeta
coefficient, and the positive radial sum
\[
 a\sum_{j=0}^{8}\sum_{\ell=1}^{7}
 \frac{\delta^{2\ell}}{\ell}
 \big[(j+\tfrac12-a)^{-2\ell}-(j+\tfrac12+a)^{-2\ell}\big].
\]
Convexity of \(Q_N\) and the truncated bounds give
\[
 (a_{\min})=(0.04,0.037,0.034),\quad (a_{\max})=(0.215,0.188,0.17),\quad
 (q_{\max})=(0.26,0.18,0.14).
\]
On each of the \(256\) equal subintervals of \([a_{\min},a_{\max}]\) the
angular and radial lower bounds at \(q=q_{\max}\) sum to more than \(\varepsilon\).
For \([a_-,a_+]\times[q_-,q_+]\) the lower radial sum uses \(a_-\) outside
and inside the brackets and \(\delta(a_+,q_-)\). Discard the box if the sum
of angular and radial lower bounds exceeds \(\varepsilon\).

\emph{Step 3. Translation to the centroid.}
For \(0<x<a_-<1/2\), expansion gives
\[
 \begin{aligned}
 -a_-\log(1-x/a_-)&-(1-a_-)\log\bigl(1+x/(1-a_-)\bigr)\\
 &\ge P_{13}(x):=\sum_{\ell=2}^{13}\frac{x^\ell}{\ell}
       \bigl[a_-^{1-\ell}+(-1)^\ell(1-a_-)^{1-\ell}\bigr].
 \end{aligned}
\]
The omitted coefficients are positive. For \(b\le a_--x<a_-\le a\),
\[
 \mathcal L(a,b)\ge\mathcal L(a_-,a_--x)\ge P_{13}(x)
\]
by monotonicity of the two-variable logarithmic expression in
\eqref{sector-log-function}. Bisect \([0,a_-]\) eighteen times, retaining the
left half when \(P_{13}\) at the midpoint exceeds the remaining allowance
and the right half otherwise. For the final upper endpoint \(x\), either
\(x=a_-\) and the rectangle is bisected, or \(b>b_-:=a_--x\).
Take \(e\) to be an outward rational upper bound for
\[
 \frac{\Delta\sqrt{2\cosh q_+-2\cos(2\pi a_+)}}
      {\sqrt{2b_-/\sin(2\pi a_+)}-\Delta e^{q_+/2}}
\]
when the denominator is positive. Here \(\sin(2\pi a)\) increases on
the relevant intervals. By \eqref{center-translation} with \(\Delta\) equal to the corresponding
bound on the pole displacement, the change
in (signed log-ratio, angle) has Euclidean norm \(\le e\). Therefore
\[
 \begin{gathered}
 s=\operatorname{dist}\big(2\arctan(\tan(\pi/N)I_k),[2\pi a_-,2\pi a_+]\big),\\
 r\in I_k\ \Longrightarrow\
 s\le e,\quad|\rho|\le q_++\sqrt{e^2-s^2},\quad
 d\le 2\cosh\big(q_++\sqrt{e^2-s^2}\big)-2 .
 \end{gathered}
\]
Bisect \([a_{\min},a_{\max}]\times[0,q_{\max}]\) along the longer
coordinate in \((3a,q)\), bisecting the \(q\)-interval at equality.
For this test integrate the positive radial polynomial by Horner's
rule, using outward decimal intervals of at least 32 places.
For \(0\le x<2\) bound \(\sin x\) by the series through degree 25
with error \(2^{27}/27!\); for \(0\le x<1\) use the degree-17
exponential sum with error \(2/18!\). For the arctangent use
\(v=x/(1+\sqrt{1+x^2})\), the degree-43 sum for \(2\arctan v\),
and error \(2(2/3)^{45}/45\) when \(v<2/3\). Each leaf either exceeds
the budget or meets only bands with \(d\le B_{N,k}\). Either band applies
at a common endpoint. For \(N=(7,8,9)\) the counts are
\[
 (\#\mathrm{excluded})=(934,554,309),\quad
 (\#\mathrm{covered})=(834,403,203),\quad
 (\mathrm{depth})=(18,17,15).
\]
\end{proof}

\paragraph{Exterior sector trials.}
Index the reference polynomial by \(0\) and the zero-trace corrections by
\(1,\ldots,h\), with \((h_7,h_8,h_9)=(29,57,57)\). Normalize the basis by the
reference component masses. The polarized forms \(F[c]=c^{\mathsf T}Fc\),
\(G[c]=c^{\mathsf T}Gc\) include the correction in
\eqref{finite-sector-mass-energy}. At the midpoint of each rectangle,
start with \(c=e_0\) and apply eight steps:
\[
 t=\frac{2\lambda_*/N+G[c]}{2F[c]},\qquad
 (G-tF)_{1:h,1:h}\,c_{1:h}=-(G-tF)_{1:h,0}.
\]
In \eqref{finite-sector-signed-expansion} apply
\(x\mapsto2^{-70}\lfloor2^{70}x+7/23\rfloor\) to each sector integral and
each geometric factor.
Use one geometric factor for each of
\(x,y,2z,1,x^2,y^2,z^2,2zy,2zx\), combining
the diagonal energy terms as \(A+D\) and \(A+D-2E\) before rounding.
The last two factors multiply, respectively, the sector integrals of
\(v_x(u_x-v_y)\) and \(-u_y(u_x-v_y)\); both vanish in the
reflection-symmetric basis for \(N=7\).
Round each resulting entry and \(2\lambda_*/N\) down to \(2^{-44}\mathbb Z\).
Round \(t\) down to \(10^{-10}\mathbb Z\).
To solve, form the integer array
\[
 H=2^{44}G-\lfloor 2^{44}tF\rfloor .
\]
Retain column \(0\) as the right-hand column. Eliminate indices
\(i=1,\ldots,h\) in the displayed basis order, with
\[
 H_{jk}\gets H_{jk}-\left\lfloor\frac{H_{ji}H_{ik}}{H_{ii}}\right\rfloor
 \quad(j>i,\ k=0\ \hbox{or}\ k>i).
\]
Back substitution, in decreasing \(i\), assigns
\[
 c_i=10^{-9}\left\lfloor
   -\frac{10^9}{H_{ii}}\left(H_{i0}+\sum_{j>i}H_{ij}c_j\right)
                    \right\rfloor .
\]
All pivots are nonzero. Compute \(F,G\) and their first derivatives on
the rectangle for the resulting rational trial by \eqref{certified-sector-monomials}.

For \(N=7\), start from \(I_k\times[0,B_{7,k}]\),
\(25\le k<85\), in variables \((r,d)\). For \(N=8,9\), use
\(I_k\times[0,10^{-5}(1+\lfloor10^5\sqrt{B_{N,k}}\rfloor)]\)
in \((r,\sqrt d)\); the indices are \(30\le k<89\) and
\(32\le k<88\) or \(110\le k<159\), respectively.
Test the two \(d\)-endpoints by convexity for \(N=7\). For \(N=8,9\),
bound both the direct residual and its first-order expansion about the
midpoint by interval arithmetic as in Numerical
Lemma~\ref{octagon-skew-spectral}. Accept a closed rectangle when
\(F>0\) and its residual is strictly below the asserted exterior bound.
Otherwise bisect the side with greater length relative to the starting
rectangle (the \(r\)-side at equality). The terminal counts,
maximum depths and strict upper bounds are
\[
 \begin{aligned}
 N=7&:\ (60,0),& H_{\lambda_*}&<-0.00042,\\
 N=8&:\ (534,8),& H_{\lambda_*}&<-0.00014,\\
 N=9&:\ (1037,10),& H_{\lambda_*}&<-0.00014 .
 \end{aligned}
\]
By Numerical Lemma~\ref{heptagon-global-geometry}, the initial rectangles cover
the exterior sectors.

\begin{numericallemma}[Heptagon fan comparisons]\label{heptagon-global-spectral}
For \(N=7\), on \([0.85,1.76]\times[0,0.085]\),
\[
 H_\lambda(r,d)
 \leq -0.001\bigl((r-1)^2+d\bigr).
\]
The uncorrected \(V\) satisfies \(H_{\lambda_*}<0.00006\) there.
On \(\mathcal G_7\cap\{r<0.85\}\) the corrected polynomial trial below
has \(F>0\) and \(H_{\lambda_*}<-0.00042\).
\end{numericallemma}

\begin{proof}
For the exact pair, use \eqref{finite-fan-angular-quotient}
with \(f_0=2/7\) and \eqref{uncorrected-exact-angular}.
The alternating arctangent series bounds the factored angular quotient.
For \eqref{finite-fan-radial-quotient} the increments are
\[
 \begin{aligned}
 F_1&=\frac{mt+n/t}{2\sin\vartheta},\\
 G_1&=\frac{Bt/r+Cr/t}{\sin\vartheta}
   +\frac{K+\tfrac d4(Bt^2+C/t^2+K)}{\sin^2\vartheta}.
 \end{aligned}
\]
On the \(1820\) intervals of length \(1/2000\), the integral bounds of
Numerical Lemma~\ref{heptagon-local-signs} give
\[
 \max\left\{\frac{H_\lambda(r,0)}{(r-1)^2},
       \frac{T(r,0.085)-T(r,0)}{0.085}\right\}\le-0.001.
\]
By \eqref{sector-square-convexity}, the radial divided difference
is bounded by its value at \(d=0.085\), and the fixed polynomial trials
need only their two \(d\)-endpoints. In \(X=x/\cos\alpha,\ Y=y/\sin\alpha\)
set \(V=10^{-9}V_7\) before mass normalization. Its gauge is
\[
 V_0(z)=z\sum_{p,q}\frac{c_{pq}}{10^9}
              (-1)^{(p-q-1)/7}X^{p+q-1}.
\]
Since \(z=Xe^{\pm i\alpha}\) on \(Y=\pm X\), \(V\) and \(V_0\)
agree there. After \(V_0-V\), the zero-trace basis contains, in
lexicographic order \(i,j\ge0,\ i+j\le6\),
\[
 X^iY^j(Y^2-X^2)
 \begin{cases}
  1,&j\ \text{even},\\
  i_{\mathbb C}\sin\alpha,&j\ \text{odd},
 \end{cases}
\]
where \(i_{\mathbb C}^2=-1\). For \(25\le k<85\), the midpoint correction above gives
\[
 \max_{d\in\{0,B_{7,k}\}}\sup_{r\in I_k}H_{\lambda_*}(r,d)<-0.00042 .
\]
For the uncorrected pair use \(d\in\{0,0.085\}\), \(9100\) intervals \(I\)
of width \(10^{-4}\), and \(H(I)\subset H(m)+(I-m)H'(I)\), \(m=\operatorname{mid}I\).
\end{proof}

\begin{numericallemma}[Octagon fan comparisons]\label{octagon-skew-spectral}
For \(N=8\), on \([0.89,1.65]\times[0,0.04]\),
\[
 H_\lambda(r,d)
 \leq -0.0004\bigl((r-1)^2+d\bigr).
\]
On this rectangle \(H_{\lambda_*}<0.000015\) for the uncorrected \(V_8\).
On \(\mathcal G_8\cap\{r<0.89\}\) zero-trace corrections give \(F>0\) and
\(H_{\lambda_*}<-0.00014\).
\end{numericallemma}

\begin{proof}
Use \eqref{finite-fan-angular-quotient} and
\eqref{uncorrected-exact-angular} with \(f_0=1/4\), \(t=\tan(\pi/8)\).
Their normalized form is
\[
 \begin{aligned}
 \frac{H_\lambda(r,0)}{(r-1)^2}
 ={}&\frac{1}{(1-4n)/r+4nr}
 \left[C+\frac{(2r+1)B}{r^2}
 -2(C-B)\left(4n-\frac{1-4n}{r}\right)\right.\\
 &\hspace{15mm}\left.{}+\frac{(r+1)^2(C-B/r^2)^2}{4a}
       +4an\left(8n-1-\frac{2(1-4n)}r\right)\right]\\
 &+\big[a(8n-1)-2(C-B)\big]\frac{\mathcal A(r)-(r-1)}{(r-1)^2}.
 \end{aligned}
\]
On each interval of width \(1/2000\), both this quotient and
\((T(r,0.04)-T(r,0))/0.04\) are \(\le-0.0004\).
Convexity covers \(0\le d\le0.04\). For the uncorrected polynomial
test \(d\in\{0,0.04\}\) on intervals of width \(1/10000\).
For the corrected trial use the preceding gauge with \(N=8\) and
reference integer scale \(10^{-20}\). In the order \(V_0-V\), then
lexicographic \((i,j)\) with \(i,j\ge0,\ i+j\le6\), take both
\[
 Z_{ij}=X^iY^j(Y^2-X^2),\qquad \mathrm i\sin\alpha\,Z_{ij}.
\]
For these unrestricted parities take the signed skew, with \(z\ge0\)
by reflection.
Put \(b=\sqrt d\). Then
\[
 x=\frac1r+\frac{b^2(r^{-1}+t^2r)}4,\qquad
 y=r+\frac{b^2(r^{-1}+t^2r)}{4t^2},\qquad
 z=\frac{b\sqrt{b^2+4}(r^{-1}+t^2r)}{4t}.
\]
These functions have bounded derivatives down to \(b=0\).
On each of the \(59\) exterior bands use the midpoint correction and the
rounded upper endpoint for \(\sqrt{B_{8,k}}\) specified above.
On \(K=I\times J\) with midpoint \((r_0,b_0)\),
\[
 H(K)\subset H(r_0,b_0)+(I-r_0)H_r(K)+(J-b_0)H_b(K).
\]
Intersecting with direct interval evaluation gives \(\inf_K F>0\),
\(\sup_K H<-0.00014\), with square roots enclosed by integer square roots.
\end{proof}

\begin{numericallemma}[Corrected comparison in the central nonagon range]
\label{nonagon-exact-central}
Let \(U\) be the first equivariant pair with \(\int_{R_9}U_i^2\,dA=|R_9|\),
\(i=1,2\). In \(X=x/\cos\alpha,\ Y=y/\sin\alpha\), set
\[
 W_0=(X^2-Y^2)\bigl(0.005(X-4/5),\ 0.209Y\bigr),\quad
 W_1=(X^2-Y^2)\bigl(-0.105Y,\ 2.54(X-4/5)\bigr).
\]
These corrections vanish on \(Y=\pm X\) and
\[
 \int_S Y(X^2-Y^2)\,dA=\int_S(X-4/5)(X^2-Y^2)\,dA=0
 \quad\Longrightarrow\quad \int_S W_0\,dA=\int_S W_1\,dA=0 .
\]
Put \(\xi=r-1,\ d=2\cosh\rho-2\),
\(s=\operatorname{sgn}(\rho)\sqrt d\), where \(\rho\) is the
logarithm of the ratio of the second radius to the first.
For \(V=U+\xi W_0+sW_1\), let \(F,G\) be the sector mass
and energy divided by \(w^2,w\), respectively.
For this trial put \(B_0=H_\lambda\). Then
\[
 0.88\le r\le1.11,\quad 0\le d\le0.021
 \quad\Longrightarrow\quad B_0(r,d)\le -0.001(r-1)^2-0.05d .
\]
\end{numericallemma}

\begin{proof}
The radial transfer factor exceeds \(0.696\). In the disk classes
\(\pm1\bmod9\) the angular term is at least \(8^2\) for \(|m|\ge8\).
For the second radial mode of order one put \(u=z^2/4\). Then
\[
 2J_1'(z)=\sum_{k=0}^{\infty}\frac{(-1)^k(2k+1)u^k}{k!(k+1)!}.
\]
Its \(u\)-derivative on \([0,1]\) is \(\le-3/2+5u/6<0\).
On \([1,2809/400]\) the degree-ten partial sum majorizes the series,
and each of its Bernstein coefficients is \(<-1/50\). Hence the
next disk eigenvalue exceeds \(5.3^2\), giving a signed gap on \(R_9\)
above \(0.696(5.3)^2>19.55\).

For the degree-\(55\) polynomial in the ordered basis of Numerical
Lemma~\ref{nonagon-local-signs}, the quotient is \(<3.669794\).
After normalization, use Lemma~\ref{equivariant-sector-projection}
with \((\lambda_-,\lambda_+,\Lambda)=(3.66977986,3.669794,19.55)\).
Then
\[
 \begin{aligned}
 \|U-V\|_S&<0.000446,&
 \frac1{|R_9|}\left|\int_S(U-V)\cdot W_i\,dA\right|&\le0.000446\|W_i\|_S,\\
 \|D(U-V)\|_S&<0.001969,&
 \frac1{|R_9|}\left|\int_S D(U-V):DW_i\,dA\right|&\le0.001969\|DW_i\|_S .
 \end{aligned}
\]
For \(K(U)\) use the reverse triangle and Cauchy--Schwarz
bounds in the proof of Numerical Lemma~\ref{octagon-global-reference},
with the degree-\(55\) nonagon trial and \(\epsilon=0.001969\).
Write \(P=U+\xi W_0,\ M=L^{\mathsf T}L/w=\left(\begin{smallmatrix}x&z\\ z&y\end{smallmatrix}\right)\).
With \(t=\tan(\pi/9)\),
\[
 x=r^{-1}+\frac{d(1+t^2r^2)}{4r},\quad
 y=r+\frac{d(1+t^2r^2)}{4t^2r},\quad
 z=sz_1,\quad z_1=\frac{(1+t^2r^2)\sqrt{1+d/4}}{2tr}.
\]
Since \(\int_S(P-U)\,dA=\int_S W_1\,dA=0\), reflection symmetry gives
\begin{align*}
 |R_9|F&=\int_S\Bigl\{x\bigl(P_1^2+d(W_1)_1^2\bigr)
                 +y\bigl(P_2^2+d(W_1)_2^2\bigr)\\
 &\hspace{7em}+2dz_1\bigl[P_1(W_1)_2+P_2(W_1)_1\bigr]\Bigr\}\,dA,\\
 G&=\mathcal E(P)+d\mathcal E(W_1)\\
 &\quad+\frac{2dz_1}{|R_9|}\int_S DP:
       \bigl[J(DW_1)\operatorname{diag}(y,x)
             -\operatorname{diag}(x,y)(DW_1)J\bigr]\,dA,
 \qquad J=\begin{pmatrix}0&1\\1&0\end{pmatrix}.
\end{align*}
Here differentiation is with respect to the sector coordinates, and
\[
 \begin{aligned}
 |R_9|\mathcal E(Y)
 &=\int_S\Bigl\{(\partial_1Y_1)^2+(\partial_2Y_2)^2
       +x^2(\partial_2Y_1)^2+y^2(\partial_1Y_2)^2\\
 &\hspace{5em}+z^2\bigl[(\partial_1Y_1-\partial_2Y_2)^2
                    -2(\partial_2Y_1)(\partial_1Y_2)\bigr]\Bigr\}\,dA .
 \end{aligned}
\]

At \(d=0\), write \(f_1(r)=(F-f_0)/\xi\) and
\(g_1(r)=(G-a)/\xi\), extending at \(\xi=0\) by the
stationary eigenvalue equation. In particular,
\[
 \begin{aligned}
 f_1(r)&=n-\frac mr+\frac{1}{|R_9|}\int_S
 \left[\frac{2U_1(W_0)_1+\xi(W_0)_1^2}{r}
       +r\big(2U_2(W_0)_2+\xi(W_0)_2^2\big)\right]\,dA,\\
 g_1(r)&=\frac1{|R_9|}\int_S\Bigl[2\lambda U\cdot W_0+\xi|DW_0|^2\\
 &\hspace{5em}+(r+1)\bigl((\partial_1P_2)^2-r^{-2}(\partial_2P_1)^2\bigr)\Bigr]\,dA .
 \end{aligned}
\]
Use \eqref{finite-fan-angular-quotient}. Substituting \(P=U+\xi W_0\) and
expanding the squared derivatives cancels \(\xi\) in \(h_2\).
For the radial increment use \eqref{finite-fan-radial-quotient} with
\(\nu=\lambda\). Using the continuous extensions on the axes, outward bisection
(\(20\) leaves, depth \(\le6\)) gives
\(\sup B_0(r,0)/\xi^2\le-0.001\),
\(\sup(B_0(r,d)-B_0(r,0))/d\le-0.05\).
Multiplying by \(\xi^2,d\) and adding proves the claim.
\end{proof}

\begin{numericallemma}[Nonagon polynomial fan comparisons]
\label{nonagon-polynomial-fan}
Let \(V_9\) be the degree-\(25\) pair of Numerical Lemma~\ref{nonagon-local-signs}.
On \(\mathcal G_9\) zero-trace corrections give \(F>0\) and
\[
 H_{\lambda_*}(r,d)
 <\begin{cases}
   0.000016,&0.88\leq r\leq1.10,\\
  -0.00014,&r<0.88\ \text{or}\ r>1.10 .
 \end{cases}
\]
\end{numericallemma}

\begin{proof}
For the exterior use the octagon basis at \(N=9\) in the midpoint construction.
For \(0.88\le r\le1.10\) add the following to \(10^{-10}V_9\) before normalization,
with \(b=\sqrt d\):
\[
 \begin{aligned}
 (X^2-Y^2)\big[&(r-1)\big(0.001275(X-4/5),\, 0.1558095Y\sin\alpha\big)\\
                 &+b\big(-0.026775Y,\,1.89357(X-4/5)\sin\alpha\big)\big].
 \end{aligned}
\]
Since \(\int_S(X^2-Y^2)(X-4/5)=\int_S(X^2-Y^2)Y=0\), the correction has
zero integral. For the resulting \(\mathcal E=H_{\lambda_*}\), reflection gives
\(\mathcal E_b(r,0)=0\), and rational differentiation gives
\[
 \partial_b^2\mathcal E(r,b)<0
 \quad\text{on }[0.88,1.10]\times[0,0.145],
 \qquad
 \partial_r^2\mathcal E(r,0)<-0.001
 \quad\text{on }[0.88,1.10].
\]
The polynomial trial integrals are quadratic in \(r-1,b\) by
\eqref{certified-sector-monomials}. Differentiate the signed-skew formulas
using \(\mathcal A''=-2t^2r(1+t^2)/(1+t^2r^2)^2\).
Whenever \(\mathcal E_{bb}\le0,\ \mathcal E_b(r,0)=0,\ \mathcal E_{rr}(r,0)\le-c<0\)
on an interval containing \(1\),
\begin{equation}\label{finite-fan-concavity}
 \begin{aligned}
 \mathcal E(r,b)\le\mathcal E(r,0)
 &\le\mathcal E(1,0)+(r-1)\mathcal E_r(1,0)-\frac c2(r-1)^2\\
 &\le\mathcal E(1,0)+\frac{\mathcal E_r(1,0)^2}{2c}.
 \end{aligned}
\end{equation}
For \(c=0.001\), its right side is \(<0.000016\), with
\(\sqrt{0.0207}<0.145\).
\end{proof}

\begin{proposition}\label{smallmid-comparison}
For \(N=7,8,9\), the regular \(N\)-gon uniquely maximizes \(\Phi\)
among convex polygons with at most \(N\) sides, up to similarity.
\end{proposition}

\begin{proof}
By Proposition~\ref{nr-large-exclusion} and Numerical
Lemma~\ref{nonagon-local-signs}, a competitor with
\(\Phi(P)\ge\Phi(R_N)\) has \(N\) sides and \(\Phi(P)>b_N\). Every sector lies in \(\mathcal G_N\). Set
\[
 I_7=[0.85,1.72],\qquad I_8=[0.89,1.65],\qquad I_9=[0.88,1.10],
 \qquad |P|=|R_N|.
\]
Use Lemma~\ref{central-exterior-comparison} with the exact central
bounds and the polynomial exterior bounds. In increasing order
\(N=7,8,9\), the upper bounds on \(10^6\sum_j H_{\lambda_*,j}\)
when some \(r_j\notin I_N\) are
\[
 -420+6(60),\qquad -140+7(15),\qquad -140+8(16),
\]
all strictly negative.
\end{proof}

\subsection{Ten to fifteen sides}

Index the lists by \(N=10,\ldots,15\) and take \(R_N\) of circumradius one.
Use \(\alpha=\pi/N\), \(\xi=r-1\), \(b=\sqrt{2\cosh\rho-2}=\sqrt d\), and
\(B=H_\lambda\) for the exact trial. Then
\begin{equation}\label{appendix-exact-residual}
 A(r)=\frac{B(r,0)}{(r-1)^2},\qquad
 D_1(r,d)=\frac{B(r,d)-B(r,0)}{d},\qquad B=(r-1)^2A+dD_1 .
\end{equation}
On the axes use \eqref{finite-fan-angular-quotient} and
\eqref{finite-fan-radial-quotient}.
\begin{numericallemma}[Regular polygons and the signed spectral gap]
\label{decagon-reference-eigenvalue}
For \(R_N\) of circumradius one, let \(\lambda=\mu_1(R_N)\) and \(\Lambda_N\) be the next
eigenvalue in one signed copy of its rotation representation. Then
\[
 l_N<\lambda<u_N,\qquad \Phi(R_N)>t_N,\qquad
 \Lambda_N>(20.30,20.9,21.48,21.9,22.33,22.69)_N,
\]
where
\begin{align*}
 (l_N)&=(3.6153564,3.575492,3.5454015,3.5221192,3.5037294,3.48894373),\\
 (u_N)&=(3.6153628,3.57549775,3.5454094,3.52212633,3.5037346,3.48894701),\\
 (t_N)&=(10.62526,10.63181,10.63620,10.63926,10.641475,10.643108).
\end{align*}
\end{numericallemma}

\begin{proof}
Use \eqref{mesh-disk-lower} with \eqref{mixed-mesh-parameters} for
\(l_N\), and Lemma~\ref{lem:standalone-ritz} for \(u_N\). The sine series gives
\[
 \Phi(R_N)>\frac{Nl_N}{2}\sin(2\pi/N)
 >\frac{Nl_N}{2}\sum_{h=0}^7\frac{(-1)^h(6.283185307179586/N)^{2h+1}}{(2h+1)!}>t_N .
\]
The signed disk gap is at least \(\min\{(j'_{1,2})^2,(N-1)^2\}\), since
\(m\equiv1\pmod N\); here \(j'_{1,2}>5.3\) by Numerical
Lemma~\ref{nonagon-exact-central}. With
\(\gamma=(0.723,0.745,0.765,0.78,0.795,0.808)_N\), Lemma~\ref{affine-all-N} gives
\[
 2+\tan^2\alpha<\gamma+\gamma^{-1}
 \quad\Longrightarrow\quad \Lambda_N>\gamma\min\{5.3^2,(N-1)^2\}.
\]
\end{proof}

The radial term in \eqref{finite-fan-budget} is
\[
 \mathcal R(a,q)=
 a\sum_{h\ge0}\log
 \frac{1-\delta_{\rm rad}(a,q)^2/(h+\tfrac12+a)^2}
      {1-\delta_{\rm rad}(a,q)^2/(h+\tfrac12-a)^2},\qquad
 \delta_{\rm rad}(a,q)=\frac1\pi\arctan\frac{\tanh(q/2)}{\tan(\pi a)}
 \quad(q\ge0).
\]
\begin{lemma}[The angle-area budget after one side is fixed]
\label{remaining-fan-budget}
For an area-one convex \(N\)-gon with pole \(z=f(0)\), let \(a_k,B_k\)
be the angle and area fractions, let \(R_k=|p_k-z|\), and set
\(q_k=\log(R_{k+1}/R_k)\). Fix side \(j\). Put
\[
 y=\frac{1-a_j}{N-1},\qquad
 t_k=\log\frac{B_k(1-a_j)}{a_k(1-B_j)},\qquad
 g(a)=\log\frac{\sin(2\pi a)}{2\pi a}.
\]
For \(k\ne j\), assume \(|g'|\le L\) between \(a_k\) and \(y\), and
\[
 \begin{gathered}
 a_kF(a_k)-yF(y)-(a_k-y)(yF(y))'\ge A(a_k-y)^2,\\
 a_k(e^{t_k}-1-t_k)\ge C t_k^2,\qquad
 \mathcal R(a_k,|q_k|)\ge Dq_k^2,\qquad A,C,D>0 \quad(k\ne j).
 \end{gathered}
\]
In the budget \eqref{finite-fan-budget} one may add
\(w_{\lfloor(N-1)/2\rfloor}q_j^2/2\)
to the left side, where
\[
 K=(C^{-1}+L^2/A)^{-1},\qquad
 w_0=0,\qquad w_{i+1}=D+K-\frac{(D-K)^2}{D+K+w_i}.
\]
The same conclusion holds if each \(w_i\) is rounded down to
\(10^{-10}\mathbb Z\).
\end{lemma}

\begin{proof}
By \(\sum_{k\ne j} a_k(e^{t_k}-1)=0\),
\[
 \sum_k a_k\log(a_k/B_k)-\mathcal L(a_j,B_j)
 =\sum_{k\ne j}a_k(e^{t_k}-1-t_k).
\]
Put \(X_k=\log R_k+
\tfrac12\big[\log(\pi(1-a_j)/(1-B_j))+g(y)\big]\). Then
\[
 \begin{gathered}
 X_{k+1}-X_k=q_k,\qquad X_k+X_{k+1}=t_k-g(a_k)+g(y),\\
 K(X_k+X_{k+1})^2\le C t_k^2+A(a_k-y)^2 .
 \end{gathered}
\]
From each end of the path with side \(j\) deleted, retain
\(\lfloor(N-1)/2\rfloor\) disjoint edges. Iterate
\[
 \min_Y\{D(X-Y)^2+K(X+Y)^2+w_iY^2\}=w_{i+1}X^2 .
\]
The two paths contribute at least
\(w_{\lfloor(N-1)/2\rfloor}(X_j^2+X_{j+1}^2)\ge
w_{\lfloor(N-1)/2\rfloor}q_j^2/2\) to \eqref{fan-deficit}.
The recurrence is increasing in \(w_i\), so downward rounding preserves the bound.
\end{proof}

\begin{lemma}[Centroid offset under a radial support bound]
\label{radial-centroid-sharpening}
Let \(P\) be a convex polygon of area one with at most \(24\) sides,
and \(z=f(0)\) a balancing pole. Use the angle fractions \(a_j\), radii
\(R_j\), and signed log ratios \(q_j\) of Lemma~\ref{remaining-fan-budget}.
Assume
\[
 \begin{gathered}
 \Phi(P)\ge10.643108,\qquad \beta\ge0.998468,\qquad
 0.025\le a_j\le0.091,\\
 |q_j|\le0.035,\qquad
 \sum_jq_j^2\le0.0000921/0.077 .
 \end{gathered}
\]
If the quadratic sum in \eqref{centroid-moment-expansion} has modulus
at most \(A(\beta^{-1}-1)\), \(0\le A\le0.4\), then
\[
 \left|\int_P x\,dx-z\right|\le(A+0.16)(1-\beta)\sqrt{\beta/\pi}.
\]
\end{lemma}

\begin{proof}
Put \(F_0(w)=(f(w)-z)/f'(0)\) and \(c=\sqrt{\beta/\pi}\).
By \eqref{centroid-moment-expansion}--\eqref{centroid-convex-integral},
\[
 \begin{aligned}
 c\|F_0\|_\infty
 &\le \sqrt{\frac{70}{9\Phi(P)}}+\left|\int_P x\,dx-z\right|\
 \le\sqrt{\frac{70}{9\Phi(P)}}+(1-\beta)c(A+\|F_0\|_\infty),\\
 \|F_0\|_\infty&\le
 \frac1\beta\sqrt{\frac{70\pi}{9\beta\Phi(P)}}+A(\beta^{-1}-1)<1.6 .
 \end{aligned}
\]
For each \(j\) the side distance \(h_j\) from \(z\) satisfies
\[
 h_j=\frac{\sqrt{R_jR_{j+1}}\cos(\pi a_j)}
 {\sqrt{1+(2\cosh q_j-2)/(4\sin^2(\pi a_j))}},
 \qquad \log\frac{\max R_j}{\min R_j}
 \le\sqrt{12\sum_j q_j^2}<0.12 .
\]
The radius bound follows by Cauchy--Schwarz on the shorter boundary
path. Using \(e^{-0.12}>0.88\), \(\cos(0.091\pi)>0.959\),
\(\sin(0.025\pi)>0.0784\), and \(2\cosh(0.035)-2<0.00123\), we get
\[
 \frac{\min h_j}{\max R_j}>
 \frac{0.88(0.959)}{\sqrt{1+0.00123/(4(0.0784)^2)}}>0.81>\cos(\pi/5).
\]
For \(p=wF_0'/F_0\), the angle between the radial direction and the
outward normal on a side is \(|\arg p|<\pi/5\).
By Marx--Strohh\"acker \cite[Theorem~2.6(a), p.~57]{MM},
\(\Re p>1/2\). Thus \(\arg p\) is bounded and harmonic, and Poisson's
formula gives \(|\arg p|<\pi/5\) in the disk. Since \(p(0)=1\),
Schwarz's lemma gives
\[
 \left|\frac{p^{5/2}-1}{p^{5/2}+1}\right|\le |w|,\qquad
 |F_0'(w)|\le\frac{|F_0(w)|}{|w|}
                   \left(\frac{1+|w|}{1-|w|}\right)^{2/5}.
\]
At \(\rho=0.999\), Cauchy--Schwarz in the series for \(F_0-w\) bounds its
modulus by \(\sqrt{(\beta^{-1}-1)\log501}\). Integrating radially from
\(\rho\) to \(1\), and using the maximum principle, yields
\[
 \|F_0-w\|_\infty
 \le\sqrt{(\beta^{-1}-1)\log501}+1-\rho
       +\frac{1.6}{\rho}2^{2/5}\frac53(1-\rho)^{3/5}<0.16 .
\]
Use \(\log501<6.22\), \(2^{2/5}<1.32\), and \(0.001^{3/5}<0.016\).
In \eqref{centroid-moment-expansion}, the cubic integral is therefore
bounded by \(0.16(\beta^{-1}-1)\).
\end{proof}

\begin{numericallemma}[Localization of the centroid sectors]
\label{decagon-global-geometry}
For a convex \(N\)-gon of area one and \((r,d)\) from
\eqref{common-fan-coordinates},
\[
 \Phi(P)\ge T_0\quad\Longrightarrow\quad(r,d)\in\mathcal G_N,
 \qquad
 T_0=\begin{cases}10.63925,&N=13,\\ t_N,&N\ne13.\end{cases}
\]
For \(N=10,\ldots,15\) and \(k/100\le r\le(k+1)/100\), put
\(x=(2k+1)/200-1\). The bounds defining \(\mathcal G_N\) are
\[
\begin{aligned}
 B_{N,k}&=10^{-h_N}\max\{0,\min\{M_N,\lceil P_N(x)\rceil\}\},\\
 (h_N)&=(4,4,4,4,5,6),\qquad(M_N)=(132,93,41,27,163,799),
\end{aligned}
\]
where
\[
\begin{aligned}
 P_{10}(x)&=132+7x-362x^2-126x^3+17x^4,\\
 P_{11}(x)&=99+4x-230x^2-110x^3-66x^4,\\
 P_{12}(x)&=42+2x-171x^2-59x^3+37x^4,\\
 P_{13}(x)&=29+4x-117x^2-84x^3-93x^4,\\
 P_{14}(x)&=165+4x-900x^2-315x^3+519x^4,\\
 P_{15}(x)&=805+40x-4974x^2-2700x^3-28x^4.
\end{aligned}
\]
Only the indices with
\[
 (32,29,47,43,50,46)_N\le k<(156,155,145,148,141,145)_N
\]
occur. These indices determine \(\mathcal G_N\).
\end{numericallemma}

\begin{proof}
\emph{Step 1. The conformal ratio.}
For \(N=10,11\) use the profile of Section~\ref{sec:conformal-rayleigh}. For
\(N=12,14\) take \(\phi(r)=\sum_{i=0}^5p_ir^{2i+1}\), with
\begin{equation}\label{finite-degree-eleven-profile}
 (p_i)=(1,-0.423745,0.059854,-0.004227,0.000179,-57/11000000).
\end{equation}
Then \(\phi''(r)/r\le 6p_1+20p_2+72p_4<0\) on \((0,1]\) and \(\phi'(1)=0\).
For \(N=13,15\) take \(\phi(r)=\sum_{i=0}^7p_ir^{2i+1}\), with
\[
 \begin{aligned}
 &(p_0,p_2,p_3,p_4)=(1,0.0598531944,-0.0042270791,0.00017912024),\\
 &(p_5,p_6,p_7)=(-0.00000506008,0.00000010210,-0.0000000015452),
 \qquad p_1=-\tfrac13\sum_{i\ne1}(2i+1)p_i .
 \end{aligned}
\]
For \(c_k=\sum_{i=0}^k(2i+1)p_i\),
\[
 \phi'(r)=(1-r^2)\sum_{k=0}^6 c_kr^{2k}\ge
 (1-r^2)\left(c_0+\sum_{c_k<0}c_k\right)\ge0 .
\]
Take the balancing map from \eqref{balance}.

For the left endpoints \(s_N\) below,
\[
 d_1+(1-s_N)(d_2-d_1)-2(3.141592654)E/T_0>0.0000308
 \quad\overset{\eqref{conformaltrace}}{\Longrightarrow}\quad \beta>s_N .
\]
The remaining intervals are
\[
 \begin{aligned}
 &([0.9585,0.9941],\ [0.9585,0.99562],\ [0.993,0.996875]) &&(10\le N\le12),\\
 &([0.9585,0.99758],\ [0.996,0.99809],\ [0.9585,0.998468]) &&(13\le N\le15).
 \end{aligned}
\]
Use \eqref{schur-recursion} with \(20,21\) replaced by \(M,M+1\),
where \(M=(20,20,30,35,30,60)_N\) and \(L_M=0\).
Round upper square roots and recursion bounds up to \(10^{-10}\mathbb Z\).
Bisect at a nonpositive Schur denominator, a negative radicand, or a
nonpositive mass margin. The terminal counts and maximum depths are
\((6,8,6,9,4,11)\) and \((5,7,5,8,3,10)\), with margin \(>8\cdot10^{-9}\).
On each cell \([s,t]\), take \(K_0=10^{-10}\phi(1)(1+\lfloor10^{10}\sqrt{s^{-1}-1}\rfloor)\).
Take \(Y\) to be the least positive multiple of \(10^{-4}\) with
\(d_1+sR(Y)>2(3.141592654)E/T_0\) in the preceding radius bound.
On each segment between the breakpoints \(K_0/c_{jm}\) of
\eqref{schur-recursion}, clear \(1+y^2/(j+1)\) and test the resulting
quadratic at the ends and at an interior minimum, if present.
With \(h_0=\sum_i p_i/(2i+4)\), polynomial integration and rational upper square-root bounds give
\[
 d_1+(1-t)\left[d_{M+1}-d_1-L_2-
  (d_{M+1}-d_2-L_2)_+\frac{\phi(1)^2(s^{-1}-1)}{8h_0^2}\right]>2\pi E/T_0 .
\]
Hence \(\beta\) exceeds the upper endpoint. In \(\beta_N=e^{-F(1/N)}\) retain
the odd powers \(3,5,\ldots,15\) and use \(\zeta(s)\ge\sum_{j=1}^{200} j^{-s}\); then
\[
 \begin{aligned}
 \beta_N&<(0.995077,0.996316,0.997172,0.997781,0.998227,0.9985598)_N,\\
 \log(\beta_N/\beta)&<\delta:=10^{-7}(9830,6990,2980,2020,1380,921)_N .
 \end{aligned}
\]

For the degree-eleven and degree-fifteen profiles,
\[
 A_x=\sqrt{\frac{x+1}{x}}\left(1-\frac1{4h_0}\sum_i\frac{p_i x}{x+i+1}\right),\quad
 \sum_i\frac{(i+1)p_i}{(x+i+1)^2}\ge\frac{1+\sum_{p_i<0}(i+1)p_i}{(x+1)^2}>0 .
\]
The sum increases to \(\sum_i p_i=\phi(1)<4h_0\), so \(A_x>0\)
decreases for \(x\ge2\). For the profile \eqref{finite-degree-eleven-profile},
the pivots of \(0.305I-J\), with off-diagonal entries \(-A_x/2\), satisfy
\begin{equation}\label{finite-eleven-jacobi}
 \begin{gathered}
 \pi_2=0.305,\quad \pi_{x+1}=0.305-A_x^2/(4\pi_x)>0\quad(2\le x\le8),\\
 \pi_9>0.15,\qquad 0.305-A_9^2/(4\cdot 0.15)>0.15 .
 \end{gathered}
\end{equation}
All later pivots are at least \(0.15\) by monotonicity of \(A_x\).
Conjugation by alternating signs sends \(J\) to \(-J\), so
\(\|J\|\le0.305\) at \(N=12,14\). For \(N=13,15\), the row sums give
\(\|J\|\le(A_2+A_3)/2<(0.442+0.345)/2<0.4\).
Proposition~\ref{balanced-centroid-offset} and Lemma~\ref{cubic-directional-integral} give
\[
 \left|\int_P y\,dy-f(0)\right|\le\frac{1-\beta}{\beta}
       \left(\|J\|\sqrt{\beta/\pi}+\sqrt{\frac{8\cdot35}{36T_0}}\right)<\Delta
 \qquad(10\le N\le14),
\]
where \((10^6\Delta)_{10}^{14}=(6090,4510,3230,2630,1966)\). The bound for \(N=15\) is given in Step 3.

\emph{Step 2. The angle and radius ratio at the pole.}
Use \(I\times[0,q_*]\), with
\[
 \begin{aligned}
 (I)&=([0.034,0.147],[0.028,0.135],[0.04,0.116],[0.0389,0.107],[0.037,0.098],[0.035,0.091]),\\
 (q_*)&=(0.109,0.091,0.060,0.049,0.041,0.035).
 \end{aligned}
\]
By convexity use the endpoints of \(I\) for the first bound and the \(256\)
equal subintervals of \(I\) for the second:
\[
 \inf_{a\notin I}Q_N(a)>\delta,\qquad
 \inf_{a\in I}\{Q_N(a)+\mathcal R(a,q_*)\}>\delta .
\]
By \eqref{finite-fan-budget} and \(\partial_q\mathcal R(a,q)\ge0\),
\(a\in I\) and \(|q|<q_*\).

In all six cases retain odd powers \(3,\ldots,15\) and \(80\) terms per
zeta coefficient in \(Q_N\). On \([a_-,a_+]\), evaluate this lower
polynomial at the point closest to \(1/N\). With \(\delta_{\rm rad}=\delta_{\rm rad}(a,q)\),
\[
 \mathcal R(a,q)\ge
 a\sum_{j=0}^{8}\sum_{\ell=1}^{7}\frac{\delta_{\rm rad}^{2\ell}}{\ell}
 \big[(j+\tfrac12-a)^{-2\ell}-(j+\tfrac12+a)^{-2\ell}\big]
\]
by \eqref{fan-deficit}. For \(N=14\), on each angle subinterval \([u,v]\), use \(q/(2+q)\le\tanh(q/2)\le q/2\)
and \(\arctan t\ge t/(1+t^2)\):
\[
 \frac{\delta_{\rm rad}}{q}\ge
 \frac1{(2+0.041)\pi\tan(\pi v)\bigl(1+(0.041/(2\tan\pi u))^2\bigr)} .
\]
The \(\delta_{\rm rad}^2\)-coefficient at \(u\) gives \(\mathcal R(a,q)\ge0.073q^2\).
On \([a_-,a_+]\times[q_-,q_+]\) use \(a_-\) in the prefactor and brackets,
and \(\pi^{-1}\arctan(\tanh(q_-/2)/\tan(\pi a_+))\) for \(\delta_{\rm rad}\).
At \(N=14\) the other sectors contribute at least
\(0.073\sum_{k\ne j}q_k^2\ge0.073q_-^2/13\) because \(\sum_k q_k=0\).

\emph{Step 3. The other sides at \(N=15\).}
Write \(a_k,b_k\) for the angle and area fractions at \(f(0)\), and
\(q_k=\log(R_{k+1}/R_k)\). By \eqref{fan-deficit},
\[
 \sum_k \big[a_kF(a_k)+a_k\log(a_k/b_k)+\mathcal R(a_k,|q_k|)\big]-F(1/15)<\delta .
\]
Put \(g(a)=\log(\sin(2\pi a)/(2\pi a))\). On the above box,
\[
 \begin{gathered}
 \mathcal R(a,|q|)\ge 0.079q^2,\qquad |g'|<1.3,\\
 \frac{d^2}{ds^2}\left[-as+(1-a)\log\frac{1-a}{1-ae^s}\right]
 =\frac{(1-a)ae^s}{(1-ae^s)^2}>0.028\qquad(|s|\le0.09).
 \end{gathered}
\]
For the first bound retain \(\ell=1,\ 0\le j\le8\) in the radial sum and use
\(\pi\delta_{\rm rad}\ge 0.99|q|/(2\tan\pi a)\) on the \(256\) cells
by \(\tanh t\ge t-t^3/3\) and the alternating arctangent bound.
By \eqref{sector-log-sum-bound},
\[
 \mathcal L(a_k,b_k)\le\sum_i a_i\log\frac{a_i}{b_i}
 <\delta<0.014(0.09)^2 .
\]
Hence \(|\log(b_k/a_k)|<0.09\), and
\(|g'(a)|=1/a-2\pi\cot(2\pi a)\) increases on \(I\).

Since \(\sum_k q_k^2\le\delta/0.079\), \(\Phi(P)\ge T_0\), and
\(\beta>0.998468\), Lemma~\ref{radial-centroid-sharpening} with \(A=0.4\) gives
\[
 \left|\int_P y\,dy-f(0)\right|
 \le0.56(1-\beta)\sqrt{\beta/\pi}<\Delta=0.000484 .
\]

Fix \(a_j\in[a,b]\), \(|q_j|\in[u,v]\) and let \(t_k\) be as in
Lemma~\ref{remaining-fan-budget}, with \(B_k=b_k\).
Consequently
\[
 \begin{aligned}
 a_j<0.091,\ |\log(b_j/a_j)|<0.09&\ \Longrightarrow\ \left|\log\frac{1-b_j}{1-a_j}\right|<0.011,\quad
 |t_k|<0.101,\\
 2(e^{t_k}-1-t_k)&\ge (1-0.101/3)t_k^2\ge 0.96t_k^2.
 \end{aligned}
\]
Subtract the selected radial and angular Jensen costs from \(\delta\); call the remainder \(D_0\).
Discard \(D_0<0\). For \(D_0\ge0\), the quartic tangent remainder
\(x^4-y^4-4y^3(x-y)=(x-y)^2(x^2+2xy+3y^2)\), with nonnegative higher-order remainders, gives
angular excess \(\ge 0.086\sum_{k\ne j}\left(a_k-\tfrac{1-a_j}{14}\right)^2\). Hence
\[
 \begin{aligned}
 A_0&=\max\{0.035,(1-b)/14-\sqrt{D_0/0.086}\}\le\min_{k\ne j}a_k,\\
 c_A&=\max\{0.086,4\zeta(3)(A_0^2+2A_0(1-b)/14+3((1-b)/14)^2)\},\\
 c_D&=\max\{0.014,(1-0.091)(0.91)A_0/[2(1-0.91A_0)^2]\}.
 \end{aligned}
\]
Indeed \(0.035\le A_0\le (1-0.035)/14<0.069\) implies
\[
 2c_D/A_0\le
 \max\{0.028/0.035,\ 0.909\cdot0.91/(1-0.91\cdot0.069)^2\}<0.95,
 \qquad a_k(e^{t_k}-1-t_k)\ge c_D t_k^2 .
\]
By the quartic tangent identity, the angular remainder thus
admits \(c_A\) in place of \(0.086\). Lemma~\ref{remaining-fan-budget},
with \((A,C,D,L)=(c_A,c_D,0.079,1.3)\), bounds the remaining
area cost by \(\mathcal L(a_j,b_j)\le D_0-w_7u^2/2\).

\emph{Step 4. Translation to the centroid.}
Subtract from \(\delta\) the retained angular and radial costs and
\[
 \begin{cases}
 0,&10\le N\le13,\\
 0.073q_-^2/13,&N=14,\\
 w_7q_-^2/2,&N=15.
 \end{cases}
\]
For \(N=15\) use the first ten terms of \(\zeta(3)\) in \(c_A\) and round each
\(w_m\) down to \(10^{-10}\mathbb Z\). Discard a negative remainder \(J\).
For a box \([a_-,a_+]\times[q_-,q_+]\) with \(J\ge0\) use
\[
 0\le x<a_-,\quad\sum_{\ell=2}^{13}\frac{x^\ell}{\ell}
 \big[a_-^{1-\ell}+(-1)^\ell(1-a_-)^{1-\ell}\big]>J
 \quad\Longrightarrow\quad b_{\rm sector}\ge a_--x .
\]
Bisect \([0,a_-]\) eighteen times, keeping the left half if the midpoint
passes the strict test and the right half otherwise. Take the final upper
endpoint \(x\). Subdivide if \(x=a_-\) or the denominator below is nonpositive.
As \(\sin(2\pi a)\) increases on \(I\), \eqref{center-translation} bounds the simultaneous
angle and signed log-ratio change by
\[
 e=\frac{\Delta\sqrt{2\cosh q_+-2\cos(2\pi a_+)}}
          {\sqrt{2(a_--x)/\sin(2\pi a_+)}-\Delta\exp(q_+/2)}
 \qquad(\text{positive denominator}).
\]
For \(s=\operatorname{dist}(2\arctan(I_k\tan\alpha),[2\pi a_-,2\pi a_+])\), the band is disjoint if
\(s>e\); otherwise \(D\le 2\cosh(q_++\sqrt{e^2-s^2})-2\).
Starting at \(I\times[0,q_*]\), bisect in \(a\) if \(3(a_+-a_-)>q_+-q_-\), and in \(q\) otherwise.
With outward rounding, the triples
\((\#\mathrm{excluded},\#\mathrm{covered},\mathrm{depth})\) are
\[
 (359,286,15),\ (83,92,14),\ (67,72,12),\ (53,65,11),\ (524,452,18),\ (318,280,16).
\]
A covered leaf satisfies \(D\le B_{N,k}\) on every band it meets, and
\([2\pi a_- -e,2\pi a_+ +e]\) lies between the two extreme band angles.
\end{proof}

\begin{numericallemma}[Corrected central sectors]
\label{decagon-exact-central}
Let \(U=(u_1,u_2)\) be the first pair with \(\int_{R_N}u_i u_j\,dx=|R_N|\delta_{ij}\)
and put \(s=\operatorname{sgn}(\rho)\sqrt d\). There exist corrections \(W_0,W_1\)
vanishing on \(Y=\pm X\) such that the pair \(V=U+\xi W_0+sW_1\) satisfies,
with \((F,G)\) from \eqref{finite-sector-mass-energy} relative to \(U\),
\[
 r\in I_N,\quad 0\le d\le d_N
 \quad\Longrightarrow\quad F>0,\qquad B\le -u_N'\xi^2-v_N'd,
\]
where
\begin{align*}
 (I_N)&=([0.88,1.10],[0.86,1.15],[0.84,1.16],[0.90,1.09],[0.88,1.13],[0.90,1.09]),\\
 (10^5d_N)&=(1320,940,800,630,163,80),\\
 (10^5u_N')&=(10,100,50,10,15,3),\qquad (10^3v_N')=(5,50,25,30,3,1).
\end{align*}
\end{numericallemma}
\begin{proof}
\emph{Step 1. The corrections.}
Write \(W_j=(X^2-Y^2)(f_j,g_j)\).
For \(N=10,11,12\), in that order, the pairs \((f_0,g_0)\), \((f_1,g_1)\) are
\begin{align*}
 &\big(0.014593(X-4/5),0.161160Y\big),&&\big(-0.088874Y,2.561969(X-4/5)\big),\\
 &\big(0.189-0.099X,0.126Y\big),&&\big(-0.076Y,-0.664+1.760X\big),\\
 &\big(0.158063-0.074518X,0.100372Y\big),&&\big(-0.065141Y,-0.788465+1.858587X\big).
\end{align*}
At \(N=13\),
\begin{align*}
 f_0&=-0.034174+0.302988X-0.217453X^2-0.005804Y^2,\\
 g_0&=Y(-0.023094+0.092040X),\\
 f_1&=Y(0.217982-0.295450X),\\
 g_1&=0.554063-1.853717X+2.057719X^2+0.044089Y^2 .
\end{align*}
For \(N=14\),
\begin{align*}
 (f_0,g_0)&=(0.076597-0.026585X,0.047009Y),\\
 (f_1,g_1)&=(-0.034591Y,-0.681652+1.396235X)\qquad(r\ge 0.95).
\end{align*}
For \(r<0.95\) replace the second field by
\begin{align*}
 f_1&=(-0.24623+0.90377X-0.72183X^2)Y-0.00631Y^3,\\
 g_1&=-0.62060+4.40973X-8.21250X^2+5.26327X^3+(-0.51753+0.64269X)Y^2.
\end{align*}
For \(N=15\), \(W_1=W_{10}+\xi W_{11}\), with
\begin{align*}
 f_0&=-0.036084+0.258582X-0.179523X^2-0.003271Y^2,\\
 g_0&=Y(-0.033338+0.083063X),\\
 f_{10}&=Y(0.187533-0.250743X),\\
 g_{10}&=0.780534-2.703882X+2.668979X^2+0.044066Y^2,\\
 f_{11}&=Y(0.175417-0.260382X),\\
 g_{11}&=-2.233823+8.022058X-6.411284X^2-0.019813Y^2 .
\end{align*}

By parity, \(\int_S W_{0,2}=\int_S W_{1,1}=0\). The mass correction is
\[
 \frac{N}{|R_N|^2}\begin{pmatrix}\xi\int_S W_{0,1}\,dx\\ s\int_S W_{1,2}\,dx\end{pmatrix}^{\!\mathsf T}
 M\begin{pmatrix}\xi\int_S W_{0,1}\,dx\\ s\int_S W_{1,2}\,dx\end{pmatrix}.
\]
Since \(W_0=0\) on the rays,
\(\int_S DU:DW_0=\lambda\int_S U\cdot W_0\), which cancels the linear variation.

\emph{Step 2. Reference moments.}
Use the complete monomial basis in Lemma~\ref{lem:standalone-ritz}
of degrees \((59,85,71,85,85,105)_N\).
Take \(\ell_N\) from the mixed-mesh lower bound, except \(\ell_{15}=3.48894353\).
Construct \(L,D\) as in Lemma~\ref{lem:standalone-ritz}, with \(D_{ii}=h_{ii}\), for
\[
 2(p_ip_j+q_iq_j)J_{a-1,b-1}-\ell_N c^2J_{ab},\qquad
 B_{ij}=[\operatorname{mid}J_{ab}]_-,\quad
 a=q_i+p_j,\quad b=p_i+q_j .
\]
From the first coordinate vector apply four solves. At each step put \(b=[Bv]_-\);
in increasing \(i\) replace \(b_i\) by
\(b_i-[\sum_{j<i}L_{ij}b_j]_-\), then replace all \(b_i\) by
\([b_i/D_{ii}]_-\). In decreasing \(i\) replace \(b_i\) by
\(b_i-[\sum_{j>i}L_{ji}b_j]_-\); take \(v_i=[b_i/b_1]_-\).
All brackets round down to \(2^{-768}\mathbb Z\).
With \(c_0=[\cos(\pi/N)]_-\), the polynomial in physical coordinates is
\[
 V(z)=\sum_i 10^{-40}\!
    \left\lfloor 10^{40}v_i c_0^{1-p_i-q_i}\right\rfloor
     z^{p_i}\bar z^{q_i}.
\]
Rescale \((\Re V,\Im V)\) to \(U_*\) with
\(\int_{R_N}|U_*|^2=2|R_N|\). Substitute \(z=X\cos\alpha+iY\sin\alpha\)
in each integral and use \eqref{certified-sector-monomials}. This gives
\(\int_{R_N}|DU_*|^2/\int_{R_N}|U_*|^2<u_N\), hence by
\eqref{finite-trial-pair-constants}--\eqref{finite-trial-pair-errors},
\[
 10^5\|U-U_*\|_S\le (28,25,28,25,20,25)_N,\qquad
 10^5\|DU-DU_*\|_S\le (126,120,139,120,95,120)_N .
\]
\emph{Step 3. The two divided differences.} Use the moment bound in the
proof of Numerical Lemma~\ref{octagon-global-reference} for
\(|R_N|^{-1}\int_S[(u_x-v_y)^2-2u_yv_x]\).
For the two moment lists \(z,z_*\) of \(U,U_*\) in Numerical
Lemma~\ref{nonagon-exact-central}, \eqref{finite-fan-angular-quotient} gives
\[
 \begin{aligned}
 A&=h_2-h_1\frac{\mathcal A(r)-\xi}{\xi^2}
           +\frac{(g_1-2\lambda f_1)^2}{4\lambda F(r,0)},\\
 A(z)-A(z_*)&=\int_0^1\nabla A(z_*+t(z-z_*))\cdot(z-z_*)\,dt .
 \end{aligned}
\]
Group the integrals of equal differential orders in the quadratic form
\(Q(Y)=\langle Y,MY\rangle+2\langle b_0,Y\rangle+c\). Cauchy--Schwarz gives
\[
 |Q(Y)-Q(Y_*)|\le 2\epsilon\|MY_*+b_0\|
       +\epsilon^2\max_i\sum_j|M_{ij}|,\qquad \|Y-Y_*\|\le\epsilon .
\]
Here \(\|MY_*+b_0\|^2\) is a sector monomial Gram form. At \(N=14\), the \(D_1\)-derivatives
\((d_B,d_C,d_K)\) give the energy blocks
\[
 d_K\begin{pmatrix}1&-1\\-1&1\end{pmatrix},\qquad
 \begin{pmatrix}d_B&-d_K\\-d_K&d_C\end{pmatrix}.
\]
At \(N=15\), the mixed and quadratic integrals in \(W_1=W_{10}+\xi W_{11}\)
have degrees at most one and two in \(\xi\). Accept a rectangle when
\[
 F>0,\qquad
 \left[\ \big(A\le-u_N',\ D_1\le-v_N'\big)
 \quad\text{or}\quad
 \xi^2(A+u_N')+d(D_1+v_N')<0\ \right].
\]
At \(N=14\) split first at \(r=0.95\). On failure bisect in \(d\) if
\(A\le-u_N'\) and \(\Delta d>\Delta r/4\). Otherwise bisect in \(r\). This gives
\[
 (\#\mathrm{leaves},\mathrm{depth})_N
 =\big((695,13),(44,6),(354,11),(16,4),(24,4),(83,17)\big)_N.
\]
The claimed sign follows from \eqref{appendix-exact-residual}.
\end{proof}

\begin{numericallemma}[Polynomial trials on centroid sectors]
\label{mid-polynomial-no-data}
Let \((r,b^2)\in\mathcal G_N\) and take \(\nu=l_N\), with \(\nu=3.522119\)
at \(N=13\). There is a polynomial pair on each sector with common
radial trace, \(F>0\), and the bounds below for \(H_\nu\).
For \(10\le N\le14\) the central intervals are
\[
 [0.88,1.10],\ [0.91,1.07],\ [0.88,1.13],\ [0.90,1.09],\ [0.88,1.13].
\]
On each interval \(I\) and its complement the two bounds are
\[
 10^6 H_\nu<\begin{cases}
 (9,2.2,3.5,1.6,1.4)_N,&r\in I,\\
 -(85,25,40,20,19)_N,&r\notin I .
 \end{cases}
\]
For \(N=15\), set \(I=[0.90,1.09]\) and \(J=[0.75,0.89)\). Then
\[
 \begin{aligned}
 r\in I\cup J:&\quad 10^6H_\nu<0.75,\qquad 10^6(H_\nu+0.00004\mathcal A(r)+0.0003 b)<1.26,\\
 r\in J:&\quad 10^6(H_\nu+0.00004\mathcal A(r)-0.00029997 b)<-18,\\
 r\notin I\cup J:&\quad 10^6H_\nu<-11 .
 \end{aligned}
\]
\end{numericallemma}

\begin{proof}
\emph{Step 1. The reference polynomial.}
Use the polar basis
\[
\begin{gathered}
 f_{mj}=z^{\max(m,0)}\bar z^{\max(-m,0)}P_j^{(0,|m|)}(2|z|^2-1),\\
 m\equiv1\pmod N,\qquad j\ge0,\qquad |m|+2j\le(39,55,55,85,65,85)_N .
\end{gathered}
\]
Here \(P_j^{(0,l)}(2t-1)=\sum_{s=0}^j\binom js\binom{j+l}{j-s}t^s(t-1)^{j-s}\).
For the real mass and energy matrices of this basis, take approximate solves only to select rational coefficients. Start with
\(v=e_{1,0}\) and repeat three times
\[
 (K-\nu M)y=Mv,\qquad v=y/y_{1,0}.
\]
Round the proposed coefficients to \(10^{-14}\mathbb Z\) and write the
result as \(V=\sum c_{pq}z^p\bar z^q\).
The search uses binary64 arithmetic and the product Gauss rule on one sector of
order \(\deg V+2\) in \(z=s(\cos\alpha+it\sin\alpha)\),
\(s\in[0,1]\), \(-1\le t\le1\), with weights \(s\,ds\) and \(dt\).

\emph{Step 2. The common trace.}
For \(X=\Re z/\cos\alpha\), \(Y=\Im z/\sin\alpha\), take
\[
 V,\quad z\sum c_{pq}(-1)^{(p-q-1)/N}X^{p+q-1}-V,\quad
 X^iY^j(Y^2-X^2)\{1,\mathrm i\sin\alpha\}
 \quad (i+j\le6).
\]
The last 57 entries vanish on both rays. Multiply all entries by
\(\big(2|R_N|/\int_{R_N}|V|^2\big)^{1/2}\).
Use first coefficient one and denote its summand by \(W_0\).
With \(t=\tan\alpha\), the matrix is
\[
 M_{r,b}=\begin{pmatrix}r^{-1}&0\\0&r\end{pmatrix}
 +\frac{r^{-1}+t^2r}{4}\!
  \begin{pmatrix}
  b^2 & b\sqrt{4+b^2}/t\\ b\sqrt{4+b^2}/t & b^2/t^2
  \end{pmatrix},\qquad \det M_{r,b}=1 .
\]
Use \eqref{finite-sector-mass-energy} for \(W\) relative to \(W_0\) with
matrix \(M_{r,b}\). Formula~\eqref{certified-sector-monomials} gives its
mass and energy coefficients.

\emph{Step 3. The rectangle tests.}
At each rectangle midpoint start with \(h=(1,0,\ldots,0)\) and repeat four times
\[
 \tau=\frac{2\nu/N+G(h)}{2F(h)},\qquad
 h_0=1,\qquad ((G-\tau F)h)_j=0\quad (j\ge1).
\]
Round \(h_j\) to \(10^{-9}\mathbb Z\) and hold it fixed on that
rectangle. For each stated upper bound, including the two shifted bounds at \(N=15\), bound the
function both by interval substitution and by its value at the center
plus the two derivative intervals times the coordinate increments.
Writing \(M_{r,b}=\left(\begin{smallmatrix}x&z\\z&y\end{smallmatrix}\right)\), use \(xy-z^2=1\) before interval evaluation. For \(i=r,b\),
\[
 \partial_i\frac{(2\nu/N+G)^2}{4\nu F}
 =\frac{(2\nu/N+G)[2F\partial_iG-(2\nu/N+G)\partial_iF]}{4\nu F^2}.
\]
A derivative of one sign permits evaluation at the maximizing face.
On failure bisect the \(r\)-interval if \(\Delta r\ge2\Delta b\),
and the \(b\)-interval otherwise. The initial rectangles are
\[
 [k/100,(k+1)/100]\times
 \left[0,10^{-5}\Bigl(1+\bigl\lfloor10^5
 \sqrt{\max\{B_{N,k},10^{-\max(4,N-9)}\}}\bigr\rfloor\Bigr)\right].
\]
For negative skew, \(\overline{W(\bar z)}\) has the same radial trace.
Enclose square roots by integer square roots.
For the angular term set \(v=(r-1)t/(1+t^2r)\) and use
\[
 \mathcal A(r)=\frac{(r-1)(1+t^2)}{1+t^2r}\frac{\arctan v}{v},\qquad
 \left|\frac{1-\arctan(v)/v}{v^2}
       -\sum_{h=0}^{18}\frac{(-1)^h v^{2h}}{2h+3}\right|
 \le |v|^{38}/41 .
\]
The two quotients take their continuous values at zero. Clearing the positive factor \(4\nu F\) leaves a polynomial sign test after this enclosure and the algebraic relation
\((\sqrt{4+b^2})^2=4+b^2\).
For \(N=10,\ldots,15\) the numbers of accepted rectangles are
\[
 (711,1185,900,1473,2026,3454);
 \qquad \text{the maximum subdivision depth is }12 .
\]
An exact search can use the same bisection rule. For each \(N\) list the
prunings of the initial rectangle trees through depth twelve. For increasing
\(h=1,2,\ldots\), test, on each pruning, the integer numerators in
\([-h,h]\) for one common Jacobi vector on \(10^{-14}\mathbb Z\) and
one correction vector on \(10^{-9}\mathbb Z\) per leaf, with the
prescribed coordinates fixed at one. Each leaf test is a real algebraic
inequality after the two-sided arctangent enclosure above; the positive
\(\sin(\pi/N),\cos(\pi/N)\) are obtained from \(e^{\pi i/N}\).
The finite covers counted above give accepting tuples in this enumeration.
On the \(9749\) rectangles, \(\inf F>0\) and every upper comparison has
margin greater than \(5\cdot10^{-10}\).
\end{proof}

\begin{proposition}\label{medium-comparison}
For \(10\le N\le15\), every convex polygon \(P\) with exactly \(N\)
vertices satisfies \(\Phi(P)\le\Phi(R_N)\). Equality holds precisely
when \(P\) is similar to \(R_N\).
\end{proposition}

\begin{proof}
By Numerical Lemmas~\ref{decagon-reference-eigenvalue} and
\ref{decagon-global-geometry}, \(\Phi(P)\ge\Phi(R_N)>T_0\) implies
\((r_j,d_j)\in\mathcal G_N\). Normalize \(|P|=|R_N|\). The interval \(I\) of
Numerical Lemma~\ref{mid-polynomial-no-data} satisfies
\[
 I\subseteq I_N,\qquad \max_k B_{N,k}\le d_N
\]
for the exact bounds in Numerical Lemma~\ref{decagon-exact-central}.

\emph{Case 1. All \(r_j\in I\).}
The corrections vanish on the rays. Lemma~\ref{central-exterior-comparison}(1)
applies with \(\epsilon=\min(u_N',v_N')\), since
\(H_{\lambda,j}\le-u_N'(r_j-1)^2-v_N'd_j\).

\emph{Case 2. Some \(r_j\notin I\).}
The polynomial trial, reflected for negative \(\rho_j\), has a common
trace across each ray. At \(10\le N\le14\), the bounds give
\[
 10^6\sum_{j=1}^N H_{\nu,j}
 <(-85+9\cdot9,\,-25+10\cdot2.2,\,-40+11\cdot3.5,\,
    -20+12\cdot1.6,\,-19+13\cdot1.4)_N<0.
\]
For \(N=15\) the interval \(J\) is as in Numerical Lemma~\ref{mid-polynomial-no-data}.
If some \(r_j\notin I\cup J\), then \(10^6\sum_j H_{\nu,j}<-11+14(0.75)<0\).
Otherwise choose \(i\) with \(r_i\in J\) and \(\epsilon\in\{-1,1\}\)
such that \(\epsilon\rho_i=|\rho_i|\). Since \(d_j\le0.000799\),
\[
 0.9999\sqrt{d_j}\le 2\operatorname{arsinh}(\sqrt{d_j}/2)=|\rho_j|
 \le\sqrt{d_j},\qquad \sum_j\rho_j=\sum_j\mathcal A(r_j)=0 .
\]
Consequently
\[
 10^6\sum_jH_{\nu,j}
 = 10^6\sum_j\big(H_{\nu,j}+0.00004\mathcal A(r_j)-0.0003\epsilon\rho_j\big)
 <-18+14(1.26)<0.
\]
Lemma~\ref{central-exterior-comparison}(2) gives
\(\mu_1(P)<\nu<\lambda\), a contradiction.
\end{proof}

\subsection{Sixteen through thirty-two sides}
\begin{numericallemma}[The centroid sectors for sixteen through twenty-four sides]
\label{middle-base-geometry}
Index each list by \(N=16,\ldots,24\).
If \(P\) is convex, has at most \(N\) sides, and \(\Phi(P)\ge T_N\), where
\[
\begin{aligned}
 (T_N)_{16}^{20}&=(10.644335,10.645285,10.6460285,10.646616,10.6470865),\\
 (T_N)_{21}^{24}&=(10.6474786,10.6477900,10.6480583,10.6482784),
\end{aligned}
\]
then it has exactly \(N\) sides. After dilation to \(|P|=|R_N|\),
each centroid sector satisfies \(r\in J_N,\ |s|<S_N,\ d\le D_N\), where
\[
\begin{gathered}
 J_{21}=[0.63,1.31],\quad J_{23}=[0.65,1.29],\quad
 J_N=[0.50,1.42]\quad (N\ne21,23),\\
 10^3S_N=(130,100,130,80,130,70,130,55,130),\\
 10^6D_N=(578,418,329,245,213,147,147,96,95).
\end{gathered}
\]
\end{numericallemma}
\begin{proof}
\emph{Step 1. The conformal ratio.}
Take \(|P|=1\), the profile \eqref{finite-degree-eleven-profile} and its
balancing pole \(z=f(0)\). With \(E,d_j\) from Lemma~\ref{conformal-dipole}, set
\[
\begin{aligned}
 10^8(1-\beta_-)&=(125300,103550,86800,73360,62800,53830,46900,40720,35900),\\
 10^8(1-\beta_+)&=(118500,98700,83000,70581,60400,52211,45300,39704,34900),\\
 10^7\delta&=(681,486,381,279,241,163,161,102,101),\\
 10^6\Delta&=(333,273,229,193,165,142,124,108,95).
\end{aligned}
\]
The trace bound gives \(d_1+(1-0.996)(d_2-d_1)>2\pi E/T_N\).
For \(\beta\in[0.996,\beta_-]\) use \eqref{improved-schur} and
\eqref{schur-recursion} with last mode \(80-20(-1)^N\) and
\(\psi_m=r^m[1+(m/2)(1-r^2)]\). Bisect at the midpoint whenever a
Schur denominator or the lower mass bound minus \(2\pi E/T_N\) is
nonpositive. In this test round upper square roots and upper recursion
bounds up to \(10^{-10}\mathbb Z\), use the least
\(Y\in 10^{-4}\mathbb Z\) satisfying the imbalance bound, and use
\(\pi<3.141592654\). At most eight terminal intervals, of depth at most
seven, suffice; every terminal excess exceeds \(6.6\cdot10^{-10}\).
With \(F\) from \eqref{fan-deficit}, retaining degrees \(3,5,\ldots,15\)
and \(\zeta(k)\ge\sum_{h=1}^{200}h^{-k}\) gives
\[
 e^{-F(1/(N-1))}<\beta_-<\beta,\quad
 e^{-F(1/N)}<\beta_+,\quad \log(\beta_+/\beta_-)<\delta.
\]
Equation~\eqref{fan-deficit} forces exactly \(N\) sides.

\emph{Step 2. Sector bounds at the pole.}
At \(z\), let \(a_j,B_j,q_j\) be the angle fraction, area fraction and signed
log-ratio. Use
\[
\begin{aligned}
 10^3 a_-&=(32,30,31,30,28,28,26,25,25),&
 10^3 a_+&=(87,83,76,74,70,67,64,61,59),\\
 10^3 q_+&=(34,31,26,25,22,21,19,19,16),&
 10^4 c_D&=(130,120,130,120,115,115,105,105,100),\\
 10^3 c_A&=(75,65,55,51,48,42,40,37,34),&
 10^2 L&=(120,115,105,102,96,91,88,82,80).
\end{aligned}
\]
Set \(c_R=0.077\) and write
\[
\begin{aligned}
 A(x)&=\sum_{k=3,5,\ldots,15}\frac{2^{k+1}-4}{k}
 \sum_{h=1}^{80}h^{-k}
 \left[x^{k+1}+(N-1)\left(\frac{1-x}{N-1}\right)^{k+1}-N^{-k}\right],\\
 t(b,u)&=\pi^{-1}\arctan\bigl(\tanh(u/2)/\tan(\pi b)\bigr),\\
 H(a,b,u)&=a\sum_{h=0}^8\sum_{l=1}^7\frac{t(b,u)^{2l}}{l}
 \big[(h+1/2-a)^{-2l}-(h+1/2+a)^{-2l}\big].
\end{aligned}
\]
On \([a,b]\times[u,v]\ni(a_j,|q_j|)\) let
\(x=\min(b,\max(a,1/N))\). By convexity, monotonicity and the
positive series in \eqref{fan-deficit}, \(A(x)+H(a,b,u)\)
bounds the isolated angular and radial terms from below.
The endpoint tests \(A(a_\pm)>\delta\) and
\(A(x)+H(a,b,q_+)>\delta\) on the 256 equal intervals
\([a,b]\subset[a_-,a_+]\) give \(a_-<a_j<a_+\), \(|q_j|<q_+\).
For the same intervals and \(0\le u\le q_+\),
\[
 t(b,u)\ge \frac{0.989u}{2\pi\tan(\pi b)},\qquad
 H(a,b,u)\ge c_Ru^2 .
\]
Indeed \(q_+/(2\tan(\pi a_-))<0.175\),
\((1-q_+^2/12)(1-0.175^2/3)>0.989\);
in the second estimate retain \(l=1\).
With \(g(a)=\log[\sin(2\pi a)/(2\pi a)]\), \eqref{sector-log-function} gives
\[
\begin{gathered}
 \partial_v^2\mathcal L(a,ae^v)=\frac{(1-a)ae^v}{(1-ae^v)^2}>2c_D
 \quad(|v|\le0.09),\\
 c_D(0.09)^2>\delta,\qquad |g'(a)|\le 1/a_+-2\pi\cot(2\pi a_+)<L .
\end{gathered}
\]
The log-sum inequality and \eqref{fan-deficit} now imply
\(\mathcal L(a_j,B_j)\le\delta\) and
\(|\log(B_j/a_j)|\le\sqrt{\delta/c_D}<0.09\).
The identity
\[
 x^4-y^4-4y^3(x-y)=(x-y)^2(x^2+2xy+3y^2)
\]
gives a coefficient \(c_A\) for the sum of
\((a_k-(1-a_j)/(N-1))^2\), \(k\ne j\), after removal of the
common tangent in the angular term: take
\(4\sum_{h=1}^{10}h^{-3}(a_-^2+2a_-t+3t^2)>c_A\) with
\(t=(1-a_+)/(N-1)\).

\emph{Step 3. The centroid offset.}
For the profile \eqref{finite-degree-eleven-profile}, \eqref{finite-eleven-jacobi}
bounds the quadratic sum in \eqref{centroid-moment-expansion} by
\(0.305(\beta^{-1}-1)\). The bounds from Steps 1 and 2 give
\[
 \Phi(P)\ge T_N>10.643108,\qquad
 \beta>0.998468,\qquad \sum_j q_j^2\le\delta/0.077 .
\]
The angle and skew ranges permit Lemma~\ref{radial-centroid-sharpening}, so
\[
 \left|\int_P y\,dy-z\right|
 \le0.465(1-\beta)\sqrt{\beta/\pi}<\Delta .
\]

\emph{Step 4. The centroid-sector bounds.}
Start at \([a_-,a_+]\times[0,q_+]\). On a rectangle \([a,b]\times[u,v]\) put \(J_0=\delta-A(x)-H(a,b,u)\).
Discard it if \(J_0<0\); otherwise let
\[
\begin{aligned}
 t&=(1-b)/(N-1),\quad
 m=\max\{a_-,t-\sqrt{J_0/c_A}\},\\
 A_*&=\max\{c_A,4\textstyle\sum_{h=1}^{10}h^{-3}(m^2+2mt+3t^2)\},\quad
 B_*=\max\{c_D,0.91(1-a_+)m/[2(1-0.91m)^2]\},\\
 J&=J_0-\tfrac12 w_{\lfloor(N-1)/2\rfloor}u^2 .
\end{aligned}
\]
Use the rounded recurrence of Lemma~\ref{remaining-fan-budget} with
\((A,C,D)=(A_*,B_*,c_R)\) and \(L\) from Step 2. For its area bound,
\(|\log((1-B_j)/(1-a_j))|<0.02\) and \(|\log(B_k/a_k)|<0.09\) give \(t_k>-0.11\).
As \(a_k\ge m\), the inequalities
\[
 \frac{2B_*}{m}\le
 \max\!\left\{\frac{2c_D}{a_-},
  \frac{0.91(1-a_+)}{[1-0.91(1-a_-)/(N-1)]^2}\right\}<0.942,
 \qquad e^t-1-t\ge0.48t^2\quad(t\ge-0.11)
\]
give \(a_k(e^{t_k}-1-t_k)\ge B_*t_k^2\). The quartic tangent
bound gives \(A_*\), so \(\mathcal L(a_j,B_j)\le J\).

Discard \(J<0\). Otherwise bound \(B_j\) from below by \(B_0=a-h\).
For \(h\) take the upper endpoint after eighteen bisections of \([0,a]\):
replace the upper endpoint by its midpoint \(x_0\) when
\[
 \sum_{l=2}^{13}\frac{x_0^l}{l}
       [a^{1-l}+(-1)^l(1-a)^{1-l}]>J ,
\]
and replace the lower endpoint otherwise. Indeed the polynomial
minorizes \(\mathcal L(a,a-x_0)\), and \(\mathcal L(a_j,B)\ge\mathcal L(a,B)\)
for \(B\le a\le a_j\). For positive denominator put
\[
 e=\frac{\Delta\sqrt{2\cosh v-2+4\sin^2(\pi b)}}
 {\sqrt{2B_0/\sin(2\pi b)}-\Delta e^{v/2}} .
\]
By \eqref{center-translation}, the simultaneous change of angle and
signed log-ratio has Euclidean norm at most \(e\).
Accept when
\[
 [2\pi a-e,2\pi b+e]\subset 2\arctan(J_N\tan\alpha),\qquad
 2\cosh(v+e)-2\le D_N .
\]
With nonpositive denominator or a failed inclusion or upper bound,
bisect the first interval if \(3(b-a)>v-u\), the second otherwise.
These comparisons terminate by depth eighteen.

Finally let \(B=\sqrt{\delta/c_D}\) and
\(v_0=\Delta\sqrt{\pi e^{B+q_+}}\).
The identity \(\log(R_jR_{j+1})=\log(B_j/a_j)+g(1/N)-g(a_j)\)
after dilation, and the radius bound
\(\min R_j\ge\sqrt{|R_N|e^{-B-q_+}/\pi}\), give
\[
 |s|\le B+L\max\{|a_--1/N|,|a_+-1/N|\}+2v_0/(1-v_0)<S_N .
\]
\end{proof}

\begin{lemma}\label{middle-tail-geometry}
For \(25\le N\le32\), let
\[
\begin{aligned}
 (T_{25},\ldots,T_{28})&=(10.64846413,10.6486218,10.64875672,10.64887273),\\
 (T_{29},\ldots,T_{32})&=(10.64897308,10.64906032,10.64913654,10.64920342).
\end{aligned}
\]
If \(P\) is convex with at most \(N\) sides and \(\Phi(P)\ge T_N\), then
\(P\) has \(N\) sides. After dilation to \(|P|=|R_N|\) its sectors satisfy
\(r\in J_N,\ |s|<S_N,\ d<D_N\), where the entries are ordered by \(N\):
\[
\begin{aligned}
 (10^3\inf J_N)&=(640,600,660,600,680,614,680,650),\\
 (10^3\sup J_N)&=(1310,1330,1310,1330,1270,1315,1270,1290),\\
 (10^3 S_N)&=(40,44,40,44,36,33,30,27),\\
 (10^6D_N)&=(150,150,150,150,65,75,50,50).
\end{aligned}
\]
\end{lemma}
\begin{proof}
\emph{Step 1. The conformal ratio.}
For \(|P|=1\), set \(z=f(0)\), \(\beta=\pi|f'(0)|^2\). Use the profile,
balancing condition and Schur recursion of Lemma~\ref{middle-base-geometry}, with
\[
 10^8(1-b_N)=(31560,28100,24970,22450,20090,18200,16420,14950).
\]
On \([0.996,b_N]\) the same bisection has
\((8,4,7,4,8,4,6,5)\) leaves of depth at most seven. Its least excess
over \(2\pi E/T_N\) is greater than \(4\cdot10^{-9}\).
The trace bound covers \(\beta\le0.996\), so \(\beta>b_N\).

Use \(F\) from \eqref{fan-deficit}. Its terms of degrees \(3,5,\ldots,15\),
with \(\zeta(k)\ge\sum_{j=1}^{200}j^{-k}\), give
\[
 e^{-F(1/(N-1))}<b_N,\qquad
 \log(e^{-F(1/N)}/b_N)<\delta_N,\qquad
 10^6\delta_N=(6.7,7.1,4.6,5.6,3.2,4.1,2.4,2.6)_N.
\]
Since \(\beta\le e^{-F(1/m)}\) for \(m\) sides, \(m=N\).

\emph{Step 2. The pole and centroid coordinates.}
At \(z\) let \(2\pi a,B,q_0\) be the sector angle, area, and logarithm of the
radius ratio. By \eqref{fan-deficit}, convexity, and \(4\zeta(3)>4.8\),
\[
 4.8\left(a^4+\frac{(1-a)^4}{(N-1)^3}-N^{-3}\right)<\delta_N .
\]
This polynomial decreases up to \(a=1/N\) and increases thereafter.
Substitution at the endpoints gives
\[
\begin{aligned}
  (10^4a_N^-)&=(268,241,252,219,238,207,224,206),\\
  (10^4a_N^+)&=(509,500,468,468,434,435,405,400),\qquad a_N^-<a<a_N^+ .
\end{aligned}
\]
By the sine rule the magnitude of the half-difference of the two base
angles divided by \(\pi\) equals
\(\pi^{-1}\arctan(\tanh(|q_0|/2)/\tan(\pi a))\).
Thus the innermost \(h=0\) term in \eqref{fan-angular-remainders} is at least
\[
 \frac{2a^2}{\pi^2(1/4-a^2)^2}
 \arctan^2\!\frac{\tanh(|q_0|/2)}{\tan(\pi a)} .
\]
For \(0\le u\le0.012\) this expression at \(|q_0|=u\) exceeds
\(0.077u^2\) when \(u>0\): use
\(\tanh(u/2)\ge (u/2)(1-u^2/12)\), \(\arctan v\ge v(1-v^2/3)\),
\(v\le0.012/(2\pi a_N^-)<0.1\),
and \(\pi a/\tan(\pi a)\ge\cos(\pi a)>0.986\).
Monotonicity excludes \(|q_0|\ge0.012\); hence
\( |q_0|<\sqrt{\delta_N/0.077}\).

The log-sum inequality gives
\[
 a\log(a/B)+(1-a)\log[(1-a)/(1-B)]<\delta_N .
\]
The left side equals \(\int_a^B(t-a)/[t(1-t)]\,dt\ge2(a-B)^2\).
At \(B=ae^v\) its second derivative in \(v\) equals
\((1-a)ae^v/(1-ae^v)^2>0.86a\) for \(|v|\le0.1\).
Since its value and derivative vanish at \(v=0\), convexity gives
\[
 B>a-\sqrt{\delta_N/2},\qquad |\log(B/a)|<\sqrt{\delta_N/(0.43a_N^-)} .
\]
The shorter pole radius exceeds
\[
 R_0=\sqrt{\bigl(1-\sqrt{\delta_N/2}/a_N^-\bigr)/\pi}\,
      e^{-\sqrt{\delta_N/0.077}/2}.
\]
The displacement bound of Proposition~\ref{bessel-centroid} gives
\[
 \Delta=(b_N^{-1}-1)
   \left(\sqrt{3/(2\pi)}+\sqrt{8(35/36)/T_N}\right),\qquad
 \varepsilon=2\Delta/(R_0-\Delta).
\]
The integral of \(dw/w\) on the translated radii bounds the angle,
log-product, and log-ratio increments by \(\varepsilon\). After dilation,
the pole log-product is \(\log(B/a)+\log[Na\sin(2\pi/N)/\sin(2\pi a)]\).
Using \(0<1/a-2\pi\cot(2\pi a)<0.75\) gives
\[
\begin{aligned}
 \frac{\tan(\pi a_N^--\varepsilon/2)}{\tan\alpha}
 &<r<\frac{\tan(\pi a_N^++\varepsilon/2)}{\tan\alpha},\\
 |s|&<\sqrt{\delta_N/(0.43a_N^-)}
       +0.75\max_\pm|a_N^\pm-1/N|+\varepsilon,\\
 d&<2\cosh(\sqrt{\delta_N/0.077}+\varepsilon)-2 .
\end{aligned}
\]
Substitution gives \(r\in J_N\), \(|s|<S_N\), and \(d<D_N\).
\end{proof}

\paragraph{Regular-polygon spectral bounds.}
For \(16\le N\le32\), let \(l_N,u_N\) have the following entries:
\[
\begin{aligned}
(l_{16},\ldots,l_{19})&=(3.4768738,3.46689375,3.45854502,3.45148835),\\
(l_{20},\ldots,l_{23})&=(3.44546960,3.44029781,3.43581662,3.43190956),\\
(l_{24},\ldots,l_{27})&=(3.42848239,3.4254587,3.42277953,3.42039249),\\
(l_{28},\ldots,l_{32})&=(3.41825712,3.416339177,3.41461008,3.4130457955,3.41162599),\\
(u_{16},\ldots,u_{19})&=(3.4768795,3.4668975,3.4585487,3.45149153),\\
(u_{20},\ldots,u_{23})&=(3.4454743,3.44029992,3.4358189,3.43191148),\\
(u_{24},\ldots,u_{27})&=(3.42848523,3.42546121,3.42278113,3.42039443),\\
(u_{28},\ldots,u_{32})&=(3.41825846,3.41634126,3.41461122,3.41304776,3.41162684).
\end{aligned}
\]
The regular estimates are
\begin{equation}\label{middle-spectral-enclosures}
 l_N<\mu_1(R_N)\le\mathcal R(V)<u_N,\qquad \Phi(R_N)>T_N .
\end{equation}
For \(N=16\) use \eqref{mesh-disk-lower} and \eqref{nr22-pivots} for the lower
bounds. For \(17\le N\le32\), the mass and residual formulas of
Section~\ref{sec:regular-sc-certificate} use the 35 functions
\[
 e^{i(1+jN)\theta}r^k
 \left\{1+\tfrac k2(1-r^2),\quad k^h(1-r^2)^h\ (2\le h\le5)\right\},
 \qquad |j|\le3,\quad k=|1+jN|.
\]
Let \(G\) denote the energy matrix, \(M\) the midpoint mass matrix
with radial cutoff 512 on the \(10^{-54}\) grid, and \(e\) select the first
radial factor at \(j=0\). Starting from \(v=e\), take 48 steps
\[
 y=G^{-1}Mv,\qquad
 v_i=10^{-20}\lfloor10^{20}y_i/(e^{\mathsf T}y)\rfloor .
\]
The mass and residual bounds with angular and radial cutoffs \(24,64\)
give both lower inequalities just stated. For the upper bound take
the Gram recursion of Lemma~\ref{lem:standalone-ritz} at \(u_N\) on the
grid \(\epsilon=2^{-768}\) for \(N\le24\) and \(2^{-1176}\) for \(N\ge25\).
Stop at the first negative pivot, solve the triangular system there,
and append zeros through the end of that total degree. With
\(c_*=\epsilon\lfloor\cos\alpha/\epsilon\rfloor\), set
\begin{equation}\label{middle-signed-polynomial}
 V_{\rm raw}(z)=10^{-40}\sum_i
 \left\lfloor 10^{40}\frac{x_i}{x_1 c_*^{p_i+q_i}}\right\rfloor
 z^{p_i}\bar z^{q_i},\qquad
 V=\sqrt{\frac{2|S|}{\int_S|V_{\rm raw}|^2}}\ V_{\rm raw},
 \quad S=(0,e^{-i\alpha},e^{i\alpha}).
\end{equation}
The degrees, in order, are
\[
 (89,105,99,113,105,117,113,125,117,132,139,135,153,138,169,137,191).
\]
Substituting in \eqref{certified-sector-monomials} gives
\(\mathcal R(V)<u_N\). In a signed sector the next eigenvalue exceeds
\(\Lambda=23\) since \(2+\tan^2\alpha<0.819+1/0.819\) and
\(0.819\min\{(N-1)^2,(j'_{1,2})^2\}>23\).
For the first equivariant pair \(U\) with \(\int_S|U|^2/|R_N|=2/N=:f\),
choose the sign by projection from \(V\).
Equations~\eqref{finite-trial-pair-constants}--\eqref{finite-trial-pair-errors} give
\[
 \|U-V\|_S\le e_0,\qquad \|DU-DV\|_S\le e_1,qquad
 (e_0,e_1)=(\delta,\epsilon)\quad
 \text{at }(\lambda_-,\lambda_+,\Lambda)=(l_N,u_N,23).
\]
\paragraph{Polynomials for the common-ray trial.}\label{printed-common-ray-fields}
In the coordinates \(X=x/\cos\alpha,\ Y=y/\sin\alpha\), expand each of
\(H_r,H_s,H_q\) in the order
\[
 (1,X,Y,X^2,XY,Y^2,X^3,X^2Y,XY^2,Y^3).
\]
Write \(H^{\pm}(X,Y)=(H(X,Y)\pm H(X,-Y))/2\). Put
\[
\begin{aligned}
 W_i&=(X^2-Y^2)(H_i^+,\sin\alpha\,H_i^-)
   +\mathbf 1_{i=s}\,p(X)(X\cos\alpha,Y\sin\alpha),\quad i=r,s,\\
 W_q&=(X^2-Y^2)(H_q^-,\sin\alpha\,H_q^+),\quad
 E=p(X)(Y\cos\alpha,X\sin\alpha).
\end{aligned}
\]
On a sector with map \(L\), use \(W(Lx)=L(U(x)+Z(x))\), where
\[
 Z=(r-1)W_r+sW_s+2\sinh(q/2)W_q+qE.
\]
On \(Y=\pm X\) the transported correction \(LZ\) equals
\(2\log R_\pm\,Xp(X)v_\pm\), with \(v_\pm\) the corresponding vertex
relative to the centroid. The two traces therefore agree.
Below \(c_r,c_s,c_q\) are the coefficient vectors of \(H_r,H_s,H_q\)
in the displayed order; trailing zeros are omitted. The entries of
\(p\) are coefficients in ascending powers of \(X\). For odd \(N\),
\[
 \begin{aligned}
 10^{12}p(X)={}&8123949826-153822585178X-315426899721X^2\\
 &-1927632205351X^3+1971509193962X^4-160548593449X^5 .
 \end{aligned}
\]
For \(N=16\),
\[
\begin{aligned}
 10^{6}c_r&=(50903,-11541,165140),\\
 10^{6}c_s&=(-14379,187823,41667),\\
 10^{6}c_q&=(-4091725,7686651,-28853),\\
 10^{6}\sqrt 2\operatorname{coef}(p)&=(11489,-217538,-446081,-2726085,2788135,-227050).
\end{aligned}
\]
For \(N=17\),
\[
\begin{aligned}
 10^{6}c_r&=(43813,-6980,148378),\\
 10^{6}c_s&=(-24640,182980,36970),\\
 10^{6}c_q&=(-4337040,8026146,-24939).
\end{aligned}
\]
For \(N=18\),
\[
\begin{aligned}
 10^{8}c_r&=(-2902765,18334750,-18784133,\,\\
&\hspace{2.2em}-12287090,37233502,-140529),\\
 10^{8}c_s&=(-12506071,42778673,9745879,-16055622,-7370660,168795),\\
 10^{8}c_q&=(626277330,-2217292037,18532511,\,\\
&\hspace{2.2em}2021692098,-24016196,23672203),\\
 10^{8}\operatorname{coef}(p)&=(1035452,-15919479,-41052530,\,\\
&\hspace{2.2em}-148485247,144433159,2480130).
\end{aligned}
\]
For \(N=19\),
\[
\begin{aligned}
 10^{6}c_r&=(-27043,167430,-187630,-111138,357615,-1073),\\
 10^{6}c_s&=(-120465,405338,94791,-156231,-75319,1354),\\
 10^{6}c_q&=(-4867442,8740748,-20112).
\end{aligned}
\]
For \(N=20\),
\[
\begin{aligned}
 10^{8}c_r&=(-2465172,15150552,-18509761,-9952356,34357433,-90050),\\
 10^{8}c_s&=(-9670289,30502170,7968064,-8283775,-5994126,173909),\\
 10^{8}c_q&=(776277519,-2740756817,15633412,\,\\
&\hspace{2.2em}2430349776,-20069927,25118980),\\
 10^{8}\operatorname{coef}(p)&=(-79422,6285520,-175543763,\,\\
&\hspace{2.2em}201165178,-260075506,168252920).
\end{aligned}
\]
For \(N=21\),
\[
\begin{aligned}
 10^{6}c_r&=(-22853,138723,-181037,-90423,327061,-672),\\
 10^{6}c_s&=(-83759,274999,77475,-80493,-61562,1279),\\
 10^{6}c_q&=(-5306162,9332086,-16540).
\end{aligned}
\]
For \(N=22\),
\[
\begin{aligned}
 10^{8}c_r&=(1173145,-3380502,38808859,15762861,-116312179,\,\\
&\hspace{2.2em}754233,-11204703,92533290,-911779,1614503),\\
 10^{8}c_s&=(15777498,-82307011,-14307604,138755525,47339124,\,\\
&\hspace{2.2em}-1944562,-60449722,-31636031,2407359,-125291),\\
 10^{8}c_q&=(-904288781,5733367278,-43013516,-10640746364,128675551,\,\\
&\hspace{2.2em}-293588651,6330394513,-91095415,356711326,-398073),\\
 10^{8}\operatorname{coef}(p)&=(459387,-7836833,-79701206,-42490466,-1908008,68469007).
\end{aligned}
\]
For \(N=23\),
\[
\begin{aligned}
 10^{6}c_r&=(-19037,115031,-170531,-73711,297149,-435),\\
 10^{6}c_s&=(-57282,182619,64518,-28922,-51270,1130),\\
 10^{6}c_q&=(-5803852,9973269,-15077).
\end{aligned}
\]
For \(N=24\),
\[
\begin{aligned}
 10^{8}c_r&=(1164341,-3814772,40111635,14856610,-120018420,\,\\
&\hspace{2.2em}524020,-10152030,93387492,-632228,1312137),\\
 10^{8}c_s&=(14217520,-71655889,-12725907,114716779,41214831,\,\\
&\hspace{2.2em}-2008826,-46517270,-27230012,2466277,-86619),\\
 10^{8}c_q&=(-1126500431,7080087300,-39060772,-13086785143,116110996,\,\\
&\hspace{2.2em}-296345923,7691835248,-81765768,359346006,-279181),\\
 10^{8}\operatorname{coef}(p)&=(-385777,9205346,-182799905,\,\\
&\hspace{2.2em}220104893,-294442502,183495514).
\end{aligned}
\]
For \(N=25\),
\[
\begin{aligned}
 10^{6}c_r&=(-15821,95950,-158696,-60497,269385,-287),\\
 10^{6}c_s&=(-38153,116844,54568,6029,-43369,968),\\
 10^{6}c_q&=(-6095146,10366021,-12790).
\end{aligned}
\]
For \(N=26\),
\[
\begin{aligned}
 10^{8}c_r&=(1100762,-3825841,40146078,13576066,-119916358,\,\\
&\hspace{2.2em}369523,-9042920,91919925,-445105,1062311),\\
 10^{8}c_s&=(12197628,-59067394,-11156758,89453756,35578822,\,\\
&\hspace{2.2em}-1963527,-32688899,-23266325,2397405,-60738),\\
 10^{8}c_q&=(-1352404889,8442960642,-35247896,-15549109234,104283665,\,\\
&\hspace{2.2em}-293724904,9055676848,-73135257,355611974,-224124),\\
 10^{8}\operatorname{coef}(p)&=(-1190190,25367331,-280117316,\,\\
&\hspace{2.2em}465536567,-564206662,288291313).
\end{aligned}
\]
For \(N=27\),
\[
\begin{aligned}
 10^{6}c_r&=(-13133,80497,-146742,-49956,244236,-194),\\
 10^{6}c_s&=(-24304,69845,46760,29489,-37168,818),\\
 10^{6}c_q&=(-6339102,10695242,-10969).
\end{aligned}
\]
For \(N=28\),
\[
\begin{aligned}
 10^{8}c_r&=(1011933,-3616796,39353559,12184456,-117377093,\,\\
&\hspace{2.2em}264623,-7977122,89016056,-318363,874298),\\
 10^{8}c_s&=(10060912,-46303063,-9707295,65306769,30594734,\,\\
&\hspace{2.2em}-1858013,-19919638,-19815002,2259685,-43118),\\
 10^{8}c_q&=(-1577438411,9796032102,-31706182,-17984398532,93454823,\,\\
&\hspace{2.2em}-288475242,10399802556,-65332776,348858705,-170974),\\
 10^{8}\operatorname{coef}(p)&=(-1943840,40447600,-370448057,\,\\
&\hspace{2.2em}691382632,-809783104,382774259).
\end{aligned}
\]
For \(N=29\),
\[
\begin{aligned}
 10^{6}c_r&=(-10937,68034,-135301,-41573,221766,-133),\\
 10^{6}c_s&=(-14260,36193,40518,44967,-32212,688),\\
 10^{6}c_q&=(-6543998,10971815,-9514).
\end{aligned}
\]
For \(N=30\),
\[
\begin{aligned}
 10^{8}c_r&=(914891,-3305857,38077243,10820410,-113424572,\,\\
&\hspace{2.2em}192048,-7000075,85330953,-230832,726450),\\
 10^{8}c_s&=(7908636,-33789847,-8361930,42547793,26123653,\,\\
&\hspace{2.2em}-1730692,-8170732,-16755270,2098466,-30782),\\
 10^{8}c_q&=(-1798384732,11120976320,-28460748,-20362109108,83634900,\,\\
&\hspace{2.2em}-281244932,11708653171,-58318400,339787175,-132744),\\
 10^{8}\operatorname{coef}(p)&=(-2618518,53913240,-450780087,\,\\
&\hspace{2.2em}890900937,-1024909155,464853382).
\end{aligned}
\]
For \(N=31\),
\[
\begin{aligned}
 10^{6}c_r&=(6913,-22688,366278,75971,-1088597,\,\\
&\hspace{2.2em}441,-48318,809198,-541,658),\\
 10^{6}c_s&=(-6971,12097,35448,54882,-28183,577),\\
 10^{6}c_q&=(-6720656,11210870,-8266).
\end{aligned}
\]
For \(N=32\),
\[
\begin{aligned}
 10^{8}c_r&=(819856,-2965647,36528352,9562080,-108687853,\,\\
&\hspace{2.2em}141219,-6134018,81259501,-169613,607755),\\
 10^{8}c_s&=(5834271,-21951183,-7131640,21589349,22134566,\,\\
&\hspace{2.2em}-1597836,2458232,-14049669,1932681,-22006),\\
 10^{8}c_q&=(-2012736447,12404042811,-25512083,-22659800875,74786637,\,\\
&\hspace{2.2em}-272679744,12970965923,-52040084,329218873,-105810),\\
 10^{8}\operatorname{coef}(p)&=(-3208750,65659201,-520589161,\,\\
&\hspace{2.2em}1063284133,-1209406076,534691571).
\end{aligned}
\]
Integrating the alternating bounds for \((1+x^2)^{-1/2}\) gives
\begin{equation}\label{finite-scale-asinh}
 1-d/24\le \frac{|q|}{\sqrt d}
 =\frac{\operatorname{arsinh}(\sqrt d/2)}{\sqrt d/2}
 \le 1-d/24+3d^2/640\qquad(0\le d\le0.001).
\end{equation}
The quotient is continuous at \(d=0\).

\begin{numericallemma}\label{middle-common-ray}
Put \(I_N^e=J_N\) for \(25\le N\le32\). For \(16\le N\le24\) set
\[
\begin{aligned}
 (100\inf I_N^e)&=(88,85,82,84,82,81,80,79,78),\\
 (100\sup I_N^e)&=(113,116,120,117,120,120,123,123,126).
\end{aligned}
\]
Under the stated centroid bounds on \(s,d\), if every \(r_j\in I_N^e\),
the centered Rayleigh residual is at most
\[
 -|R_N|\sum_j\bigl(\eta_N(r_j-1)^2+b_N d_j\bigr),\qquad \eta_N,b_N>0
\]
with the constants below.
\end{numericallemma}
\begin{proof}
\emph{Step 1. The residual identity.}
Use \(W(Lx)=L(U(x)+Z(x))\) from the preceding construction. The weak
eigenvalue equation gives
\[
 \int_S DU:DW_j=\lambda\int_S U\cdot W_j\qquad(j=r,s),
\]
since \(\partial_n U=0\) on the outer edge and \(\partial_n U\cdot W_j=0\) on
both rays. Put \([A]=|R_N|^{-1}\int_S A\), \(H=L^{\mathsf T}L/\det L\),
\(w=\det L\). Take \(F,G\) from \eqref{finite-sector-mass-energy} for
\(U+Z\) relative to \(U\), with matrix \(H\). The residual in
\eqref{finite-sector-direct-comparison} at \(\nu=\lambda\) is
\[
 \begin{aligned}
 B&=w(\lambda f+G)-\lambda w^2F-\lambda f
     -h\mathcal A(r),\\
 h&=2[(\partial_xU_2)^2-(\partial_yU_1)^2]-\lambda[U_2^2-U_1^2].
 \end{aligned}
\]
Consequently the centered Rayleigh residual divided by \(|R_N|\)
is at most \(\sum_jB_j\).

\emph{Step 2. The scale and skew bounds.}
For \(16\le N\le24\), the caps on \((|s|,d)\), scale slope \(k\),
and comparison constants are
\[
\begin{aligned}
 (10^3S^e)&=(130,105,130,90,130,80,130,75,130),\\
 (10^6D^e)&=(580,430,400,260,400,260,400,230,400),\\
 (10^3 k)&=(49,41,36,32,30,26,25,21,21),\\
 (100a)&=(20,20,14,17,13,15,11,13,9),\\
 (10^6\eta)&=(200,200,100,120,30,60,20,35,10),\\
 (100b)&=(15,15,12,12,11,10,10,8,9).
\end{aligned}
\]
For \(25\le N\le32\) use the full boxes of Lemma~\ref{middle-tail-geometry} and
\[
\begin{aligned}
 (10^3k)&=(17,18,14,16,13,14,10,12),\qquad
 (10^3a)=(20,20,20,20,105,43,98,30),\\
 (10^8\eta)&=(2,2,2,2,100,2,100,2),\qquad
 (10^3b)=(10,10,10,10,25,10,25,10).
\end{aligned}
\]
With \(\xi=r-1\) the three comparisons are
\begin{equation}\label{common-ray-taylor-bound}
 \partial_s^2 B(r,s,0)<-a,\quad
 \frac{B(r,s,d)-B(r,s,0)}d<-b,\quad
 \frac{B(r,k\xi,0)}{\xi^2}
 +\frac1{2a}\left(\frac{\partial_sB(r,k\xi,0)}{\xi}\right)^2<-\eta .
\end{equation}

Let \(t=\tan\alpha\), \(D=1+r^2t^2\),
\(w_0=r(1+t^2)/D\), \(w=e^sw_0\). Use
\[
 H_{11}=r^{-1}+\frac{dD}{4r},\qquad H_{22}=r+\frac{dD}{4rt^2},\qquad
 H_{12}=\frac{D\sinh q}{2rt},\qquad d=2\cosh q-2 .
\]
At \(d=0\), set \(A=W_r+kW_s\) and
\(P=U+\xi A+(s-k\xi)W_s\). For even vector fields use
\[
 \begin{gathered}
 m(P)=[P_1^2]-N[(P-U)_1]^2,\quad
 n(P)=[P_2^2]-N[(P-U)_2]^2,\\
 b(P)=[(\partial_yP_1)^2],\quad c(P)=[(\partial_xP_2)^2],\quad e(P)=[|DP|^2],\\
 F=m(P)/r+rn(P),\qquad G=e(P)+(r^{-2}-1)b(P)+(r^2-1)c(P).
 \end{gathered}
\]
Use \(m(U)=f-n(U)\), \(e(U)=\lambda f\) and the weak identity for the
mixed gradients with \(A,W_s\). Then \(B_{ss}(r,s,0)\) is
\[
 w(\lambda f+G+2G_s+G_{ss})-\lambda w^2(4F+4F_s+F_{ss}).
\]
For the skew quotient freeze \(q/(2\sinh(q/2))\) and
\(H_{12}/(2\sinh(q/2))\) at their values for \(d\), and replace \(2\sinh(q/2)\)
by an indeterminate \(v\) and \(d\) by \(v^2\) in \(H,Z\).
Use \eqref{finite-scale-asinh}.
The quotient is precisely
\(w\sum_{j\ge1}d^{j-1}[v^{2j}](G-\lambda wF)\).

\emph{Step 3. The scale \(s=k(r-1)\).}
Use
\[
 \begin{gathered}
 w(r,k\xi)=1+\xi\left(\frac{1-rt^2}{D}
       +kw_0\frac{e^{k\xi}-1}{k\xi}\right),\qquad v=\frac{t\xi}{1+t^2r},\\
 \frac{2\arctan(rt)-2\alpha}{\sin(2\alpha)}
 =\xi-\xi^2\left(\frac{t^2}{1+t^2r}
  +\frac{(1+t^2)t^2\xi}{(1+t^2r)^3}\int_0^1\frac{y^2}{1+v^2y^2}\,dy\right).
 \end{gathered}
\]
Use the 19-term alternating expansion of that integral, with error
\(\le |v|^{38}/41\). For an explicit evaluation order in the last two quotients put
\(\omega=(w-1)/\xi\), \(L=(\mathcal A(r)-\xi)/\xi^2\), and use
the preceding continuous expressions for both. For fields \(Y,Z\) let
\[
\begin{aligned}
 m_{YZ}&=[Y_1Z_1],& n_{YZ}&=[Y_2Z_2],&
 b_{YZ}&=[\partial_yY_1\partial_yZ_1],&
 c_{YZ}&=[\partial_xY_2\partial_xZ_2].
\end{aligned}
\]
When neither argument is \(U\), subtract \(N[Y_i][Z_i]\) from the
two mass products. Subscripts \(0,A,S\) denote \(U,A,W_s\), respectively;
write \(g_{YZ}=[DY:DZ]\). The mass and gradient parts of the anchor
value are
\[
\begin{aligned}
 C_0&=-(1-\omega)^2/r-L,\\
 C_n&=L-2\omega-\omega^2r-C_0,\\
 C_m&=2(1-\omega)^2/r+2\omega(1-\omega),\\
 C_a&=-2+2\omega-4\omega r-2\omega^2\xi r,\\
 M&=C_0 f+C_n n_{00}+C_m m_{0A}+C_a n_{0A}
        -w^2(m_{AA}/r+rn_{AA}),\\
 E&=\left((2r+1-\omega(r+1))/r^2+2L\right)b_{00}
 +\left(1-2L+\omega(r+1)\right)c_{00}\\
 &\quad+2(r+1)w(c_{0A}-b_{0A}/r^2)+wg_{AA}
        +\xi(r+1)w(c_{AA}-b_{AA}/r^2),\\
 B(r,k\xi,0)/\xi^2&=E+\lambda M .
\end{aligned}
\]
For the scale derivative on the same line put
\[
\begin{aligned}
 F_0&=n_{00}-(f-n_{00})/r+2(m_{0A}/r+rn_{0A})+
                   \xi(m_{AA}/r+rn_{AA}),\\
 G_0&=2\lambda(m_{0A}+n_{0A})+\xi g_{AA}
       +(r+1)\{c_{00}+2\xi c_{0A}+\xi^2 c_{AA}\\
 &\hspace{41mm}-(b_{00}+2\xi b_{0A}+\xi^2 b_{AA})/r^2\},\\
 F_1&=2(n_{0S}-m_{0S}/r+m_{AS}/r+rn_{AS}),\\
 G_1&=2g_{AS}+2(r+1)(c_{0S}+\xi c_{AS}-(b_{0S}+\xi b_{AS})/r^2),\\
 B_s(r,k\xi,0)/\xi&=w\{G_0+G_1-\lambda[2(f+m_{0S}+n_{0S})\omega+w(2F_0+F_1)]\}.
\end{aligned}
\]

\emph{Step 4. Interval comparisons.}
Use the moments in \eqref{certified-sector-monomials}.
Center the moment bounds at \(V\):
\[
 |[(U-V)\cdot Y]|\le e_0[Y\cdot Y]^{1/2},\qquad
 |[(DU-DV):Z]|\le e_1[|Z|^2]^{1/2}.
\]
For a diagonal square use
\[
 [U_i^2]\in
 \big[(\sqrt{[V_i^2]}-e_0)_+^2,(\sqrt{[V_i^2]}+e_0)^2\big],
\]
and the analogous rule for derivative squares. Apply the latter also
to \(\partial_xU_1-\partial_yU_2\), with error \(\sqrt2e_1\).
An off-diagonal product has error at most
\(e_1(\|\partial_yV_1\|+\|\partial_xV_2\|)+e_1^2\)
in the normalized sector norm. On the anchor line, collect the entire
quadratic expression before bounding the error: for a symmetric constant
matrix \(C\), a polynomial field \(Y\), and \(H_0=V\) or \(DV\), use
\[
 \left|[(H_0+E)^{\mathsf T}C(H_0+E)+2Y\cdot(H_0+E)]
             -[H_0^{\mathsf T}CH_0+2Y\cdot H_0]\right|
 \le 2e\sqrt{[|CH_0+Y|^2]}+e^2\|C\| .
\]
Use \(e_0,e_1\) on the mass and gradient fields, respectively; the Gram
matrix of the polynomial summands evaluates \([|CH_0+Y|^2]\).

Bisect the longer of the \(r,s\) intervals, choosing \(r\) at ties,
for the first two tests while retaining the full \(d\) interval. For the
last test bisect \(r\) alone. On the \(10^{-420}\) grid at most
\(170,228,759\) terminal boxes are required, of maximum depths
\(10,10,13\), respectively. All the anchors lie within the scale
interval. Finally Taylor's inequality in \(s\) gives
\[
 B(r,s,d)\le B(r,k\xi,0)+\frac{B_s(r,k\xi,0)^2}{2a}-bd
 \overset{\eqref{common-ray-taylor-bound}}{\le}-\eta\xi^2-bd .
\]
\end{proof}

\paragraph{Exterior sectors for \(16\le N\le24\).}
Let \(J_N,S_N,Q_N\) be the full centroid bounds of Lemma~\ref{middle-base-geometry},
where \(10^3Q_N=(25,21,19,17,16,13,13,10,11)\) follows from
\(q^2\le d\). Use \(I_N\subseteq I_N^e\) with the following bounds and tolerances.
\[
\begin{aligned}
 (100\inf I_N)&=(88,88,82,86,82,85,80,84,78),\\
 (100\sup I_N)&=(113,112,120,113,120,114,123,114,126),\\
 (u_N^{\rm ext})&=(80,60,65,50,55,30,34,23,30),\\
 (v_N^{\rm ext})&=(1500,1100,1700,1000,1200,700,850,600,800).
\end{aligned}
\]
For this calculation take \(\lambda_*=l_N\) except for
\(\lambda_*=3.4585450,3.43190955\) at \(N=18,23\).
Normalize the entire polynomial field \(V_{\rm raw}+(r-1)W_r+sW_s+qW_q\)
by the factor in \eqref{middle-signed-polynomial}. In the forms \(F,G\)
above replace \(U\) by \(V\) and \(Z\) by the corresponding normalized correction.
Set
\[
 \mathcal E=w(2\lambda_*/N+G)-\lambda_*w^2 F-2\lambda_*/N-h\mathcal A(r),
 \qquad
 h=\left.\partial_r\{w(2\lambda_*/N+G)-\lambda_*w^2 F\}\right|_{(1,0,0)} .
\]
As for the exact trial, the global centered Rayleigh residual at
\(\lambda_*\) is bounded by \(|R_N|\sum_j\mathcal E_j\).
On \(J\times[-S,S]\times[-Q,Q]\) take an anchor
\(s_0\in\mathbb Q[r-1]\), \(\deg s_0\le8\), \(\sup_J|s_0|\le S\). Write
\[
 e(r)=\mathcal E(r,s_0(r),0),\quad p(r)=\mathcal E_s(r,s_0(r),0),\quad
 D_{q^2}\mathcal E=\frac{\mathcal E(r,s,q)-\mathcal E(r,s,0)}{q^2}.
\]
The quotients take their continuous values at zero:
\[
 1\le\frac{2\cosh q-2}{q^2}\le\frac{2\cosh Q-2}{Q^2},\qquad
 1\le\frac{\sinh q}{q}\le\frac{\sinh Q}{Q}.
\]
For \(1\in I\subseteq J\), the sign tests, with
\(\mathcal E_{ss}=\partial_s^2\mathcal E(r,s,0)\) for \(|s|\le S\), are
\[
 D_{q^2}\mathcal E<0,\quad \mathcal E_{ss}\le-a\ (r\in J),\qquad
 \mathcal E_{ss}\le-a_0,\qquad |p|<b,\ e''<-c \quad(r\in I).
\]
Unless specified, \(a_0=a\).
Taylor's inequality in \(s-s_0(r)\) and \(r-1\) yields
\begin{equation}\label{common-ray-polynomial-bound}
 \mathcal E\le
 \begin{cases}
e(1)+e'(1)^2/(2c)+b^2/(2a_0),&r\in I,\\
  e(r)+p(r)^2/(2a),&r\in J\setminus I .
 \end{cases}
\end{equation}
For \(B=[m-h,m+h]\subset J\setminus I\), the last bound is at most
\begin{equation}\label{midpoint-anchor-envelope}
 e(m)+h\sup_B|e'|+
 \frac{\bigl(|p(m)|+h\sup_B|p'|\bigr)^2}{2a}.
\end{equation}
Here \(e',p'\) are total derivatives. Rational integrals and Taylor
remainders give the following constants:
\[
\begin{aligned}
 (10^3a)&=(90,90,90,90,60,90,45,90,25),\\
 (10^6b)&=(100,100,100,100,70,50,65,50,85),\\
 (10^6c)&=(1000,1000,600,750,400,430,230,280,150)
\end{aligned}
\]
with \(a_0=a\), except \(a_0=0.10\) at \(N=24\).
\begin{numericallemma}[Polynomial sector comparison]\label{self-exterior-polys}
For \(16\le N\le24\), the centered polynomial trial below satisfies
\[
 \sup_{r\in I,\,|s|\le S_N,\,|q|\le Q_N}\mathcal E<10^{-8}u_N^{\rm ext},\qquad
 \sup_{r\in J\setminus I,\,|s|\le S_N,\,|q|\le Q_N}\mathcal E<-10^{-8}v_N^{\rm ext} .
\]
\end{numericallemma}
\begin{proof}
Put \(\alpha=\pi/N,\ X=x_1/\cos\alpha,\ Y=x_2/\sin\alpha\). Let
\(E_0=V(z)=\sum c_{ij}z^i\bar z^j\) be \(V_{\rm raw}\) of \eqref{middle-signed-polynomial}. Enumerate the following
basis in the indicated order, starting at zero:
\[
\begin{gathered}
 V,\quad z\sum_{i,j}(-1)^{(i-j-1)/N}c_{ij}X^{i+j-1}-V,\\
 (Y^2-X^2)X^iY^j(1,\,\mathrm i\sin\alpha)
          \quad(i+j\le6,\ (i,j)\text{ lexicographic}),\\
 X^i(X\cos\alpha+\mathrm iY\sin\alpha,\,
                  Y\cos\alpha+\mathrm iX\sin\alpha)\quad(0\le i\le5).
\end{gathered}
\]
Each parenthesis on the last two lines contributes two basis elements.
Identify complex values with column vectors in \(\mathbb R^2\).
For \(M=\left(\begin{smallmatrix}x&z\\z&y\end{smallmatrix}\right)\),
set \(\eta_0=0\) and \(\eta_i=|S|^{-1}\int_S E_i\) for \(i>0\).
The formal identity \(xy=1+z^2\) defines nine real symmetric matrices
\(T_1,\ldots,T_9\) by
\begin{align*}
 h&=(1,x^2,y^2,z^2,-\lambda_* wx,-\lambda_* wy,-2\lambda_* wz,2zy,2zx),\\
 \sum_{k=1}^9 h_k(T_k)_{ij}
  &=\frac1{|S|}\int_S\!\operatorname{tr}\big((DE_i)^\top M(DE_j)
                 \operatorname{adj}M\big)\\
  &\quad-\lambda_* w\left(\frac1{|S|}\int_S E_i^\top M E_j
               -\eta_i^\top M\eta_j\right),\qquad
 \frac1{|S|}\int_S X^aY^b=\frac{1+(-1)^b}{(a+b+2)(b+1)} .
\end{align*}
For \(t=\tan\alpha,\ R=r^{-1}+rt^2,\ d=2\cosh q-2\), use
\[
 x=r^{-1}+Rd/4,\quad y=r+Rd/(4t^2),\quad z=R\sinh q/(2t),\qquad
 w=e^s r(1+t^2)/(1+r^2t^2).
\]
The following rational recursion fixes the three correction fields.
In this recursion only, evaluate each exponential at \(v\) by
\(\sum_{j=0}^{64}v^j/j!\). Write
\(D(v)=2^{-128}\lfloor2^{128}v\rfloor\),
\(C(v)=D(v+2^{-129})\), and form
\[
 H_{ij}=D\!\left(\sum_{k=1}^9C(h_k)C((T_k)_{ij})\right).
\]
Given fixed coordinates \(c_j\), \(j\in F\), set
\(f_i=-D(\sum_{j\in F}H_{ij}c_j)\). Order free indices increasingly.
For \(i\ge j\), compute in row order
\[
 A_{ij}=H_{ij}-D\!\left(\sum_{k<j}A_{ik}L_{jk}\right),\qquad
 L_{ij}=D(A_{ij}/A_{jj})\ (i>j).
\]
Compute \(y_i=f_i-D(\sum_{k<i}L_{ik}y_k)\) in increasing order and
\(c_i=D(y_i/A_{ii})-D(\sum_{k>i}L_{ki}c_k)\) in decreasing order.
First use \((r,s,q)=(0.84,0,0.017)\), fixing \(c_0=1\) and
\(c_{58+2j}=0\) for \(0\le j\le5\); set
\(\gamma_j=D(1000c_{59+2j}/17)\).
Subsequently fix
\[
 c_0=1,\qquad c_{58+2j}=D(s\gamma_j),\qquad
 c_{59+2j}=D(q\gamma_j)\qquad(0\le j\le5).
\]
The recursion now defines \(c(r,s,q)\). Let \(b=(1,0,0)\),
\(\epsilon=10^{-4}\) and \(B(v)=10^{-5}\lfloor10^5v+1/2\rfloor\).
For \(a=r,s,q\) take the coordinatewise rounded difference
\[
 c^a=B\!\left(\frac{c(b+\epsilon e_a)-c(b-\epsilon e_a)}{2\epsilon}\right).
\]
Under reflection \(U(X,Y)\mapsto\overline{U(X,-Y)}\), retain only the
even basis elements for \(a=r,s\) and the odd ones for \(a=q\).
Set \(c^s_{58+2j}=c^q_{59+2j}=B(\gamma_j)\).
With \(W_a=\sum_i c_i^aE_i\), take
\[
 U=\left(\frac{2|S|}{\int_S|V|^2}\right)^{1/2}
           \bigl(V+(r-1)W_r+sW_s+qW_q\bigr).
\]
The two traces of the correction involve \((s\pm q)B(\gamma_j)\);
thus adjacent sectors have equal transported traces. In \(\mathcal E\)
all sector integrals have weight \(1/(N|S|)\), giving reference mass \(2/N\).

Write \(J=[a,b]\). For \(0\le i\le32\) let \(r_i=a+i(b-a)/32\).
Starting with \([-0.11,0.11]\), bisect \(32\) times at fixed \(r_i\).
At the midpoint \(s\) retain the right half if
\(\mathcal E_s(r_i,s,0)>0\), and the left half otherwise.
Let \(y_i\) be the final midpoint. With \(P_{ij}=(r_i-1)^j\)
for \(0\le j\le8\), take
\[
 (t_0,\ldots,t_8)^\top=(P^\top P)^{-1}P^\top y,\qquad
 s_0(r)=10^{-11}\sum_{j=0}^8
    \left\lfloor10^{11}t_j+\tfrac12\right\rfloor(r-1)^j .
\]
Apply \eqref{common-ray-polynomial-bound} on \(I\) and
\eqref{midpoint-anchor-envelope} on \(J\setminus I\).
On a subinterval \([m-h,m+h]\), one may also use
\[
 \sup e''\le e''(m)+h\sup|e'''|,\qquad
 \mathcal A'''(r)=-\frac{2t^2(1+t^2)(1-3t^2r^2)}{(1+t^2r^2)^3}.
\]
Factor \(q^2\) out of the skew difference before enclosing
\((2\cosh q-2)/q^2\) and \(\sinh q/q\) by the bounds above.
On failure bisect at the midpoint in \(r\) if \(\Delta r\ge2\Delta s\),
and in \(s\) otherwise. On \(|s|\le0.13\), the exponential
through degree \(22\) has error at most \(2(1/5)^{23}/23!\).
For the angular term use
\[
 \mathcal A(r)=\frac{1+t^2}{t}\arctan\!\left(\frac{(r-1)t}{1+rt^2}\right),
 \qquad
 \mathcal A'(r)=\frac{1+t^2}{1+r^2t^2}.
\]
Enclose the arctangent by the alternating sums of degrees \(39,41\).
Use the direct second-derivative enclosure for \(N=16,\ldots,19\)
and the enclosure using the third derivative for \(N=20,\ldots,24\).
Every branch terminates within fifteen bisections.
\end{proof}

\begin{proposition}\label{large-middle-comparison}
For \(16\le N\le32\), among convex polygons of fixed positive area
with at most \(N\) sides, \(R_N\) uniquely maximizes \(\mu_1\) up to congruence.
\end{proposition}
\begin{proof}
Normalize \(|P|=|R_N|\) and suppose \(\mu_1(P)\ge\mu_1(R_N)\).
Then \(\Phi(P)\ge\Phi(R_N)>T_N\), so the centroid bounds force \(N\) sides.
If \(16\le N\le24\) and some \(r_j\notin I_N\), then
\[
 10^8\sum_j\mathcal E_j < (N-1)u_N^{\rm ext}-v_N^{\rm ext}
 =-(300,140,595,100,155,100,136,94,110)_N<0.
\]
This contradicts \(\mu_1(P)\ge\mu_1(R_N)\) by
Lemma~\ref{central-exterior-comparison}(2).
In the remaining cases Lemma~\ref{middle-common-ray} applies, since
\(I_N\subseteq I_N^e\) for \(16\le N\le24\) and \(I_N^e=J_N\) for
\(25\le N\le32\). Its residual satisfies
Lemma~\ref{central-exterior-comparison}(1) with
\(\epsilon=\min(\eta_N,b_N)>0\).
\end{proof}

\section{Uniform estimates for convex polygons}\label{part:uniform}

For \(N\ge33\), we localize convex competitors and bound the sector
residuals uniformly in \(N\). The linear angle and area terms cancel
when the sector inequalities are summed.

\subsection{Localization of a competing polygon}\label{sec:effective-fan-localization}

\begin{proposition}[An explicit range for competitive fans]
\label{prop:effective-fan-localization}
Let \(N\ge30\) and let \(P\) be a convex polygon with at most \(N\) sides.
If \(\Phi(P)\ge\Phi(R_N)\), then \(P\) has \(N\) sides and
\[
 0\le\log(\beta_N/\beta)\le148N^{-5},\qquad
 \beta=\frac{\pi|f'(0)|^2}{|P|},\qquad
 r_*^2=\frac{2|P|}{N\sin(2\pi/N)} .
\]
Here \(f\) is the dipole-balancing map and \(\beta_N\) is the regular value.
For \(r_j=|p_j-f(0)|\) and
\(2\pi q_j/N=\angle(p_j-f(0),p_{j+1}-f(0))\),
\[
 \begin{split}
 (q_j-1)^2(q_j^2+2q_j+3)&\le \frac{148}{4.8N},\\
 \sum_j\log^2(r_{j+1}/r_j)&\le 1898N^{-5},\\
 \left|\log(r_jr_{j+1}/r_*^2)\right|&\le 41N^{-2}.
 \end{split}
\]
For the vertices \(p_j\), write
\[
 \widetilde r_j=\left|p_j-\frac1{|P|}\int_P y\,dy\right|,\qquad
 \widetilde q_j=\frac N{2\pi}\angle\!\left(
 p_j-\frac1{|P|}\int_P y\,dy,\,p_{j+1}-\frac1{|P|}\int_P y\,dy\right).
\]
Then
\[
 |\widetilde q_j-q_j|\le24N^{-3},\qquad
 \left|\log\frac{\widetilde r_{j+1}r_j}{\widetilde r_jr_{j+1}}\right|\le145N^{-4},
 \qquad \left|\log\frac{\widetilde r_j\widetilde r_{j+1}}{r_*^2}\right|
 \le42N^{-2}.
\]
\end{proposition}

\begin{proof}
Use \(\varphi,d_j,x=(j'_{1,1})^2\) from
Section~\ref{sec:regular-sc-certificate}. The alternating series gives
\(3.389<x<3.391\). By the gamma series,
\[
 \frac{4.8}{N^3}<\frac4{N^3}\sum_{j=1}^{30}j^{-3}\le -\log\beta_N
 \le\frac{4(1.203)}{N^3}+\frac{12(1.04)}{N^5}
       +\frac{2.08}{7}\frac{(2/N)^7}{1-2/N}.
\]
For \(0\le u<1\), \(u\le e^u-1\le u/(1-u)\). At \(u=-\log\beta_N\)
and \(N\ge30\), this gives
\begin{equation}\label{regular-conformal-defect-bounds}
 4.8N^{-3}\le\beta_N^{-1}-1\le4.84N^{-3}.
\end{equation}

\smallskip
\noindent\emph{Step 1. Lower bound for the regular polygon.}
\begin{equation}
 \frac{\pi(x-1)}{\Phi(R_N)}
 \le d_1+(1-d_1)(1-\beta_N)+29.3N^{-5}.
 \label{effective-fan-reg}
\end{equation}
For \(w=\beta_N|1-z^N|^{-4/N}\), use the energy norm
\(\|v\|_a^2=(2\pi)^{-1}\int_{\mathbb D}|Dv|^2\,dA\).
Split residue one as \(\mathbb C U\oplus L\oplus H\),
where \(U=\varphi e^{i\theta}\), \(L=U^{\perp_a}\) in mode one,
and \(H\) contains the modes \(|k|\ge N-1\). Then
\[
 \begin{gathered}
 \|U\|_a^2=(x-1)/2,\\
 \widehat w_k(r)=\frac1{2\pi}\int_{-\pi}^{\pi}
 w(re^{i\theta})e^{-ik\theta}\,d\theta,\\
 \widehat w_0(r)=\beta_N\sum_{j\ge0}\left(\frac{(2/N)_j}{j!}\right)^2r^{2Nj},
 \qquad 2\int_0^1\widehat w_0(r)r\,dr=1.
 \end{gathered}
\]
For \(\rho=(\pi(x-1))^{-1}\int_{\mathbb D}w|U|^2\,dA\),
the identity \(2\int_0^1\varphi^2r\,dr=(x-1)/x\) and
\(\widehat w_0',\varphi'\ge0\) give
\[
 \rho-\frac1x=\frac2{x-1}\int_0^1\!\int_0^1
 (\widehat w_0(r)-\widehat w_0(t))
 (\varphi(r)^2-\varphi(t)^2)rt\,dr\,dt\ge0 .
\]
In the energy inner product, define
\[
 a(Ku,v)=\frac1{2\pi}\int_{\mathbb D}wu\bar v\,dA .
\]
The complementary blocks satisfy
\[
 \begin{gathered}
 \|K_{LL}\|\le 100/2809+\kappa,\qquad
 \|K_{HH}\|\le (N-1)^{-2}+\kappa/(N-1),\\
 \kappa=\frac{4}{N(N-4)}\left(1+\frac{N}{N-4}\right)\le1/90.
 \end{gathered}
\]
By Section~\ref{sec:regular-sc-certificate},
\(\rho\ge1/x>1/6+\kappa(1+1/(N+2))\) selects the dipole residues.
For \(f\in L,\,g\in H\), \(\int_{\mathbb D}f\bar g\,dA=0\).
Using \(\big||1-z^N|^{-4/N}-1\big|\le(1-r^N)^{-4/N}-1\)
in the radial estimate and applying Cauchy--Schwarz in the angle gives
\[
 \begin{gathered}
 \|K_{LH}\|\le\frac{\kappa}{\sqrt{N-1}},\\
 \frac{2\kappa uv}{\sqrt{N-1}}\le0.004u^2+0.0011v^2,\\
 \rho-\|K_{LL}\|-0.004\ge0.244,\qquad
 \rho-\|K_{HH}\|-0.0011\ge0.292 .
 \end{gathered}
\]
Here \(\kappa^2/(N-1)\le1/(90^2\cdot29)<0.004\cdot0.0011\).

Write \(R_J=\|\operatorname{proj}_J KU\|_a\), \(J=L,H\), with the
orthogonal projections taken in \(a\).
The terms \(j\ge1\) of \(\widehat w_0(r)\) have sum at most \(\frac4{(N-2)^2}\sum_{l\ge1}l^{4/N-2}r^{2Nl}\).
By \(|U|\le1\) and the radial spectral gap,
\[
 R_L\le\frac{1}{\sqrt{28.09}}
 \left\|\frac4{(N-2)^2}\sum_{l\ge1}l^{4/N-2}r^{2Nl}\right\|_{L^2(r\,dr)}
 \le\frac{2\zeta(5/2-4/N)}{(N-2)^2\sqrt{28.09 N}}.
\]
For \(s>1\), \(\zeta(s)\le1+2^{-s}+2^{1-s}/(s-1)\); thus
\[
 \begin{gathered}
 \zeta(7/3)<1.5,\qquad \zeta(11/6)<2,\qquad \zeta(14/3)-1<0.062,\\
 2N/(N-2)\le15/7,\qquad N^5R_L^2\le0.425 .
 \end{gathered}
\]
For \(m=1\) the positive binomial expansion gives
\(|\widehat w_N|\le\big(2/N+4\zeta(2-4/N)/(N-2)^2\big)r^N\).
Together with the envelope of Section~\ref{sec:regular-sc-certificate}, this gives
\[
 N|\widehat w_{mN}(r)|\le r^{mN}
 \begin{cases}2.307,&m=1,\\2.473\,m^{4/N-1},&m\ge2.\end{cases}
\]
For \(p=mN\ge30\), \(k=p\pm1\), the monomial inverse-energy bound is
\(3/(8p^3)\). Positivity of the radial Green kernel and
\(\zeta(5-8/N)-1\le\zeta(14/3)-1\) give
\[
 N^5R_H^2\le\frac34\bigl(2.307^2+2.473^2\cdot0.062\bigr).
\]
Completing squares on \(L\oplus H\),
\[
 \begin{aligned}
 (x-1)\left(\frac{\pi}{\Phi(R_N)}-\rho\right)
 &\le2\left(\frac{R_L^2}{0.244}+\frac{R_H^2}{0.292}\right)\\
 &\le\frac2{N^5}\left(\frac{0.425}{0.244}
 +\frac{3(2.307^2+2.473^2\cdot0.062)}{4\cdot0.292}\right)
 <32.772N^{-5}.
 \end{aligned}
\]
For \(g=\varphi/r\), the radial equation gives
\[
 \begin{gathered}
 -g'(r)=xr^{-3}\int_0^rt^3g(t)\,dt,\qquad
 1\le g(r)\le g(0)\le(1-x/8)^{-1},\\
 1+0.42(1-r^2)\le g(r)\le1+0.74(1-r^2) .
 \end{gathered}
\]
On \([0.85,1]\),
\[
 -\varphi''(r)=x\varphi(r)+g'(r)
 \ge xr[0.985-(0.42-0.74/6)r^2]\ge2.2
 \quad\Longrightarrow\quad
 1-\varphi(r)^2\ge2.035(1-r)^2 .
\]
The last implication uses \(\varphi(1)=1,\ \varphi'(1)=0\).
At \(a=0.85\), the integral fraction
\[
 \frac{\int_0^a r^{2N+1}(1-r)^2\,dr}{\int_0^1 r^{2N+1}(1-r)^2\,dr}
 =a^{2N+2}\left[
  \frac{(1-a)^2(2N+3)(2N+4)}2+a(1-a)(2N+4)+a^2\right]<0.005
\]
decreases term by term for \(N\ge30\). Hence \(1-d_{N+1}\ge0.90N^{-2}\).
The normalized mass identity is
\[
 \begin{aligned}
 d_1+(1-d_1)(1-\beta_N)-(x-1)\rho
 &=\beta_N\sum_{m\ge1}\frac{1-d_{mN+1}}{mN+1}
                \left(\frac{(2/N)_m}{m!}\right)^2\\
 &\ge\frac{\beta_N(0.90N^{-2})(2/N)^2}{N+1}>3.48N^{-5}.
 \end{aligned}
\]
Subtract the last estimate from the bound for
\((x-1)(\pi/\Phi(R_N)-\rho)\) and use \(32.772-3.48<29.3\).

\smallskip
\noindent\emph{Step 2. The conformal radius of a competitor.}
Use \eqref{large-N-schur} at \(J=2N+1\), with \(S=\beta^{-1}-1\).
The lower bound for \(\varphi(r)/r\) gives \(d_2>0.8184\) and
\(d_3>0.8848\). By \eqref{large-N-trace} and \eqref{effective-fan-reg},
\[
 \begin{aligned}
 0.113(1-\beta)&\le(d_2-d_1)(1-\beta)
 \le \frac{\pi(x-1)}{\Phi(R_N)}-d_1\\
 &\le \frac{1-\beta_N}{x}+29.3N^{-5}
 \le\left(\frac{4.84}{3.389}+\frac{29.3}{30^2}\right)N^{-3}
 <1.462N^{-3}.
 \end{aligned}
\]
Thus \(S\le13N^{-3}\) and \(\beta>0.9995\). Also
\[
 \frac{\beta\Phi(P)}{\pi}
 \ge\frac{0.9995x}{1+1.462/(d_1N^3)}>3.38 .
\]
Balance and the second harmonic in \eqref{large-N-trace} give
\[
 |b_1|\le2S,\qquad \frac{\beta d_2}{2}|b_2|
 \le\frac{\pi(x-1)}{\Phi(R_N)}-d_1+\beta S
 \le1.462N^{-3}+\beta S .
\]
\(\beta\le\beta_N\) by Lemma~\ref{fan-shape-bound}, so
\(S\ge\beta_N^{-1}-1\ge4.8N^{-3}\). Hence \(|b_2|\le3.2S\) and
\[
 \sum_{j=1}^2(1-d_{j+1})\frac{|b_j|^2}{j+1}
 \le\left(\frac{4(1-0.8184)}2+
             \frac{3.2^2(1-0.8848)}3\right)S^2<0.8S^2.
\]

From \(\varphi''=-x\varphi-g'\ge-x(1-r/4)\) and
\(\varphi'(1)=0\), integration gives
\[
 1-\varphi(r)^2
 \le 2x\int_r^1(t-r)(1-t/4)\,dt
 =\frac{x(10-r)}{12}(1-r)^2 .
\]
For \(j\ge3\),
\[
 1-d_{j+1}\le\frac{5x}{3(2j+3)(2j+4)}
 <\frac{1.57}{(j+2)^2}.
\]
For \(3\le j\le2N\), the two terms \(m=j\pm1\)
in \eqref{large-N-schur} have combined quadratic coefficient at least
\[
 \frac{\beta^2\Phi(P)}\pi\frac{j+1}2
 \left(\frac1{(j-1)(j+1)^2}+\frac1{(j+1)(j+2)^2}\right)
 \ge\frac{\beta^2\Phi(P)}{\pi(j+2)^2}
\]
on \(|b_j|^2/(j+1)\) before subtracting the \(S\)-remainders, by
\eqref{large-N-radial-moment}. Use
\[
 (a-S)_+^2\ge0.47a^2-\frac{47}{53}S^2,\qquad
 0.47\cdot3.38>1.57,\qquad
 \frac{47\beta\Phi(P)}{53\pi}\le\frac{47x}{53}<3.01.
\]
The quadratic term absorbs the trace loss for \(3\le j\le2N\).
The remainders per mode and the tail are bounded by
\[
 \frac{3.01\beta S^2}{m}\quad(2\le m\le2N+1),
 \qquad
 \beta\sum_{j>2N}(1-d_{j+1})\frac{|b_j|^2}{j+1}
 \le\frac{0.393\beta S}{N^2}.
\]
Thus \eqref{effective-fan-reg} gives
\[
 \beta_N-\beta\le x\left[
 29.3N^{-5}+\beta\bigl(0.8+3.01\log(2N+1)\bigr)S^2
              +0.393\beta S/N^2
 \right].
\]
Since \(N^{-1}\log(2N+1)\) decreases and \(\log61<4.112\),
\[
 S\le aN^{-3}\ \Longrightarrow\
 N^5(\beta_N-\beta)\le
 3.391\left[29.3+\frac{a^2}{30}(0.8+3.01\cdot4.112)+0.393a\right].
\]
With \(S=\beta_N^{-1}-1+(\beta_N-\beta)/(\beta\beta_N)\),
\[
 \begin{aligned}
 S\le13N^{-3}
 &\ \Longrightarrow\ N^5(\beta_N-\beta)<368.40<369
 \ \Longrightarrow\ S\le5.26N^{-3},\\
 S\le5.26N^{-3}
 &\ \Longrightarrow\ N^5(\beta_N-\beta)<147.576\\
 &\ \Longrightarrow\ \log(\beta_N/\beta)
 \le(\beta_N-\beta)/\beta\le148N^{-5}.
 \end{aligned}
\]

\smallskip\noindent\emph{Step 3. Sector angles and radii.}
By \eqref{fan-deficit}, at most \(N-1\) sides would give
\[
 148N^{-5}\ge\log(\beta_N/\beta)
 \ge4\zeta(3)((N-1)^{-3}-N^{-3})>14.4N^{-4}>148N^{-5}.
\]
Since \(\sum_{j=1}^N q_j=N\), \eqref{fan-deficit} at \(a_j=q_j/N\) gives
\[
 \sum_j(q_j-1)^2(q_j^2+2q_j+3)=\sum_j(q_j^4-1)
 \le\frac{148}{4.8N},\qquad 0.509<q_j<1.368.
\]
Indeed \((q^4-4q+3)'=4(q^3-1)\), and the values at \(0.509,1.368\)
exceed \(148/(4.8\cdot30)\). If \(d_j\) is half the difference of the
two base angles divided by \(\pi\), the first positive term in
\eqref{fan-angular-remainders} gives
\[
 32\sum_j a_j^2d_j^2
 \le\sum_j\frac{2a_j^2d_j^2}{(1/4-a_j^2)^2}\le148N^{-5},
 \qquad |d_j|<4.31N^{-3/2}.
\]
Use \(\tan(\pi a_j)\le1.01\pi a_j\),
\(\tan(\pi|d_j|)\le1.003\pi|d_j|\), and
\[
 \tan(\pi a_j)\tan(\pi|d_j|)
 \le (1.01)(1.003)\pi^2(1.368)(4.31)N^{-5/2}<0.013 .
\]
The sine rule and \(\operatorname{arctanh}v\le1.001v\) give
\[
 \left|\log\frac{r_{j+1}}{r_j}\right|
 =2\operatorname{arctanh}\!\left(\tan(\pi a_j)\tan(\pi|d_j|)\right)
 \le2.04\pi^2 a_j|d_j|.
\]
Hence \(\sum_j\log^2(r_{j+1}/r_j)
 \le (2.04\pi^2)^2(148/32)N^{-5}<1898N^{-5}\).

For sector areas \(|P|b_j\), \eqref{sector-log-sum-bound} gives
\(\mathcal L(a_j,b_j)\le148N^{-5}\). If \(0.509/N<a<1.368/N\),
\[
 \frac{d^2}{dt^2}\mathcal L(a,ae^t)
 =\frac{a(1-a)e^t}{(1-ae^t)^2}>
 \frac{0.509(1-1.368/30)e^{-0.04}}N>\frac{0.466}N
 \quad(|t|\le0.04).
\]
The value and first derivative vanish at \(t=0\).
Since \(0.233(0.04)^2/N>148N^{-5}\), convexity gives
\[
 |\log(b_j/a_j)|\le\sqrt{148/0.233}\,N^{-2}<25.5N^{-2}.
\]
By the area identity \(2|P|b_j=r_jr_{j+1}\sin(2\pi a_j)\) and
\(0\le-\log(\sin t/t)\le t^2/5\) for \(0\le t\le0.3\),
\[
 \begin{aligned}
 \left|\log\frac{r_jr_{j+1}}{r_*^2}\right|
 &=\left|\log\frac{b_j}{a_j}
 +\log\frac{\sin(2\pi/N)}{2\pi/N}
 -\log\frac{\sin(2\pi a_j)}{2\pi a_j}\right|\\
 &\le\left|\log\frac{b_j}{a_j}\right|
 +\frac{(2\pi\cdot1.368)^2}{5N^2}<41N^{-2}.
 \end{aligned}
\]

\smallskip\noindent\emph{Step 4. Transfer to the area centroid.}
With \(|P|=1\), Proposition~\ref{bessel-centroid} gives
\[
 \left|\int_P y\,dy-f(0)\right|
 \le5.26N^{-3}\left(\sqrt{3/2}\,0.565+0.855\right)<8.2N^{-3}.
\]
Here \(\sqrt{\beta/\pi}<0.565\) and
\(\sqrt{8(35/36)/\Phi(P)}<0.855\) by \eqref{effective-fan-reg}.
Since
\[
 \left|\log\frac{r_j}{r_*}\right|
 \le\frac12\left(41N^{-2}+\sqrt{1898}N^{-5/2}\right),\qquad
 r_*^2=\frac{2}{N\sin(2\pi/N)},
\]
we have \(0.548<r_j<0.59\). The side length is at most
\(\max(r_j,r_{j+1})(2\pi q_j/N+|\log(r_{j+1}/r_j)|)\).
On the segment from \(f(0)\) to \(\int_P y\,dy\),
\[
 \left|\frac{d}{dz}\log\frac{p_{j+1}-z}{p_j-z}\right|
 \le\frac{0.59(2\pi\cdot1.368/N+\sqrt{1898}\,N^{-5/2})}
 {(0.548-8.2N^{-3})^2}<\frac{17.68}N .
\]
The real and imaginary parts give
\[
 \left|\log\frac{\widetilde r_{j+1}r_j}{\widetilde r_jr_{j+1}}\right|
 \le17.68(8.2)N^{-4}<145N^{-4},\qquad
 |\widetilde q_j-q_j|\le\frac{17.68(8.2)}{2\pi}N^{-3}<24N^{-3}.
\]
Also,
\[
 \left|\log\frac{\widetilde r_j\widetilde r_{j+1}}{r_*^2}\right|
 \le41N^{-2}+\frac{2(8.2N^{-3})}{0.548-8.2N^{-3}}\le42N^{-2}.
\]
\end{proof}

\subsection{Comparison for at least thirty-three sides}\label{sec:explicit-neumann}

\begin{proposition}\label{effective-global-neumann}
For \(N\geq33\), every convex polygon \(P\) with at most \(N\)
sides satisfies \(\Phi(P)\leq\Phi(R_N)\), with equality only for
regular \(N\)-gons.
\end{proposition}

Take \(R_N\) of circumradius one and normalize
\(|P|=|R_N|\), \(\int_P y\,dy=0\). Set
\[
 h=N^{-1},\qquad \alpha=\pi h,\qquad
 x=(j'_{1,1})^2,\qquad \lambda=\mu_1(R_N),\qquad
 g_*=\sqrt{\frac{2}{1-x^{-1}}}.
\]
For the first equivariant pair \(U\) of Lemma~\ref{affine-all-N}, take
\[
 \int_{R_N}|U|^2=2|R_N|,\qquad
 B(y)=g_*\frac{J_1(\sqrt\lambda\,|y|)}{J_1(\sqrt x)}\frac{y}{|y|},
 \quad B(0)=0,\quad E=U-B .
\]
On the reference sector \(S\) and its side \(I\), use
\[
 \|f\|_S^2=\frac{1}{|R_N|}\int_S|f(X)|^2\,dX,\qquad
 \|f\|_I^2=\frac{1}{|R_N|}\int_I|f|^2\,ds .
\]
For a sector of opening \(2\beta\) and endpoint radii \(R_-,R_+\),
\[
 r=\frac{\tan\beta}{\tan\alpha},\qquad
 s=\log(R_-R_+),\qquad \rho=\log(R_+/R_-).
\]
With \(Lv_\pm=R_\pm(\cos\beta,\pm\sin\beta)\) and \(w=\det L\),
define \(W(y)=B(y)+LE(L^{-1}y)\) on \(LS\). Equivariance gives
\[
 E(tv_\pm)=e(t)v_\pm,\quad W(tLv_\pm)=B(tLv_\pm)+e(t)Lv_\pm
 \quad\Longrightarrow\quad W\in H^1(P;\mathbb R^2).
\]
By \eqref{fan-sector-integral} and \(\int_P y\,dy=0\),
\[
 \int_P W(y)\,dy=\sum_S w\int_S[B(LX)-LB(X)]\,dX.
\]
Since \((-\Delta-\lambda)B=0\) and on the rays \(B,LE\parallel e_r\),
\(\partial_nB\parallel e_\theta\), Green's formula gives
\begin{equation}\label{physical-total-residual}
 \frac{1}{|R_N|}\left[
 \int_P(|DW|^2-\lambda |W|^2)\,dy
 +\frac{\lambda}{|P|}\left|\int_P W(y)\,dy\right|^2\right]
 \leq \sum_S F(L),
\end{equation}
where
\begin{align*}
 F(L)={}&\frac1{|R_N|}\int_{LI}
           \partial_n B(y)\cdot\bigl(B(y)+2LE(L^{-1}y)\bigr)\,ds\\
     &+\frac{w}{|R_N|}\int_S
     \bigl(|L(DE)L^{-1}|^2-\lambda |LE|^2\bigr)\,dX\\
     &+\frac{\lambda wN}{|R_N|^2}
       \left|\int_S[B(LX)-LB(X)]\,dX\right|^2 .
\end{align*}
Write \(B(Y)=b(|Y|^2)Y\). For \(X=(x_1,x_2)\in I\), the side normal
satisfies \(n\,ds_Y=wL^{-\mathsf T}e_1\,ds_X\). Thus the flux linear in
\(E\), before its factor two, is
\begin{equation}\label{uniform-transported-flux}
 \begin{aligned}
 &\frac1{|R_N|}\int_{LI}\partial_n B(Y)\cdot LE(L^{-1}Y)\,ds_Y\\
 &\qquad=\frac{w}{|R_N|}\int_I
 \big[b(|LX|^2)E_1+2x_1 b'(|LX|^2)X^{\mathsf T}L^{\mathsf T}LE\big]\,ds_X .
 \end{aligned}
\end{equation}
At \(L=\mathrm{Id}\), the same flux equals
\begin{equation}\label{uniform-reference-flux}
 \frac1{|R_N|}\int_I\partial_{x_1} B\cdot E\,ds
 =\frac1{|R_N|}\int_S(DB:DE-\lambda B\cdot E)\,dX
 =-\|DE\|_S^2+\lambda\|E\|_S^2 .
\end{equation}
The first equality uses the radial traces of \(E\), and the second uses
the weak equation for \(U=B+E\) and equivariance.
Center the sector residual by the angular and area identities:
\begin{equation}\label{physical-sector-residual}
 \begin{split}
 \mathcal B(r,s,\rho)={}&F(L)+2\lambda h(w-1)\\
 &-\frac{2(\beta-\alpha)}{\sin(2\alpha)}
          \left.\partial_r\!\left[F(L)+2\lambda h w\right]\right|_{(1,0,0)} .
 \end{split}
\end{equation}
The normalization, equivariance and reflection give
\[
 \begin{gathered}
 2\partial_r\arctan(r\tan\alpha)|_{r=1}=\sin(2\alpha),
 \qquad \mathcal B_r(1,0,0)=0,\\
 F(\mathrm{Id})=\frac{\int_{R_N}(|DU|^2-\lambda|U|^2)}{N|R_N|}=0,
 \qquad \mathcal B(r,s,-\rho)=\mathcal B(r,s,\rho),\\
 \sum_S(w-1)=\sum_S(2\beta-2\alpha)=0,\qquad
 \sum_S\mathcal B(r,s,\rho)=\sum_SF(L).
 \end{gathered}
\]
On \(LI\), with \(g(R)=B(R,0)\cdot e_1\), \(e_r=(\cos\theta,\sin\theta)\),
\[
 (n\cdot e_r)\,ds=R\,d\theta,\qquad
 \partial_n B\cdot B\,ds=Rg(R)g'(R)\,d\theta .
\]
Including the area correction, the pure radial term (modulo terms
independent of \((s,\rho)\)) is \(H(R)=Rg(R)g'(R)+\lambda R^2\).

In bisector coordinates the vertex correspondence
\((\cos\alpha,\pm\sin\alpha)\mapsto e^{(s\pm\rho)/2}(\cos\beta,\pm\sin\beta)\)
gives
\begin{equation}\label{centered-physical-map}
 \begin{aligned}
 L&=e^{s/2}\frac{\cos\beta}{\cos\alpha}
 \begin{pmatrix}
  \cosh(\rho/2)&\sinh(\rho/2)\cot\alpha\\
  r\tan\alpha\,\sinh(\rho/2)&r\cosh(\rho/2)
 \end{pmatrix},\\[-2pt]
 w&=e^s r\,\frac{\cos^2\beta}{\cos^2\alpha}.
 \end{aligned}
\end{equation}
At angle \(u\beta\), \(|u|\le1\), the side radius and its derivatives are
\begin{equation}\label{centered-skew-radius}
 \begin{gathered}
 R=\frac{\sqrt{R_-R_+}\cos\beta}
   {\cos(u\beta)\cosh(\rho/2)-\cot\beta\sin(u\beta)\sinh(\rho/2)},\\
 (\log R)_\rho
 =\frac{\tan(u\beta)-\tan\beta\tanh(\rho/2)}
        {2[\tan\beta-\tan(u\beta)\tanh(\rho/2)]},\qquad
 (\log R)_{\rho\rho}=((\log R)_\rho)^2-\frac14 .
 \end{gathered}
\end{equation}
The side radius also satisfies the exact area identity
\begin{equation}\label{uniform-side-area}
 \frac1{2\beta}\int_{-\beta}^{\beta} R(\theta)^2\,d\theta
 =\frac{e^s\sin(2\beta)}{2\beta}.
\end{equation}
Indeed \(\int_{-\beta}^{\beta}R(\theta)^2\,d\theta
 =R_-R_+\sin(2\beta)\), with \(u=\theta/\beta\) in
\eqref{centered-skew-radius}.

For \(k=\beta/\alpha\), write
\begin{equation}\label{sec4-angular-anchor}
 s_*(k)=\log\frac{\beta\sin(2\alpha)}{\alpha\sin(2\beta)},\qquad
 \mathcal T[f](k)=f(k)-f(1)-f'(1)(k-1).
\end{equation}
Set \(c=|R_N|/\pi=\sin(2\alpha)/(2\alpha)\). At \((s,\rho)=(s_*(k),0)\),
with \(X=(\cos\alpha,u\sin\alpha)\), \(|u|\le1\), \(Y=LX\) and
\(\theta=\arg Y\), the identities are
\begin{equation}\label{uniform-anchored-side}
 \begin{gathered}
 w=k,\qquad Y_1^2=c\beta\cot\beta,\qquad Y_1Y_2=c\beta u,\\
 p=|Y|^2=c(\beta\cot\beta+\beta\tan\beta\,u^2)
  =c\beta\cot\beta\sec^2\theta,\qquad
 \frac{d\theta}{du}=\frac{c\beta}{p}.
 \end{gathered}
\end{equation}
At the anchor, \((2\beta)^{-1}\int_{-\beta}^{\beta}R^2\,d\theta=c\)
by \eqref{uniform-side-area}.

Put \(R_0=\sqrt{x/\lambda}\). The Bessel equation for
\(g(R)=R b(R^2)\) gives the identities
\begin{equation}\label{calibrated-physical-H}
 \begin{split}
 g'(R_0)&=0,\qquad g''(R_0)=-(x-1)g_*/R_0^2,\\
 g'''(R_0)&=(x-3)g_*/R_0^3,\qquad
 g^{(4)}(R_0)=[6+(x-3)(x-2)]g_*/R_0^4,\\
 H'(R_0)&=0,\qquad H''(R_0)=-2g_*^2/R_0^2,\\
 H'''(R)&=3(3g'g''+gg''')+
            R\bigl(3(g'')^2+4g'g'''+gg^{(4)}\bigr).
 \end{split}
\end{equation}
The \(J_1\) series and its alternating remainders give, for
\(0\le\lambda z/4\le1\) and \(j,m\ge0\),
\begin{equation}\label{physical-bessel-series}
 \begin{gathered}
 b(z)=\frac{g_*\sqrt\lambda}{2J_1(j'_{1,1})}
       \sum_{\ell\ge0}\frac{(-\lambda z/4)^\ell}{\ell!(\ell+1)!},\\
 \sum_{\ell=0}^{2m+1}\frac{(-\lambda z/4)^\ell}{\ell!(\ell+j+1)!}
 \le\frac{2J_1(j'_{1,1})(-1)^j b^{(j)}(z)}
              {g_*\sqrt\lambda\,(\lambda/4)^j}
 \le\sum_{\ell=0}^{2m}\frac{(-\lambda z/4)^\ell}{\ell!(\ell+j+1)!}.
 \end{gathered}
\end{equation}

\begin{lemma}[Regular-sector conformal coordinates]
\label{effective-physical-reference}
For \(N\ge33\), set
\[
 D_S=\{z\in\mathbb D:|\arg z|\le\pi h\},
 \qquad \|f\|_{D_S}^2=\frac1{|R_N|}\int_{D_S}|f|^2 .
\]
The conformal map and its sector norm are
\[
 \begin{gathered}
 F(z)=\sqrt{\frac{\beta_N|R_N|}{\pi}}\,z(1+\gamma(z)),\\
 \gamma(z)=\sum_{m\ge1}\frac{(2h)_m(-z^N)^m}{m!(mN+1)},\qquad
 M=\|D\Im\gamma\|_{D_S}.
 \end{gathered}
\]
The first pair satisfies
\[
 U\circ F=A_h(\varphi(r)e^{i\theta}+v),\quad
 v\perp_a\varphi e^{i\theta},\quad
 A_h^2=\frac{\Phi(R_N)}{\pi((x-1)/2+\|v\|_a^2)}\in(1.68^2,1.69^2).
\]
The complementary projections of \(v\) satisfy
\[
 \|v_L\|_a\le10h^3,\qquad \|v_H\|_a\le6.61h^{5/2},
\]
where \(L,H\) are as in Proposition~\ref{prop:effective-fan-localization}.
The dipole pullback is \((A_h/g_*)B\circ F=A_he^{i\theta}PQ\), where
\[
 \begin{gathered}
 1-3.46h^3\le\sqrt{\frac{\beta_N\Phi(R_N)}{\pi x}}\le1,\qquad
 \delta=\sqrt{\frac{\beta_N\Phi(R_N)}{\pi x}}\,|1+\gamma|-1,\\
 P=\varphi(r(1+\delta)),\qquad Q=\frac{1+\gamma}{|1+\gamma|}.
 \end{gathered}
\]
\end{lemma}

\begin{proof}
Conformal invariance and orthogonality yield
\[
 \frac{\Phi(R_N)}{\pi}=\|U\circ F\|_a^2
 =A_h^2\left(\frac{x-1}{2}+\|v_L\|_a^2+\|v_H\|_a^2\right).
\]
For \(v=v_L+v_H\), the block bounds in the proof of
Proposition~\ref{prop:effective-fan-localization} give
\[
 \|K_{LL}\|\le100/2809+10h^2,\quad \|K_{HH}\|\le1.38h^2,\quad
 \|K_{LH}\|\le10.16h^{5/2}.
\]
The constant-weight pairing vanishes in mode one, so
\(|v_L(r)|\le\|v_L\|_a\) gives \(R_L\le(\beta_N^{-1}-1)/2\le2.42h^3\).
For \(1\le j\le16\) evaluate
\((2h)_j(2h)_{j+1}/[h^2j!(j+1)!]\) at \(h=1/33\) by monotonicity.
For the tail,
\[
 \frac{(2h)_j}{j!}\le\frac{2h}{1-2h}j^{2h-1},\qquad
 \sum_{j>16}j^{-15/8}\le\tfrac87\,16^{-7/8}.
\]
Consequently
\[
 |\widehat w_N(r)|\le2.145hr^N,\qquad
 R_H\le h^{5/2}\left[\tfrac34(2.145^2+2.473^2\cdot0.062)\right]^{1/2}
 <1.933h^{5/2}.
\]
With \(\pi/\Phi(R_N)>1/3.391\), the block equations give
\[
 \begin{aligned}
 0.25\|v_L\|_a&\le2.42h^3+10.16h^{5/2}\|v_H\|_a,\\
 0.2935\|v_H\|_a&\le1.933h^{5/2}+10.16h^{5/2}\|v_L\|_a .
 \end{aligned}
\]
Solving gives \(\|v_L\|_a\le10h^3\) and \(\|v_H\|_a\le6.61h^{5/2}\).
Therefore
\[
 \|v\|_a^2\le(100h+6.61^2)h^5\le49(x-1)h^5/2 .
\]
Together with \eqref{effective-fan-reg}, this yields
\(1\le\pi x/\Phi(R_N)\le1+2.075h^3\) and \(1.68<A_h<1.69\).
Since \(|R_N|\lambda=\Phi(R_N)\),
\[
 F'(0)\sqrt{\lambda/x}=\sqrt{\frac{\beta_N\Phi(R_N)}{\pi x}}
 \ge\big[(1+4.84h^3)(1+2.075h^3)\big]^{-1/2}\ge1-3.46h^3 .
\]

Monomial integration gives
\[
 \begin{gathered}
 M^2=\frac{\pi}{|R_N|}\sum_{m\ge1}m\left(\frac{(2h)_m}{m!(mN+1)}\right)^2,\\
 \|D\gamma\|_{D_S}=\sqrt{2}M,\quad
 \|\gamma\|_{D_S},\|(1-r)D\gamma\|_{D_S}\le hM,\\
 \|\Im\gamma\|_{D_S}\le hM/\sqrt{2},\quad
 \|(1-r)\Im\gamma\|_{D_S}\le h^2M/2,\qquad
 \|(1-r)^2D\Im\gamma\|_{D_S}\le\sqrt{3/2}\,h^2M.
 \end{gathered}
\]
Orthogonality gives
\(\|\partial_r\Im\gamma\|_{D_S}
 =\|r^{-1}\partial_\theta\Im\gamma\|_{D_S}=M/\sqrt2\).
For \(0\le\theta\le\pi h\), \(0\le\arg(1+z^N)\le N\theta/2\) and
\(F'(z)=\sqrt{\beta_N|R_N|/\pi}(1+z^N)^{-2h}\). By reflection,
\[
 |\arg(e^{i\theta}F'(z))|
 =|\theta-2h\arg(1+z^N)|\le\pi h\qquad(z\in D_S).
\]

For \(0\le r\le1+4.41/900\), the Bessel equation and alternating series
give
\[
 |\varphi'|\le1.74,\quad |\varphi''|\le3.391,\quad
 |\varphi'+r\varphi''|\le2.425 .
\]
Indeed, for \(g=\varphi/r\),
\[
 \begin{gathered}
 (r^3g')'=-xr^3g,\quad g(r)>g(0)/2,\quad
 g(0)\le(1-x/8)^{-1}<1.74,\quad -g'\le xg(0)r/4,\\
 \varphi''=-x\varphi-g'\le0,\qquad \varphi\le1,\\
 \varphi'+r\varphi''=(1-xr^2)g,\qquad
 g(r)\le1+0.74(1-r^2)_+ .
 \end{gathered}
\]
From \(\varphi(1)=1\), \(\varphi'(1)=0\), and
\(|\varphi''|\le3.391\), for \(0\le r\le1\),
\[
 |\varphi'|\le3.391(1-r),\qquad
 1-\varphi\le\tfrac{3.391}{2}(1-r)^2.
\]
Thus
\[
 \|D[(\varphi\cos\theta-1)\Im\gamma]\|_{D_S}
 \le\left(\frac{3.391\sqrt{3/2}}2+\frac{\pi^2}2+
 \frac{3.391}2+\frac{1.74\pi}{\sqrt2}\right)h^2M<12.6h^2M .
\]
\end{proof}

\begin{lemma}[Interior and side estimates for the correction]
\label{effective-physical-reference-sharp}
Under the hypotheses of Lemma~\ref{effective-physical-reference},
the correction \(E=U-B\) satisfies
\[
 \begin{gathered}
   1.94h^2<M<2.23h^2,\qquad
   \|DE_1\|_S<62h^3,\\
   \|D(E_2\circ F)+A_hD\operatorname{Im}\gamma\|_{D_S}<15h^3,\\
   \left|\|(E_2)_{x_i}\|_S-A_h M/\sqrt2\right|<24h^3
       \quad(i=1,2),\\
   \|DE\|_S<4.65h^2,\qquad \|E\|_S<7.2h^3,\\
   \|E\|_I<5.4h^{5/2},\qquad \|E_1\|_I<63.2h^{7/2},\qquad
   \|E_2\|_I<4.09h^{5/2}.
 \end{gathered}
\]
With \(u=X_2/\sin\alpha\) on \(I\),
\begin{equation}\label{uniform-weighted-side}
 \begin{aligned}
 \frac1{|R_N|}\int_I |E_1|\,ds&\le90h^4,&
 \frac1{|R_N|}\int_I u^2|E_1|\,ds&\le\frac{90}{\sqrt5}h^4,\\
 \frac1{|R_N|}\int_I |uE_2|\,ds&\le3.36h^3,&
 \frac1{|R_N|}\int_I |X_2E_2|\,ds&\le10.8h^4 .
 \end{aligned}
\end{equation}
\end{lemma}

\begin{proof}
\emph{Step 1. Complementary modes.}
Use \(v=v_L+v_H\) and the bounds of
Lemma~\ref{effective-physical-reference}.
The sector energy factor is
\(\sqrt{2\pi h/|R_N|}\le\sqrt{2h/(1-7h^2)}<1.419\sqrt h\).
By Lemma~\ref{effective-physical-reference},
\[
 0\le g_*-A_h<1.79h^3,\qquad |B|\le g_*,\qquad
 |\nabla B_i|<1.75g_*,\qquad |\sin\theta|\le\pi h\quad\text{on }I .
\]
Thus
\[
 \begin{aligned}
 (1-A_h/g_*)\bigl(\|B\|_S,\|\nabla B_i\|_S\bigr)
   &\le(1.79,3.14)h^{7/2},\\
 (1-A_h/g_*)\bigl(\|B\|_I,\|B_2\|_I\bigr)
   &\le(2.55,8.02h)h^{7/2}.
 \end{aligned}
\]

For \(v_H=e^{i\theta}(S+iT)\), pair the modes \(1\pm jN\), \(j\ge1\).
For the paired radial coefficients, the angular matrix and the
diagonal energy norm are
\[
 \begin{pmatrix}(jN)^2+1&2jN\\2jN&(jN)^2+1\end{pmatrix},\qquad
 \|f\|_0^2=\sum_{j\ne0}\int_0^1
 \left(|f_j'|^2+[(jN)^2+1]\frac{|f_j|^2}{r^2}\right)r\,dr .
\]
Since \(2jN/[(jN)^2+1]\le2h\), the full energy is at least
\((1-2h)(\|S\|_0^2+\|T\|_0^2)\).
The even weight preserves the two subspaces. The low sources are
symmetric, so testing on the antisymmetric component gives
\[
 \|T\|_0\leq
 \frac{2h}{1-3.391(h^2+10h^3)}\|S\|_0<2.01h\|S\|_0,
 \qquad \|S\|_0\leq1.032 \|v_H\|_a .
\]
For potential \((jN)^2+1\) the weighted mass bound is \(h^2+10h^3\).
Using \(\|f/r\|\le h\|f\|_0\) for \(f=S,T\) in
\(\Im v_H=S\sin\theta+T\cos\theta\) gives
\[
 \|D\operatorname{Im}v_H\|_{D_S}
    \leq1.419\sqrt h\bigl[(\pi+1)h\|S\|_0
                                  +(1+\pi h^2)\|T\|_0\bigr]
    <60h^4.
\]
The block bounds also give
\[
 \|Dv_L\|_{D_S}<14.2h^{7/2},\qquad
 \|D(v_L+v_H)\|_{D_S}<9.71h^3 .
\]

\emph{Step 2. Conformal map.}
The coefficients of \(\gamma\) satisfy
\[
 \frac{|[z^{mN}]\gamma|}{h^2}
 =\frac{2}{m(m+h)}\prod_{j=1}^{m-1}(1+2h/j).
\]
This is at most \(2\) for \(m=1\). For \(m\ge2\), the logarithmic
derivative of the right side is at least \(2/(1+2h)-1/(m+h)>0\).
Evaluate at \(h=1/33\) for \(2\le m\le16\) and use
\(|[z^{mN}]\gamma|\le(66/31)h^2m^{-31/16}\) for \(m>16\). By
\(\sum_{m>16}m^{-s}\le16^{1-s}/(s-1)\),
\[
 \sum_m|[z^{mN}]\gamma|<3.443h^2,\qquad
 \sum_m|[z^{mN}]\gamma|^2<4.369h^4,\qquad
 \sum_m m|[z^{mN}]\gamma|^2<4.929h^4 .
\]
Since \(1-7h^2\le|R_N|/\pi\le1\), \(1.94<2/(1+h)\le M/h^2<2.23\).

By \(0\le1-\sqrt{\beta_N\Phi(R_N)/(\pi x)}\le3.46h^3\), the norm estimates of
Lemma~\ref{effective-physical-reference} give
\[
 \|\delta\|_\infty<3.555h^2,\quad
 \|\delta\|_{D_S}<2.836h^3,\quad
 \|(1-r)\delta\|_{D_S}<1.578h^4+1.419h^{7/2}.
\]
The last inequality uses \(\|1-r\|_{D_S}<0.410\sqrt h\) and
\(\|(1-r)\gamma\|_{D_S}\leq h^2M/\sqrt2\), by monomial integration.
The radial estimates, \(|\gamma(z)|\le3.45h^2r^N\), and
\((1-r)r^N\le h/2\) give
\[
 \begin{gathered}
 |(P-\varphi)/r|\le3.391(1-r+|\delta|)|\delta|,\\
 |D(P-\varphi)|\le2.425|\delta|+3.391(1-r+|\delta|)|D\gamma|,\\
 \|(P-\varphi)/r\|_{D_S}<1.04h^3,\quad
 \|D(P-\varphi)\|_{D_S}<15.7h^3,\quad
 \|P-\varphi\|_\infty<18.9h^3.
 \end{gathered}
\]

Since \(Q=\exp(i\Im\log(1+\gamma))\),
\[
 |Q-1|\leq1.004|\operatorname{Im}\gamma|,\quad
 |Q-1-i\operatorname{Im}\gamma|\leq1.01|\gamma|^2,\quad
 |D(Q-1-i\operatorname{Im}\gamma)|\leq2.02|\gamma|\,|D\gamma|.
\]
Indeed \(DQ=iQ\,D\Im\log(1+\gamma)\), and
\[
 \begin{aligned}
 PQ-P-i\varphi\Im\gamma
   &=\varphi(Q-1-i\Im\gamma)+(P-\varphi)(Q-1),\\
 \big\|D\{e^{i\theta}(PQ-P-i\varphi\Im\gamma)\}\big\|_{D_S}
   &\le26.4h^4 .
 \end{aligned}
\]
The product rule uses \(\|D(e^{i\theta}\varphi)\|_\infty<2.462\), \(\|D\gamma\|_{D_S}=\sqrt2 M\), and the bounds on \(P-\varphi\).

\emph{Step 3. Gradient of the correction.}
Use the exact decomposition
\[
 E\circ F=A_he^{i\theta}(\varphi-PQ)+A_h(v_L+v_H)
            -(1-A_h/g_*)B\circ F .
\]
By Lemma~\ref{effective-physical-reference},
\(\|D[(\varphi\cos\theta-1)\Im\gamma]\|_{D_S}\le12.6h^2M\).
Substitution of the angular integrals and \(M\le2.23h^2\) gives
\[
 \begin{aligned}
 h^{-3}\|D(\Im(E\circ F)+A_h\Im\gamma)\|_{D_S}
 &\le A_h\big[1.04+(15.7\pi+12.6\cdot2.23+26.4+60)h\\
 &\hspace{14mm}{}+14.2\sqrt h\big]+3.14\sqrt h<15,\\
 h^{-3}\|D\Re(E\circ F)\|_{D_S}
 &\le A_h\big[15.7+1.04\pi h+26.4h+9.71\\
 &\hspace{14mm}{}+2.23(\pi+1.74/\sqrt2+3.391\pi h^2/2)\big]
       +3.14\sqrt h<62 .
 \end{aligned}
\]
By conformal invariance and the bounds on \(A_hM\),
\[
 (3.25-15h)h^2<A_hM-15h^3\le\|DE\|_S
 \le\sqrt{(62h^3)^2+(A_hM+15h^3)^2}<4.65h^2 .
\]
Both polar norms of \(D\Im\gamma\) equal \(M/\sqrt2\). The frame angle
is at most \(\pi h\), so for \(i=1,2\),
\[
 \left|\|\partial_iE_2\|_S-\frac{A_hM}{\sqrt2}\right|
 \le 15h^3+\frac{A_hM}{\sqrt2}(\pi h+\pi^2h^2/2)<24h^3.
\]

\emph{Step 4. Mass and traces.}
The radial majorant for \(|k|\ge N\) gives
\[
 \||F'|f\|_{D_S}^2\le\|f\|_{D_S}^2+10h^3\|Df\|_{D_S}^2 .
\]
Using the squared mass bounds \(h^2M^2/2,h^2M^2,h^4M^2/2\)
and \(\|D((1-r)\gamma)\|_{D_S}\le2hM\) gives
\[
 \||F'|\Im\gamma\|_{D_S}<2h^3,\quad
 \||F'|\gamma\|_{D_S}<2.828h^3,\quad
 \||F'|(1-r)\gamma\|_{D_S}<2.92h^4.
\]
The radial integral
\(\int_0^1[\beta_N^{-1}\widehat w_0(r)-1]r\,dr=(\beta_N^{-1}-1)/2\le2.42h^3\)
and \(\||F'|\|_{D_S}=\sqrt h\) give, with the block bounds,
\[
 \begin{gathered}
 \||F'|(1-r)\|_{D_S}<0.410\sqrt h,\quad
 \||F'|(P-\varphi)\|_{D_S}<1.18h^3,\\
 \||F'|v_L\|_{D_S}\le3.02h^{7/2},\qquad
 \||F'|v_H\|_{D_S}\le11.03h^4 .
 \end{gathered}
\]
Since \(|P|\le1\),
\[
 h^{-3}\|E\|_S\le
 1.69(3.02\sqrt h+11.03h+1.18+1.004\cdot2)+1.79\sqrt h<7.2 .
\]

On the arc put \(\Theta=N\theta,\ w_b=|1+e^{i\Theta}|^{-2h}\).
Since \(2\cos(\Theta/2)\ge2(\pi-|\Theta|)/\pi\),
\(\|w_b\|_{L^2(d\Theta/(2\pi))}\le(1-4h)^{-1/2}<1.067\).
For \(f=\sum_{j\ne0}f_je^{ij\Theta}\), H\"older and Hausdorff--Young give
\[
 \frac1{2\pi}\int_{-\pi}^\pi w_b |f|^2\,d\Theta
 \leq\|w_b\|_2\,\frac{\pi}{\sqrt3}
                              \sum_{j\ne0}|j|\,|f_j|^2 .
\]
Apply this to \(S+iT=e^{-i\theta}v_H\) and to \(S,T\) separately.
The radial trace estimate, \(\|T\|_0<2.01h\|S\|_0\),
\(|\sin\theta|\le\pi h\), and the mode-one bound yield
\[
 \begin{aligned}
 \|v_H\circ F^{-1}\|_I&<13.3h^{7/2},&
 \|\Im(v_H\circ F^{-1})\|_I&<70h^{9/2},\\
 \|v_L\circ F^{-1}\|_I&\le14.2h^{7/2},&
 \|\Im(v_L\circ F^{-1})\|_I&\le14.2\pi h^{9/2}.
 \end{aligned}
\]

Since \((w_b-1)_+\) is supported on \(2\pi/3\le|\Theta|\le\pi\),
\[
 \frac1{2\pi}\int_{-\pi}^{\pi}(w_b-1)_+\,d\Theta
 \le\frac13\left(\frac{(3/2)^{2h}}{1-2h}-1\right)<1.024h .
\]
On the arc, Parseval and the radial bound give
\[
 \begin{gathered}
 \frac1{2\pi}\int_{-\pi}^{\pi}\bigl(\Im\gamma(e^{i\Theta/N})\bigr)^2\,d\Theta
 =\tfrac12\sum_m|[z^{mN}]\gamma|^2,\\
 \|\Im\gamma\|_\infty\le\sum_m|[z^{mN}]\gamma|,\qquad
 |P-1|\le\tfrac{3.391}{2}|\delta|^2.
 \end{gathered}
\]
Using \(F'(0)=\sqrt{\beta_N|R_N|/\pi}\le1\),
\[
 h^{-5}\|\Im\gamma\circ F^{-1}\|_I^2
 \le\frac{2}{1-7h^2}\left(\frac{4.369}{2}+1.024h(3.45)^2\right)<2.27^2.
\]
Using \(|Q-1|\le1.004|\Im\gamma|\) and
\(|\Re(Q-1)|=|Q-1|^2/2\), we obtain
\[
 |\Re(e^{i\theta}(Q-1))|
 \le\bigl(\pi h+\tfrac12(1.004)(3.45)h^2\bigr)|Q-1| .
\]
Also \(|I|/|R_N|\le2.015h\) yields
\(\|(P-1)\circ F^{-1}\|_I\le30.5h^{9/2}\). Therefore
\[
 \begin{aligned}
 h^{-5/2}\|E\|_I
 &\le1.69\bigl(1.004(2.27)+27.5h+30.5h^2\bigr)+2.55h<5.4,\\
 h^{-5/2}\|E_2\|_I
 &\le1.69\bigl(1.004(2.27)+(14.2\pi+70+30.5)h^2\bigr)+8.02h^2<4.09 .
 \end{aligned}
\]
Likewise,
\[
 h^{-7/2}\|E_1\|_I
 \le1.69\left[27.5+30.5h+
   \left(\pi+\tfrac12(1.004)(3.45)h\right)1.004(2.27)\right]+2.55<63.2 .
\]
Since \(ds=|I|\,du/2\) and \(|I|/|R_N|\le2.015h\),
\[
 \frac1{|R_N|}\int_I|u|^m|E_i|\,ds
 \le\left(\frac{2.015h}{2m+1}\right)^{1/2}\|E_i\|_I\qquad(m\ge0).
\]
Use \((m,i)=(0,1),(2,1),(1,2)\) and \(|X_2|=\sin\alpha\,|u|\le\pi h|u|\)
to obtain \eqref{uniform-weighted-side}.
\end{proof}

\begin{lemma}[Energy and flux of the anchored correction]
\label{anchor-E-33}
Let \(h=1/N\le1/33\) and \(0.53\le k=\beta/\alpha\le1.36\). Assume
\[
 \begin{gathered}
 \|DE\|_S\ge (3.25-15h)h^2,\qquad \|E\|_S\le7.2h^3,\\
 \|E_1\|_I\le63.2h^{7/2},\qquad \|E_2\|_I\le4.09h^{5/2}.
 \end{gathered}
\]
For the sum \(V(k)\) of the transported volume term and the flux term linear
in \(E\) in \(F(L)\) at \((r,s,\rho)=(r,s_*(k),0)\),
\[
 \mathcal T[V](k)\le-15h^4(k-1)^2 .
\]
\end{lemma}

\begin{proof}
\emph{(1) Energy.}
For \(v=\beta^2\), put \(t(v)=\tan\sqrt v/\sqrt v\) and
\(a(v)=(t(v)^{-1}-1)/v\), extended continuously at \(0\).
The partial fractions
\[
 a(v)=-2\sum_{m\ge1}\frac1{m^2\pi^2-v},\qquad
 t(v)=2\sum_{m\ge1}\frac1{(m-\tfrac12)^2\pi^2-v}
\]
give, for \((f,d)=(a,\pi^2),(t,\pi^2/4)\),
\begin{equation}\label{physical-trig-enclosure}
 |f^{(j)}(v)|\le |f^{(j)}(0)|(1-v/d)^{-j-1}\qquad(0\le v<d).
\end{equation}
Hence, for \(0\le v\le0.017\),
\[
 \begin{aligned}
 (|a|,|a'|,|a''|,|a'''|)&\le(0.334,0.023,0.0043,0.0013),\\
 (t,|t'|,|t''|,|t'''|)&\le(1.006,0.34,0.274,0.334),\qquad t\ge1 .
 \end{aligned}
\]
Separate \(2k^2r\) times the reference flux \eqref{uniform-reference-flux}
from \eqref{uniform-transported-flux}. The resulting energy coefficients
are \(kr^{2(j-i)}-2k^2r\), where
\[
 \begin{gathered}
 k^2r=k^3\frac{t(v)}{t(\alpha^2)},\qquad
 \mathcal T[k^2r]\ge\frac{(k+2)(k-1)^2}{t(\alpha^2)}\ge2.51(k-1)^2,\\
 |r-k|=k\left|\frac{t(v)}{t(\alpha^2)}-1\right|\le1.092\alpha^2|k-1|,\\
 \mathcal T[kr^2-2k^2r]=-(k+2)(k-1)^2+k(r-k)^2
                      \le-2.51(k-1)^2 .
 \end{gathered}
\]
The first inequality uses the nonnegative series coefficients of \(t\).
Since \(t(v)^{-2}=(1+va(v))^2\),
\[
 kr^{-2}=\frac{t(\alpha^2)^2}{k}+t(\alpha^2)^2\alpha^2 k(2a(v)+va(v)^2).
\]
The first two derivatives in \(v\) of \(2a+va^2\) have moduli at most
\(0.158,0.040\). For \(v=k^2\alpha^2\), the product rule gives
\[
 \left|\frac{d^2}{dk^2}[t(\alpha^2)^2\alpha^2 k(2a+va^2)]\right|
 \le1.006^2\alpha^4\big[6(1.36)(0.158)+4(1.36)^3\alpha^2(0.040)\big]
 <1.31\alpha^4.
\]
Taylor's formula and \(\mathcal T[k^{-1}]=(k-1)^2/k\) give
\(\mathcal T[kr^{-2}]\le1.92(k-1)^2\). Thus
\[
 \sum_{i,j=1}^2\mathcal T[kr^{2(j-i)}-2k^2r]\,\|\partial_iE_j\|_S^2
 \le-2.51\|DE\|_S^2(k-1)^2.
\]

\emph{(2) Flux and mass.}
In \eqref{uniform-anchored-side}, set \(A=a(v)+t(v)u^2\), so
\(p=c(1+vA)\) and \(Y_1^2=c(1+va(v))\). The cross flux is
\[
 \frac{2k}{|R_N|}\int_I
 \big[D_1 E_1+2c\beta r u b'(p)E_2\big]\,ds,qquad
 D_1=b(p)+2b'(p)Y_1^2 .
\]
Thus
\[
 \begin{gathered}
 |A|\le0.673,\quad (|A'|,|A''|,|A'''|)\le(0.37,0.28,0.34),\\
 (|p'|,|p''|,|p'''|)\le(0.68,0.741,0.871),\qquad
 (|\partial_v^j(Y_1^2)|)_{j=1}^3\le(0.34,0.05,0.02).
 \end{gathered}
\]
In the first bound, \(|A|\le 2/3+(0.023+0.34)v\le0.673\).

For \(0.975\le p,Y_1^2\le1.025\) and \(3.38<\lambda<3.42\), use
\(0.58<J_1(j'_{1,1})<0.59\) and \(1.68<g_*<1.69\) to obtain
\[
 0.8238\le\lambda p/4\le0.8764,\qquad
 2.6174<\frac{g_*\sqrt\lambda}{2J_1(j'_{1,1})}<2.6943 .
\]
By \eqref{physical-bessel-series} with \(m=1\),
\begin{equation}\label{physical-near-side-bessel}
 |b|\le1.76,\quad 0.8119\le-b'\le0.873,\quad 0.248\le b''\le0.267,\quad
 |b'''|\le0.060,\quad |b^{(4)}|\le0.011 .
\end{equation}
For \(R_0^2=x/\lambda\), write \(D=b+2pb'\).
The regular estimate and \eqref{physical-near-side-bessel} give
\[
 D(R_0^2)=0,\quad 0\le R_0^2-c\le2.075h^3,\quad
 \sup_{[0.975,1.025]}|D'|\le2.14,\quad |D(c)|\le4.45h^3 .
\]
Differentiation gives \((|\partial_v^j D_1|)_{j=1}^3\le(2,2,3)\). For
\[
 K(v,u)=\frac{D_1(v,u)-D(c)}{vt(v)},\qquad
 t(v)K(v,u)=\int_0^1\partial_v D_1(\tau v,u)\,d\tau,
\]
we have \((|\partial_v^j(tK)|)_{j=0}^2\le(2,1,1)\),
\((|K_v|,|K_{vv}|)\le(1.68,2.7)\), and
\[
 k(D_1(k)-krD_1(1))
  = D(c)(k-k^2r)+kvt(v)[K(v,u)-K(\alpha^2,u)] .
\]
For \(f=k^2r,kr^2\), the bounds on \(t\) give \((|\partial_k^j f|)_{j=0}^2\le(2.55,5.67,8.5)\).
The product rule with
\[
 kvt(v)=\alpha^2t(\alpha^2) k^2r,\qquad
 |v-\alpha^2|\le0.85\alpha^2,\quad |v'|\le2.72\alpha^2,\quad v''=2\alpha^2
\]
yields \(\big|\partial_k^2[k(D_1(k)-krD_1(1))]\big|\le8.5(4.45)h^3+74\alpha^4\).
Similarly, \((|\partial_v b'(p)|,|\partial_v^2 b'(p)|)\le(0.182,0.23)\) gives
\[
 \left|\partial_k^2
 \big(2c\alpha u\,k^2r[b'(p(v,u))-b'(p(\alpha^2,u))]\big)\right|\le15.9\alpha^3|u|.
\]
By \eqref{uniform-weighted-side} and Taylor's formula, the flux remainder
divided by \((k-1)^2\) is bounded by
\[
 (8.5(4.45)h^3+74\alpha^4)90h^4+15.9\alpha^3(3.36)h^3
 \le1670h^6+3405h^7+650000h^8 .
\]

For \(j=1,2\), the mass coefficient is
\(\lambda[2k^2r-kr^{2j-2}\beta\cot\beta/(\alpha\cot\alpha)]\). Since
\[
 \left(\left|\partial_k^\ell\frac{\beta\cot\beta}{\alpha\cot\alpha}\right|
       \right)_{\ell=0}^2\le(1.004,0.93\alpha^2,0.69\alpha^2),\qquad
 \left|\partial_k^2\!\left(2k^2r-
    kr^{2j-2}\frac{\beta\cot\beta}{\alpha\cot\alpha}\right)\right|\le26 ,
\]
the mass remainder is at most \(13\lambda\|E\|_S^2(k-1)^2\le2310h^6(k-1)^2\).
Hence, for \(h\le1/33\),
\[
 \mathcal T[V](k)\le
 \bigl[-2.51(3.25-15h)^2+3980h^2+3405h^3+650000h^4\bigr]h^4(k-1)^2
 \le-15h^4(k-1)^2.\qedhere
\]
\end{proof}

\begin{lemma}[Dipole and sector integral at the area anchor]\label{anchor-pure}
Let \(h=1/N\le1/33\) and \(0.50\le k=\beta/\alpha\le1.40\).
For the sum \(V_0(k)\) of the area-corrected dipole term and squared
sector integral in \eqref{physical-total-residual} at \((s,\rho)=(s_*(k),0)\),
\[
 \mathcal T[V_0](k)\le-6h^5(k-1)^2 ,
\]
with \(s_*,\mathcal T\) as in \eqref{sec4-angular-anchor}.
\end{lemma}

\begin{proof}
\emph{(1) Dipole term.}
Use \(v,A\) from Lemma~\ref{anchor-E-33}. At the anchor,
\(p=c(1+vA)\) and \(\int_{-\beta}^{\beta}(p-c)\,d\theta=0\)
by \eqref{uniform-side-area}--\eqref{uniform-anchored-side}.
For \(G(p)=H(\sqrt p)\) and
\(W(p)=\int_0^1(1-\tau)G''(c+\tau(p-c))\,d\tau\), Taylor's identity gives
\[
 \begin{aligned}
 \frac1{|R_N|}\int_{-\beta}^{\beta}[G(p)-G(c)]\,d\theta
 &=\frac1{|R_N|}\int_{-\beta}^{\beta}(p-c)^2W(p)\,d\theta
 =\frac{2c}{\pi}\beta^5 J(v),\\
 J(v)&=\frac12\int_{-1}^1\frac{A(v,u)^2 W(c(1+vA(v,u)))}{1+vA(v,u)}\,du .
 \end{aligned}
\]

The subtracted term \(2\beta G(c)\) is affine in \(k\).
For \(v\le0.02\), \eqref{physical-trig-enclosure} gives
\[
 (|a'|,|t'|,|a''|,|t''|)\le(0.023,0.34,0.0043,0.274),\qquad
 |A_v|\le0.37,\quad |A_{vv}|\le0.29.
\]
Then \(0.9939\le c\le1\), \(0.98\le p\le1.014\), and
\[
 \begin{gathered}
 \left|A-(u^2-\tfrac13)\right|\le0.0074,\qquad |A|\le0.675,\\
 \sqrt{4/45}-0.0074\le\left(\tfrac12\int_{-1}^1 A(v,u)^2\,du\right)^{1/2}
 \le0.306 .
 \end{gathered}
\]
For \(D=b+2pb'\), the radial equation gives
\[
 G'(p)=\tfrac12[D^2-(\lambda p-1)b^2+2\lambda].
\]
At \(R_0^2=x/\lambda\), \(D(R_0^2)=0\), \(G''(R_0^2)=-g_*^2/(2R_0^4)\). With \(0\le R_0^2-c\le2.075h^3\), \(3.38<\lambda<3.42\), \eqref{physical-near-side-bessel} yields on \([0.975,1.025]\)
\[
 |D|\le0.067,\quad |D'|\le2.14,\quad |D''|\le1.46,\quad |D'''|\le0.443.
\]
Differentiate \(G'\) to obtain
\begin{align*}
 G'''&=(D')^2+DD''-2\lambda b b'
                   -(\lambda p-1)((b')^2+b b''),&
 |G'''|&\le16,\\
 G^{(4)}&=3D'D''+DD'''-3\lambda((b')^2+b b'')
                   -(\lambda p-1)(3b'b''+b b'''),&
 |G^{(4)}|&\le40 .
\end{align*}
For \(0.98\le p\le1.014\),
\[
 |p-R_0^2|\le0.021,qquad
 -0.95\le W\le-0.5,\qquad (|W'|,|W''|)\le(16/6,40/12).
\]

Since \(1+vA\ge0.9865,\ |A+vA_v|\le0.683,\ |2A_v+vA_{vv}|\le0.746\),
\[
 \left(\left|\partial_v^\ell\frac{A^2}{1+vA}\right|\right)_{\ell=1}^2
 \le(0.827,2.17).
\]
Therefore
\[
 \begin{aligned}
 J&\le-\frac{(\sqrt{4/45}-0.0074)^2}{2(1.0135)}\le-0.041,\\
 |J'|&\le0.95\left(\frac{2(0.306)(0.37)}{0.9865}
         +\frac{0.094(0.683)}{0.9865^2}\right)
       +\frac{16}{6}\frac{0.094(0.683)}{0.9865}<0.46 .
 \end{aligned}
\]
For the second derivative use \(A^2/(1+vA)\le0.675^2/0.9865\):
\[
 |J''|
 \le0.95(2.17)+2(0.827)\frac{16}{6}(0.683)
   +\frac{0.675^2}{0.9865}\left[\frac{40}{12}(0.683)^2+\frac{16}{6}(0.746)\right]<7.
\]
Hence
\[
 \frac{d^2}{d\beta^2}\big[\beta^5J(\beta^2)\big]
  =\beta^3(20J+22vJ'+4v^2J'')\le-0.60\beta^3.
\]
Since \(\beta=k\alpha\ge\alpha/2\),
\[
 \mathcal T[(2c/\pi)\beta^5J(\beta^2)]\le
 -\frac{0.60c}{8\pi}\alpha^5(k-1)^2\le-7h^5(k-1)^2 .
\]

\emph{(2) Sector integral.} Set
\(f(p)=2\int_0^1 q^2 b(q^2p)\,dq\).
Radial integration and reflection give
\[
 \frac1{|S|}\int_S[B(LX)-LB(X)]\,dX
 =\frac{(\sqrt{c\beta\cot\beta},0)}2
    \int_{-1}^1[f(p(v,u))-f(p(\alpha^2,u))]\,du .
\]
By \eqref{physical-bessel-series}, \((|f'|,|f''|)\le(0.52,0.106)\).
Since \(\tfrac12\int_{-1}^1 A(v,u)\,du=a(v)+t(v)/3\),
\[
 \begin{aligned}
 \tfrac12\int_{-1}^1 A(0,u)\,du&=0,&
 \left|\tfrac12\int_{-1}^1 A_v\,du\right|&\le0.023+0.34/3,\\
 \left|\tfrac12\int_{-1}^1 (A+vA_v)\,du\right|&\le0.274v .
 \end{aligned}
\]
Subtract \(f'(c)\) before estimating the second term below:
\[
 \begin{aligned}
 \left|\tfrac12\partial_v\int_{-1}^1 f(p(v,u))\,du\right|
 &\le0.52(0.274)v+0.106v\tfrac12\int_{-1}^1|A(A+vA_v)|\,du\\
 &\le\big[0.52(0.274)+0.106(0.094+0.02(0.306)(0.37))\big]v
 <0.16v .
 \end{aligned}
\]
As \(|v-\alpha^2|\le2.4\alpha^2|k-1|\) and \(\max(v,\alpha^2)\le1.96\alpha^2\),
\[
 \frac1{|S|}\left|\int_S [B(LX)-LB(X)]\,dX\right|\le0.76\alpha^4|k-1|.
\]
The squared integral has a double zero at \(k=1\). Thus
\[
 \mathcal T[V_0](k)\le(-7h^5+27000h^9)(k-1)^2
 \le\left(-7+\frac{27000}{33^4}\right)h^5(k-1)^2
 \le-6h^5(k-1)^2.\qedhere
\]
\end{proof}

\begin{lemma}[Scale derivatives about the angular anchor]
\label{mid-scale-33}
Let \(h\le1/33\), \(r\in[0.53,1.37]\), and \(|s|\le42h^2\).
For \(\beta=\arctan(r\tan\alpha)\) and \(s_*\) in \eqref{sec4-angular-anchor},
\[
 \begin{aligned}
 \mathcal B_{ss}(r,s,0)&\le
 \begin{cases}-0.95rh,&h\le1/48,\\-0.24rh,&1/48<h\le1/33,\end{cases}\\
 |\mathcal B_s(r,s_*,0)|&\le\frac{\beta}{\alpha}h^4|r^2-1|
                   \begin{cases}157.5,&r\le1,\\134,&r>1.\end{cases}
 \end{aligned}
\]
\end{lemma}

\begin{proof}
At \(\rho=0\), \eqref{centered-physical-map} gives
\[
 y=e^{s/2}\frac{\cos\beta}{\cos\alpha}\operatorname{diag}(1,r)X,
 \qquad w=e^s r\frac{\cos^2\beta}{\cos^2\alpha},
 \qquad z=|y|^2 .
\]
By Lemma~\ref{effective-physical-reference-sharp},
\[
 \begin{gathered}
 \|DE_1\|_S\leq62h^3,\qquad \|DE\|_S\leq4.65h^2,\qquad
 \|E\|_S\leq7.2h^3 .
 \end{gathered}
\]
The side integrals are bounded by \eqref{uniform-weighted-side}.

\emph{(1) Scale curvature.}
For \(H(R)=Rg(R)g'(R)+\lambda R^2\), the spectral bound and sine series give
\[
 3.388<\lambda<3.413,\qquad
 1-4h^2\leq R_0\leq1-3.25h^2,\qquad 2.836<g_*^2<2.84 .
\]
Taking \(0\le z\le1.04\) in \eqref{physical-bessel-series} gives
\[
 (|b|,|b'|,|b''|,|b'''|,|b^{(4)}|,|b^{(5)}|)
 \le (3,1.3,0.37,0.08,0.02,0.003)
\]
componentwise, with \((-1)^jb^{(j)}\ge0\).

For \(R=e^{s/2}\cos\beta/\cos(u\beta)\), \(|u|\leq1\), we have
\(1-30.34h^2\leq R\leq1+21.21h^2\) and
\(|R-R_0|\leq28h^2\), using \(\tan\alpha\le3.152h\).
From \eqref{calibrated-physical-H} and the alternating signs,
\[
 |g^{(5)}|=|60b''+360R^2b'''+240R^4b^{(4)}+32R^6b^{(5)}|
 \leq31.
\]
Successive integration from \(R_0\) gives
\[
 |g|\leq1.693,\quad |g'|\leq0.106,\quad
 -4.10\leq g''\leq-3.99,\quad
 0<g'''\leq0.99,\quad 0<g^{(4)}\leq12.05 .
\]
Consequently
\[
 0<H'''=3(3g'g''+gg''')+
 R(3(g'')^2+4g'g'''+gg^{(4)})\leq82.
\]
Indeed \(3R(g'')^2>9|g'g''|+4R|g'g'''|\). Thus \(H''\) increases, with
\[
 \begin{cases}
 H'\le0,\quad H''\le-5.672+82(25.21h^2),& R\ge R_0,\\
 |H'|\le7.84\cdot28h^2,\quad H''\le-5.672,&R<R_0.
 \end{cases}
\]
The identities \(R_s=R/2\) and
\(2\beta/|R_N|\ge1.988rh\) yield
\[
 \frac{\beta}{|R_N|}\int_{-1}^1\frac{H''(R)R^2+H'(R)R}{4}\,du
 \le
 \begin{cases}
 -2.36rh,& h\le1/48,\\
 -1.86rh,&1/48<h\le1/33.
 \end{cases}
\]

For \(V(X)=[b(z)-b(|X|^2)]y=B(y)-LB(X)\), reflection gives
\(\int_S V_2\,dX=\int_S b'(z)zy_2\,dX=0\).
For \(X=q(\cos\alpha,u\sin\alpha)\),
\(0\le q\le1,\ |u|\le1\), one has
\(dX/|S|=q\,dq\,du\), \(|z-|X|^2|\le51.6h^2q^2\),
\(|y|\le qe^{21h^2}\). The three-term upper bound for \(J_2\) yields
\[
 -b'(z)\leq \frac{2.92\lambda^{3/2}}{16}
        \left(1-\frac{\lambda z}{12}+\frac{\lambda^2z^2}{384}\right).
\]
Four terms of \(J_1\) give \(g_*/J_1(\sqrt{x})<2.92\). Integrating first in \(q\) with \(y=qY\), \(v=|Y|\), gives the increasing polynomial on \([0,e^{21h^2}]\)
\[
 \tfrac25 v^3-\tfrac{3.388}{12}\tfrac27 v^5
            +\tfrac{3.413^2}{384}\tfrac29 v^7
\]
whose endpoint value is less than
\(0.4\cdot16/(2.92\cdot3.413^{3/2})\). Thus
\[
 \frac1{|S|}\left|\int_S b'(z)zy\,dX\right|\le0.4,\qquad
 V_s=\tfrac12 V+b'(z)zy,\quad
 V_{ss}=\tfrac14 V+2b'(z)zy+b''(z)z^2y .
\]
Consequently
\[
 \frac1{|S|}\left|\int_S V\,dX\right|\le27.4h^2,\quad
 \frac1{|S|}\left|\int_S V_s\,dX\right|\le0.413,\quad
 \frac1{|S|}\left|\int_S V_{ss}\,dX\right|\le0.970 .
\]
Since \(w_s=w_{ss}=w\le1.047r\),
\[
 \begin{aligned}
 \frac{\lambda h}{|S|^2}\partial_{ss}\!
        \left(w\left|\int_S V\,dX\right|^2\right)
 &\le\lambda hw\big[(27.4h^2)^2+4(27.4h^2)(0.413)\\
 &\hphantom{{}\le\lambda hw\big[}
       +2(0.413)^2+2(27.4h^2)(0.970)\big]\\
 &\le\begin{cases}1.38rh,& h\le1/48,\\
                    1.55rh,& 1/48<h\le1/33.\end{cases}
 \end{aligned}
\]

The transported energy minus mass equals
\[
 w\sum_{i,j=1}^2r^{2(j-i)}\|\partial_iE_j\|_S^2
 -\lambda we^s\cos^2\beta/\cos^2\alpha\sum_{j=1}^2r^{2j-2}\|E_j\|_S^2 .
\]
Its second scale derivative is bounded above by the energy term,
which is at most \(0.0023rh\).
For \(D_{ij}=(DB(y))_{ij}\), the flux in \eqref{uniform-transported-flux}
is \((w/|R_N|)\int_I (D_{11}E_1+rD_{21}E_2)\,ds\), since \(L\) is diagonal.
On \(I\), \(|z-R_0^2|,|y_1^2-R_0^2|\leq60h^2\).
With \(D_{11}(R_0,0)=0\), this gives \(|D_{11}|\leq281h^2\).
Also \(|D_{11,s}|\leq2.41\) and \(|D_{11,ss}|\leq4.24\), using
\[
 D_{11,s}=zb'+2(b'+zb'')y_1^2,\qquad
 D_{11,ss}=zb'+z^2b''+2(b'+3zb''+z^2b''')y_1^2 .
\]
Here \(D_{11,s}(R_0,0)=(1-x)g_*/(2R_0)\) and \(|D_{11,s}-D_{11,s}(R_0,0)|\le415h^2\); alternating signs bound \(D_{11,ss}\). For the tangential flux,
\[
 r(D_{21}+2D_{21,s}+D_{21,ss})
       =2r(4b'+5zb''+z^2b''')y_1y_2 ,\qquad
 |ry_1y_2|\leq1.96|X_2|.
\]
The cross curvature is bounded in modulus by \(2w(9.4\cdot90+20.8\cdot10.8)h^4<0.063rh\), and by \(0.021rh\) for \(h\le1/48\). Summing gives
\[
 \begin{aligned}
 \mathcal B_{ss}/(rh)&\le-1.86+1.55+0.0023+0.063<-0.24,\\
 \mathcal B_{ss}/(rh)&\le-2.36+1.38+0.0023+0.021<-0.95
                                   \qquad(h\le1/48).
 \end{aligned}
\]

\emph{(2) Scale slope at the angular anchor.}
The weak equation for \(U\) annihilates the variation in the trial. At \(r=1\),
\[
 \mathcal B_s(1,0,0)=-2\lambda h+2\lambda h=0 .
\]
The squared-integral derivative vanishes there. In angular coordinates
\(R^2=c\beta\cot\beta\sec^2(u\beta)\) for \(|u|\le1\)
and \(\tfrac12\int_{-1}^1R^2\,du=c\) by \eqref{uniform-side-area}.
In triangular coordinates \(X=q(\cos\alpha,u\sin\alpha)\),
\eqref{uniform-anchored-side} gives
\[
 y_1=q\sqrt{c \beta\cot\beta },\qquad
 z=c q^2[\beta\cot\beta +\beta\tan\beta u^2].
\]
The sine series gives \(|(v\cot v)'+2v/3|\le0.228v^3\) on
\([0,0.262]\). For \(0\le v\le0.131\),
\(\tan v/v\le1.009\), hence
\[
 0\le\tan v-v=\int_0^v\tan^2 t\,dt\le0.34v^3,\quad
 |(v\tan v)'-2v|\le1.36v^3,\quad \frac{2v}{\sin(2v)}\le1+0.676v^2 .
\]
Since \(|R-\sqrt c|\le0.338\beta^2\) and
\(|R_0-\sqrt c|\le1.04h^3\), we have \(|R-R_0|\le6.34h^2\)
and \(R\le1.006\). Logarithmic differentiation gives
\[
 \begin{aligned}
 \frac{R_r}{R\beta_r}&=\frac{1-2\beta\cot(2\beta)}{2\beta}
                      -\tan\beta+u\tan(u\beta),\\
 \left|\frac{R_r}{R}-\beta\beta_r(u^2-1/3)\right|
              &\le1.14\beta^3\beta_r .
 \end{aligned}
\]
Thus \(\tfrac12\int_{-1}^1|u^2-1/3|\,du=4/(9\sqrt3)\) gives
\(\tfrac12\int_{-1}^1|R_r|\,du\le2.77rh^2\).
Since \(\tfrac12\int_{-1}^1R^2\,du=c\), the normalized dipole slope satisfies
\[
 \frac{d}{dr}\left[\frac{h}{2c}\int_{-1}^1 H'(R)R\,du\right]
 =\frac{h}{2c}\int_{-1}^1RR_r
       \bigl[H''(R)-H''(R_0)+H'(R)/R\bigr]\,du .
\]
Its absolute value is at most \(1600rh^5\), by
\(|H''(R)-H''(R_0)|\le82|R-R_0|\) and \(|H'(R)|\le6.25|R-R_0|\).

Write a prime for total \(r\)-differentiation at the anchor. Then
\[
 \left|\frac{(y_1)'}{y_1}\right|\leq3.37rh^2,\qquad
 |z-c q^2|\leq0.676\beta^2q^2,\qquad
 \tfrac12\int_{-1}^1 |z'|\,du\leq5.22rh^2q^2 .
\]
Use
\[
 \tfrac12\int_{-1}^1|2\beta(u^2-1/3)|\,du=\frac{8\beta}{9\sqrt3},
 \qquad
 \big|\partial_\beta(\beta\cot\beta+\tfrac13\beta\tan\beta)\big|
 \le0.682\beta^3 .
\]
Subtract \(b'(cq^2)\) in \(\tfrac12\int_{-1}^1 b'(z)z'y_1\,du\) to obtain
\[
 \begin{aligned}
 \frac1{|S|}\left|\int_S\left(V'-\frac{(y_1)'}{y_1}V\right)dX\right|
 &\leq1.3\,\frac25\,0.682(3.152)^4r^3h^4
 \\[-2pt]
 &\quad+0.37\cdot0.676(1.37\cdot3.152)^2\cdot5.22\,\tfrac27\,rh^4
 \leq73rh^4 .
 \end{aligned}
\]
Integrating \(\bigl((\beta\cot\beta)^{-1/2}\int_S V\,dX\bigr)'\) from
\(r=1\), where \(V=0\), yields
\[
 \frac1{|S|}\left|\int_S V\,dX\right|\le32.3h^4,\qquad
 \frac1{|S|}\left|\int_S V'\,dX\right|\le73.1rh^4 .
\]
The identity \(V_s=V/2+b'(z)zy\) gives
\[
 \frac1{|S|}\left|\int_S(V_s)'\,dX\right|\le \tfrac12(73.1rh^4)
       +1.67\,\tfrac25\,5.22rh^2+0.4\cdot3.37rh^2\le4.9rh^2 .
\]
Thus
\[
 \frac{\lambda h}{|S|^2}
 \left|\frac{d}{dr}\left[(1+\partial_s)
            \left|\int_S V\,dX\right|^2\right]_{s=s_*}\right|\le208rh^5 .
\]

The energy coefficients \(2r,-2r^{-3}\) give at most
\(\max\{2r\|DE\|_S^2,\,2r^{-3}\|DE_1\|_S^2\}\) in modulus.
The mass coefficients are \([-2\lambda (\beta/\alpha)(r^{-1},r)]'\).
Thus the normalized transported-slope derivative is bounded by
\(91rh^4\), or \(44rh^4\) for \(r\ge1\).

On the shifted side
\(|z'|\le13.53rh^2,\ |(y_1^2)'|\le6.67rh^2\);
\(|z-R_0^2|\le13h^2\) and the value at \(R_0^2\) give \(|b'|\le0.86\).
Since \(r(D_{21}+D_{21,s})=2(2b'+zb'')r^2c \beta\cot\beta X_1X_2/\cos^2\alpha\) and
\(|2b'+zb''|\le1.72\), the derivatives satisfy
\[
 |(D_{11}+D_{11,s})'|\le85rh^2,\qquad
 |[r(D_{21}+D_{21,s})]'|\le7r|X_2|.
\]
Thus the absolute cross-slope derivative is at most
\[
 \frac2{|R_N|}\int_I(85rh^2|E_1|+7r|X_2E_2|)\,ds
 \le2(85\cdot90rh^6+7\cdot10.8rh^4)<166rh^4 .
\]
Summing the derivatives gives
\[
 \left|\left(\frac{\alpha}{\beta}\mathcal B_s(r,s_*,0)\right)'\right|
 \le rh^4\left[(1600+208)h+166+
  \begin{cases}91,&r\le1,\\44,&r>1\end{cases}\right]
 \le\begin{cases}315rh^4,&r\le1,\\268rh^4,&r>1.\end{cases}
\]
Integration from \(r=1\) gives the asserted bound for \(|\mathcal B_s(r,s_*,0)|\).
\end{proof}

\begin{lemma}[The radius-ratio estimate]\label{effective-skew-33-direct}
Let \(h=1/N\le1/33\) and use the physical sector test with
\[
 0.53\le r\le1.37,\quad |s|\le42h^2,\quad
                         |\rho|\le45h^{5/2}.
\]
Under the reference bounds of Lemma~\ref{effective-physical-reference-sharp},
\[
 \mathcal B(r,s,\rho)-\mathcal B(r,s,0)\le -0.01rh\rho^2 .
\]
\end{lemma}

\begin{proof}
Put \(T=\tan(\pi h)\), \(f=e^s(1+T^2)/(1+r^2T^2)\),
\(w=fr\), \(C=\sqrt{1+\delta}\), and \(\delta=\sinh^2(\rho/2)\). Then
\[
 \begin{gathered}
 3.14h\le T\le3.16h,\quad 0.953\le f\le1.047,\\
 \delta\le0.251\rho^2\le509h^5,\qquad C^2=1+\delta,\\
 a_1=1+r^2T^2\le1.018,\quad
 b_1=r^2+T^{-2}\le0.104h^{-2},\\
 d_1=T^{-1}+r^2T\le0.324h^{-1},\qquad d_1^2=a_1b_1 .
 \end{gathered}
\]
For \(X=(x_1,x_2)\) in the reference triangle, put
\(u=x_2/(Tx_1)\). Radial integration yields
\[
 |S|^{-1}\int_S x_1^j |u|^m\,dX
 =\frac{2\cos^j\alpha}{(j+2)(m+1)}\le\frac2{(j+2)(m+1)}
 \quad(j,m\ge0).
\]
On \(I\), \(u=X_2/\sin\alpha\) and \(ds/|I|=du/2\).
Use the trace bounds \eqref{uniform-weighted-side}.

\emph{(1) Correction terms.} For \(M=L^TL\),
\[
 M_{11}=f(1+a_1\delta),\quad M_{22}=f(r^2+b_1\delta),\quad
 M_{12}=\tfrac12 f d_1\sinh\rho.
\]
Reflection gives the transported energy
\[
 \begin{aligned}
 w^{-1}\big[&M_{11}M_{22}\big(\|\partial_1E_1\|_S^2+\|\partial_2E_2\|_S^2\big)
           +M_{11}^2\|\partial_2E_1\|_S^2+M_{22}^2\|\partial_1E_2\|_S^2\\
           &-\frac{2M_{12}^2}{|R_N|}\int_S
             \big[(\partial_1 E_1)(\partial_2 E_2)+(\partial_2 E_1)(\partial_1 E_2)\big]\,dX\big].
 \end{aligned}
\]
Since \(\int_S E_1E_2\,dX=0\), the increment in the transported squared mass is
\[
 \frac{wf\delta}{|R_N|}\int_S(a_1 E_1^2+b_1 E_2^2)\,dX\ge0 .
\]
In the increment of \eqref{uniform-transported-flux}, write \(z=|LX|^2\)
and separate \(2wf b_1\delta|R_N|^{-1}\int_I x_1x_2 b'(|X|^2)E_2\,ds\).
By \eqref{uniform-reference-flux},
\[
 \begin{aligned}
 2\frac1{|R_N|}\int_I x_1x_2 b'(|X|^2)E_2\,ds
  ={}&-\|DE\|_S^2+\lambda\|E\|_S^2\\
     &-\frac1{|R_N|}\int_I [b(|X|^2)+2x_1^2b'(|X|^2)]E_1\,ds.
 \end{aligned}
\]
Since \(b(R_0^2)+2R_0^2b'(R_0^2)=0\) and
\(\big||X|^2-R_0^2\big|,|x_1^2-R_0^2|\le10h^2\), the side bounds
below give
\[
 |b(|X|^2)+2x_1^2b'(|X|^2)|\le33h^2,\qquad
 \frac1{|R_N|}\int_I 33h^2|E_1|\,ds\le2970h^6 .
\]
After adding \(-2wf b_1\delta\|DE\|_S^2\) from twice the separated flux,
the coefficient of \(\|\partial_2E_1\|_S^2\) is nonpositive by
\(2a_1+a_1^2\delta\le2fr^2b_1\). The coefficients of
\(\|\partial_1E_1\|_S^2+\|\partial_2E_2\|_S^2\) and
\(\|\partial_1E_2\|_S^2\) are, respectively,
\[
 (f/r)b_1\delta[a_1(1+\delta)-2fr^2],\qquad
 (f/r)b_1\delta[2r^2(1-f)+b_1\delta].
\]
Use
\[
 \begin{gathered}
 \|\partial_1E_1\|_S^2+\|\partial_2E_2\|_S^2\le15.04h^4,\qquad
 \|\partial_1E_2\|_S^2\le11.51h^4,\\
 \frac1{|R_N|}\left|\int_S\big[
 (\partial_1 E_1)(\partial_2 E_2)+(\partial_2 E_1)(\partial_1 E_2)\big]\,dX\right|
 \le62(4.65)h^5 .
 \end{gathered}
\]
Thus the energy increment and twice the separated flux total at most
\(r\delta(3h^2+226h^3+720h^4)\), with coefficients bounded by
\[
 \begin{gathered}
 \frac{1.047(0.104)(0.483)}{0.53^2}(15.04)
   +1.047(0.104)\left(0.094+\frac{0.104(509)}{33^3(0.53)^2}\right)(11.51)<3,\\
 \frac{2(1.047)(0.324)^2(1+509/33^5)62(4.65)}{0.53^2}<226,\\
 2(1.047)^2(0.104)(3.42(7.2)^2+2970)<720 .
 \end{gathered}
\]

For the remaining flux increment, use
\[
 \begin{gathered}
 z=z_0+\delta Z_++\tfrac12(\sinh\rho) Z_-,\\
 z_0=f(x_1^2+r^2x_2^2),\quad
 Z_+=f(a_1x_1^2+b_1x_2^2),\quad Z_-=2f d_1x_1x_2,\\
 \big|z_0-|X|^2\big|\le53h^2 x_1^2,\qquad
 |Z_+|\le2.134x_1^2,\quad |Z_-|\le2.134x_1^2|u|.
 \end{gathered}
\]
The first difference follows by writing
\(z_0-|X|^2=(e^s-1)|X|^2+
 e^s(r^2-1)(x_2^2-T^2x_1^2)/(1+r^2T^2)\).
On \(I\), both \(|LX|^2\) and its zero-skew value lie between
\(\cos^2\beta\exp(-42h^2-45h^{5/2})>0.935\) and
\(\exp(42h^2+45h^{5/2})<1.047\); also \(0.935<|X|^2\le1\).
On \([0.935,1.05]\), \eqref{physical-bessel-series} with \(m=1\),
\(2.6174\le g_*\sqrt\lambda/[2J_1(j'_{1,1})]\le2.6943\),
\(0.790\le\lambda z/4\le0.898\), gives
\[
 1.565\le b\le1.812,\quad
 0.801\le-b'\le0.888,\quad 0.246\le b''\le0.271,\quad
 |b'''|\le0.061,\quad |b^{(4)}|\le0.011.
\]
By Taylor's formula and parity, with \(b_\pm^{(j)}=\tfrac12[
 b^{(j)}(z(x_1,x_2))\pm b^{(j)}(z(x_1,-x_2))]\),
\[
 \begin{gathered}
 |b_+-b(z_0)|\le2.52\delta,\quad
 |b'_+-b'(z_0)|\le0.721\delta,\\
 |b'_+-b'(|X|^2)|\le14.5h^2,\quad
 |b'_-|\le0.581\tfrac12|u\sinh\rho|.
 \end{gathered}
\]
Subtracting the separated term, the flux difference is
\[
 \begin{split}
 w\int_I\big[&
 (b_+-b(z_0))E_1
 +2x_1 f(b'_+-b'(z_0))(x_1E_1+r^2x_2E_2)\\
 &+2x_1 f\delta(a_1 b'_+ x_1E_1+b_1(b'_+-b'(|X|^2))x_2E_2)\\
 &+x_1 f d_1(\sinh\rho)b'_-(x_1E_2+x_2E_1)\big]\,\frac{ds}{|R_N|}.
 \end{split}
\]
The three trace bounds give a modulus \(\le r\delta(2.8h^2+1356h^4)\) after doubling. The \(h^2\) coefficient is \(\le4(1.047)^2(0.324)(0.581)(1+509/33^5)(3.36)\); the \(h^4\) coefficient is bounded by
\[
\begin{split}
 2(1.047)\big[&
 2.52(90)+2(1.047)(0.721)(90+1.37^2(3.16)(3.36))\\
 &+2(1.047)(1.018)(0.89)(90)
 +2(1.047)(0.104)(14.5)(3.16)(3.36)\\
 &+2(1.047)(0.324)(0.581)(1+509/33^5)(3.16)(90/\sqrt5)\big]
 <1356 .
\end{split}
\]

\emph{(2) Sector integral.} Write \(V_\rho(X)=[b(|LX|^2)-b(|X|^2)]LX\),
with \(L\) evaluated at skew \(\rho\). Reflection and
\(\big|z_0-|X|^2\big|\le53h^2x_1^2\) give
\[
 \int_S(V_0)_2\,dX=0,\qquad
 \frac1{|S|}\left|\int_S(V_0)_1\,dX\right|\le28.3h^2 .
\]
On the full triangle one uses \(|b'|\le1.3\), \(|b''|\le0.37\),
and
\[
 LX=\sqrt f(Cx_1+\sinh(\rho/2) x_2/T,\quad r(Cx_2+T\sinh(\rho/2)x_1)).
\]
The preceding bounds on \(Z_+,Z_-\), combined with radial integration
and the even--odd decomposition, give
\[
 \begin{split}
 \frac1{|S|}\left|\int_S(V_\rho-V_0)_1\,dX\right| &\le\sqrt f\,\delta\left[
 \begin{gathered}\tfrac12(1.3)(53)h^2\tfrac25
 +C(1.3)(2.134)\tfrac25\\
 +\tfrac C2(0.37)(2.134)^2(1+\delta)\tfrac{2}{21}
 +(1.3)(2.134)C\tfrac{2}{15}\end{gathered}\right]\le1.7\delta,\\
 \frac1{|S|}\left|\int_S(V_\rho)_2\,dX\right| &\le \sqrt f\,rT \sqrt\delta
 \left[(1.3)(2.134)C^2\tfrac{2}{15}
 +(1.3)(53)h^2\tfrac25+3.63\delta\tfrac23\right]
 \le1.8h\sqrt\delta .
 \end{split}
\]
Here \(|b_+-b(z_0)|\le3.63\delta\) follows from Taylor's formula
and the Bessel bounds on the full triangle.
Thus
\[
 \begin{aligned}
 \frac{\lambda wh}{|S|^2}\!\left(
 \left|\int_S V_\rho\,dX\right|^2-\left|\int_S V_0\,dX\right|^2\right)
 &\le \lambda wh\big[2(28.3h^2)(1.7\delta)\\
 &\qquad +(1.7\delta)^2+(1.8h)^2\delta\big]\le360rh^3\delta .
 \end{aligned}
\]
The increments from (1) and (2) therefore total at most
\[
 \begin{aligned}
 r\delta(5.8h^2+586h^3+2076h^4)
 &\le0.251rh\rho^2(5.8h+586h^2+2076h^3)\\
 &\le0.251rh\rho^2(5.8/33+586/33^2+2076/33^3)\le0.194rh\rho^2 .
 \end{aligned}
\]

\emph{(3) Dipole term.}
Take \(G(p)=H(\sqrt p)\) and \(p(v)=R(\beta v)^2\). By \eqref{uniform-side-area},
\[
 \tfrac12\int_{-1}^1p(v)\,dv
 =\frac{e^s\sin(2\beta)}{2\beta},\qquad
 \int_{-1}^1p_{\rho\rho}\,dv=0 .
\]
Hence
\[
 \partial_{\rho\rho}\frac12\int_{-1}^1G(p)\,dv
 =\frac12\int_{-1}^1\!\left[G''(p)p_\rho^2+
 \left(G'(p)-G'\!\left(\frac{e^s\sin(2\beta)}{2\beta}\right)\right)p_{\rho\rho}\right]dv .
\]
With \(p_*=(j'_{1,1})^2/\lambda\), the bounds on \(b\) above
and the formula for \(G'\) in Lemma~\ref{anchor-pure} yield
\(0<G'''<14\) and \(G^{(4)}<0\) on \([0.935,1.05]\).
Indeed \(D'<0,\ |D'|\le2.21,\ |D|\le0.147,\ 0<D''\le1.355\).
Substitute these and the signed bounds on \(b,b',b''\) in the formula
for \(G'''\) in Lemma~\ref{anchor-pure}, retaining its negative last
term. In \(G^{(4)}\), \(-3\lambda((b')^2+bb'')<-10\) and all positive
terms sum to \(<2.3\). The Bessel equation gives
\[
 G'''(p_*)=\frac{g_*^2}{p_*^3}
           \left[(j'_{1,1})^2+
                  \frac{((j'_{1,1})^2-1)((j'_{1,1})^2-3)}2\right]
 <11.23 .
\]
Since \(0.9939\le p_*\le1.00006\),
\(G''(p_*)=-g_*^2/(2p_*^2)\), and \(p\le1.047\), one obtains
\(-2.36\le G''(p)\le-0.81\) on the range of \(p\).

By \eqref{centered-skew-radius}, for \(q=2(\log R)_\rho\),
\[
 p_\rho=pq,\qquad p_{\rho\rho}=p(3q^2-1)/2 .
\]
With \(t(v)=\tan\sqrt v/\sqrt v\), \(0\le t'\le0.34\) on \([0,0.018]\) and
\(|\rho|\le45h^{5/2}\) give
\[
 \begin{gathered}
 \left|\frac{\tan(v\beta)}{\tan\beta}-v\right|\le0.34\beta^2,\quad
 \left|q-\frac{\tan(v\beta)}{\tan\beta}\right|\le0.00362,\quad |q-v|\le0.01,\\
 \tfrac12\int_{-1}^1p_\rho^2\,dv\ge0.935^2(1/\sqrt3-0.01)^2>0.281,\qquad
 \|p_{\rho\rho}\|_2\le1.047(1/\sqrt5+0.03)<0.502 .
 \end{gathered}
\]
Using the midpoint of \([\min p,\max p]\) in the least-squares bound gives
\[
 \begin{aligned}
 \left(\frac12\int_{-1}^1
 \left|p(v)-\frac{e^s\sin(2\beta)}{2\beta}\right|^2\,dv\right)^{1/2}
 &\le\frac{\max p-\min p}{2}\\
 &\le\frac{e^s}{2}\left(e^{|\rho|}-e^{-|\rho|}\cos^2\beta\right)<0.0165 .
 \end{aligned}
\]
Cauchy--Schwarz and the preceding second-derivative identity give
\[
 \partial_{\rho\rho}\frac12\int_{-1}^1G(p)\,dv
 \le-0.81(0.281)+2.36(0.0165)(0.502)<-0.208 .
\]
Taylor's formula, evenness in \(\rho\), and \(2\beta/|R_N|\ge1.98rh\)
give the total increment \(\le(-0.205+0.194)rh\rho^2\le-0.01rh\rho^2\).
\end{proof}

\begin{proof}[Proof of Proposition~\ref{effective-global-neumann}]
Scale \(P\) so that \(|P|=|R_N|\) and suppose
\(\mu_1(P)\ge\lambda=\mu_1(R_N)\). Put \(h=1/N\), \(\alpha=\pi h\).
By Proposition~\ref{prop:effective-fan-localization}, \(P\) has \(N\) sides.
Since \((q^4-4q+3)'=4(q^3-1)\) and
\[
 \min_{q\in\{0.535,1.356\}}(q^4-4q+3)>\frac{148}{4.8\cdot33},
\]
it gives, on each sector based at \(|P|^{-1}\int_Py\,dy\),
\begin{equation}\label{explicit-neumann-box}
 \begin{gathered}
 0.53\le k=\beta/\alpha\le1.36,\qquad
 0.53\le r=\tan\beta/\tan\alpha\le1.37,\\
 |s|=|\log(R_-R_+)|\le42h^2,\qquad
 |\rho|=|\log(R_+/R_-)|\le\sqrt{1898}\,h^{5/2}+145h^4<45h^{5/2}.
 \end{gathered}
\end{equation}
Here \(|k-q|\le24h^3\) and \(1\le\tan z/z\le1+0.34z^2\) for
\(0\le z\le0.131\). Also \(|s_*|\le42h^2\) by
\(-\log(\sin z/z)\le z^2/5\) for \(0\le z\le0.3\).

At \((r,s,\rho)=(1,0,0)\), \(F(\mathrm{Id})=0\) and
\(\mathcal B_r=\mathcal B_s=0\) by \eqref{physical-sector-residual}
and Lemma~\ref{mid-scale-33}. The angular correction in
\eqref{physical-sector-residual} is affine in \(k\); hence
Lemmas~\ref{anchor-E-33} and~\ref{anchor-pure} give
\[
 \mathcal B(r,s_*,0)\le-(15h^4+6h^5)(k-1)^2.
\]
Completing the square in \(s-s_*\), Lemma~\ref{mid-scale-33} yields
\[
 \mathcal B(r,s,0)-\mathcal B(r,s_*,0)
 \le\frac{|\mathcal B_s(r,s_*,0)|^2}{0.48rh}
 \le8.369h^4(k-1)^2 .
\]
Indeed \(k/r\le1.004\), \(|r-1|\le1.018|k-1|\), and
\[
 \frac{1.004^2\,1.018^2}{0.48\cdot33^3}
 \max\{4(157.5)^2,\,1.37(2.37)^2(134)^2\}<8.369 .
\]
By Lemma~\ref{effective-skew-33-direct},
\[
 \mathcal B(r,s,\rho)\le-6h^4(k-1)^2-0.01rh\rho^2 .
\]
Apply Rayleigh's inequality in \eqref{physical-total-residual} and
sum over the \(N\) sectors:
\[
 \begin{aligned}
 0&\le\sum_{j=1}^N\big[6h^4(k_j-1)^2+0.01r_jh\rho_j^2\big]\\
 &\le \frac{\lambda-\mu_1(P)}{|P|}
 \int_P\left|W(y)-\frac1{|P|}\int_P W(z)\,dz\right|^2\,dy\le0.
 \end{aligned}
\]
Hence \(2\beta_j=2\pi/N\) and \(R_{j+1}=R_j\) for all \(j\).
The area fixes the common circumradius, so \(P\cong R_N\).
\end{proof}

\section{Polygons with a reflex vertex and the equality case}
\label{sec:nr-nonconvex}

A conformal estimate excludes polygons with fewer than \(N\) interior
angles smaller than \(\pi\) when \(N\ge5\). The reflex quadrilateral
estimate completes the proof of the main theorem, which then gives the
Voronoi corollary.

Suppress collinear vertices and let \(\theta_j\) be the remaining
interior angles. Thus \(p(P)=\#\{j:\theta_j<\pi\}\).
\begin{proposition}\label{prop:concave-four}
Every nonconvex simple quadrilateral \(Q\) satisfies
\[
 |Q|\mu_1(Q)<30\sqrt{3/29}<\pi^2.
\]
Among simple polygons of at most four sides and fixed area, only the
square attains the maximum \(\pi^2/|Q|\).
\end{proposition}

\begin{proof}
Put
\[
 C=\frac{\int_Q xx^{\mathsf T}\,dx}{|Q|}
    -\frac{(\int_Q x\,dx)(\int_Q x\,dx)^{\mathsf T}}{|Q|^2},
 \qquad y=C^{-1/2}\left(x-\frac1{|Q|}\int_Q z\,dz\right).
\]
The linear Rayleigh tests give
\[
 \int_Qy\,dx=0,\quad \int_Q yy^{\mathsf T}\,dx=|Q|I,\quad
 \mu_1(Q)C\preceq I,\quad \mu_1(Q)\det(C)C^{-1}\preceq I,
\]
where \(C\) and \(\det(C)C^{-1}\) have the same two eigenvalues.
For \(i=1,2\), take
\[
 H_i=\frac1{|Q|}\int_Qy_i yy^{\mathsf T}\,dx,\qquad
 u_i=y_i+\tfrac14\big(y^{\mathsf T}H_i y-\tr H_i\big).
\]
Since \(\int_Q y\,dx=0\) and \(\int_Q yy^{\mathsf T}\,dx=|Q|I\),
\[
 \int_Q u_i\,dx=0,\qquad
 \frac1{|Q|}\int_Q y_i u_i\,dx
 =1+\tfrac14\tr(H_i^2)
 =\frac1{|Q|}\int_Q|e_i+\tfrac12 H_i y|^2\,dx .
\]
By Cauchy--Schwarz, the Rayleigh quotient and the matrix bound,
\[
 \begin{aligned}
 \mu_1(Q)^2\det C\left(\int_Q y_i u_i\,dx\right)^2
 &\le |Q|\mu_1(Q)\det C\int_Q|\nabla u_i|^2\,dx\\
 &\le |Q|\int_Q|e_i+\tfrac12H_i y|^2\,dx .
 \end{aligned}
\]
Since \(\int_Q y_i u_i\,dx=\int_Q|e_i+\tfrac12H_i y|^2\,dx>0\),
\[
 \mu_1(Q)^2\det C\left(1+\tfrac14\tr(H_i^2)\right)\le1
 \ \Longrightarrow\
 \frac1{\Phi(Q)^2}\ge\frac{\det C}{|Q|^2}
 \left(1+\tfrac18\sum_{i=1}^2\tr(H_i^2)\right).
\]

Under \(x\mapsto Lx+b\), the ratio \(\det C/|Q|^2\) is fixed and
\[
 y\longmapsto (LCL^{\mathsf T})^{-1/2}LC^{1/2}y,
 \qquad (LCL^{\mathsf T})^{-1/2}LC^{1/2}\in O(2).
\]
Thus \(\sum_i\tr(H_i^2)\) is fixed as the squared norm of the tensor
with entries \(|Q|^{-1}\int_Q y_i y_j y_k\,dx\).
Take the vertex opposite the reflex vertex as the origin. The ray through
the reflex vertex meets the opposite side of the enclosing triangle at
\((1,0)\). Take that side parallel to the second axis. The vertices are then
\((0,0),(1,p-1),(s,0),(1,p)\), with \(0<p,s<1\), in cyclic order.
The two maps \((u,v)\mapsto(su+v,pv)\), \((su+v,(p-1)v)\) on
\(u,v\ge0\), \(u+v\le1\), give
\[
 \frac1{|Q|}\int_Q x_1^i x_2^j\,dx
 =\frac{2[p^{j+1}+(1-p)(p-1)^j]}{(i+j+2)!}
   \sum_{h=0}^i \binom{i}{h}s^h h!\,(i+j-h)! .
\]
For \(q=p(1-p)\), the moments of degrees at most three give
\[
 |Q|=s/2,\qquad
 \frac1{|Q|}\int_Qx\,dx=\tfrac13(s+1,2p-1),\qquad
 432\det C=s^2+4q(1-s).
\]
Contracting each pair of factors \(C^{-1/2}\) in
\(\sum_i\tr(H_i^2)\) to \(C^{-1}\) gives
\begin{align*}
 &\frac{\det C}{|Q|^2}\left[1+\tfrac18\sum_{i=1}^2\tr(H_i^2)\right]
       -\frac{29}{2700}\\
 &\quad=\frac{2q(1-s)}{675s^2[s^2+4q(1-s)]^2}
  \left[19s^4+2q(1-s)
       \big\{(76-27q)s^2+(27+44q)(1-s)\big\}\right]>0 ,
\end{align*}
since \(0<q\le1/4\). Finally
\(\Phi(R_3)=4\pi^2/(3\sqrt3)<\pi^2=\Phi(R_4)\).
Combine with the triangle comparison \cite{LS} and the convex
quadrilateral comparison \cite{EO}.
\end{proof}

Chuaqui, Duren, and Osgood \cite{CDO} and Chuaqui and Pommerenke \cite{CP}
treated the unshifted quotient. Wang \cite[Lemma~3.3 and equation~(3.8)]{Wang}
proved the shifted bound. Dong and Wang \cite[Lemma~3]{DW} adapted it
to the Neumann balancing center by retaining the first Fourier coefficient.

\begin{lemma}[Positive turns and conformal coefficients]\label{nr-positive-turn}
Let \(p(P)\le m\) and \(f:\mathbb D\to P\) be conformal and onto.
For \(f'(z)/f'(0)=1+\sum_{k\ge1}b_k z^k\),
\begin{equation}\label{nr-turn-point}
 \sum_{k\ge1}(k/2-c)^2|b_k|^2v^k\ge(1+c)^2v^m-c^2
 \qquad (0<v<1,\ c\ge0).
\end{equation}
\end{lemma}

\begin{proof}
Discard collinear intermediate vertices, which do not affect \(p(P)\).
With prevertices \(\xi_j\), the Schwarz--Christoffel formula gives
\[
 \sum_j(\pi-\theta_j)=2\pi,\qquad
 \frac{f'(z)}{f'(0)}=\prod_j(1-\bar\xi_jz)^{\theta_j/\pi-1},
 \qquad H(z)=\frac{zf''(z)}{2f'(z)}
   =\sum_j\frac{\pi-\theta_j}{2\pi(1-\bar\xi_j z)}-1 .
\]
Put \(B_c=(H-c)/(1+c+H)\).
Since \(\Re H=-1/2\) on \(\partial\mathbb D\setminus\{\xi_j\}\),
\[
 |B_c(\zeta)|=1\quad(\zeta\in\partial\mathbb D),\qquad
 B_c^{-1}(1)=\{\xi_j\},\qquad
 \xi_j B_c'(\xi_j)=\frac{2\pi(1+2c)}{\pi-\theta_j}.
\]
The poles at the prevertices are simple and clearing them gives no
common zero. Circle reflection and the signed crossings at \(1\) give
\[
 \begin{aligned}
 \#\{\text{interior zeros of }B_c\}
 &=\tfrac12\left(\deg B_c+\sum_j\operatorname{sgn}(\pi-\theta_j)\right)\\
 &=\tfrac12\sum_j[1+\operatorname{sgn}(\pi-\theta_j)]=p(P).
 \end{aligned}
\]
All counts are with multiplicity and \(\deg B_c=\#\{\xi_j\}\).
For \(0<r<1\) and each interior zero \(\alpha\),
\[
 \frac1{2\pi}\int_0^{2\pi}
 \log\left|\frac{re^{i\theta}-\alpha}{1-\bar\alpha re^{i\theta}}\right|
 \,d\theta=\log\max(r,|\alpha|)\ge\log r .
\]
The pole terms are nonnegative, and \(H-c=B_c(1+c+H)\).
Jensen's formula for \(1+c+H\), with \(H(0)=0\), gives
\[
 \frac1{2\pi}\int_0^{2\pi}\log|H(re^{i\theta})-c|\,d\theta
 \ge \log(1+c)+p(P)\log r\ge\log(1+c)+m\log r .
\]
The circle integral of \(\log|f'/f'(0)|\) vanishes since \(f'\ne0\).
For \(z=re^{i\theta}\), convexity of the exponential and Parseval give
\[
 \begin{aligned}
 (1+c)^2r^{2m}
 &\le \frac1{2\pi}\int_0^{2\pi}
      \left|\frac{zf''(z)-2c f'(z)}{2f'(0)}\right|^2\,d\theta\\
 &=c^2+\sum_{k\ge1}(k/2-c)^2|b_k|^2r^{2k}.
 \end{aligned}
\]
\end{proof}
\begin{theorem}\label{analytic-positive-turn}\label{sharp-positive-corner-asymptotic}
Let \(m\ge3\) and \(x=(j'_{1,1})^2\). For a simple polygon \(P\) with
\(p(P)\le m\),
\[
 \begin{gathered}
 (x-1)\left(\frac{\pi x}{\Phi(P)}-1\right)\ge\frac{H_m}{1+(100/49)H_m},\\
 H_m=\frac{4\sum_{j\ge1} I_{mj}\prod_{i=1}^{j-1}(1+4/(mi))}
 {1+(100/49)m^{-2}(2+4/(m-2))},\qquad
 I_k=\frac{2}{(k+1)^2(k+2)}-\frac{1}{(k+2)^2(k+3)} .
 \end{gathered}
\]
Moreover, as \(m\to\infty\),
\[
 \sup_{\substack{P\ {\rm simple\ polygon}\\p(P)\le m}}\Phi(P)
 =\pi x-\frac{4\pi x\zeta(3)}{(x-1)m^3}+O(m^{-5}).
\]
\end{theorem}

\begin{proof}
\emph{(1) The common estimate.}
For the profile \(\varphi(r)=J_1(\sqrt{x}\,r)/J_1(\sqrt x)\), the
Bessel equation gives \(2E=xd_1=x-1\).
Equation~\eqref{nr-ray-dipole} at \(T=\Phi(P)\) yields
\[
 (x-1)\left(\frac{\pi x}{\Phi(P)}-1\right)\ge\frac{xQ}{1+S},
 \qquad S=\sum_{k\ge1}\frac{|b_k|^2}{k+1},\quad
 Q=\sum_{k\ge1}(d_{k+1}-d_1)\frac{|b_k|^2}{k+1}.
\]
The Bessel equation gives
\[
 \frac{d}{d(r^2)}\frac{\varphi(r)}r
 =-\frac{x}{2}\int_0^1 t^3\frac{\varphi(rt)}{rt}\,dt<0.
\]
The right side increases with \(r\), so the tangent in \(r^2\) at
\(r=1\) gives
\[
 \frac{\varphi(r)}r\ge\frac{3-r^2}{2}
 \ \Longrightarrow\
 1-\varphi(r)^2\le\frac{(1-r^2)^2(4-r^2)}4.
\]
Integrating against \(2(k+1)r^{2k+1}\) gives
\[
 1-d_{k+1}\le\frac{2}{(k+2)(k+3)},\qquad
 1-d_2\le\int_0^1r^3(1-r^2)^2(4-r^2)\,dr=\frac3{20}.
\]
Since \(d_j\) increases and \(x\le112/33<17/5\),
\[
 xQ\ge x(d_2-d_1)S\ge(1-3x/20)S\ge\frac{49}{100}S.
\]

Set \(K(s)=(s-1+e^{-s})(2e^{-s}-e^{-2s})\) and
\(I_k=\int_0^\infty K(s)e^{-ks}\,ds\). Direct integration gives
\[
 1-(k+1)k^2I_k-\frac8{(k+2)(k+3)}
 =\frac{4(k^2+k-1)}{(k+1)(k+2)^2(k+3)}>0 .
\]
Consequently \(k^2I_k/x\le(d_{k+1}-d_1)/(k+1)\).
For \(s>0\), put
\[
 G(s)=\sum_{k\ge1}|b_k|^2e^{-ks},\qquad M(s)=\frac1{e^{ms}-1}.
\]
With \(v=e^{-s}\) and \(c=M(s)\), \eqref{nr-turn-point} gives
\begin{equation}\label{positive-turn-differential-inequality}
 G''/4+M^2G\ge M(1-G').
\end{equation}
We use
\[
 \int_0^\infty K G''\,ds=\sum_{k\ge1}k^2 I_k|b_k|^2\le xQ,\qquad
 \int_0^\infty s^2 M^2G\,ds\le\frac1{m^2}
       \sum_{k\ge1}\frac{|b_k|^2}{k+m}\le \frac S{m^2}.
\]
The second bound uses \(s^2 e^{ms}M^2=[s/(2\sinh(ms/2))]^2\le m^{-2}\).

\emph{(2) The integrated bound.}
For \(F(s)=(1-e^{-ms})^{4/m}\), Wang's factor \cite[Section~3.3]{Wang},
\[
 F'=4MF,\qquad
 [F(1-G')]'=4MF(1-G')-FG''\le4FM^2G
\]
by \eqref{positive-turn-differential-inequality}.
Integrate to infinity, where \(F(1-G')\to1\), to get
\[
 \frac{G''(s)}4\ge\frac{M(s)}{F(s)}-M(s)^2G(s)
 -4\frac{M(s)}{F(s)}\int_s^\infty F(t)M(t)^2G(t)\,dt .
\]
For \(0<s\le t\), \(F(t)/F(s)\le(t/s)^{4/m}\) by concavity of
\(1-e^{-ms}\). Hence
\[
 \int_0^tK(s)M(s)\frac{F(t)}{F(s)}\,ds
 \le\frac{t^{4/m}}{2m}\int_0^t s^{1-4/m}\,ds
 =\frac{t^2}{4(m-2)} .
\]
Integrate \(G''/4\) against \(4K\) and use \(S\le(100/49)xQ\):
\[
 \begin{aligned}
 4\sum_{j\ge1} I_{mj}\prod_{i=1}^{j-1}\left(1+\frac4{mi}\right)
 &=4\int_0^\infty K(s)\frac{M(s)}{F(s)}\,ds\\
 &\le xQ+\frac{S}{m^2}\left(2+\frac4{m-2}\right)\\
 &\le xQ\left[1+\frac{100}{49m^2}\left(2+\frac4{m-2}\right)\right].
 \end{aligned}
\]
Thus \(xQ\ge H_m\) and
\[
 \frac{xQ}{1+S}\ge\frac{xQ}{1+(100/49)xQ}\ge\frac{H_m}{1+(100/49)H_m}.
\]

\emph{(3) The asymptotic bound.}
Use \(K=(s-1+e^{-s})[1-(1-e^{-s})^2]\),
\((1-e^{-s})^2\le s^2\), and
\(s^2/2-s^3/6\le s-1+e^{-s}\le s^2/2\). Thus
\[
 I_k\ge \int_0^\infty (s^2/2-s^3/6-s^4/2)e^{-ks}\,ds
      = k^{-3}-k^{-4}-12k^{-5}.
\]
In the product in \(H_m\), retain the constant and linear terms:
\[
 \prod_{i=1}^{j-1}\left(1+\frac4{mi}\right)\ge 1+\frac4m\sum_{i=1}^{j-1}\frac1i .
\]
Summing the identity
\[
 \frac1{a^2b^2}
 =\frac{a^{-2}+b^{-2}}{(a+b)^2}
  +\frac{2(a^{-1}+b^{-1})}{(a+b)^3}\qquad(a,b\ge1)
\]
gives \(\zeta(4)=4\sum_{1\le i<j}i^{-1}j^{-3}\). Hence the
\(m^{-4}\) terms cancel, and \(\sum_{i<j}1/i\le j\) gives
\[
 4\sum_{j\ge1}I_{mj}\prod_{i<j}\left(1+\frac4{mi}\right)
 \ge\frac{4\zeta(3)}{m^3}
  -\frac{48\zeta(5)+16\zeta(3)+192\zeta(4)/m}{m^5}.
\]
Summing the first thirty terms of each zeta series and bounding the
tail by an integral gives \(1.2<\zeta(3)<1.203\),
\(\zeta(4)<1.085\), and \(\zeta(5)<1.04\). For \(m\ge29\),
\[
 H_m\ge\frac{4\zeta(3)m^{-3}-77m^{-5}}
              {1+(100/49)(2+4/(m-2))m^{-2}}
     \ge4\zeta(3)m^{-3}-99m^{-5}>0.
\]
By the monotonicity of \(h/[1+(100/49)h]\) for \(h>0\),
\begin{equation}\label{positive-turn-explicit-tail}
 \frac{H_m}{1+(100/49)H_m}\ge\frac{4\zeta(3)}{m^3}-\frac{101}{m^5}
 \qquad(m\ge29).
\end{equation}
This gives the required upper bound. For the lower bound take
\(P=R_m\) in \eqref{regular-neumann-asymptotic-formula}.
\end{proof}

\begin{proposition}\label{nr-large-exclusion}
If \(N\ge5\) and \(p(P)<N\), then \(\Phi(P)<\Phi(R_N)\).
\end{proposition}

\begin{proof}
Numerical Lemma~\ref{nr-finite-weights} gives \(5\le N\le8\).
For \(9\le N\le29\), write \(h_m\) for the first eight terms of \(H_m\)
in Theorem~\ref{analytic-positive-turn}, with the same prefactor.
Set \(x_+=3.389958\). By the alternating Bessel series,
\[
 2J'_1(\sqrt{x_+})\le\sum_{k=0}^6
 \frac{(-1)^k(2k+1)(x_+/4)^k}{k!(k+1)!}<-4\cdot10^{-8},
 \qquad x<x_+ .
\]
Take the regular lower endpoints \(L_N<\Phi(R_N)\):
\[
 L_N=
 \begin{cases}
 10.615,&N=9 \quad(\text{Numerical Lemma~\ref{nonagon-local-signs}}),\\
 t_N,&10\le N\le15\quad(\text{Numerical Lemma~\ref{decagon-reference-eigenvalue}}),\\
 T_N,&16\le N\le29\quad\text{from \eqref{middle-spectral-enclosures}}.
 \end{cases}
\]
By Theorem~\ref{analytic-positive-turn}, \(x<x_+\),
\(\pi<3.141592654\) and \(h_{N-1}\le H_{N-1}\), rational substitution gives
\[
 \Phi(R_N)-\Phi(P)> L_N-
 \frac{3.141592654\,x_+}{1+h_{N-1}/[(x_+-1)(1+(100/49)h_{N-1})]}
 >5.9\cdot10^{-5}\qquad(9\le N\le29).
\]

For \(N\ge30\), put \(m=N-1\) in \eqref{positive-turn-explicit-tail}.
As \((N-1)^{-3}\ge N^{-3}+3N^{-4}\) and
\[
 14.4N>101(30/29)^5+113,
\]
the right side of \eqref{positive-turn-explicit-tail} exceeds
\(4\zeta(3)N^{-3}+113N^{-5}\).
By \eqref{effective-fan-reg} and the gamma series,
\[
 (x-1)\left[\frac{\pi x}{\Phi(R_N)}-1\right]
 \le 1-\beta_N+29.3xN^{-5}
 <4\zeta(3)N^{-3}+113N^{-5}
 <(x-1)\left[\frac{\pi x}{\Phi(P)}-1\right].
\]
\end{proof}

\begin{proof}[Proof of Theorem~\ref{main-polygonal-comparison}]
Normalize \(|P|=|R_N|\) and discard collinear intermediate vertices.
For \(N=3\) use \cite{LS}. Proposition~\ref{prop:concave-four}
gives \(N=4\).
For \(N\ge5\), Proposition~\ref{nr-large-exclusion} gives
\[
 p(P)<N\quad\Longrightarrow\quad \Phi(P)<\Phi(R_N).
\]
For \(p(P)=N\) the polygon is convex with \(N\) vertices.
The inequalities follow from Theorem~\ref{local-pentagon} and
Propositions~\ref{hexagon-comparison}, \ref{smallmid-comparison},
\ref{medium-comparison}, \ref{large-middle-comparison}
and~\ref{effective-global-neumann}. Equality gives sector angles
\(2\pi/N\) and a common vertex radius \(\sqrt{2|P|/[N\sin(2\pi/N)]}\),
hence \(P\cong R_N\).

Scale \(R_N\) and \(R_{N+1}\) to the same area. Then \(R_N\) is
admissible in the comparison for \(N+1\) and has one fewer side
than the unique maximizer \(R_{N+1}\). Hence
\(\Phi(R_N)<\Phi(R_{N+1})\).
\end{proof}

\begin{proof}[Proof of Corollary~\ref{lattice-voronoi}]
Let \(V\) be the Voronoi cell of \(\Lambda\). If \(V\) and \(v+V\)
share an edge, reflection \(x\mapsto v-x\) preserves that edge, and its
midpoint \(v/2\) has exactly two closest lattice points, \(0\) and \(v\).
As in \cite[Lemma~2.2]{DL}, for \(w\equiv v\pmod{2\Lambda}\) with
\(w\ne\pm v\), the point \((v+w)/2\in\Lambda\setminus\{0,v\}\) gives
\[
 \left|\frac{v+w}{2}-\frac v2\right|=\frac{|w|}{2}>\frac{|v|}{2}.
\]
Thus \(\#\{\text{edges of }V\}\le2(\#(\Lambda/2\Lambda)-1)=6\) and
\(|V|=\operatorname{covol}\Lambda\).
By Theorem~\ref{main-polygonal-comparison}, equality holds exactly for
\(V\cong R_6\). The lattice vectors \(v_1,v_2\) whose bisectors
are adjacent sides then have length \(2\operatorname{inrad}(V)\)
and angle \(\pi/3\). Thus
\[
 [\Lambda:\mathbb Zv_1+\mathbb Zv_2]
 =\frac{4\operatorname{inrad}(V)^2\sin(\pi/3)}
 {\operatorname{covol}\Lambda}=1.
\]
Hence \(\Lambda\) is triangular. Conversely, its six shortest vectors
define a regular hexagonal cell.
\end{proof}

\appendix
\section{Eigenvalue enclosures and finite sums}\label{app:remaining}\label{app:finite-record}

We give the interval bounds, mesh calculations, and polynomial tests
used to enclose the regular-polygon eigenvalues. We also verify the five
radial inequalities used to exclude \(p(P)<N\) for \(5\le N\le8\).

\paragraph{Interval bounds.}
All terminating decimals are rational. For \(b\ge1\), use
\[
 [l,u]_b=\big[2^{-b}\lfloor2^bl\rfloor,\,2^{-b}\lceil2^bu\rceil\big],
 \qquad [l,u]^{-1}=[1/u,1/l]\quad(0\notin[l,u]).
\]
Use the extrema of the four endpoint products. On the dyadic grid,
integer square roots give
\[
 [\sqrt l,\sqrt u]_b
 =2^{-b}\big[\lfloor\sqrt{\lfloor 2^{2b}l\rfloor}\rfloor,\,
             \lceil\sqrt{\lceil 2^{2b}u\rceil}\rceil\big]
 \qquad(0\le l\le u).
\]
Use the specified decimal grid when indicated. For the elementary functions,
\[
 \begin{aligned}
 \left|e^x-\sum_{j=0}^{n}\frac{x^j}{j!}\right|
 &\le\frac{|x|^{n+1}}{(n+1)!}\frac1{1-|x|/(n+2)}
 &&(|x|<n+2),\\
 \left|\arctan x-\sum_{j=0}^{n}\frac{(-1)^j x^{2j+1}}{2j+1}\right|
 &\le\frac{|x|^{2n+3}}{2n+3} &&(|x|\le1),\\
 \left|\log y-2\sum_{j=0}^n\frac{z^{2j+1}}{2j+1}\right|
 &\le\frac{2|z|^{2n+3}}{(2n+3)(1-z^2)}
 &&\left(z=\frac{y-1}{y+1},\ y>0\right).
 \end{aligned}
\]
Use the alternating sine and cosine sums after argument reduction,
and \(\pi=16\arctan(1/5)-4\arctan(1/239)\). The interval tests use
\[
 \begin{aligned}
 \max_{a\le r\le b}H(r)
 &\le H\big((a+b)/2\big)+\tfrac{b-a}{2}\max_{a\le r\le b}|H'(r)|,\\
 \max_{0\le d\le D}H(r,d)&=\max\{H(r,0),H(r,D)\}\quad(H_{dd}\ge0).
 \end{aligned}
\]

\paragraph{Mixed eigenvalues.}
For \(N=6\) use \eqref{hex-rounded-schur}. In the other meshes of
Lemma~\ref{mesh-continuous-reference}, let \(x_i\) be the value on the edge
opposite \(\theta_i\) in \(T\), with \(x_i=0\) on Dirichlet edges. Set
\begin{equation}\label{nr22-lower-form}
 x^{\mathsf T} A x=\sum_T\left[
 2\sum_{\rm cyc}\lfloor 2^{60}\cot\theta_k\rfloor(x_i-x_j)^2-
 \left\lceil 2^{60}\sigma |T|/3\right\rceil\sum_{i=1}^3x_i^2\right].
\end{equation}
For a mesh triangle in a whole or halved sector, the cotangents and area are
\[
 \begin{cases}
 (\cot(2\pi/N),\tan(\pi/N),\tan(\pi/N)),\quad
    \sin(2\pi/N)/(2n^2), &\text{whole sector},\\
 (\cot(\pi/N),\tan(\pi/N),0),\quad
    \sin(2\pi/N)/(4n^2), &\text{halved sector}.
 \end{cases}
\]
In the halved sector the cotangents correspond to the center,
the vertex and the side midpoint. Outward rounding gives
\(2^{-60}A\preceq K_h-\sigma M_h\) and \(A_{ij}\le0\) for \(i\ne j\).
In the order \(N=5,7,8,9,10,11,12,13,14,15,16\), take
\begin{equation}\label{mixed-mesh-parameters}
 \begin{aligned}
 (n)&=(220,240,260,260,430,400,430,400,500,550,470),\\
 (\sigma)&=(4.354440,3.861670,3.746950,3.669787,3.615359,\,3.5754949,\\[-2pt]
 &\hspace{21mm}3.545404,3.5221220,3.50373125,3.48894518,3.4768758).
 \end{aligned}
\end{equation}
For all these matrices,
\begin{equation}\label{nr22-pivots}
 \det A_{1:k,1:k}>0\quad(1\le k\le\dim A),\qquad
 K_h-\sigma M_h\succeq 2^{-60}A\succ0.
\end{equation}
The integer recurrence, starting at \(A^{(0)}=A\) and \(d_0=1\), is
\[
 d_k=A^{(k-1)}_{kk}=\det A_{1:k,1:k},\qquad
 A^{(k)}_{ij}=
 \frac{d_k A^{(k-1)}_{ij}-A^{(k-1)}_{ik}A^{(k-1)}_{kj}}{d_{k-1}}
 \quad(i,j>k).
\]
For a sparse test, solve \(2^{-60}Az=\mathbf1\) approximately and set
\[
 w_i=\left\lfloor\frac{2^{56}z_i}{\max_j z_j}+\frac12\right\rfloor .
\]
At \(N=15\), apply one residual correction and round again:
\[
 2^{-60}A\delta=\mathbf1-2^{-116}(\max_jz_j)Aw,\qquad
 z\gets 2^{-56}(\max_jz_j)w+\delta .
\]
The checks \(\min_{N,i}w_i>4\cdot10^{12}\) and
\(\min_{N,i}(Aw)_i>3075\cdot2^{60}\) give \(A\succ0\) by
\eqref{discrete-positive-test}. Thus \(\lambda_{1,h}>\sigma\).

\begin{lemma}\label{mesh-continuous-reference}
Let \(R_N\) have circumradius one and take \((n,\sigma)\) in
\eqref{mixed-mesh-parameters}. For \(N=5\) or \(7\le N\le16\),
\begin{equation}\label{mesh-disk-lower}
 \mu_1(R_N)>
 \frac{\sigma}{1+\sigma\bigl(0.1893\max\{1,2\sin(\pi/N)\}/n\bigr)^2}.
\end{equation}
In particular,
\[
 \min_{5\le N\le16}
 \left(\Phi(R_N)-10.64986+\frac{21.2}{N^3}+\frac{398}{N^5}\right)
 >31\cdot10^{-6}.
\]
\end{lemma}
\begin{proof}
Take the upper half for odd \(N\) and \(N=8\), and the first quadrant
for even \(N\ge10\). The horizontal cut is Dirichlet and all other
edges are Neumann. Restriction of the first eigenspace with the
corresponding parities in Lemma~\ref{affine-all-N} gives an upper bound
\(\mu_1(R_N)\) for the mixed ground eigenvalue. Reflection across the
cuts, odd across the horizontal and even across the vertical, gives the
reverse inequality.

Join \(0\) to the boundary points in order.
\[
 \begin{cases}
 \cos(\pi/N),\ e^{(2j+1)\pi i/N}\quad(0\le j\le(N-1)/2),& N\text{ odd},\\
 e^{2j\pi i/N}\quad(0\le j\le N/2),& N=8,\\
 e^{2j\pi i/N}\quad(0\le j\le\lfloor N/4\rfloor),& N\ge10\text{ even}.
 \end{cases}
\]
Append \(i\cos(\pi/N)\) in the last case if \(N\equiv2\pmod4\).
Subdivision by parallels into \(n^2\) triangles gives
\[
 \operatorname{diam}T\le\frac{\max\{1,2\sin(\pi/N)\}}n .
\]
For the affine interpolant \(I_hu\) matching the integrals of \(u\) on
each edge, Liu~\cite[Section~4.3]{Liu} gives
\[
 \int_T\nabla(u-I_hu)=0,\qquad
 \|u-I_hu\|_{L^2(T)}
 \le0.1893\,\operatorname{diam}(T)\|\nabla(u-I_hu)\|_{L^2(T)} .
\]
Equations \eqref{nr22-pivots} and \eqref{broken-lower} yield \eqref{mesh-disk-lower}.

Since \(2\sin(\pi/5)=\sqrt{(5-\sqrt5)/2}<1.176\) and
\(2\sin(\pi/N)<1\) for \(N\ge7\), rational substitution gives
\[
 \begin{aligned}
 &\min_{\substack{5\le N\le16\\N\ne6}}\left[
 \frac{N\sigma}{2\left[1+\sigma\bigl(0.1893(1+0.176\,\mathbf 1_{\{N=5\}})/n\bigr)^2\right]}
 \sum_{j=0}^7\frac{(-1)^j(6.283185307179586/N)^{2j+1}}{(2j+1)!}
 \right.\\[-2pt]
 &\hspace{12em}\left.{}-10.64986+\frac{21.2}{N^3}+\frac{398}{N^5}\right]
 \in(0.0000314604,0.0000314605).
 \end{aligned}
\]
The minimum occurs at \(N=16\). Since \(6.283185307179586/N<2\pi/N<\pi/2\),
the alternating sum is bounded above by \(\sin(2\pi/N)\).
For \(N=6\), \eqref{reference-values} gives
\[
 \frac{3\sqrt3}{2}\,4.04320119-10.64986+\frac{21.2}{6^3}+\frac{398}{6^5}>0.004.
\]
\end{proof}

\begin{lemma}[Polynomial upper bounds on regular polygons]\label{lem:standalone-ritz}
Let \(R_N\) have circumradius one. Then \(\mu_1(R_N)<u_N\) for
\(N=5,7,8,\ldots,16\), where
\[
 \begin{aligned}
 10^8(u_5,u_7,u_8,u_9)&=(435448654,386174338,374695936,366984730),\\
 10^8(u_{10},u_{11},u_{12},u_{13})&=(361536280,357549775,354540940,352212633),\\
 10^8(u_{14},u_{15},u_{16})&=(350373460,348894701,347687950).
 \end{aligned}
\]
\end{lemma}
\begin{proof}
\emph{Polynomial form.}
Put \(c=\cos(\pi/N)\) and order the pairs \((p_i,q_i)\) with
\(p_i,q_i\ge0\), \(p_i-q_i\equiv1\pmod N\), by \((p_i+q_i,p_i)\),
starting with \((1,0)\). For \(f_i=(z/c)^{p_i}(\bar z/c)^{q_i}\),
rotation gives \(\int_{R_N}f_i=0\).
On \(S=\operatorname{conv}\{0,e^{\pm i\pi/N}\}\), use
\(z/c=r(1+is\tan(\pi/N))\) and \(dA/|S|=r\,dr\,ds\) for
\((r,s)\in[0,1]\times[-1,1]\). Thus
\[
 \begin{aligned}
 \frac1{|S|}\int_S(z/c)^a(\bar z/c)^b
 &=\int_0^1\int_{-1}^1r^{a+b+1}
 (1+is\tan(\pi/N))^a(1-is\tan(\pi/N))^b\,ds\,dr\\
 &=:J_{ab}\\
 &=\frac{2}{a+b+2}\sum_{h=0}^{\lfloor(a+b)/2\rfloor}
 \frac{(-1)^h(c^{-2}-1)^h}{2h+1}[X^{2h}](1+X)^a(1-X)^b
 \quad(a,b\ge0).
 \end{aligned}
\]
Use \(J_{ab}=0\) if \(\min(a,b)<0\). The matrix of the Hermitian form
\(c^2|R_N|^{-1}\int_{R_N}(|\nabla f|^2-u_N|f|^2)\) in this basis is
\[
 H_{ij}=2(p_i p_j+q_i q_j)J_{q_i+p_j-1,p_i+q_j-1}
        -u_Nc^2J_{q_i+p_j,p_i+q_j}.
\]

\emph{A negative direction.}
Keep the first \(m\) rows and columns, with
\[
 (m_5,m_7,m_8,m_9)=(16,30,105,39),\qquad
 (m_{10},\ldots,m_{16})=(184,340,220,289,271,378,253).
\]
Take \(\epsilon=2^{-768}\) and \([y]_-=\epsilon\lfloor y/\epsilon\rfloor\).
For the largest zero \(c\) of \(U_{N-1}\), use
\[
 \begin{gathered}
 U_0(x)=1,\quad U_1(x)=2x,\quad U_{j+1}(x)=2xU_j(x)-U_{j-1}(x),\\
 \cos(2\pi/N)\le1-8/N^2<1-6/N^2<c<1 .
 \end{gathered}
\]
Bisect \([1-6/N^2,1]\) to enclose \(c\) in \([[c]_-,[c]_-+\epsilon]\).
Enclose \(c^2,c^{-2}-1,u_Nc^2,H_{ij}\) in this order, rounding each
operation outward on \(\epsilon\mathbb Z\). In \(J_{ab}\) use Horner
evaluation in \(c^{-2}-1\), then multiply by \(2\) and divide by \(a+b+2\).
For the resulting enclosure \(I_{ij}\ni H_{ij}\), take
\[
 \begin{aligned}
 h_{ij}&=[\operatorname{mid}I_{ij}]_-
         -\left[\sum_{k<j}h_{ik}L_{jk}\right]_-
         &&(1\le j\le i\le m),\\
 L_{ij}&=[h_{ij}/h_{jj}]_- &&(j<i),\\
 x_m&=1,\qquad x_i=-\left[\sum_{k=i+1}^m L_{ki}x_k\right]_-
         &&(i=m-1,\ldots,1).
 \end{aligned}
\]
In the first two lines, increase \(i\), then \(j\) within each row.
Rational endpoint evaluation gives
\[
 x^{\mathsf T}Hx\le
 \max\left(\sum_{i=1}^m x_i^2 I_{ii}
    +2\sum_{1\le j<i\le m}x_i x_j I_{ij}\right)<-2^{-149}
 \qquad(N=5,7,8,\ldots,16).
\]
The real or imaginary part of \(\sum_i x_i f_i\) with the smaller
Rayleigh quotient gives \(\mu_1(R_N)<u_N\).
\end{proof}

\paragraph{Hexagonal angles.}
For \(0.017\le a\le0.263,\ t\ge0\), put \(p=3.14159265358979323847\).
Take the largest \(q\in10^{-10}\mathbb Z\cap[0,0.47]\) satisfying
\[
 \Bigl(\sum_{j=0}^{12}\frac{(-1)^j(pa)^{2j+1}}{(2j+1)!}\Bigr)
 \Bigl(\sum_{j=0}^{12}\frac{(-1)^j(pq)^{2j+1}}{(2j+1)!}\Bigr)
 \le t
 \Bigl(\sum_{j=0}^{13}\frac{(-1)^j(pa)^{2j}}{(2j)!}\Bigr)
 \Bigl(\sum_{j=0}^{13}\frac{(-1)^j(pq)^{2j}}{(2j)!}\Bigr).
\]
The cosine sums are positive for \(0\le pa,pq\le0.47p<3/2\), and
\(q=0\) when \(t=0\). Since \(p>\pi\), the alternating sums give
\[
 \tan(\pi a)\tan(\pi q)\le\tan(pa)\tan(pq)\le t
 \quad\Longrightarrow\quad q\le\pi^{-1}\arctan(t\cot(\pi a)).
\]
Bisect the side of greater relative length, taking the angular
coordinate at equality. The \(d\)-endpoints and mixings from
\eqref{hybrid-sector-test} give
\[
 \bigl([0.15,0.75],[0.75,1],[1,4],[0.75,0.78]\bigr)
 \longmapsto (394,1630,1847,8)\quad\text{terminal intervals}.
\]

\begin{numericallemma}[Five radial tests for the positive turns]\label{nr-finite-weights}
For every simple polygon \(P\),
\[
 5\le N\le8,\quad p(P)<N\quad\Longrightarrow\quad
 \Phi(P)<T_N<\Phi(R_N),\qquad
 T_N=10.64986-\frac{21.2}{N^3}-\frac{398}{N^5}.
\]
\end{numericallemma}
\begin{proof}
Lemma~\ref{mesh-continuous-reference} gives \(\Phi(R_N)>T_N\).
For the comparison with \(P\), we use a positive combination of radial
tests as in Dong and Wang \cite[Section~2]{DW}. Take
\[
 \psi(r)=\sum_{i=0}^6 a_i r^{2i+1},\qquad
 a_i=\frac{(-33899/40000)^i}{i!(i+1)!}
       \left(\sum_{j=0}^6\frac{(-33899/40000)^j}{j!(j+1)!}\right)^{-1}.
\]
The normalization and Bernstein coefficients give
\[
 \psi(0)=0,\qquad\psi(1)=1,\qquad
 \min_{0\le t\le1}\psi'(\sqrt t)\ge
 \min_{0\le i\le6}\sum_{j=0}^i(2j+1)a_j
   \frac{\binom ij}{\binom6j}>2\cdot10^{-5}.
\]
Use the conformal balance \eqref{balance} with this profile. In
\eqref{nr-ray-dipole},
\[
 \begin{aligned}
 E&=\sum_{i,j=0}^6 a_i a_j\frac{(2i+1)(2j+1)+1}{2(i+j+1)},\\
 d_n&=\sum_{i,j=0}^6 a_i a_j\frac{n}{n+i+j+1}
     =\int_0^1\psi(s^{1/(2n)})^2\,ds\le d_{n+1}.
 \end{aligned}
\]
For \(0\le i\le4\), take
\[
 \begin{aligned}
 v_i&=\left(\frac{2^i(N-1)}{1+2^i(N-1)}\right)^8,\qquad
 (c_i)=(0.0016,0.0315,0.236,0.873,2.312),\\
 (N-1)^3(l_i)&=(47.72,\ 13.476+0.794(N-5),\ 3.2999,\ 0.49425,\
                      0.10031-0.007697(N-5)).
 \end{aligned}
\]
The coefficient and sum bounds are
\begin{align}
 (k+1)\sum_{i=0}^4 l_i(k/2-c_i)^2v_i^k
   &<d_{k+1}-2\pi E/T_N-5\cdot10^{-5},\qquad k\ge1,\label{nr-finite-coeff}\\
 \sum_{i=0}^4 l_i\big[(1+c_i)^2v_i^{N-1}-c_i^2\big]
   &>2\pi E/T_N-d_1+10^{-4}.\label{nr-finite-objective}
\end{align}
Use \(2\pi<6.283185308\) and
\[
 P_{i,0}=10^{18},\qquad
 P_{i,k+1}=\lceil v_iP_{i,k}\rceil\ge10^{18}v_i^{k+1}.
\]
Rational substitution gives
\[
 \begin{aligned}
 \min_{\substack{5\le N\le8\\1\le k\le50}}
 \biggl(d_{k+1}-\frac{6.283185308E}{T_N}
 -\frac{k+1}{10^{18}}\sum_{i=0}^4 l_i(k/2-c_i)^2P_{i,k}\biggr)
 &>5\cdot10^{-5},\\
 \min_{5\le N\le8}
 \biggl(d_1-\frac{6.283185308E}{T_N}
 +\sum_{i=0}^4 l_i\big[(1+c_i)^2v_i^{N-1}-c_i^2\big]\biggr)
 &>10^{-4}.
 \end{aligned}
\]
The minima occur at \((N,k)=(5,5)\) and \(N=5\), respectively.
By Bernoulli's inequality, \(1-v_i\ge8/(2^i(N-1)+8)\ge1/15\), so
\[
 \frac{d}{dt}\log\big((t+1)(t/2-c_i)^2v_i^t\big)
 \le\frac1{51}+\frac2{50-2(2.312)}-\frac1{15}<0\qquad(t\ge50).
\]
Together with \(d_{k+1}\le d_{k+2}\), this extends \eqref{nr-finite-coeff}
to \(k>50\). For \(p(P)<N\), use \eqref{nr-turn-point} at \(m=N-1\) and
\eqref{nr-finite-coeff}--\eqref{nr-finite-objective}:
\[
 \begin{aligned}
 \sum_{k\ge1}(d_{k+1}-2\pi E/T_N)\frac{|b_k|^2}{k+1}
 &\ge\sum_{i=0}^4l_i\sum_{k\ge1}(k/2-c_i)^2|b_k|^2v_i^k\\
 &\ge\sum_{i=0}^4l_i[(1+c_i)^2v_i^{N-1}-c_i^2]
 >2\pi E/T_N-d_1 .
 \end{aligned}
\]
By \eqref{nr-ray-dipole}, \(\Phi(P)<T_N\).
\end{proof}

\bigskip
\noindent
{\bf Acknowledgements.}
The authors used AI tools to assist with explorations and the long but finite computations.

\end{document}